\documentclass[11pt,a4paper]{article}

\usepackage{amsmath,amsfonts}
\usepackage{algorithmic}
\usepackage{algorithm}
\usepackage{array}
\usepackage[caption=false,font=normalsize,labelfont=sf,textfont=sf]{subfig}
\usepackage{textcomp}
\usepackage{stfloats}
\usepackage{url}
\usepackage{verbatim}
\usepackage{graphicx}
\usepackage{cite}
\usepackage{amsmath,amssymb,amsthm,cite, graphicx,pdfpages,todonotes}

\usepackage[utf8]{inputenc}
\usepackage{amsmath}
\usepackage{float}
\usepackage{array}
\usepackage{tabularx}
\usepackage{color}
\usepackage{caption}
\usepackage{multirow}
\usepackage{amssymb}

\begin{document}
\begin{center}
		{\large\bf{A comparison of rational and neural network based approximations}}
	\end{center}
	
	\begin{center}
		{Vinesha Peiris, Reinier D\'iaz Mill\'an, Nadezda Sukhorukova, Julien Ugon}
	\end{center}

\begin{abstract}
Rational and neural network based approximations are efficient tools in modern approximation. These approaches are able to produce accurate approximations to nonsmooth and non-Lipschitz functions, including  multivariate domain functions. In this paper we compare the efficiency of function approximation using rational approximation, neural network and their combinations. It was found that rational approximation is superior to neural network based approaches with the same number of decision variables. Our numerical experiments demonstrate the efficiency of rational approximation, even when the number of approximation parameters (that is, the dimension of the corresponding optimisation problems) is small. Another important contribution of this paper lies in the improvement of rational approximation algorithms. Namely, the optimisation based algorithms for rational approximation can be adjusted to in such a way that the conditioning number of the constraint matrices are controlled. This simple adjustment enables us to work with high dimension optimisation problems and improve the design of the neural network.  The main strength of neural networks is in their ability to handle models with a large number of variables: complex models are decomposed in several simple optimisation problems.  Therefore the the large number of decision variables is in the nature of neural networks.
\end{abstract}



\section{Introduction}
At the core of the efficiency and popularity of neural networks is the universal approximation theorem~\cite{LeshnoPinkus1993}, which states that deep neural networks can approximate any continuous function to any accuracy, providing that the activation functions are nonpolynomials. This powerful theoretical result  paves the way to investigating which activation functions work best in practice.

Of particular interest to this paper, \cite{boulle2020rational} suggests rational activation functions. Rational functions are attractive because their parameters can themselves be learnt by the network, and it was shown by~\cite{boulle2020rational} that good approximation results can be achieved. A comprehensive literature review on uniform approximation by rational function approximation, its geometrical properties and its connection with polynomial approximation can be found in~\cite{Alimov}. 

In the present paper, we build on this work and make some additional observations. First, rational functions themselves can approximate any function to any accuracy. Consequently, it is possible to replace any other activation function with a rational approximation of it. This leads to the question of whether it makes sense to do so, that is, whether the networks built on an approximate rational activation work similarly to networks built on the actual activation function. In turn, as shown by~\cite{boulle2020rational}, the rational coefficients can themselves be learnt, which would make the networks with rational approximations fundamentally more powerful than other networks.

The second observation is that the family of rational functions is closed under addition, multiplication and composition. A neural network with a rational activation function is itself a rational function. This leads to the question of how they compare to a classical approach to rational approximation. The \emph{theoretical} answer to this is obviously that classical rational approximation is better because the search space is a strict superset of the neural network approach. It remains to be investigated whether the practical difference is noticeable. This is the second purpose of this paper.

To achieve these objectives, in this paper, we compare two techniques for function approximations: rational approximation and neural network based approximation. We also consider some combinations of these two approaches in order to improve the performance of the learning system. In particular, we propose a new method for optimising the parameters of a neural network with rational activation functions using an alternating approach where one set of parameters is optimised while the other parameters are fixed and vice versa.

The neural network we are using in this study has a specific structure and is also known as deep learning. 
Deep learning is one of the key tools in the modern area of Artificial Intelligence. There are many practical applications for deep learning,  including data analysis, signal and image processing and many others~\cite{Goodfellow2016, Sun2020OptimDeepLearning}. Despite all these practical applications, deep learning is just a specific type of function approximation, where the function to be approximated is only known by its values in a finite number of points (training data), while the approximation prototype (class of approximations) is a composition of affine mappings and the so called activation functions (special type of univariate function).  A very comprehensive and thorough textbook on the modern view of deep learning can be found in~\cite{Goodfellow2016}. 

Our second main approach is rational approximation. In this case, the approximations are simply the ratios of polynomial functions. The choice of rational functions is quite natural: they provide a flexible approximation to a wide range of functions, including non-smooth and non-Lipschitz function~\cite{loeb1957rational, Rivlin1962, Ralston1965Reme}. This flexibility is even comparable with free-knot piecewise polynomial approximation~\cite{shevshevPopov1987}. At the same time, the corresponding optimisation problems are quasiconvex~\cite{loeb1960, SL, AMCPeirisSukhSharonUgon} and there are a number of efficient computational methods to tackle them~\cite{MLegazquasiconvexduality, JPCrouzeix1980quasi, DaCruzAlgorithmsQuasiconvex} just to name a few. 

Our numerical experiments highlight differences between the two methods. The direct rational approximation approach offers greater control over the rational functions, the objective and therefore ultimately produces more accurate approximations than neural networks with rational activation functions. 

Apart from comparing different approximation techniques, we also implemented an additional improvement of the rational approximation algorithm. The idea of the improvement is simple and straightforward to implement as part of our approach and increases the stability of the results. The idea is to restrict the condition number of matrices appearing in the auxiliary linear programming problems within rational approximation. This extra restriction is a simple linear constraint, which can be naturally added to the constraints of the linear programs. A similar approach was proposed in~\cite{SharonPeirisSukhorukovaUgon} for lifting matrix functions. In this paper, we applied this idea to improve rational approximation algorithms, especially in the case of multi-dimensional approximation. The results of the numerical experiments are very promising.


The paper is organised as follows. In Section~\ref{sec:preliminaries} we provide the background of the approximation methods we use in this paper: neural network, rational approximation and their combinations. In Section~\ref{sec:models}  we explain the approximation models we are using in this paper. Then in Sections~\ref{sec:resultsNN}-\ref{sec:results2D} we compare the models. Finally, in Section~\ref{sec:conclusions} we draw the conclusions and highlight possible future research directions. Appendix~\ref{appendix:resultsNN} contains supplementary results, including the results with only two nodes in the hidden layer. The approximations with only two nodes in the hidden layer are not accurate, but we include these results for consistency and comparison.

\section{Preliminaries}\label{sec:preliminaries}
\subsection{Deep learning}

Deep learning is a powerful tool for data and function approximation, which can also be used for data analysis and data classification. The power of deep learning is in its structure: the subproblems that the system has to solve are very simple from the point of view of modern optimisation theory and therefore the system can handle problems with hundreds and even thousands of decision variables very efficiently.

The origin of deep learning is mathematical in nature.  Essentially, the objective of deep learning is to solve an approximation problem: optimise the weights (parameters) of the network. These weights are the decision variables of certain optimisation problems, whose objective functions represent inaccuracy of approximation.  Therefore, it is natural to approach this problem using modern optimisation tools~\cite{Sun2020OptimDeepLearning, VidalHaeffele17, Goodfellow2016}.

The solid mathematical background of deep learning was established in~\cite{Cybenko, Hornik1991, LeshnoPinkus1993, pinkus_1999}. These works  rely on the results of the celebrated Kolmogorov-Arnold Theorem~\cite{Kol57,Arnold57}. The Kolmogorov–Arnold representation theorem states that every multivariate continuous function can be represented as a composition of continuous univariate functions over the binary operation of addition. In general, there is no algorithm for constructing these composition functions. Instead of constructing this representation,  modern deep learning techniques approximate this composition function by a composition of affine transformations and the so-called activation functions. The activation functions are univariate non-polynomial functions. The most commonly used activation functions are sigmoid functions and Rectified Linear Unit (ReLU) functions (monotonic piecewise linear function $\varphi(x)=\max\{0, x\}$). 

\cite{boulle2020rational} shows an equivalence between neural networks using rational approximation functions and ReLU networks, see also~\cite{telgarsky2017neural}. Specifically they show that any ReLU network can be approximated to within $\varepsilon>0$ by a rational network of size $\mathcal{O}(\log(\log(1/\varepsilon))$ and any rational network can be approximated to within $\varepsilon$ of size $\mathcal{O}(\log(1/\varepsilon))^3$. They also implemented numerical experiments to compare the two types of rational activation. An interesting characteristic of their approach is that the parameters of the activation function can themselves be optimised as part of the learning process. Their work opens the question of how using a direct approach to rational approximation compares with neural network-based approximation.

Most deep learning algorithms rely on least squares based measures of inaccuracy (loss function). Since loss functions are measures of inaccuracy, the goal of optimisation is to optimise the weights by minimising the corresponding loss function. Therefore, least squares-based models are minimising a smooth quadratic function, and therefore fast and simple optimisation techniques (for example, gradient descent) are applicable. In some cases, however, other loss functions are more efficient: uniform based, etc.

The goal of this paper is to compare the approximation results obtained by deep learning and other approximation techniques. In particular, we are looking at rational approximation due to its approximation power~\cite{shevshevPopov1987}. More details on rational approximation will be provided in the next section. On the other hand, in~\cite{telgarsky2017neural}, the author demonstrates that rational functions are as powerful as neural networks with standard ReLU activation functions. This result encouraged many researchers in the deep learning community to combine neural networks, rational functions and rational approximation techniques together to enhance the performance of each other~\cite{boulle2020rational, molinapade, vinesha2022}.

\subsection{Rational approximation}

Rational approximation in Chebyshev (uniform) sense was a very popular research direction in the 1950s-70s~\cite{Achieser1965,Boehm1964, cheney1964generalized, Meinardus1967rational, Ralston1965Reme, Rivlin1962} and many others.  

There are two main groups of methods to approach rational approximation. The first group~\cite{Trefethen2018} is dedicated to ``nearly optimal'' solutions. This approach (also known as the AAA approach) is very efficient and therefore very popular, but it is ``nearly optimal''. The  extension of AAA to multivariate cases is still open. Moreover, this method is only designed for unconstraint optimisation and, as we will see later in this paper, this will limit our ability to work with ill-conditioned matrices.   The second group of methods is based on modern optimisation techniques:  the corresponding optimisation problems are quasiconvex and can be solved using a general quasiconvex optimisation method. 

There are many methods for rational approximation~\cite{lee1973, DiffCorrection1972, osborne1969algorithm, Ralston1965Reme} (just to name a few), but the implementation of most of them relies on solving linear programming problems. Therefore, additional linear constraints can be easily added to these implementations without making the problems complex. Currently, the most popular optimisation methods for rational approximation are the bisection method for quasiconvex functions and the differential correction method. The advantage of these two methods is that these methods can be easily extended to the case of non-monomial cases (generalised rational approximation)~\cite{cheney1964generalized, AMCPeirisSukhSharonUgon} and to multivariate settings~\cite{AghiliSukhorukovaUgon}. Overall, the differential correction method has a quadratic convergence~\cite{DiffCorrection1972}, while the bisection method converges linearly~\cite{SL}. Therefore, in our experiments with univariate rational approximation, we use the differential correction method. At the same time, the bisection method is still a good choice due to its simplicity. Hence, we use the bisection method in our multivariate rational approximation experiments. 



\section{Computational models}\label{sec:models}
We compare six different approaches for  approximating a given continuous function. Four of them are based on neural network (NN) and the remaining two are purely rational approximation based approaches. The approaches are as follows:
\begin{enumerate}
    \item NN with ReLU activation function,
    \item NN with rational approximation to ReLU activation (we apply the differential correction method to approximate ReLU),
    \item NN with rational activation function, the coefficients of the rational activation function are learnt from NN,
    \item NN with rational activation, the coefficients of the rational activation function and the parameters of the network are learnt with split method, where the training process is done in three steps: between the input layer and the hidden layer and then between the hidden layer and the output layer, and finally, the coefficients of the rational activation, 
    \item Rational approximation (Differential correction), this is the direct rational approximation between the input and output layer,
    \item Rational approximation (AAA), this is the direct ``near-optimal'' rational approximation between the input and output layer.
\end{enumerate}

We use neural networks with three layers: input layer, hidden layer and output layer. We have the following options for the number of nodes in the hidden layer:
\begin{itemize}
    \item a network with 2 nodes in the hidden layer (see Appendix~\ref{appendix:resultsNN}),
    \item a network with 10 nodes in the hidden layer.
\end{itemize}
In our experiments, we use different numbers of epochs: 50, 100 and 200. We record the training time per epoch for each network training procedure.

For NN methods, we the MSE-based loss function (MSE stands for Mean Squared Error and simply means ``least squares'') and  uniform loss. For the MSE-based loss we use ADAM, while for the uniform loss we use ADAMAX (a version of ADAM specially designed for uniform loss). All the rational based approximation models are designed for the uniform loss, which is a stronger criterion than least squares and clearly more appropriate if the goal is to approximate a function. 

In our experiments, we consider two different settings.
\begin{itemize}
\item The domain is a segment (univariate function).
We approximate the nonsmooth function 
$$f(x) = \sqrt{|x-0.25|}, \quad x \in [-1,1],$$
by using the above approaches with different conditions (number of nodes in the hidden layer, number of epochs, etc.). The choice of the function $f(x)$ is due to its difficulty for most approximation techniques: this function is nonsmooth and non-Lipschitz at $x=0.25$.
\item The domain is bi-variate. It is also possible to extend the results to higher dimensions, but all these extensions are out of the scope of this paper.
\end{itemize}



In the next section, we only present the important results of the numerical experiments. In particular, we report the results related to neural networks with 10 nodes in the hidden layer where the loss function is in the form of uniform norm. The Appendix~\ref{sec:Appendix} contains a very thorough discussion on the remaining results. 

\section{Results: Neural network-based approximation.}\label{sec:resultsNN}

The findings in this section are from a neural network with 10 nodes in the hidden layer. The loss function is based on the uniform norm and we use ADAMAX as the optimiser.

\subsection{Neural Network with ReLU activation}\label{subsec:results1:Set4}

We start our experiments with the standard ReLU activation function. Table~\ref{tab:Results: experiments set 4 of NN with ReLU activation} demonstrates the loss function value and computational time for different epochs. 

One can see that the value of the loss function does not improve significantly when the number of epochs is increasing. Moreover, it appeared that the minimum reported value of the loss function may take place before the final epoch. For example, in the case of~200 epochs, the minimal value (0.088032) was obtained at epoch 187, while the reported value (after 200 epochs)
is 0.138172. The last column of this table reports the running time per epoch and there is no significant difference between the computational time for each epoch.

\begin{table}
    \centering
    \begin{tabular}{|c|c|c|c|}
    \hline
    Epoch & Final loss & Minimum loss & Run time (per epoch) \\
    \hline
    50  & 0.130828 &          & 2.52s $\pm$ 79.3ms\\
    45  &          & 0.117747 & \\
    \hline
    100 & 0.139669 &          & 2.51s $\pm$ 96.8ms\\
    97  &          & 0.111249 & \\
    \hline
    200 & 0.138172 &          & 2.49s $\pm$ 59.2ms\\
    187 &          & 0.088032 & \\
    \hline
    \end{tabular}
    \caption{Results: experiments of NN with ReLU activation}
    \label{tab:Results: experiments set 4 of NN with ReLU activation}
\end{table}



Figure~\ref{fig:Relu activation - set 4, epoch is 50,100,200} shows the approximations computed by the network with ReLU activation. The accuracy of the approximations is significantly better, compared to the networks with fewer nodes in the hidden layer. It also appeared that the approximation around the ``difficult point'' $x=0.25$ is better in the case of the uniform loss than it is in the case of MSE (Figures for MSE can be found in the Appendix~\ref{sec:Appendix}). Therefore, the uniform loss may be a better measure when it is essential to obtain accurate approximations around such points. 

\begin{figure}
    \centering
    \includegraphics[width=40mm]{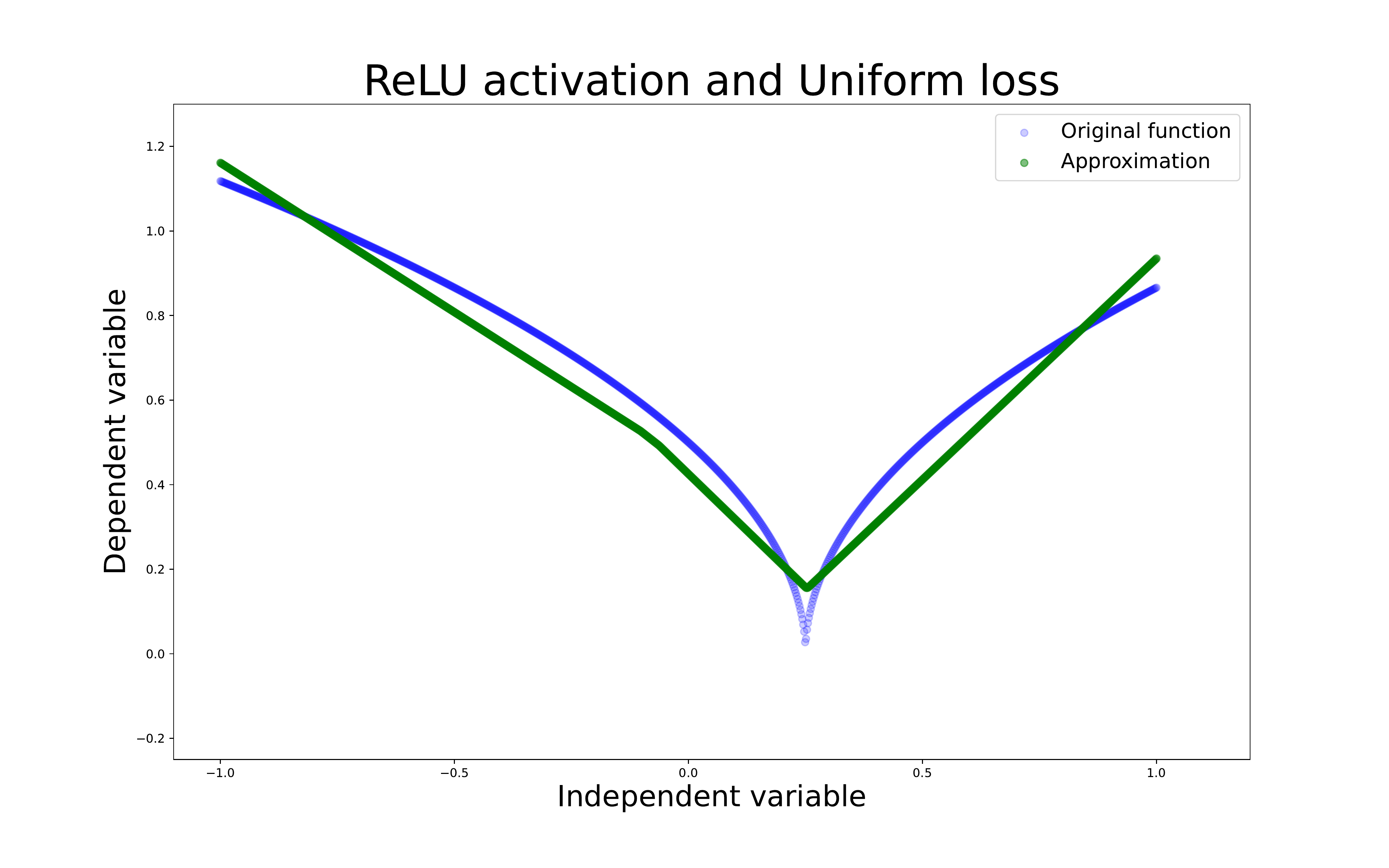}
    \includegraphics[width=40mm]{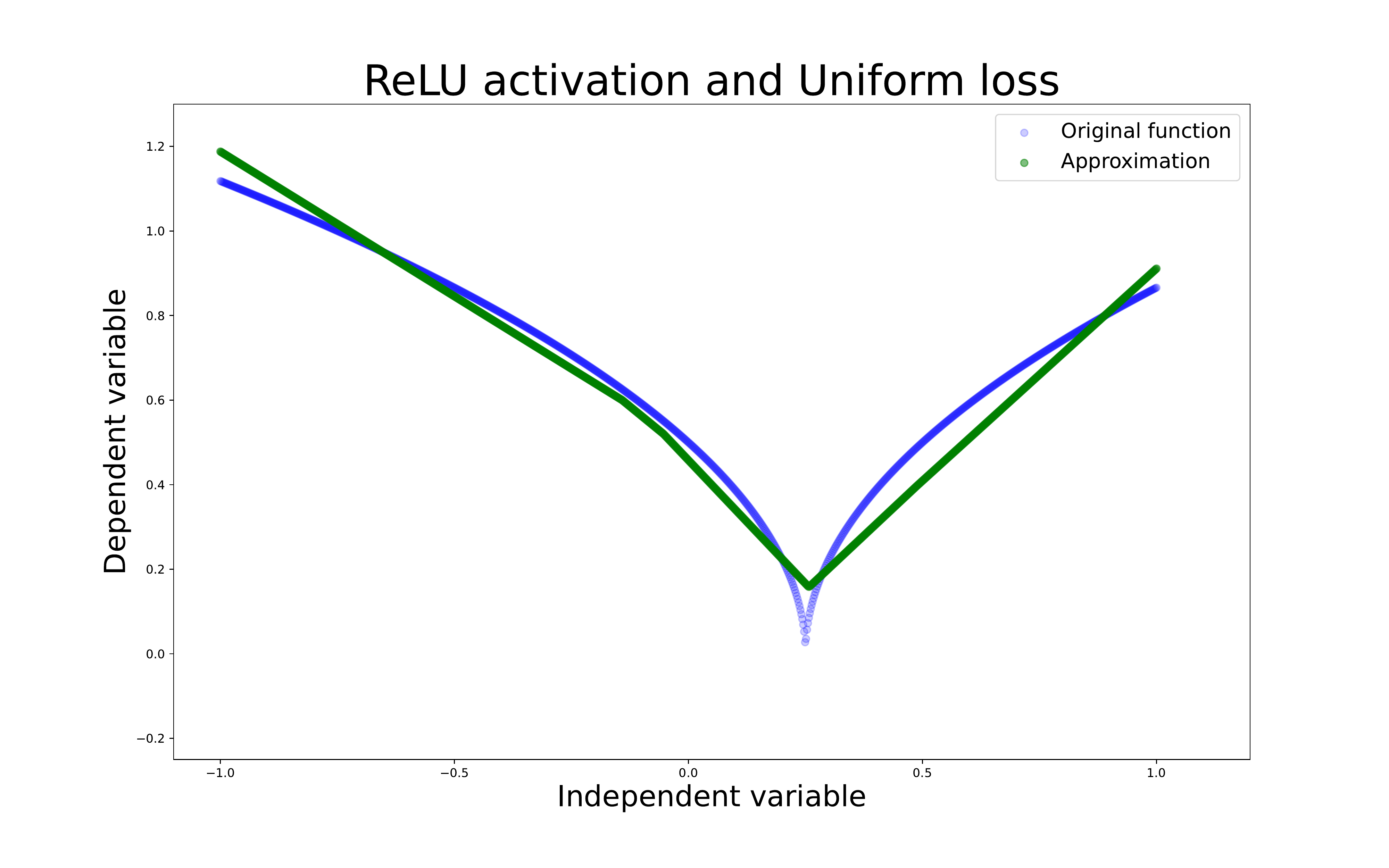}
    \includegraphics[width=40mm]{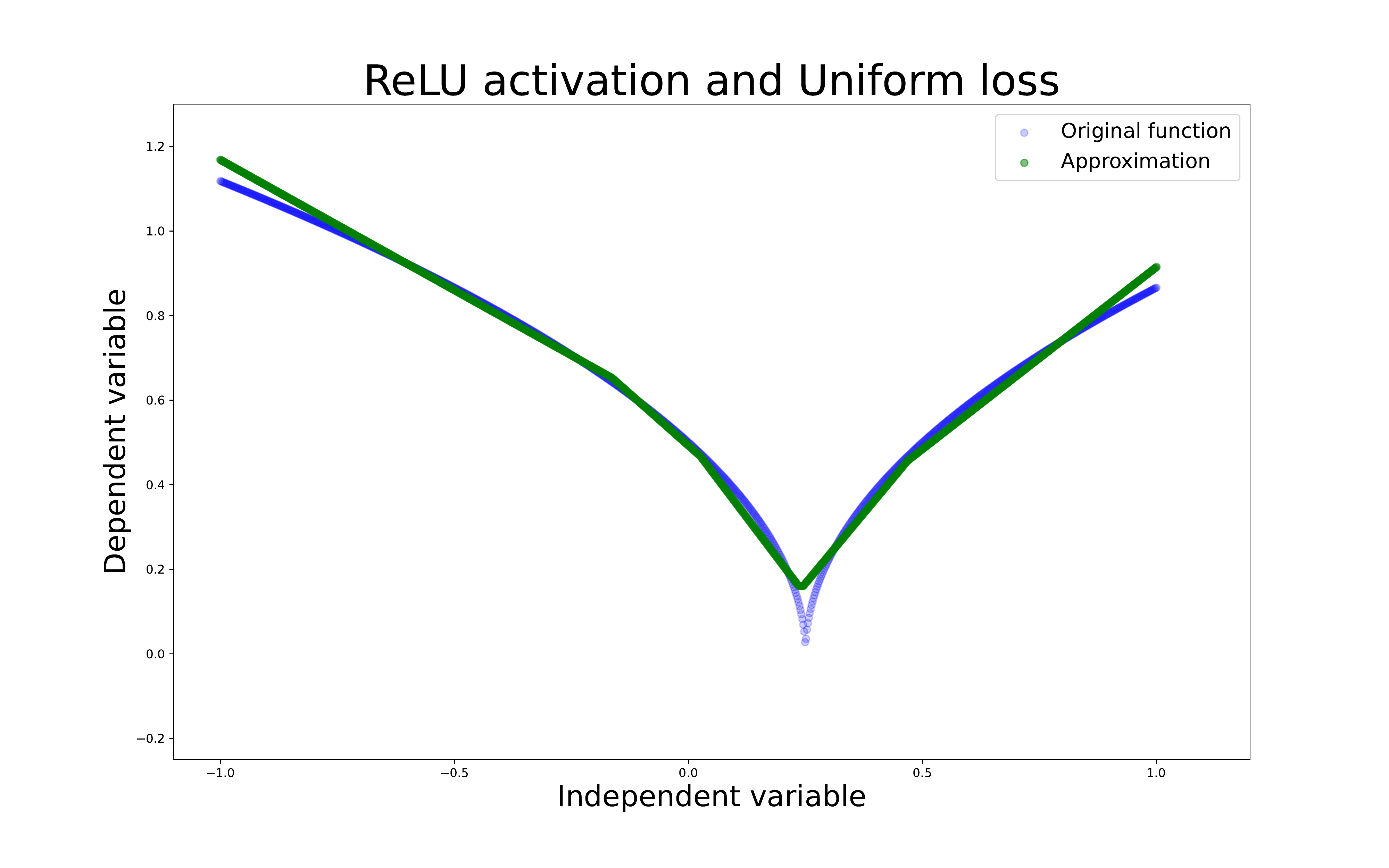}
    \caption{Approximation is computed by a neural network with ReLU activation,  50, 100 and 200~epochs.}
    \label{fig:Relu activation - set 4, epoch is 50,100,200}
\end{figure}

\subsection{Neural Network with rational approximation to ReLU activation}\label{subsec:results2:Set4}

In this section, we present the results of numerical experiments, where the experiment settings are similar to those in Section~\ref{subsec:results1:Set4}, but the activation function is the rational approximation to ReLU. This approximation was found by the differential correction method. 


Our activation function is a rational function of degree $(3,2)$: the rational approximation is the ratio of two polynomials, the degree of the numerator is~3 and the degree of the denominator is~2.  

The coefficients of the rational activation function come from the best rational $(3,2)$ approximation to the ReLU function. These coefficients are fixed throughout the whole training procedure. We do not learn the coefficients with the rest of the parameters of the network.

Table~\ref{tab:Results: experiments set 4 of NN with rational approximation to ReLU activation} shows the  improvement in the accuracy compared to the standard ReLU activation in Section~\ref{subsec:results1:Set4}. This observation is especially prominent when the number of epochs is~100 or~200. 

\begin{table}
    \centering
    \begin{tabular}{|c|c|c|c|}
    \hline
    Epoch & Final loss & Minimum loss & Run time (per epoch) \\
    \hline
    50  & 0.175439 &          & 2.82s $\pm$ 90.5ms\\
    46  &          & 0.125264 & \\
    \hline
    100 & 0.105346 &          & 2.77s $\pm$ 79.5ms\\
    99  &          & 0.071646 & \\
    \hline
    200 & 0.084491 &          & 2.83s $\pm$ 382ms\\
    195 &          & 0.053663 & \\
    \hline
    \end{tabular}
    \caption{Results: experiments of NN with rational approximation to ReLU}
    \label{tab:Results: experiments set 4 of NN with rational approximation to ReLU activation}
\end{table}



Figure~\ref{fig:Rational activation with ReLU coefficients - set 4, epoch is 50,100,200} shows that the accuracy of the approximation is improving, even for the ``difficult point'', especially when the number of epochs is~100 or more.

\begin{figure}
    \centering
    \includegraphics[width=40mm]{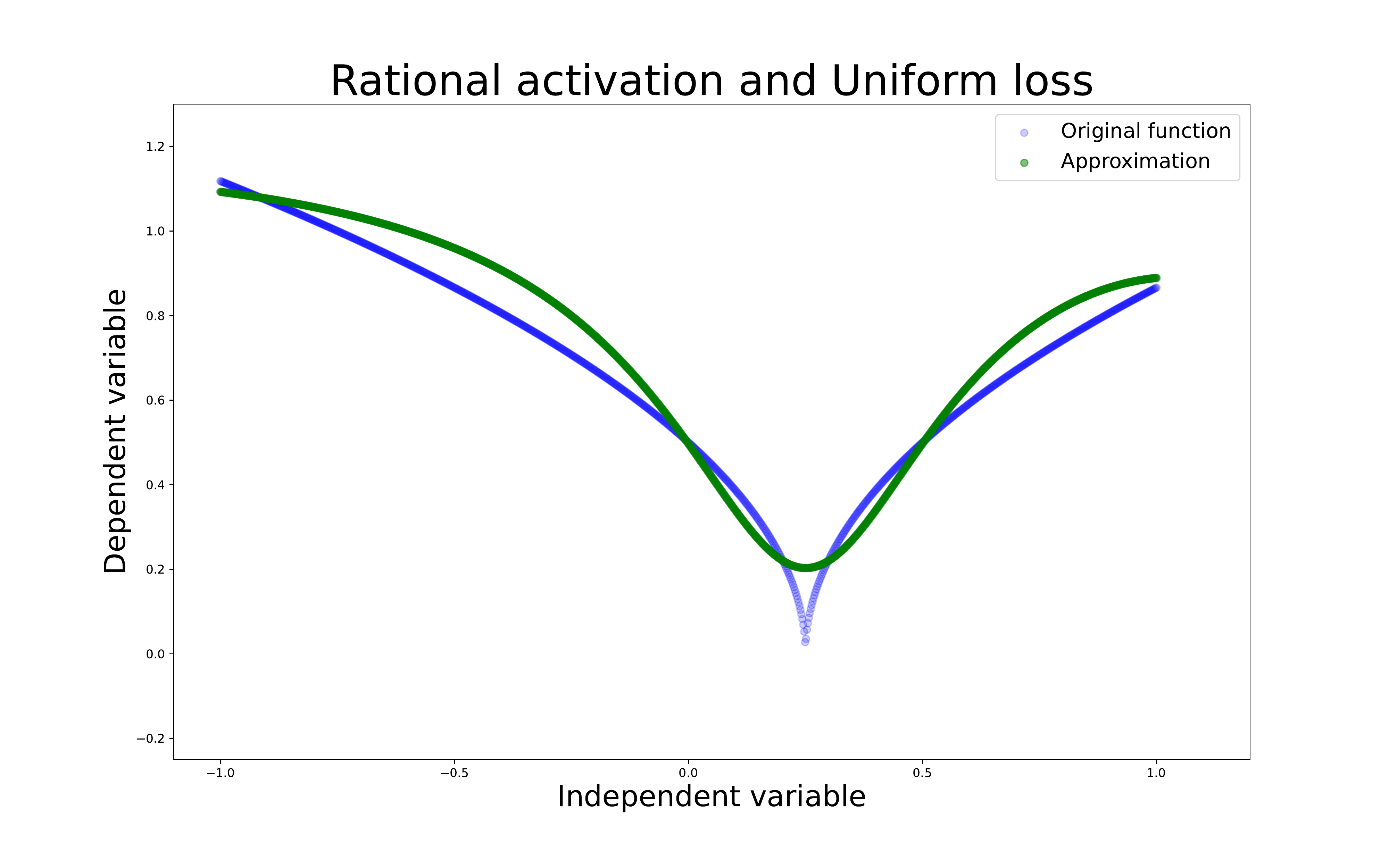}
     \includegraphics[width=40mm]{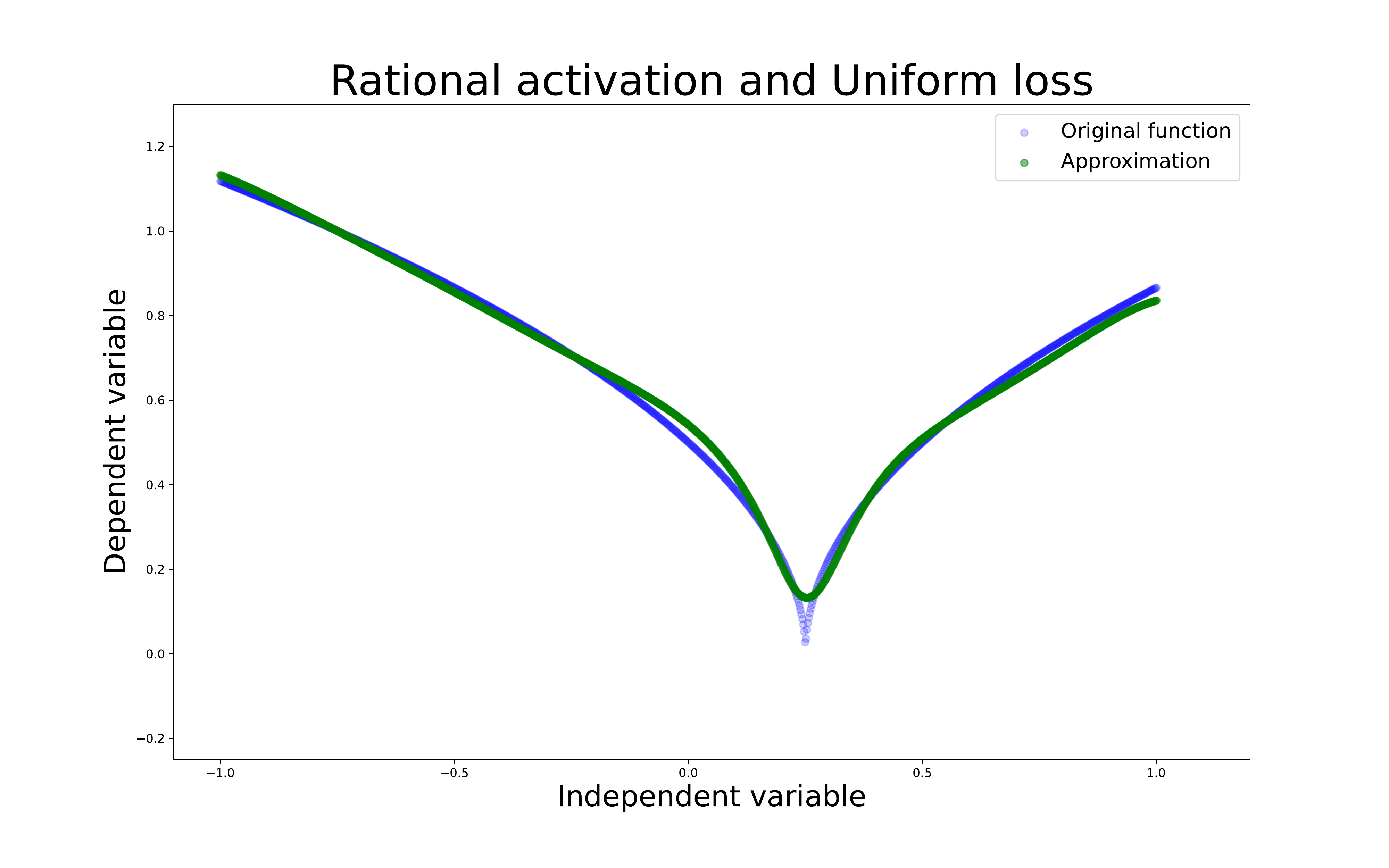}
      \includegraphics[width=40mm]{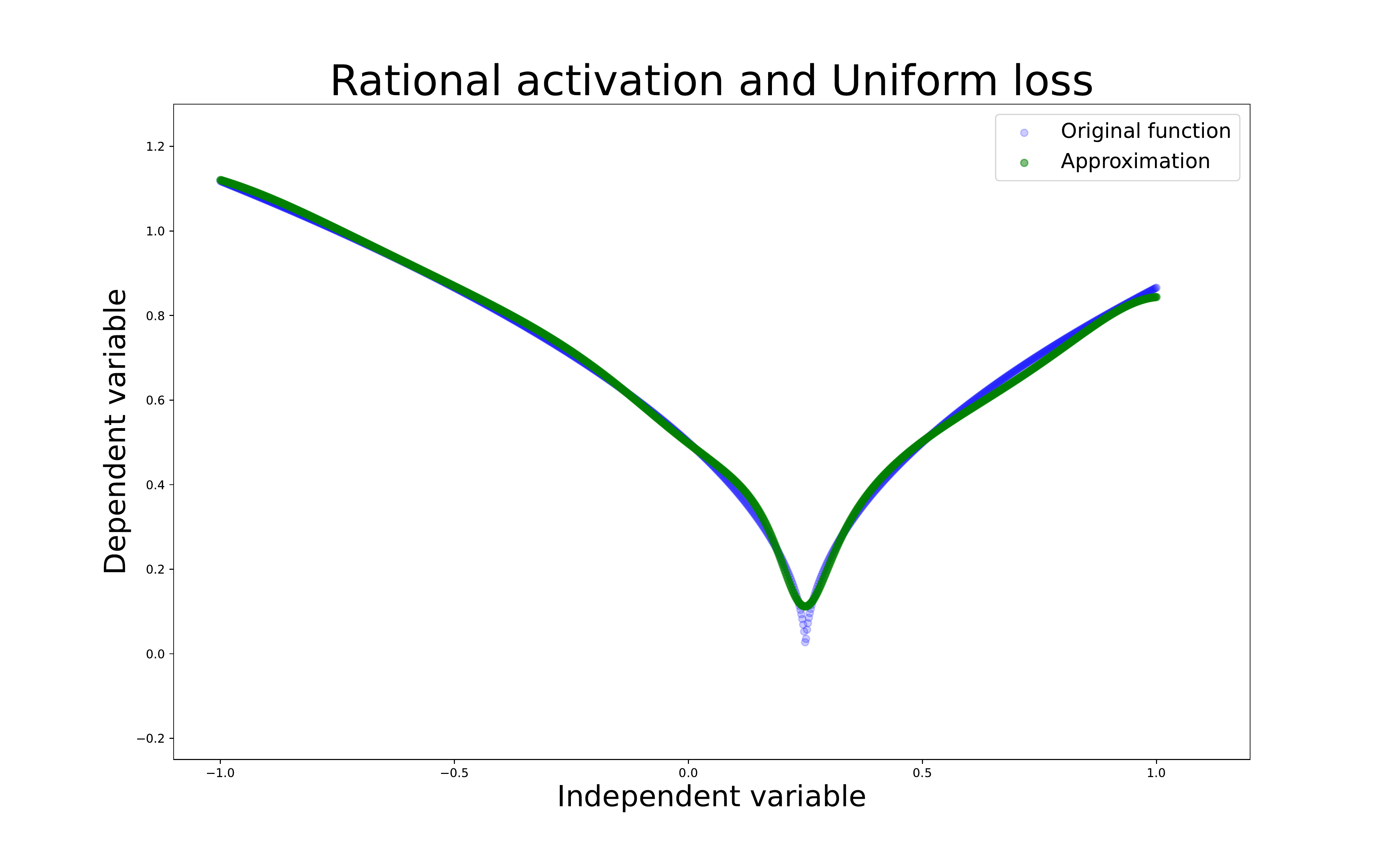}
    \caption{Approximation is computed by a neural network with rational activation to ReLU, 50, 100 and 200 epochs.}
    \label{fig:Rational activation with ReLU coefficients - set 4, epoch is 50,100,200}
\end{figure}

Overall conclusion: the approximation results are more accurate when the rational approximation to ReLU is used as the activation function rather than  ReLU itself. This observation is valid for both MSE and uniform loss. Results related to MSE can be found in the Appendix~\ref{sec:Appendix}


\subsection{Neural Network with rational activation}\label{subsec:results3:Set4}


In this case, our activation function is a rational function of degree $(3,2)$. 
The coefficients of the rational activation function are now a part of the parameter set. We learn these coefficients as we learn other parameters during the training procedure. This type of network with rational activation function whose coefficients are learnable parameters are called `rational neural networks' and more details can be found in~\cite{boulle2020rational}. The Python code for this section is also from~\cite{boulle2020rational}.

Table~\ref{tab:Results: experiments set 4 of NN with rational approximation and usual training} shows that when the number of nodes in the hidden layer is~10, the optimal loss function values are very close to the case when ReLU was approximated by the rational function. Figure~\ref{fig:Rational activation with usual training - set 4, epoch is 50,100,200} confirms these results. 

\begin{table}
    \centering
    \begin{tabular}{|c|c|c|c|}
    \hline
    Epoch & Final loss & Minimum loss & Run time (per epoch) \\
    \hline
    50  & 0.131586 &          & 2.75s $\pm$ 77.9ms\\
    44  &          & 0.114195 & \\
    \hline
    100 & 0.104754 &          & 2.64s $\pm$ 62.1ms\\
    96  &          & 0.097785 & \\
    \hline
    200 & 0.078658 &          & 2.76s $\pm$ 423ms\\
    177 &          & 0.055616 & \\
    \hline
    \end{tabular}
    \caption{Results: experiments with rational activation function}
    \label{tab:Results: experiments set 4 of NN with rational approximation and usual training}
\end{table}



The approximation appears to be accurate even around the ``difficult point'' when the number of epochs is~200 (similar to MSE).

\begin{figure}
    \centering
    \includegraphics[width=40mm]{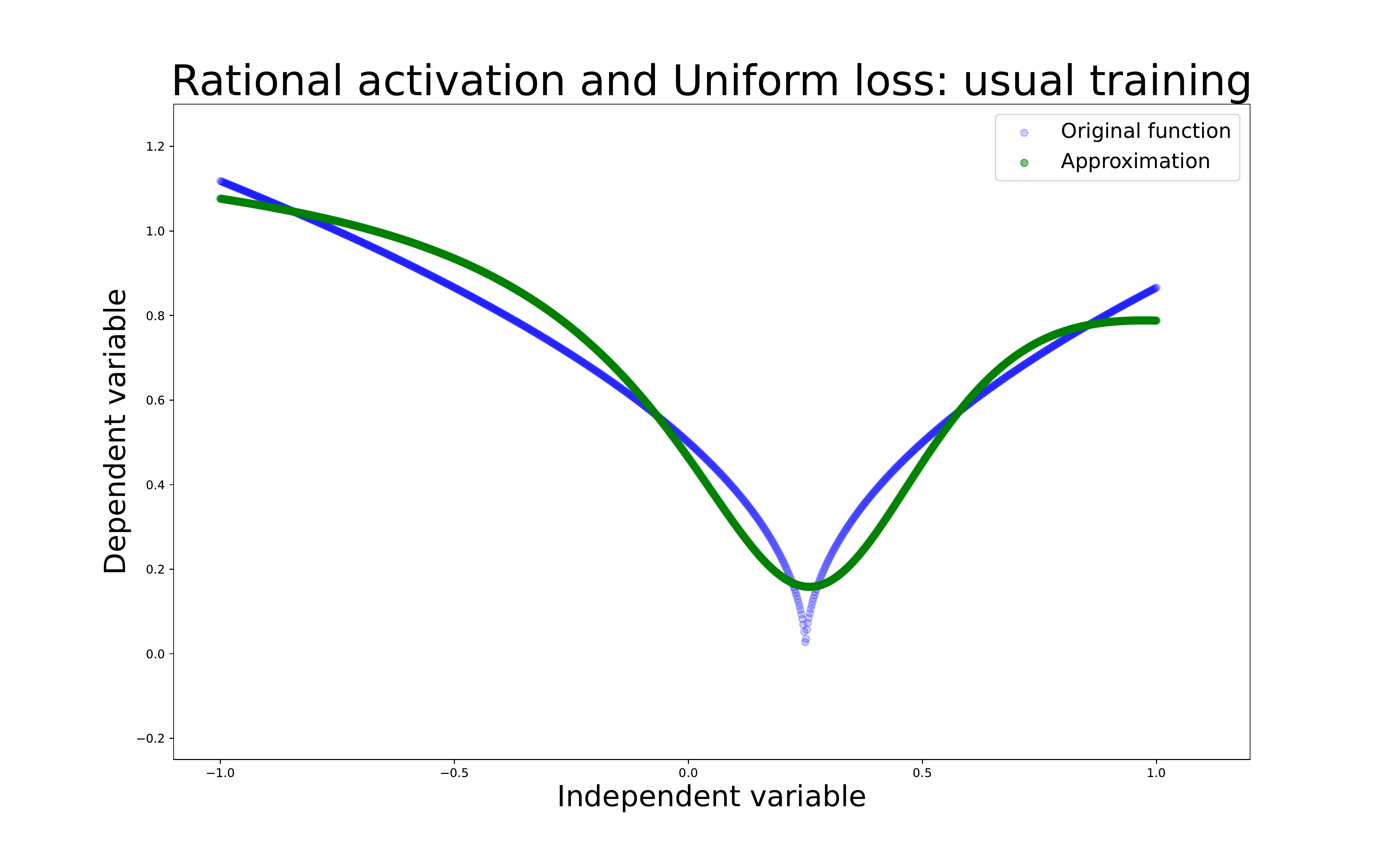}
        \includegraphics[width=40mm]{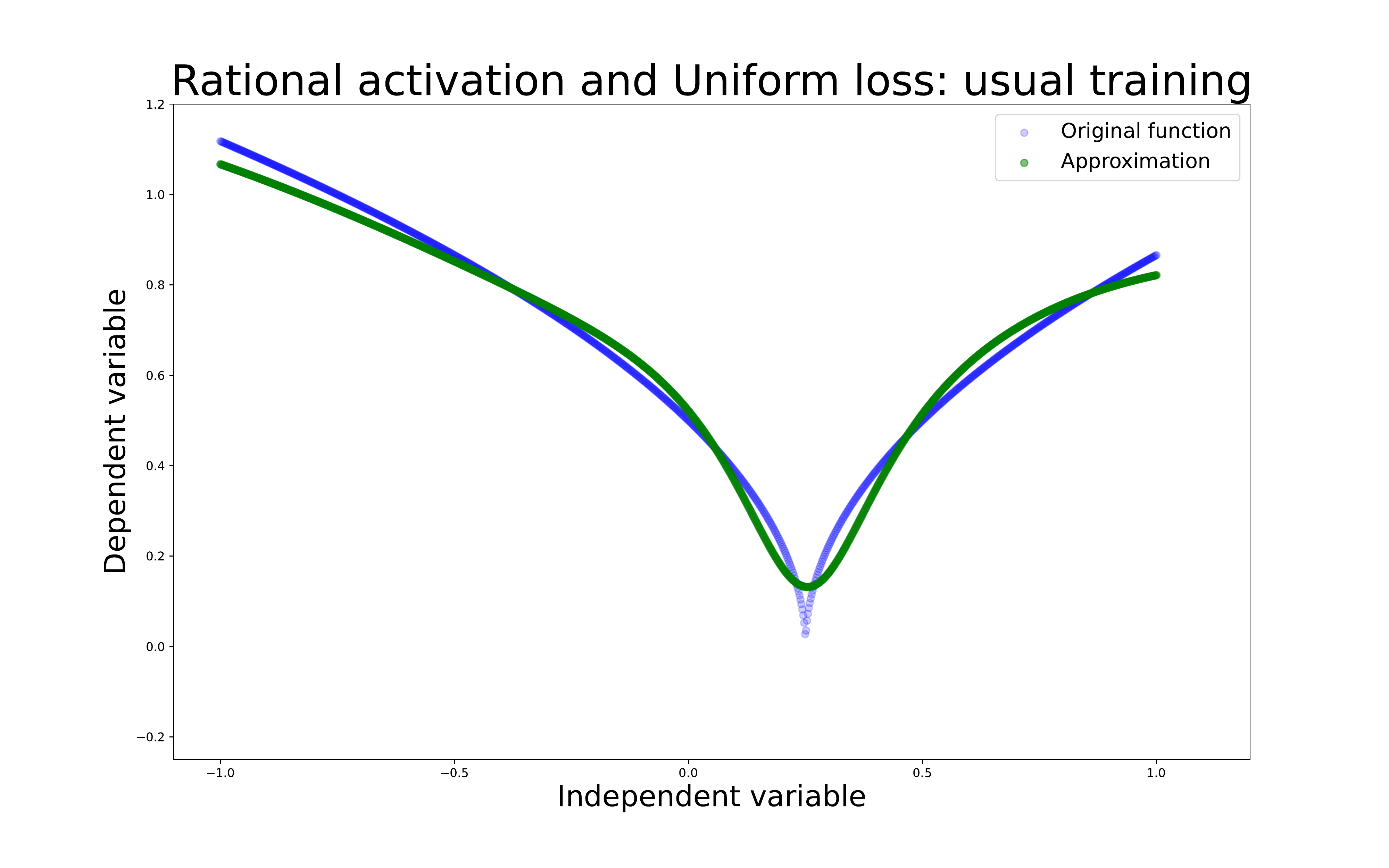}
            \includegraphics[width=40mm]{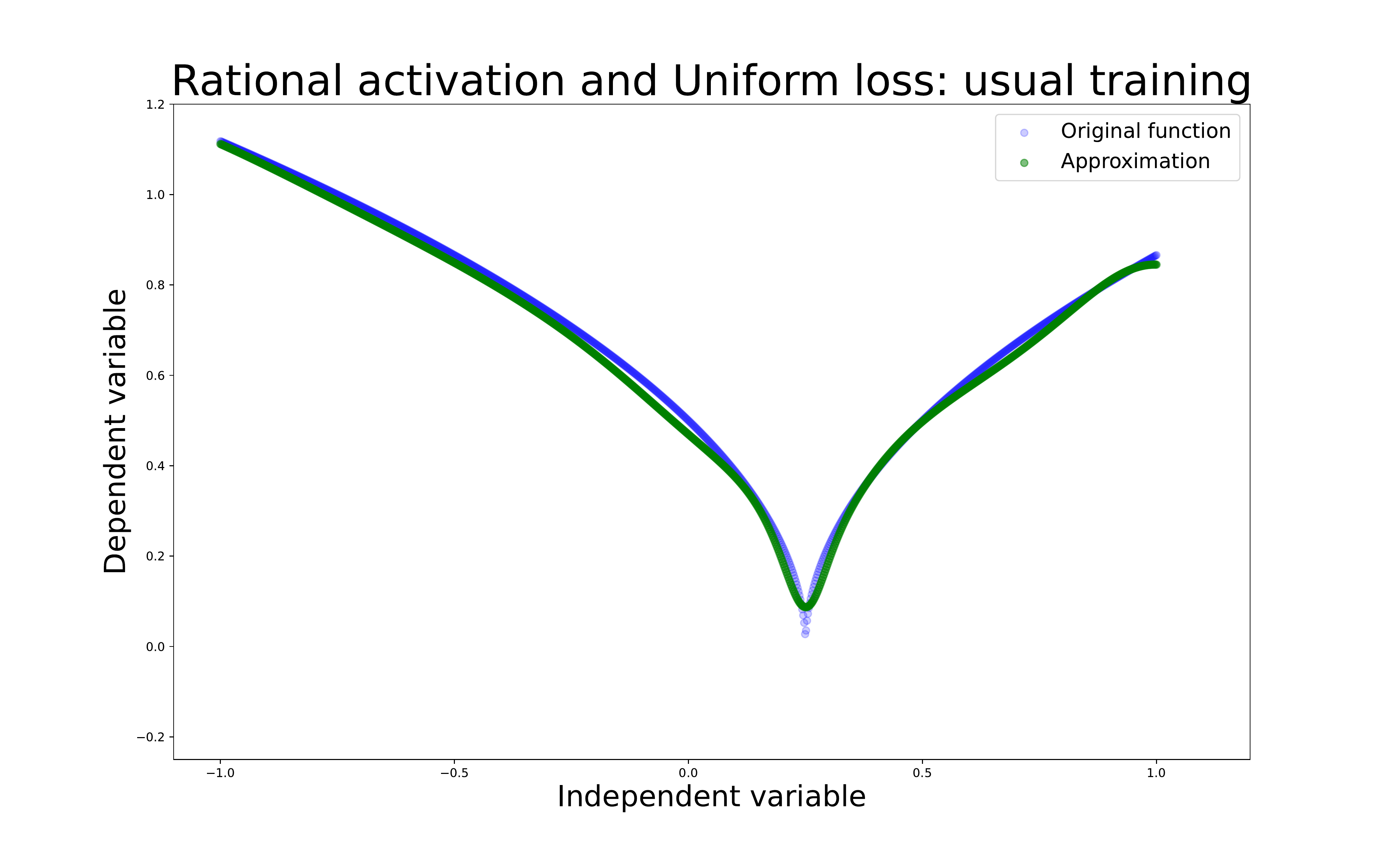}
    \caption{Approximation is computed by the usual training process, epoch is 50, 100, 200.}
    \label{fig:Rational activation with usual training - set 4, epoch is 50,100,200}
\end{figure}

\subsection{Neural Network with rational activation, learn with split training method} \label{subsec:results4:Set4}


In this section, our activation function is a rational function of degree $(3,2)$. 
The coefficients are now a part of the parameter set. Similar to the network in Section~\ref{subsec:results4:Set4}, we learn these coefficients as we learn other parameters during the training procedure. The training process is done in three steps: 
\begin{enumerate}
\item 
weights between the input and hidden layer, 
\item 
weights between the hidden layer and the output;
\item rational coefficients.
\end{enumerate}

This approach is based on a block coordinate-wise optimisation; therefore, even when it produces a sequence which decreases the  functional values, the limit point is not guaranteed to be globally optimal.

Table~\ref{tab:Results: experiments set 4 of NN with rational approximation and split training} shows that the results are comparable for this method and for the cases of rational activation and rational approximation to the ReLU activation function (in terms of the optimal values of the loss function). 

\begin{table}
    \centering
    \begin{tabular}{|c|c|c|c|}
    \hline
    Epoch & Final loss & Minimum loss & Run time (per epoch) \\
    \hline
    50  & 0.168947 &          & 2.65s $\pm$ 72.2ms\\
    49  &          & 0.113316 & \\
    \hline
    100 & 0.093841 &          & 2.71s $\pm$ 662ms\\
    75  &          & 0.087595 & \\
    \hline
    200 & 0.061735 &          & 2.82s $\pm$ 384ms\\
    190 &          & 0.053655 & \\
    \hline
    \end{tabular}
    \caption{Results: experiments with split method}
    \label{tab:Results: experiments set 4 of NN with rational approximation and split training}
\end{table}



\begin{figure}
    \centering
    \includegraphics[width=40mm]{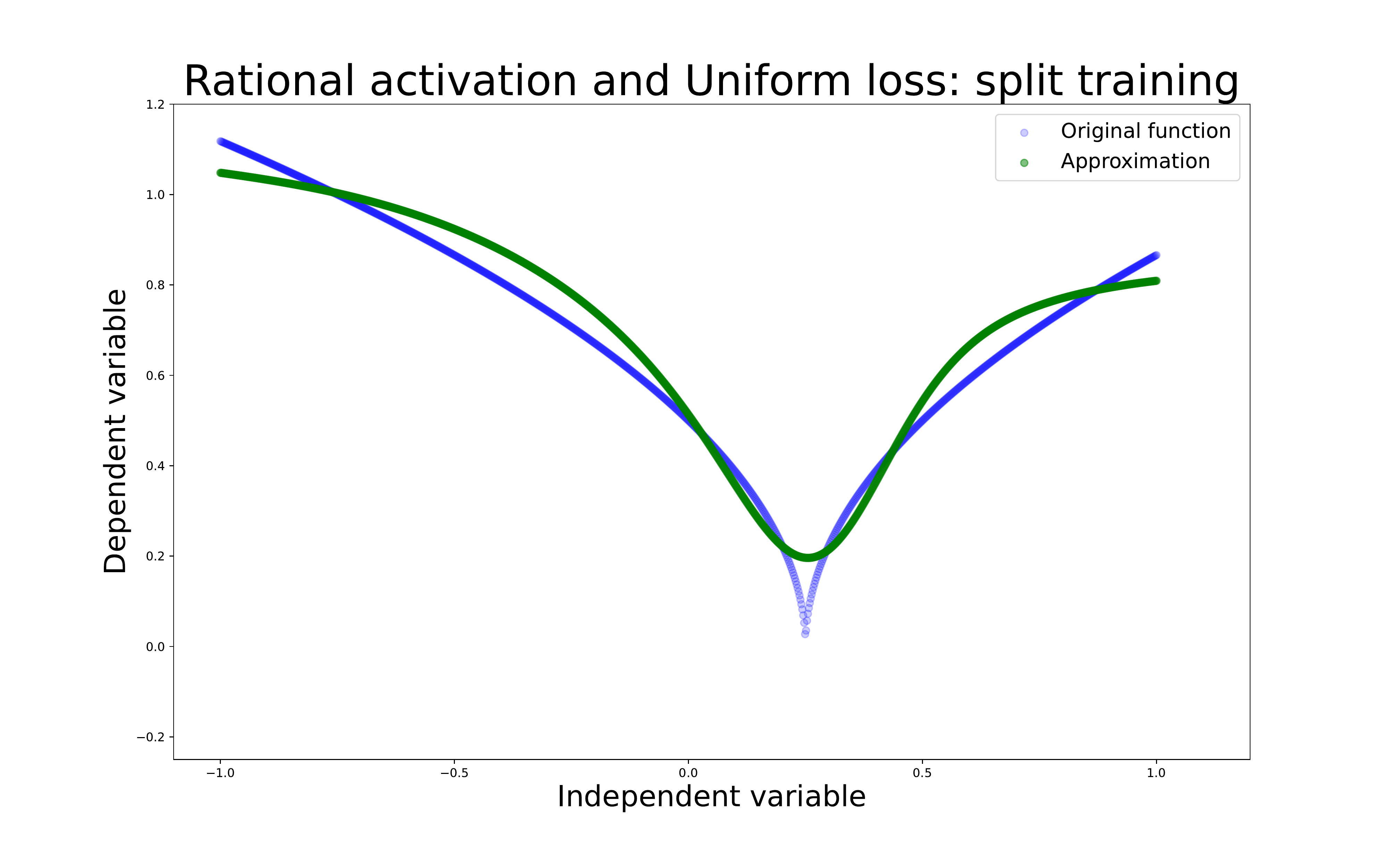}
    \includegraphics[width=40mm]{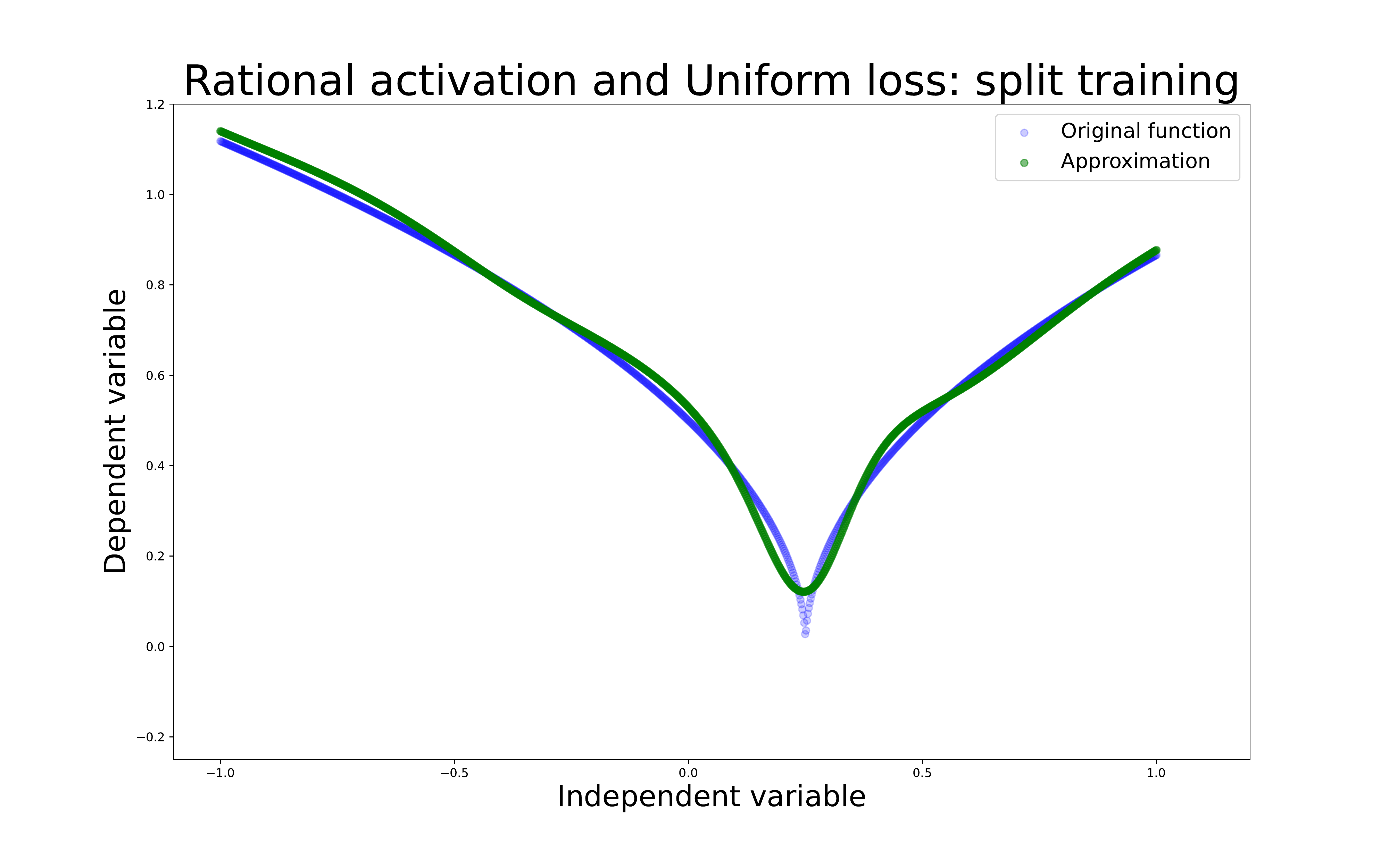}
    \includegraphics[width=40mm]{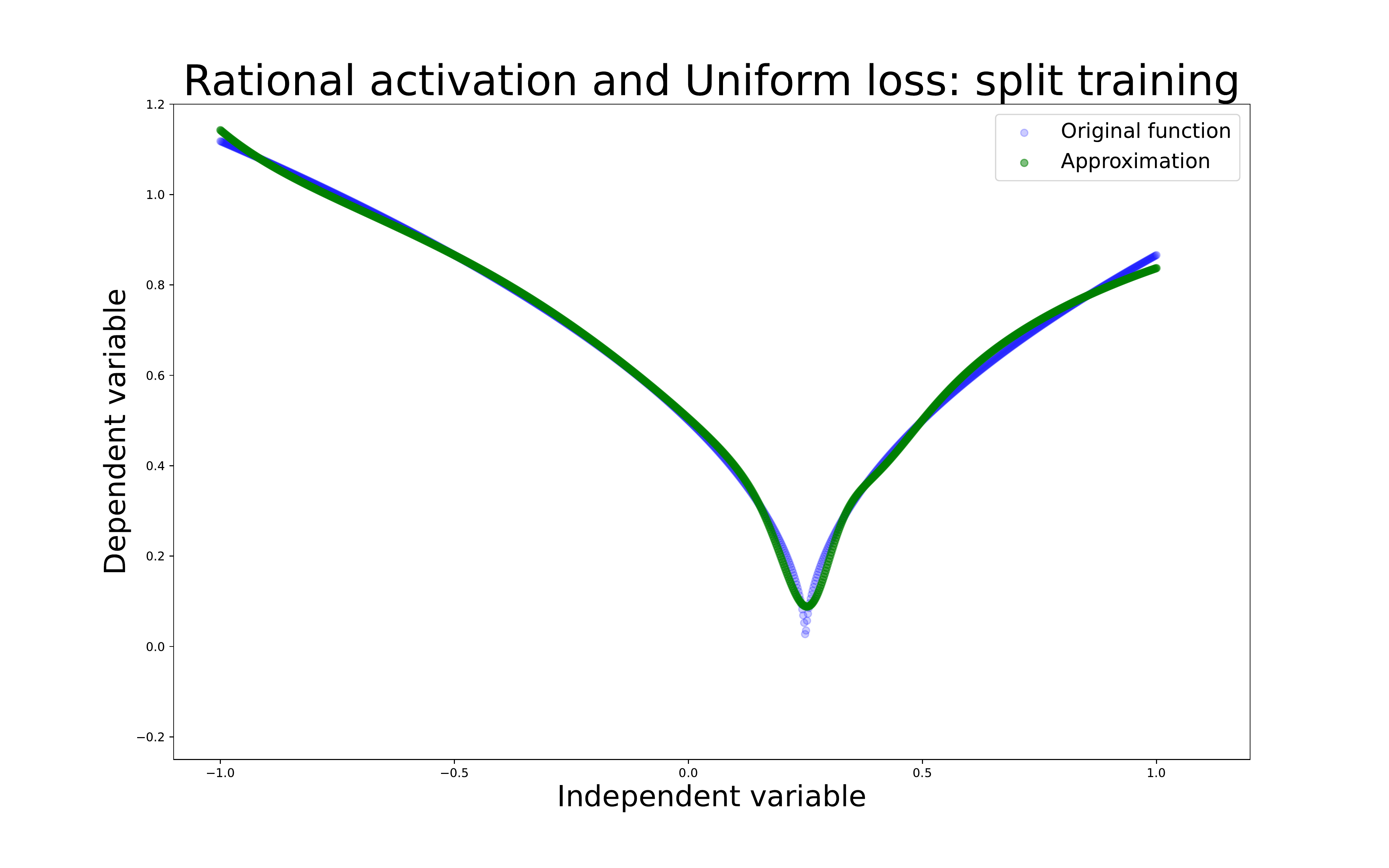}
    \caption{Approximation is computed by the split training process,  50, 100, 200~epochs.}
    \label{fig:Rational activation with split training - set 4, epoch is 50,100,200}
\end{figure}

Figure~\ref{fig:Rational activation with split training - set 4, epoch is 50,100,200} shows with the increase of the number of epochs the approximation is improving even around the ``difficult point''.

Overall conclusion: the approximation results are more accurate when the rational activation function is used with the split training method rather than the usual training approach.

\section{Direct rational approximation approach: the differential correction and AAA method.}

In this section, all the experiments correspond to uniform loss. The decision variables are the coefficients of the rational function. We compare two main methods for constructing rational approximations: the differential correction method (converges to a global optimal solution) and the AAA method (practical method, converges to ``near optimal'' solution). The codes for both methods are implemented in Python.

\subsection{The differential correction method}


\subsubsection{Rational approximation of degree (21,20), the differential correction method.}

When one considers a neural network with three layers and just 10~nodes in the hidden layer and the rational activation function of degree~$(3,2)$, the output is a rational function of degree $(21,20)$. Therefore, we can approximate the function $f(x) = \sqrt{|x-0.25|}$ by a rational function of degree~$(21,20)$  using the differential correction method and compare this approximation with the results obtained for uniform loss function and 10~nodes in the hidden layer. Results related to the case where the number of nodes on the hidden layer is~2, can be found in the Appendix~\ref{sec:Appendix}.

In theory, the direct approximation by a rational function of degree~$(21,20)$ is more flexible. We use Python  in our experiments.

We also compute the rational approximations of degree~$(20,20)$ for comparison. We will also use these additional results to compare with the AAA method, which is expecting the same degree in the numerator and the denominator. Table~\ref{tab:Results: experiments set 1 of rational approximation through differential correction method} summarises the computational time and the optimal loss. The results demonstrate that the computational time is very small for all three settings and corresponding errors are very similar. Interestingly, the increase in the degree does not always improve the optimal loss function value. This can be explained as a result of computational instability when the dimension of the corresponding optimisation problems is increasing.

\begin{table}
    \centering
    \begin{tabular}{|c|c|c|}
    \hline
    Degrees & Error (uniform norm) & Run time \\
    \hline
    (21,20) & 0.04376536257200335 & 10.9s\\
    \hline
    (20,20) & 0.04841949102524923 & 11.5s\\
    \hline
    (21,21) & 0.04411128396473672 & 12.4s\\
    \hline
    \end{tabular}
    \caption{Direct rational approximation for degrees $(21,20)$, $(20,20)$ and $(21,21)$ by the differential correction method.}
    \label{tab:Results: experiments set 2 of rational approximation through differential correction method}
\end{table}
Table~\ref{tab:Results: experiments set 2 of rational approximation through differential correction method} summarises the computational time and the optimal loss values for the degrees~(21,20),~(20,20) and (21,21). These results are very similar.

Figure~\ref{fig:Rational approximation of degree (21,20), (20,20)  and (21,21) via the differential correction method} depicts the approximation of degree~(21,20). The approximation is accurate, but there is a strong oscillation of the approximation, which is not surprising for high degree rational functions.
\begin{figure}
    \centering
    \includegraphics[width=40mm]{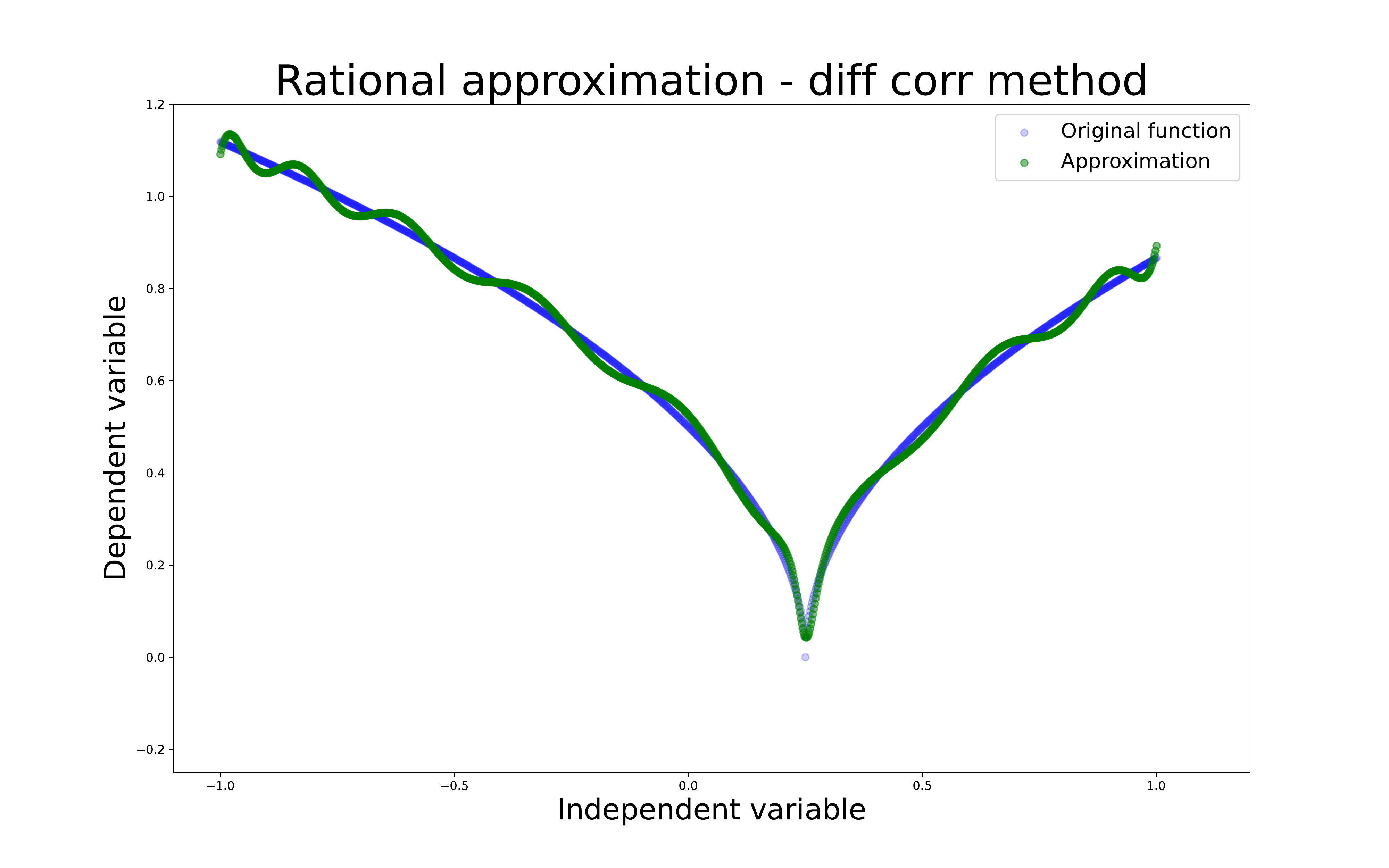}
    \includegraphics[width=40mm]{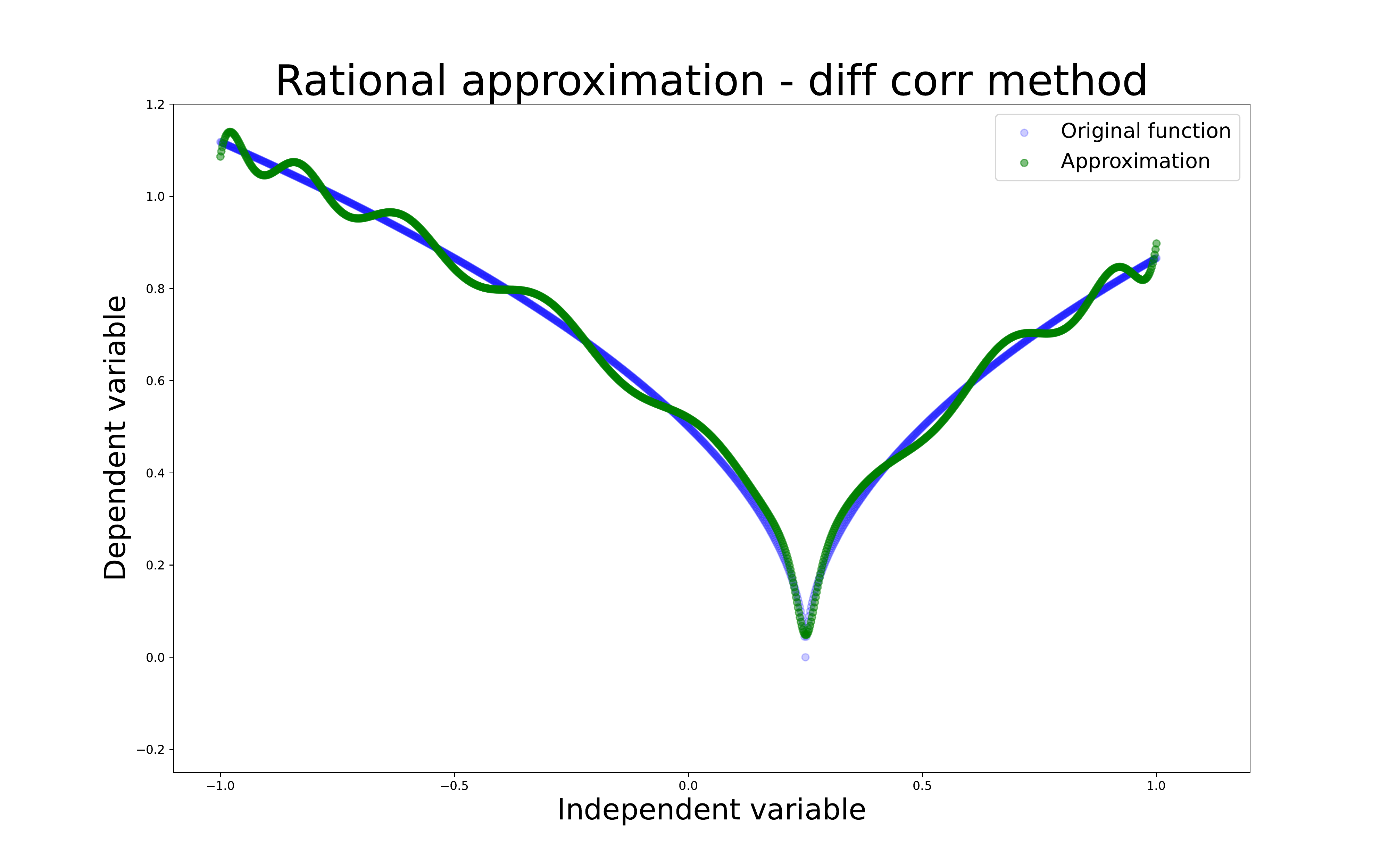}
    \includegraphics[width=40mm]{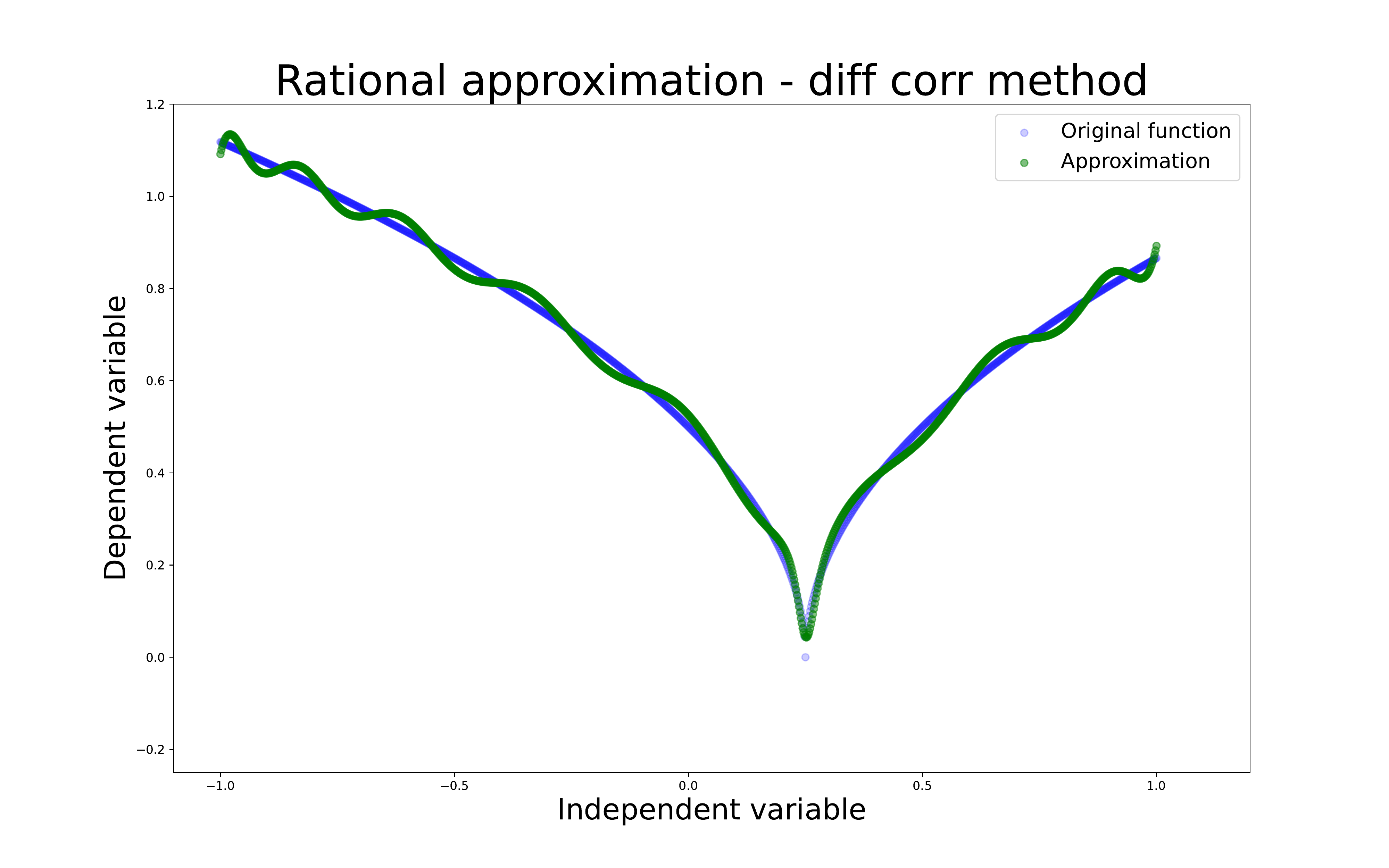}
    \caption{Approximation by the differential correction method, rational approximation degree is $(21,20)$, $(20,20)$ and $(21,21)$ respectively.}
    \label{fig:Rational approximation of degree (21,20), (20,20)  and (21,21) via the differential correction method}
\end{figure}


\subsection{The AAA method}

\subsubsection{Rational approximation of degree (21,20) with the AAA method}
In order to compare the differential correction method with AAA, we apply the AAA method in the same experimental settings as they are in the case of the differential correction method. Due to the limitations of the AAA method, we cannot use different degrees for the numerator and the denominator. Hence, for these experiments, we use the rational approximation of degree (20,20) and the rational approximation of degree (21,21). The codes are developed in Python by C. Hofreither~\cite{Hofreither21}.

\begin{table}
    \centering
     \begin{tabular}{|c|c|c|c|}
    \hline
    Degrees & Error type & Error & Run time \\
    \hline
    \multirow{2}{4em}{(20,20)} & Uniform error &
    0.00026068574873364114   & \multirow{2}{4em}{450ms}\\
    \cline{2-3}
     & MSE & 0.0005497563210130126 & \\
    \hline
    \multirow{2}{4em}{(21,21)} & Uniform error &
    3.67961692604446e-06   & \multirow{2}{4em}{496ms}\\
    \cline{2-3}
     & MSE & 4.843449639625575e-05 & \\
    \hline
    \end{tabular}
    \caption{Results: AAA method, degree (20,20) and (21,21)}
    \label{tab:Results: experiments set 2 of rational approximation through AAA method}
\end{table}


\begin{figure}
    \centering
    \includegraphics[width=50mm]{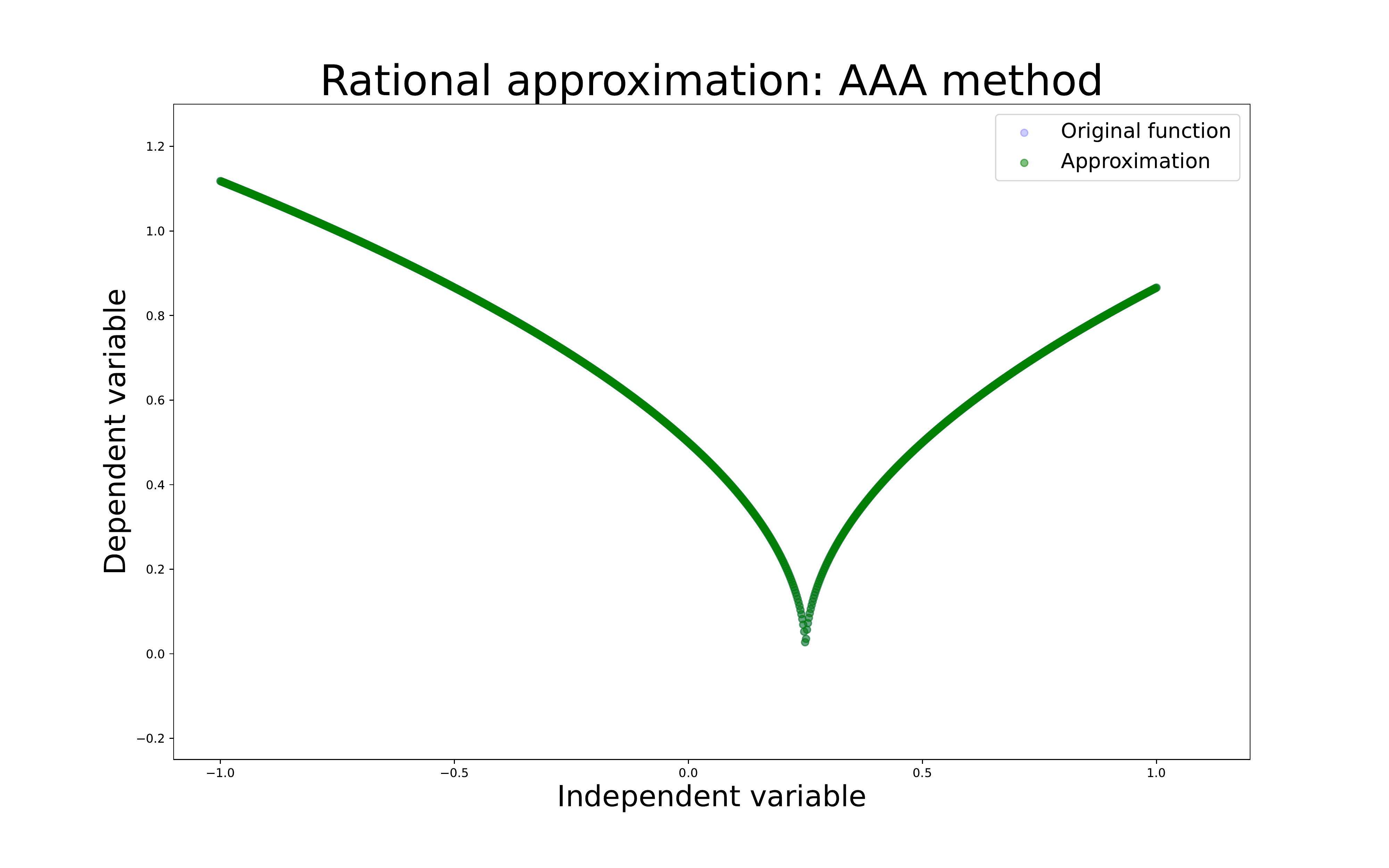}
    \includegraphics[width=50mm]{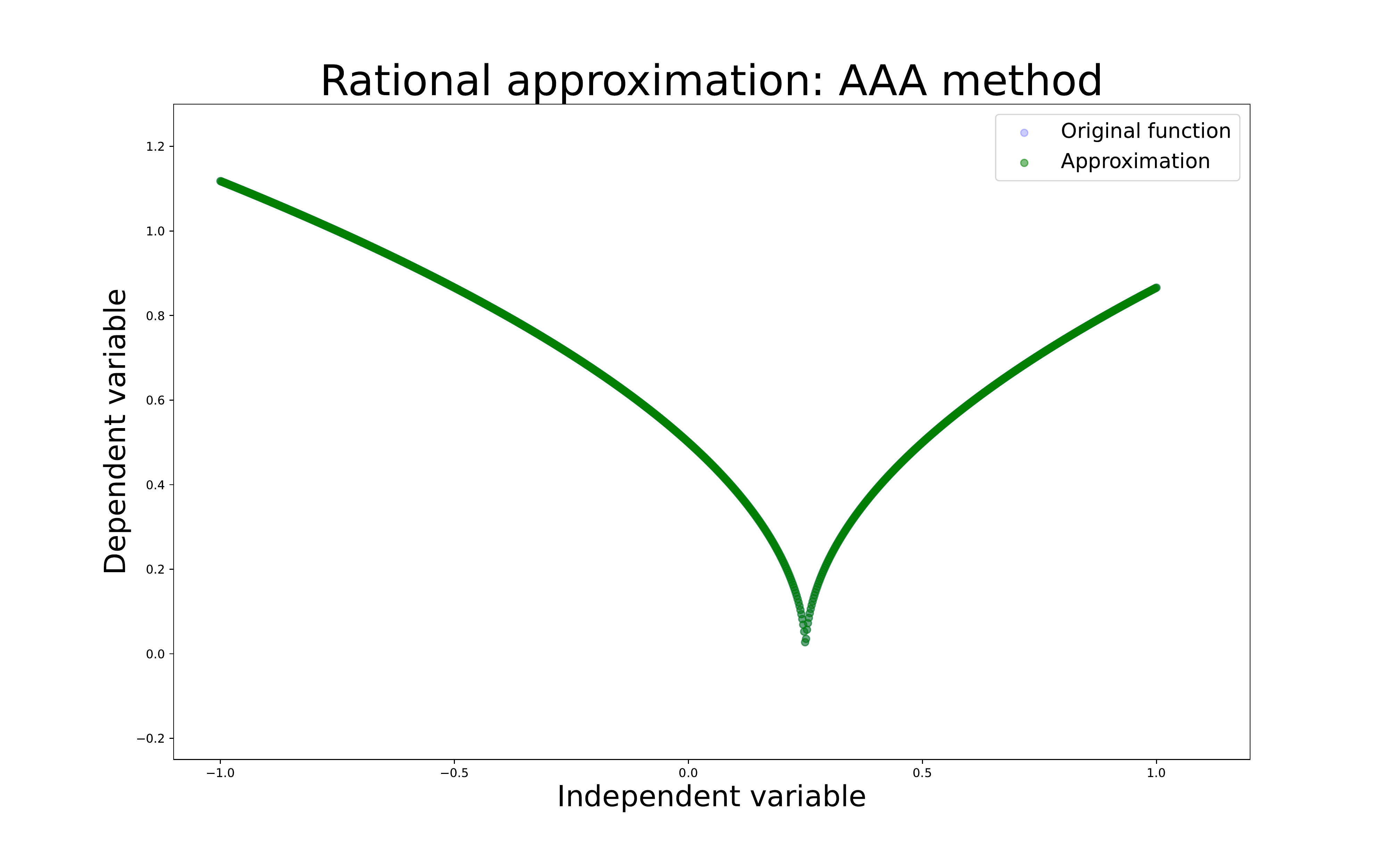}
    \caption{Approximation is computed by the AAA method, degree (20,20) and (21,21).}
    \label{fig:Rational approximation of degree (20,20) and (21,21) via AAA method}
\end{figure}

Table~\ref{tab:Results: experiments set 2 of rational approximation through AAA method} and Figures~\ref{fig:Rational approximation of degree (20,20) and (21,21) via AAA method} shows that the approximations are very accurate and the computational time is very low. These results suggest that the AAA method is very fast and accurate, even when the degree is high. The function and the approximation are indistinguishable from the picture.

\section{Results for multivariate domains}\label{sec:results2D}

In this section, we recreate the experiments presented in section~4.1 of~\cite{boulle2020rational}, where the authors approximate a solution to the KdV (Korteweg-de Vries) equation:
$$u_t = -uu_x - u_{xxx}, \quad u(x,0) = -\sin(\pi x /20)$$ by using four neural networks, which contain ReLU, sinusoid, rational, and polynomial activation functions. In this section, our goal is to approximate $u_i$ on $(x_i,t_i)$ by using the multivariate bisection method and the algorithm was implemented in Matlab. The specifications of the dataset is as follows: 
\begin{itemize}
    \item $(u_i)_{1 \leq i \leq N}$ - observation data (known function values) of size $512 \times 201.$
    \item $(x_i , t_i)_{1 \leq i \leq N}$ - spatio tempural points (known domain values). 
    \item $x_i \in [0,40]$ with the step size of 1.
    \item $t_i \in [-20, 19.84375]$ with the step size of 0.390625.
\end{itemize}
In~\cite{boulle2020rational}, authors used Euclidean norm to compute the errors but here, we compute the error of the approximation in the uniform norm. Therefore, there are not many possibilities for comparison with the findings in~\cite{boulle2020rational}. All of our experiments are constructed by using a certain set of domain points (using every $k^{\text{th}}$ point of each dimension and the corresponding function values). The approximation and the uniform error term are computed on the same reduced domain. This extra step has been taken to ensure that big data volume does not lead MATLAB programs to crash.

\subsection{The original dataset} \label{ssec:original}
\begin{figure}
    \centering
    \includegraphics[width=60mm]{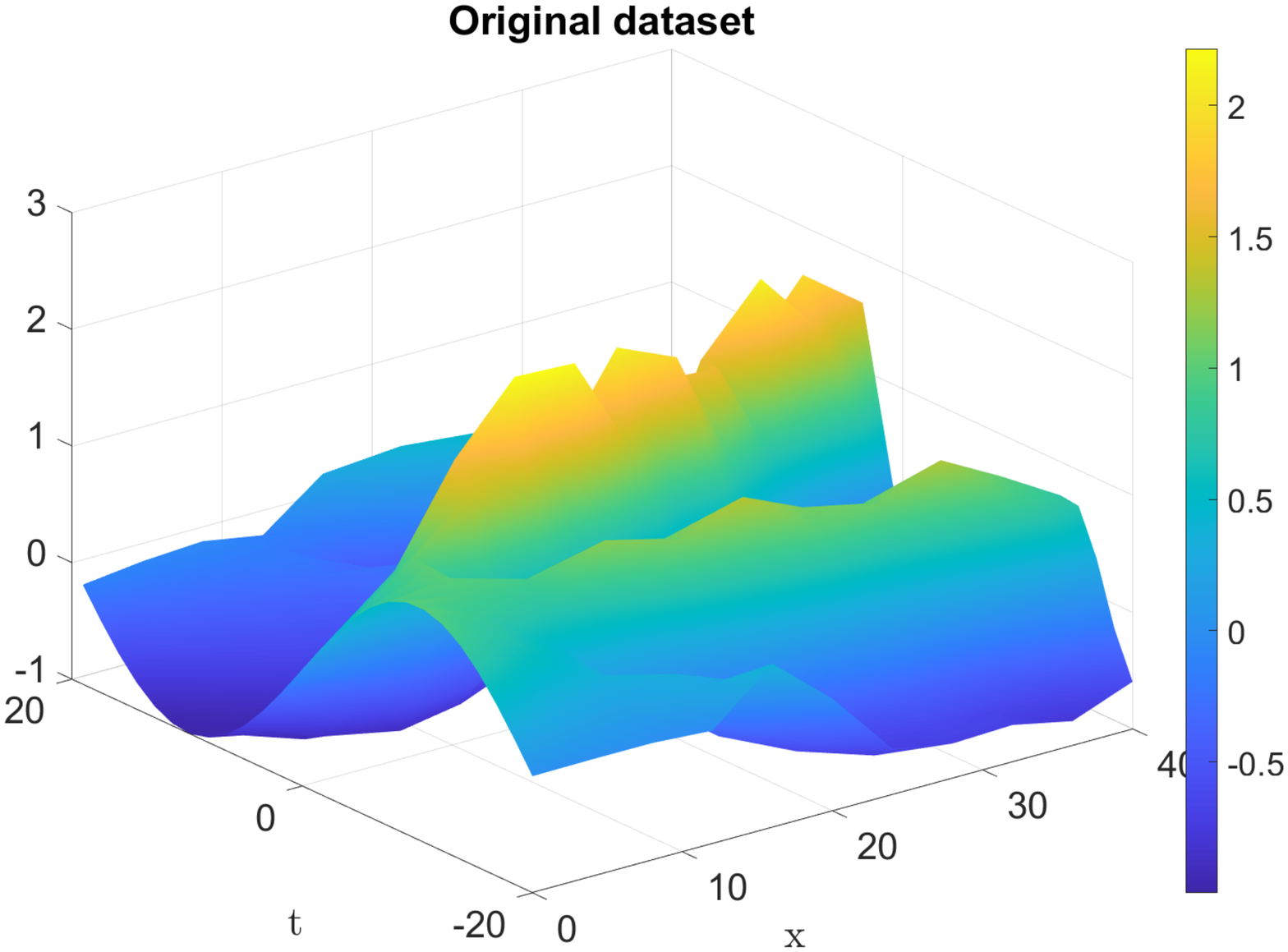}
    \includegraphics[width=60mm]{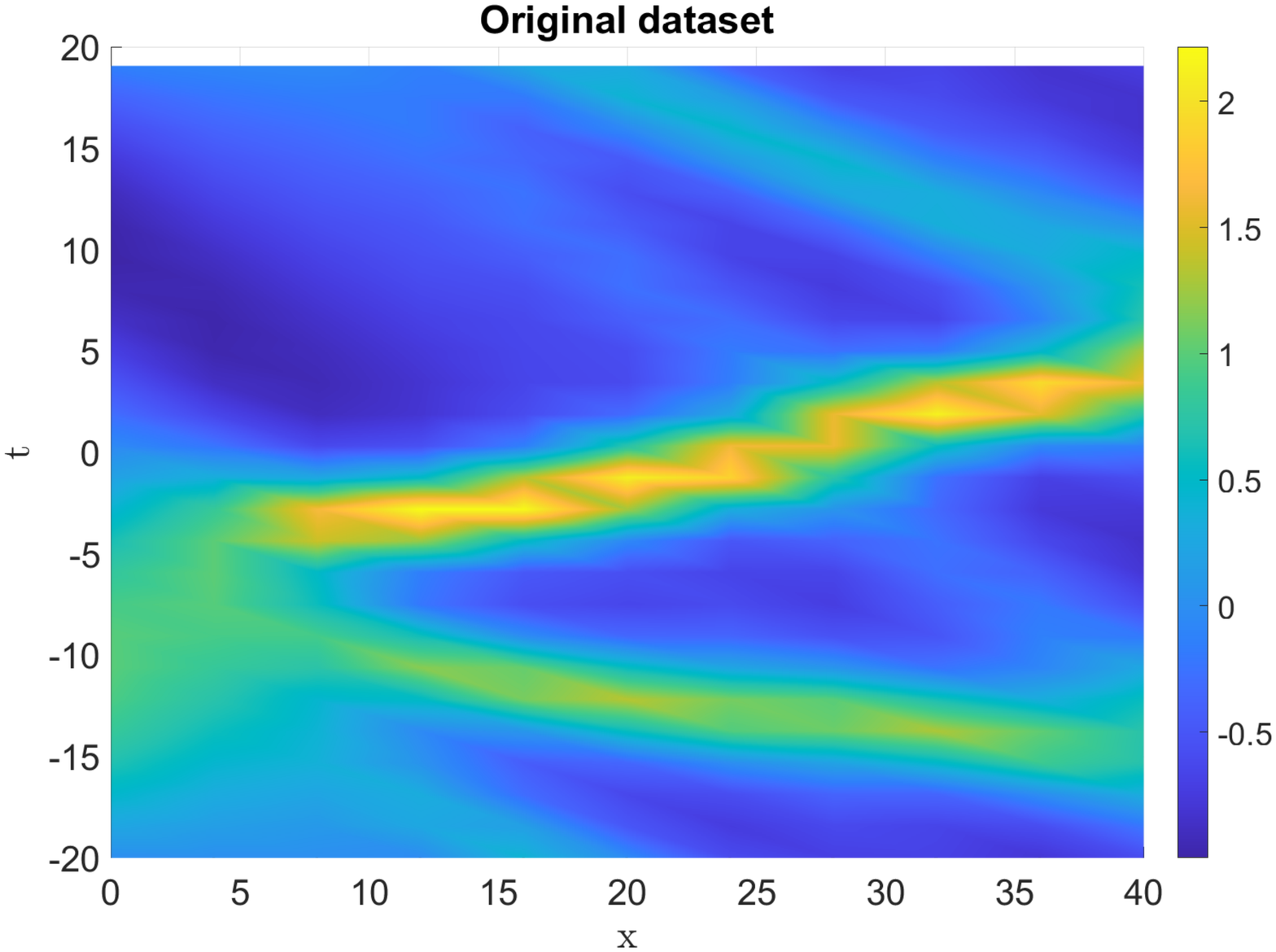}
    \caption{Pictures constructed for the original dataset taking the pairs of every $20^{\text{th}}$ point of the domain: 3D and 2D view of the original dataset}
    \label{fig:original dataset with pairs of every 20th point of the domain}
\end{figure}

\begin{figure}
    \centering
    \includegraphics[width=60mm]{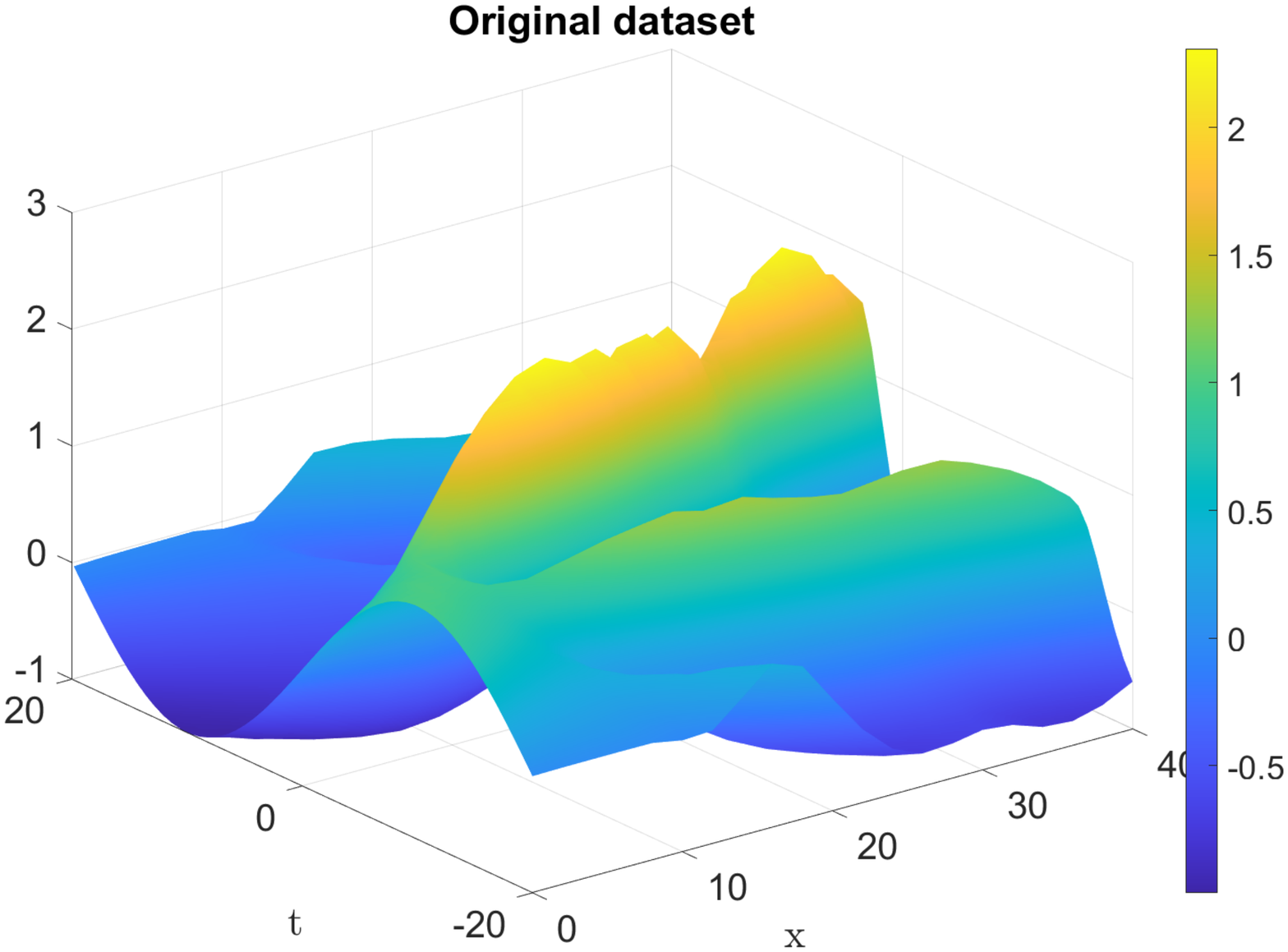}
    \includegraphics[width=60mm]{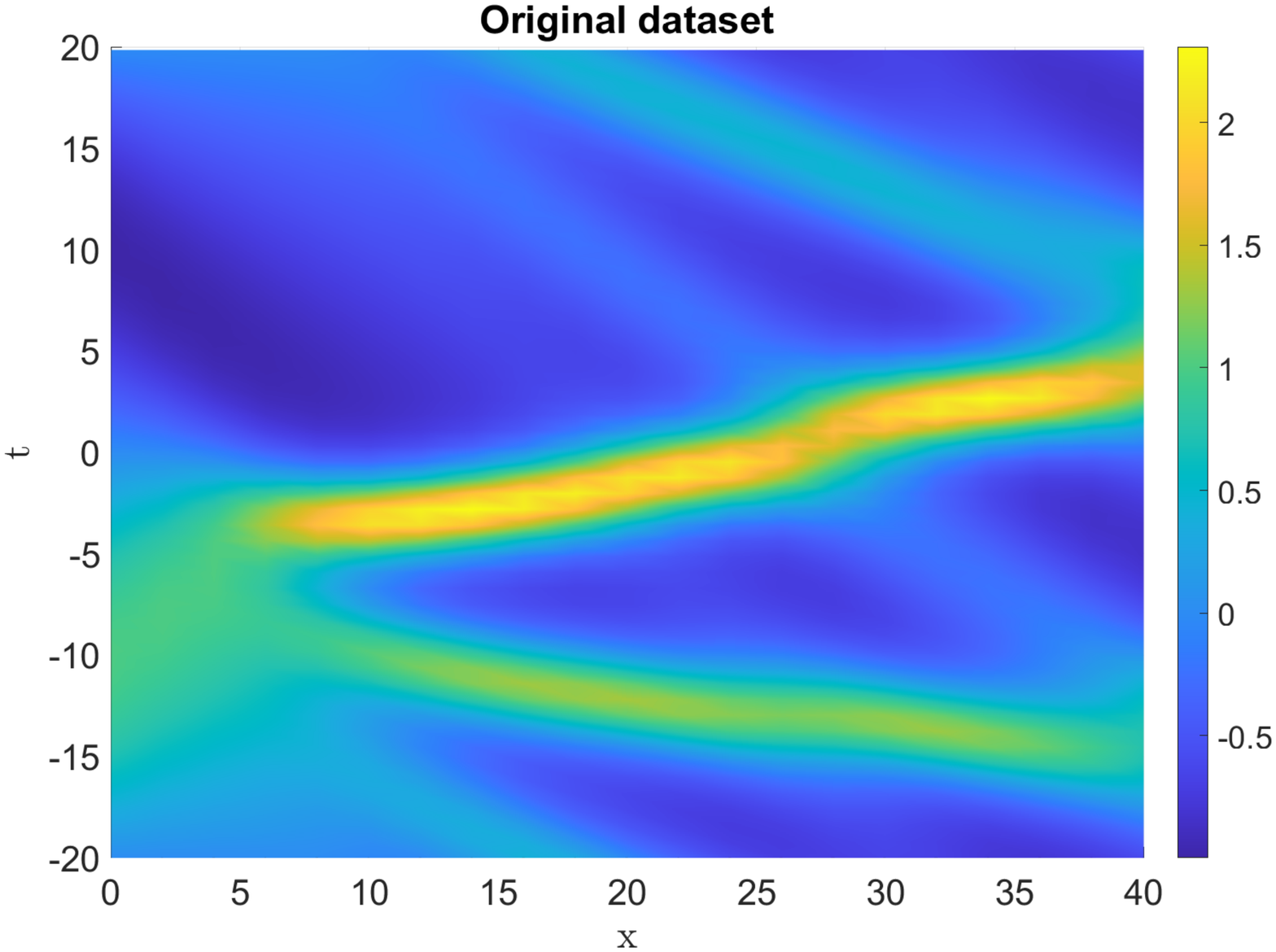}
    \caption{Pictures constructed for the original dataset taking the pairs of every $10^{\text{th}}$ point of the domain: 3D and 2D view}
    \label{fig:original dataset with pairs of every 10th point of the domain}
\end{figure}

\begin{figure}
    \centering
    \includegraphics[width=60mm]{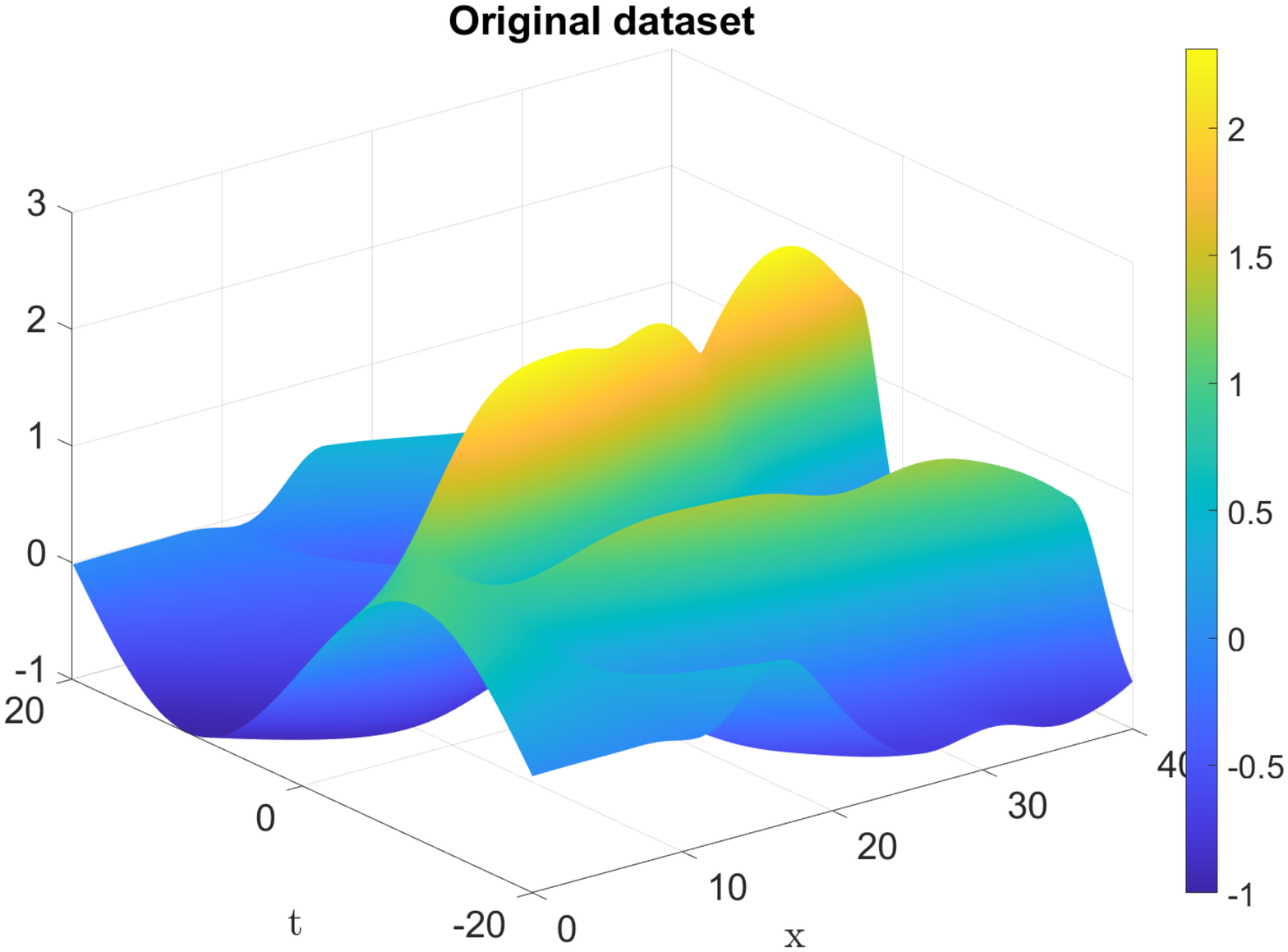}
    \includegraphics[width=60mm]{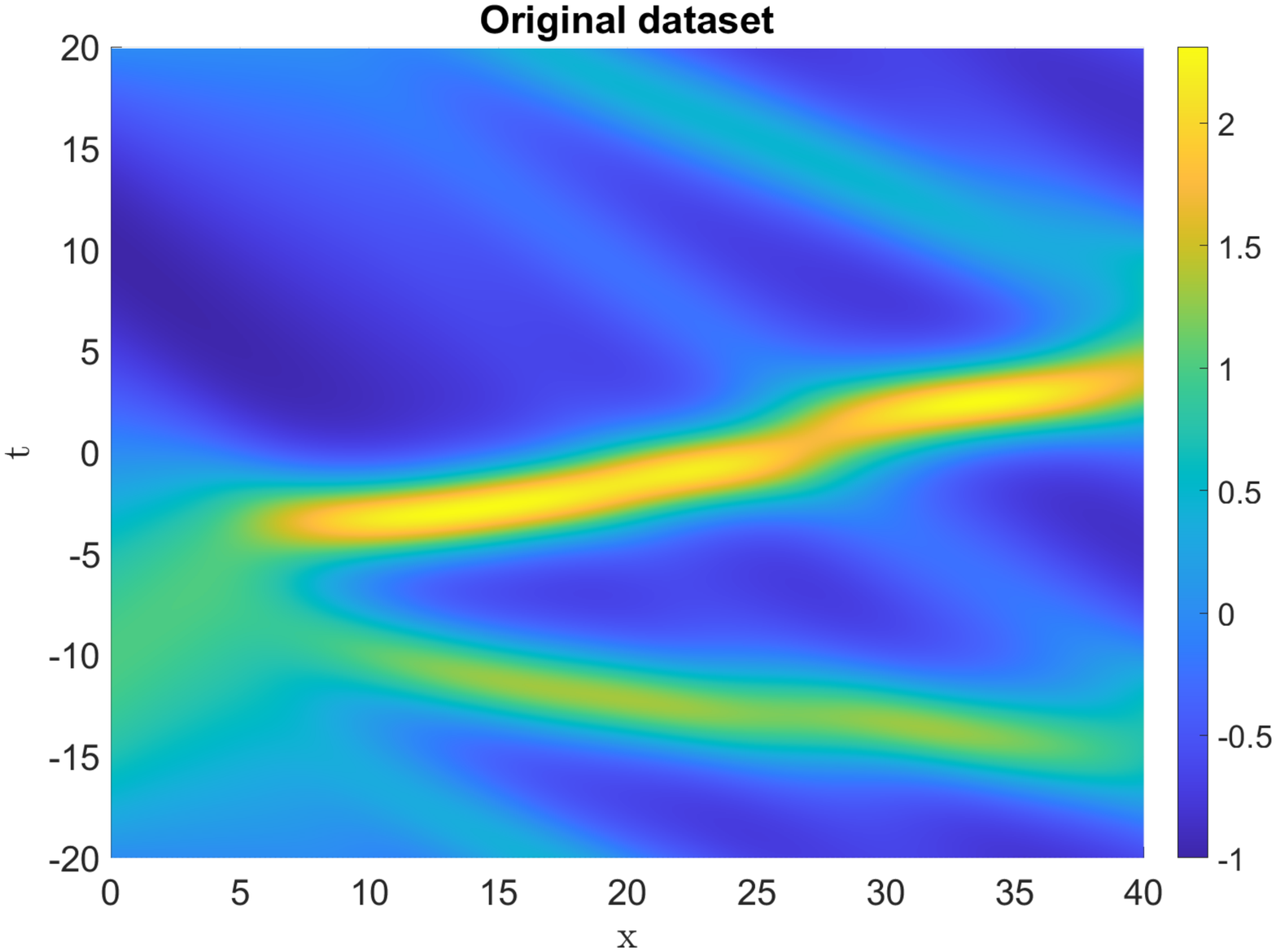}
    \caption{Pictures constructed for the original dataset taking all the points of the domain: 3D and 2D view}
    \label{fig:original dataset with all the pairs of points of the domain}
\end{figure}

The following pictures (Figure~\ref{fig:original dataset with pairs of every 20th point of the domain}, Figure~\ref{fig:original dataset with pairs of every 10th point of the domain}, Figure~\ref{fig:original dataset with all the pairs of points of the domain}) are of the original dataset. Some images are constructed by using a certain set of domain points. In particular, we use the pairs of every $20^{\text{th}}$ point of the domain (along each of the two dimensions) and pairs of every $10^{\text{th}}$ point of the domain in order to reduce the dimension of the original problem. For simplicity, we sometimes refer to these pairs as ``Every $10^{\text{th}}$ point of the domain'' or ``Every $20^{\text{th}}$ point of the domain''.

Now, we compute rational approximations of different degrees and report the uniform approximation error and the computational time for each degree when $k=20$ and $k=10$. We start our experiments with rational approximation of degree $(2,2)$ and keep increasing the degree up to $(20,20)$. The classical monomials in the rational function were replaced by Chebyshev monomials for better approximations.

\subsection{Rational approximation of degree (2,2)}
We report 3D and 2D pictures of the rational approximation of degree (2,2) along with its error curve. 

\begin{figure}
	\centering
		\includegraphics[width=60mm]{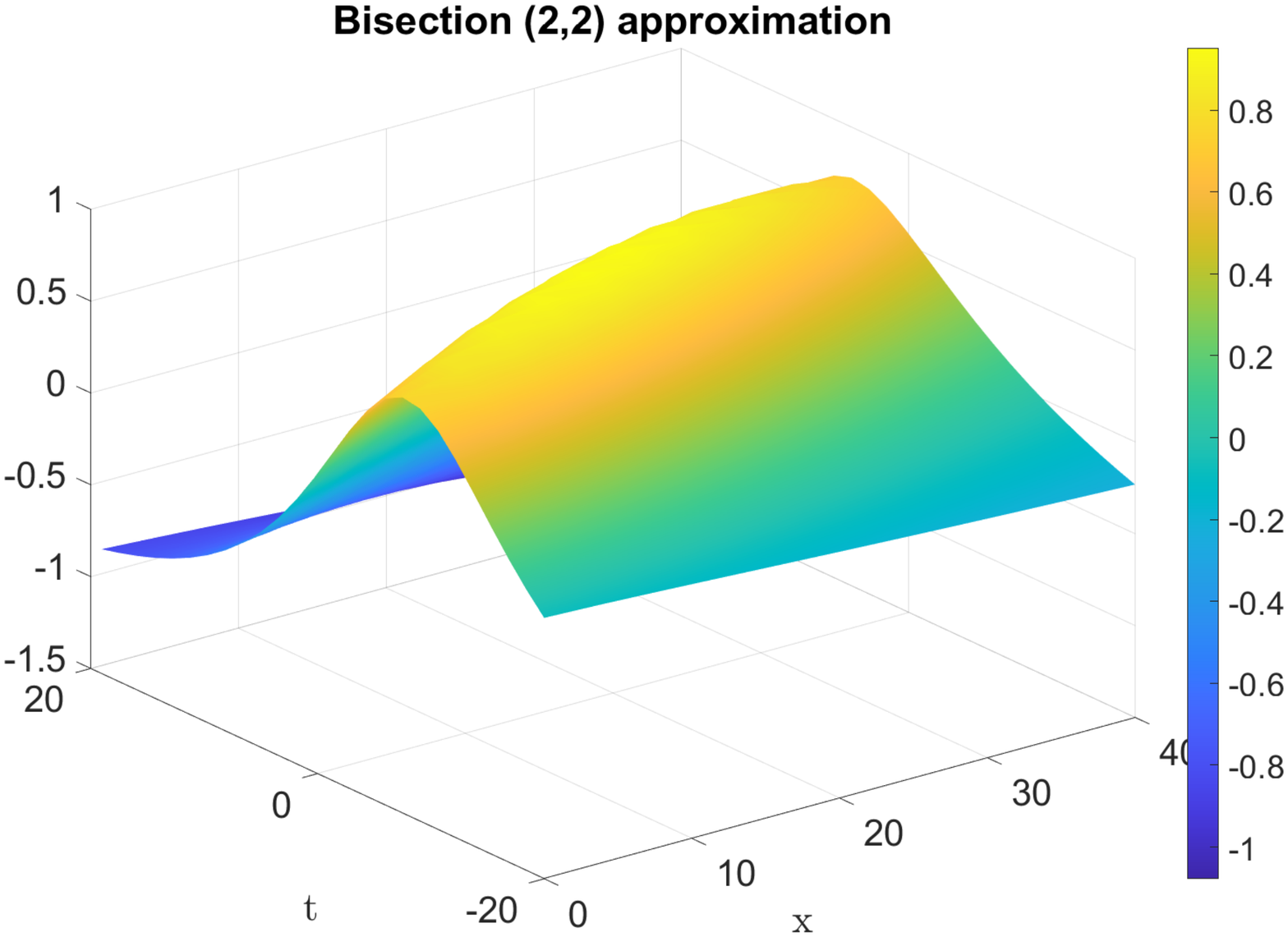}
		\includegraphics[width=60mm]{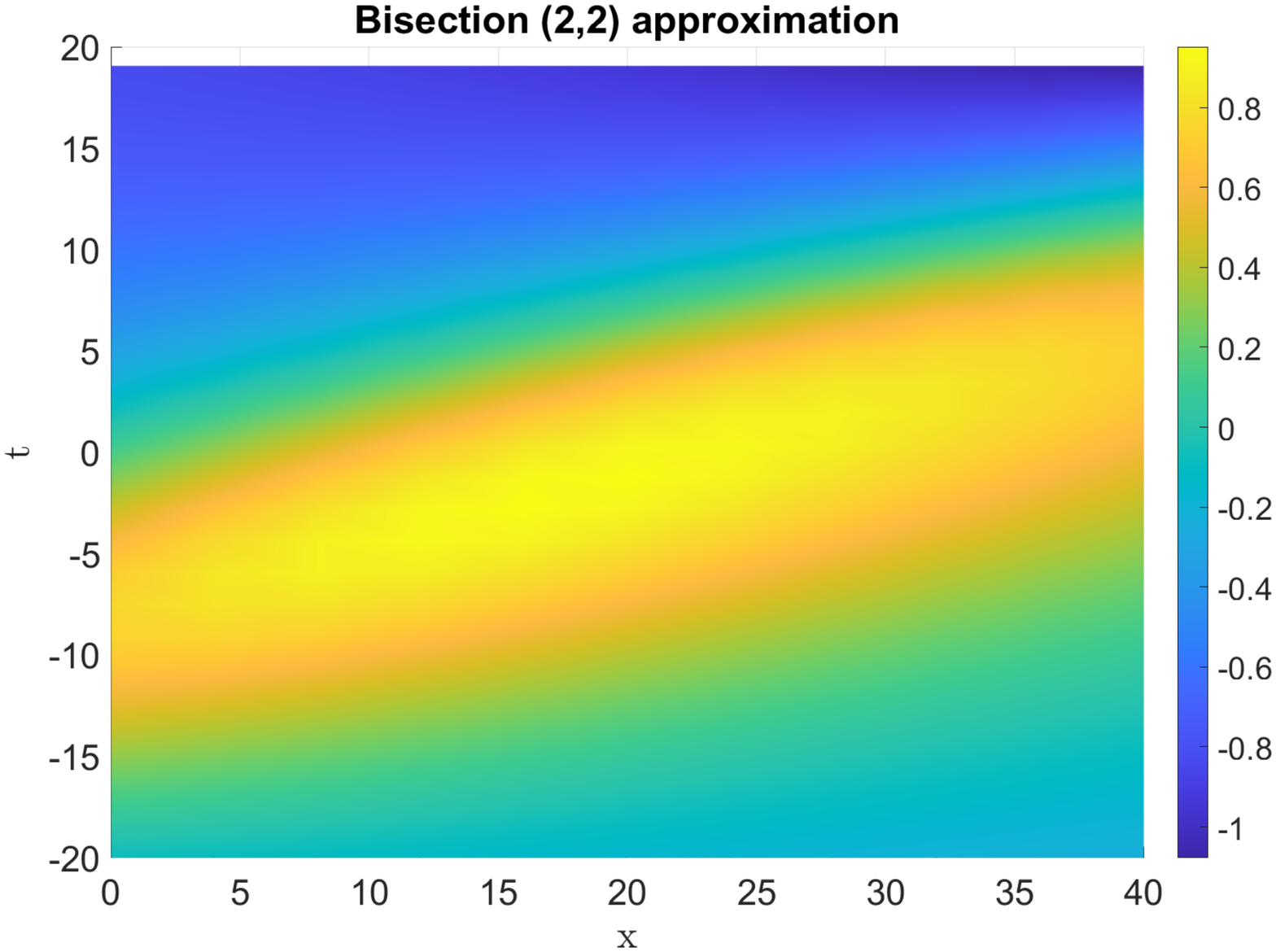}
		\includegraphics[width=60mm]{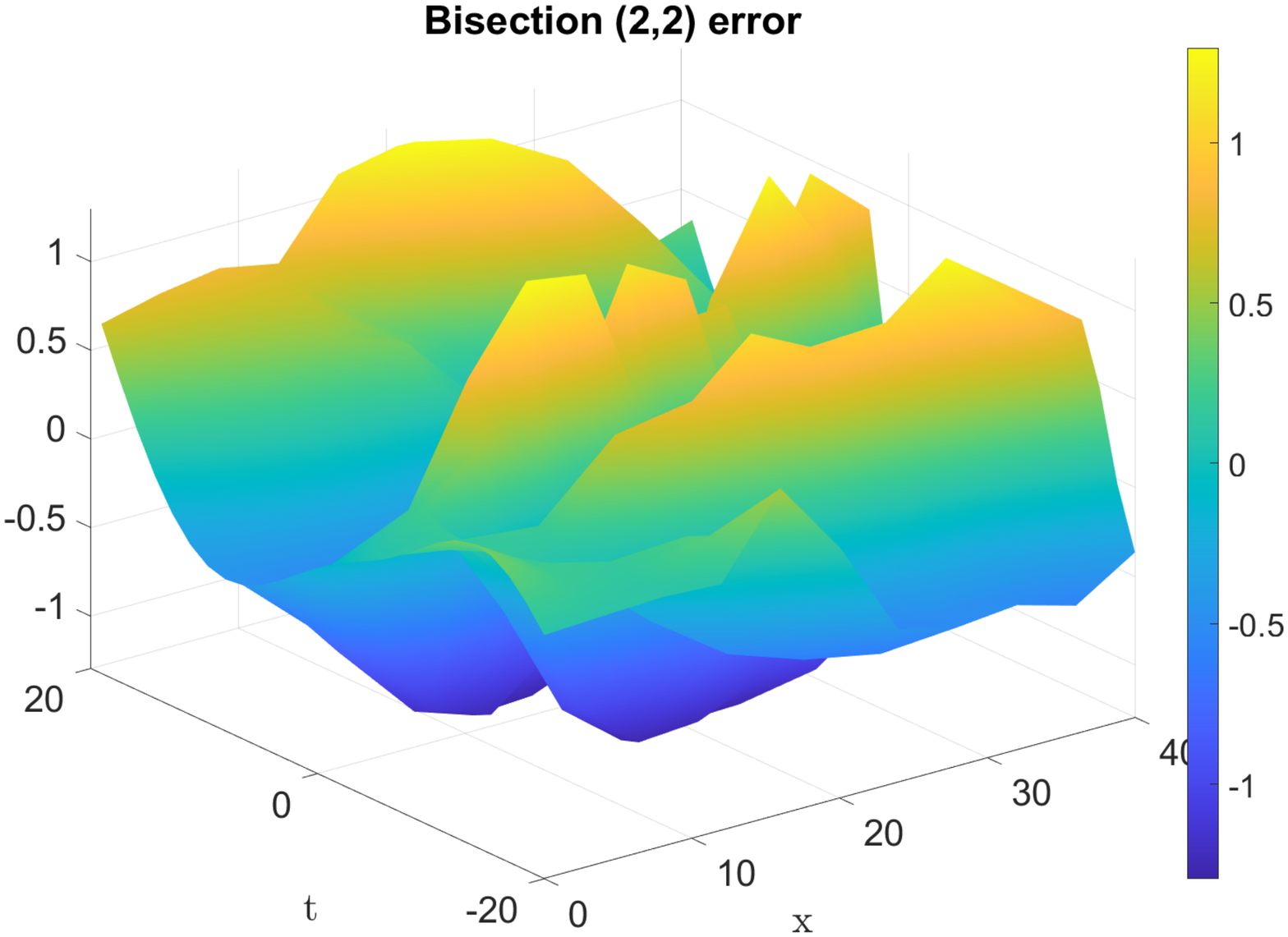}
		\includegraphics[width=60mm]{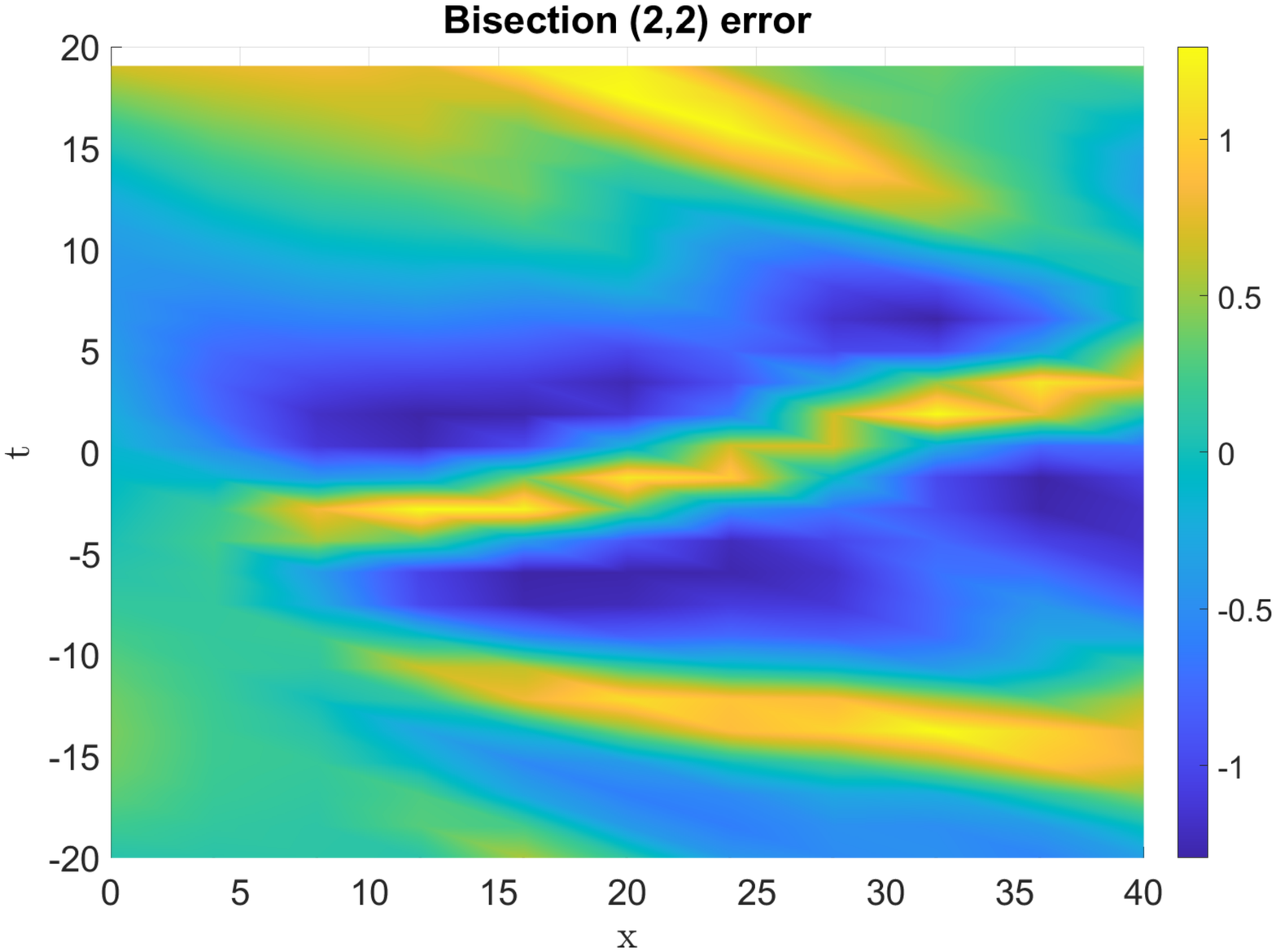}
	\caption{Approximations (above) and the error curves (below) are constructed on the domain which consists of pairs of every 20th point of the original domain: 3D and 2D view}
        \label{fig:Rational approximation of degree (2,2) multivariate domain of pairs of every 20th point}
\end{figure}

\begin{figure}
	\centering
		\includegraphics[width=60mm]{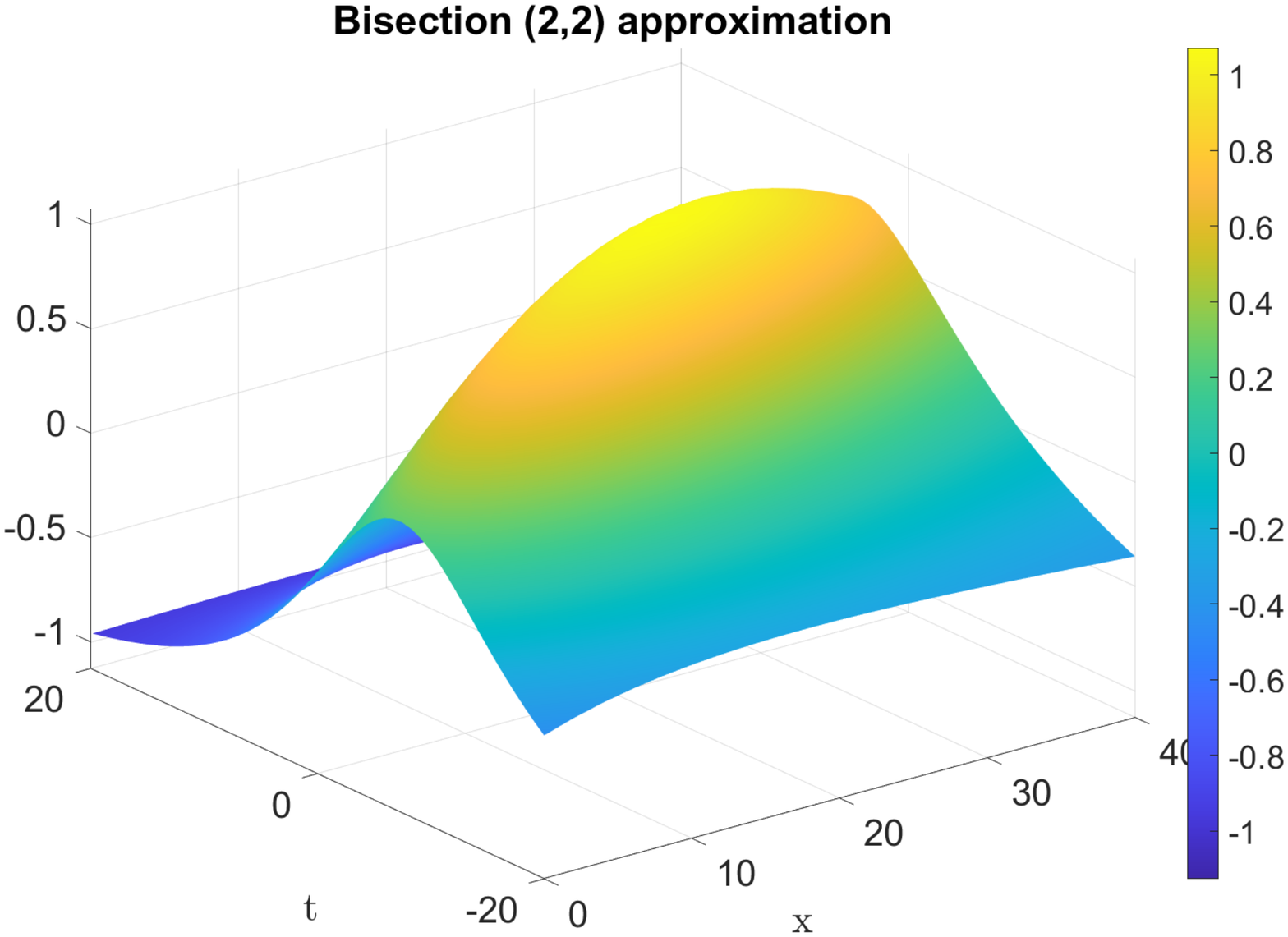}
		\includegraphics[width=60mm]{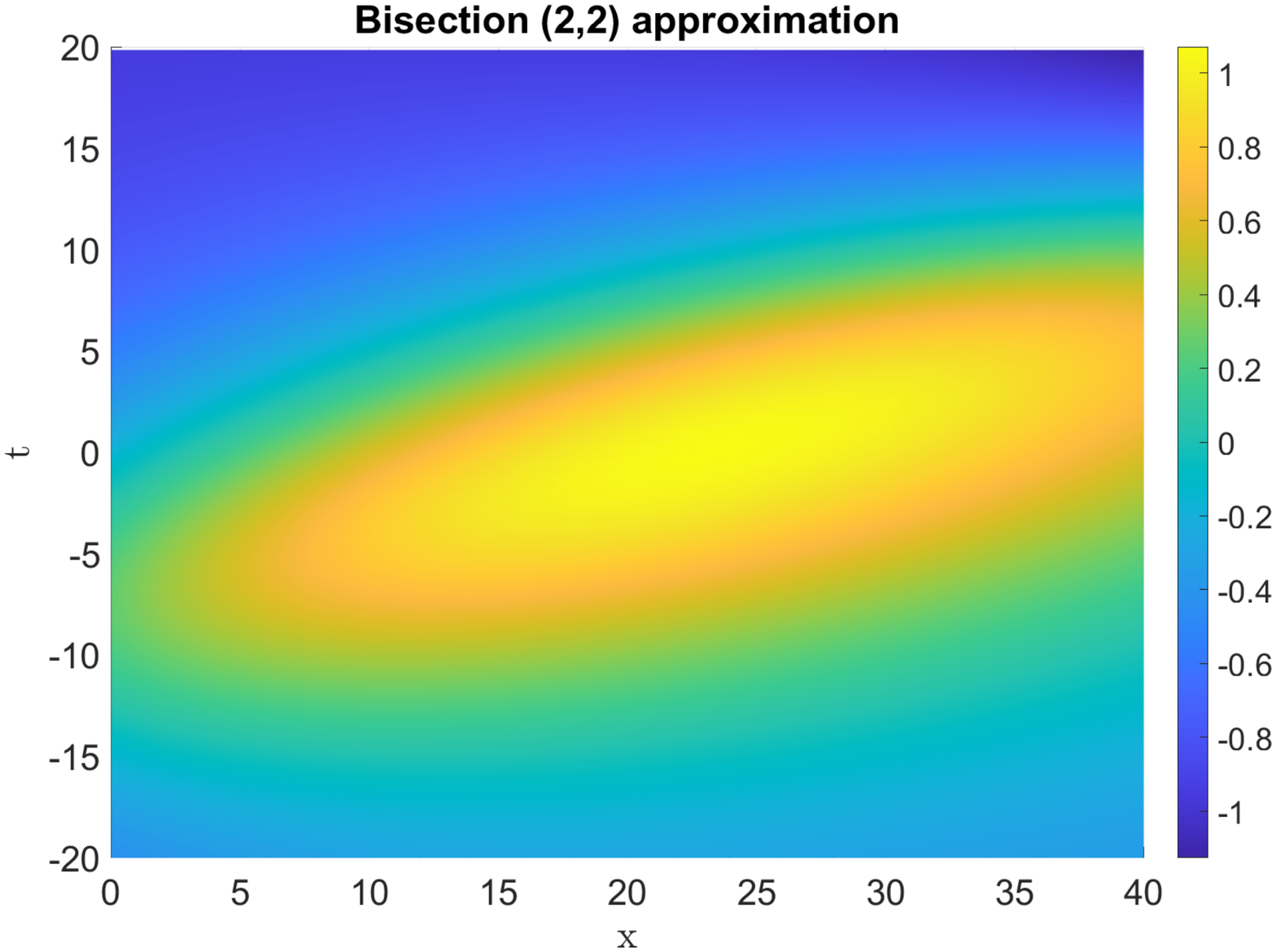}
		\includegraphics[width=60mm]{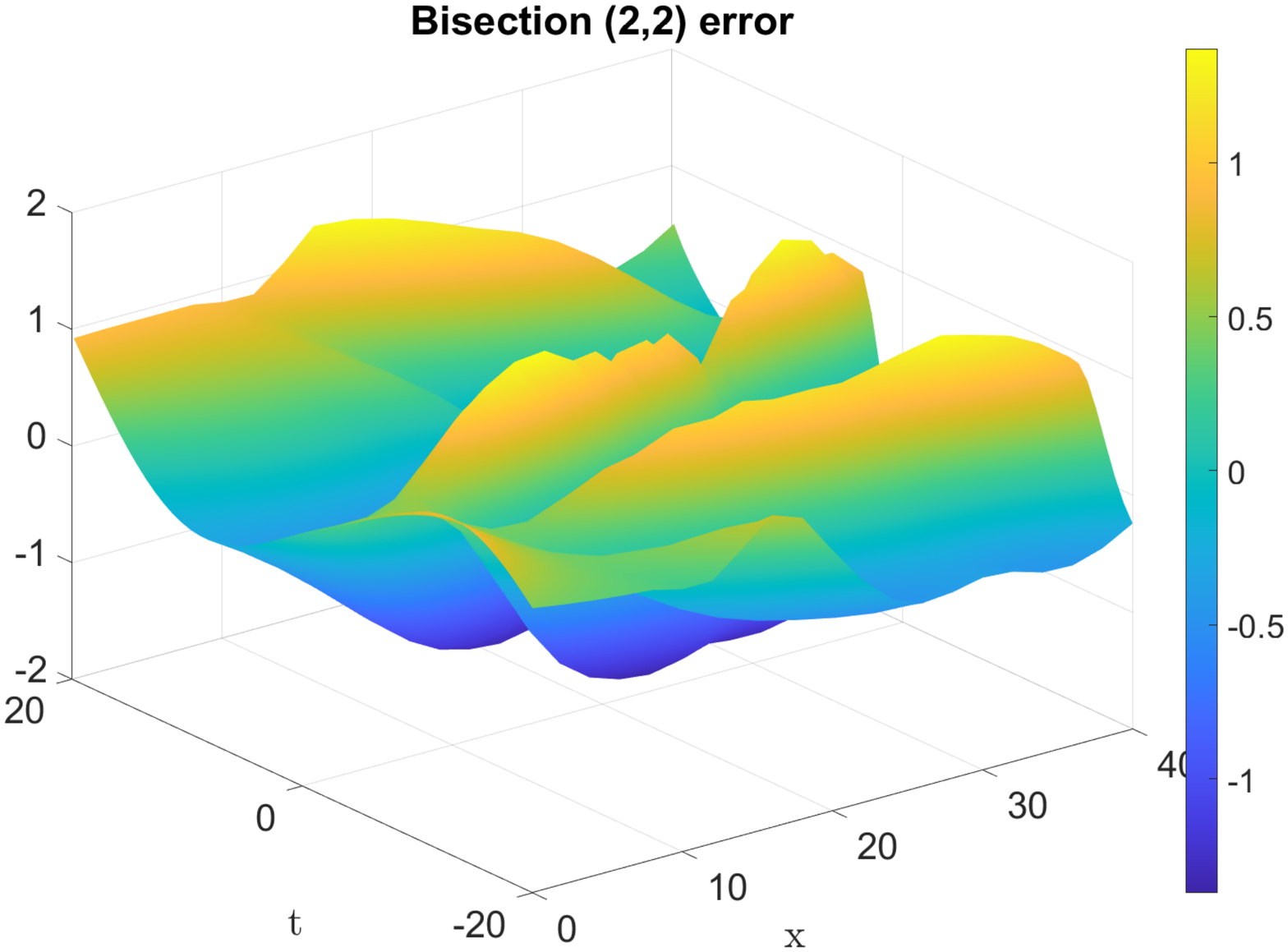}
		\includegraphics[width=60mm]{Error_2_2_20th_2D_DNB.pdf}
	\caption{Approximations (above) and the error curves (below) are constructed on the domain which consists of pairs of every $10^{\text{th}}$ point of the original domain: 3D and 2D view}
        \label{fig:Rational approximation of degree (2,2) multivariate domain of pairs of every 10th point}
\end{figure}

In Figure~\ref{fig:Rational approximation of degree (2,2) multivariate domain of pairs of every 20th point}, we compute the approximation by only considering the pairs of every $20^{\text{th}}$ point of the domain and in Figure~\ref{fig:Rational approximation of degree (2,2) multivariate domain of pairs of every 10th point}, the approximation is computed on the domain which contains pairs of every $10^{\text{th}}$ point of the domain. The computational time and the uniform error is presented in Table~\ref{tab:Aproximation of degree (2,2)}.

\begin{table}
	\centering
	\begin{tabular}{|c|c|c|}
		\cline{2-3}
		\multicolumn{1}{c|}{} &   Uniform error & Time (sec.)\\
		\hline
		Every $20^{\text{th}}$ point of the domain & 1.29449744612994 & 2.049840 \\
		\hline
		Every $10^{\text{th}}$ point of the domain & 1.36991465992617 & 2.504411 \\
		\hline
	\end{tabular}
	\caption{Uniform error and computational time for the degree (2,2) approximation}
        \label{tab:Aproximation of degree (2,2)}
\end{table}

One can clearly see that the algorithm runs faster when the number of discretised points in the domain is smaller.

\subsection{Rational approximation of degree (5,5)}
Here, we now increase the degree of the approximation to degree (5,5).

\begin{figure}
	\centering
		\includegraphics[width=60mm]{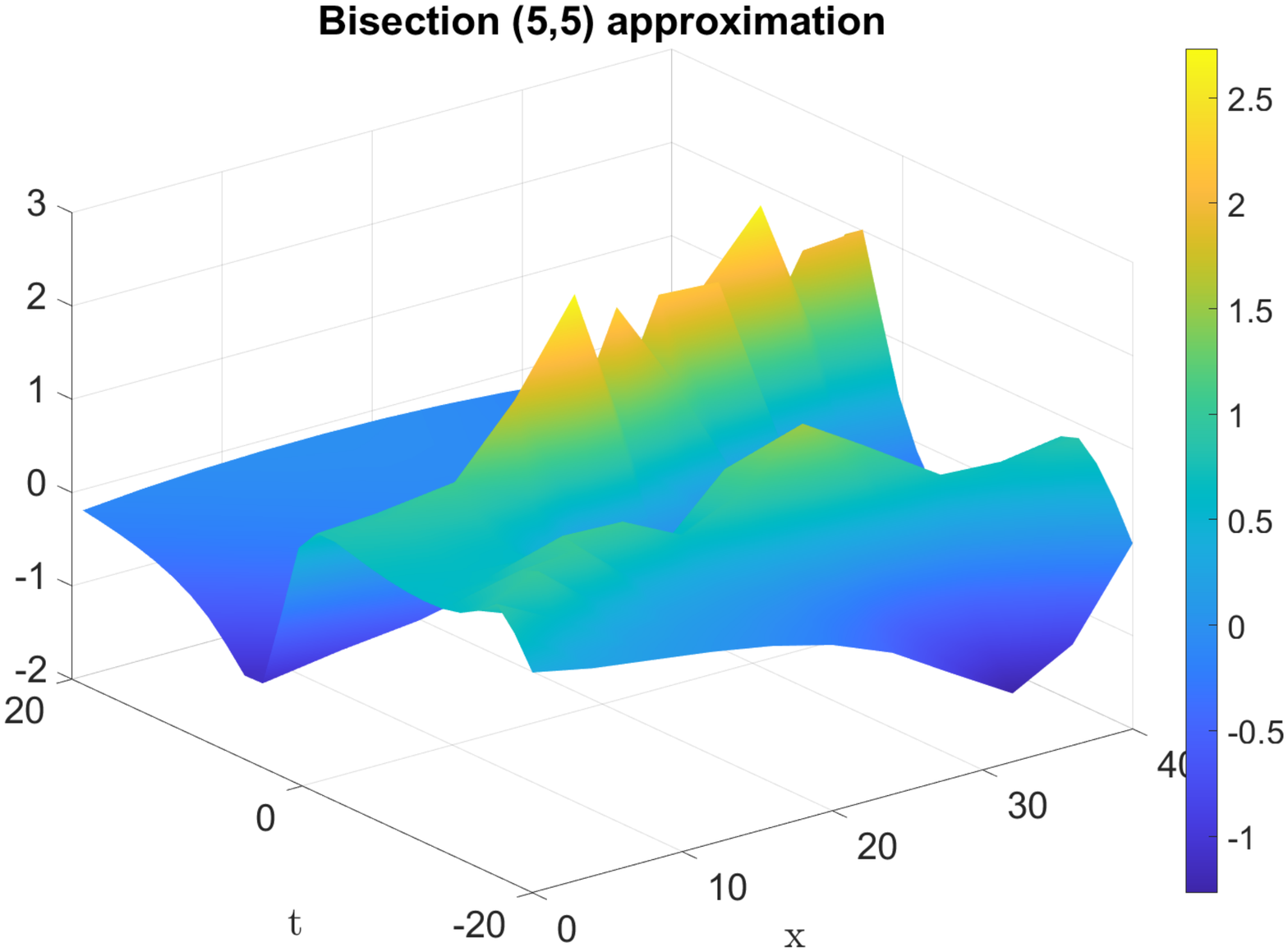}
		\includegraphics[width=60mm]{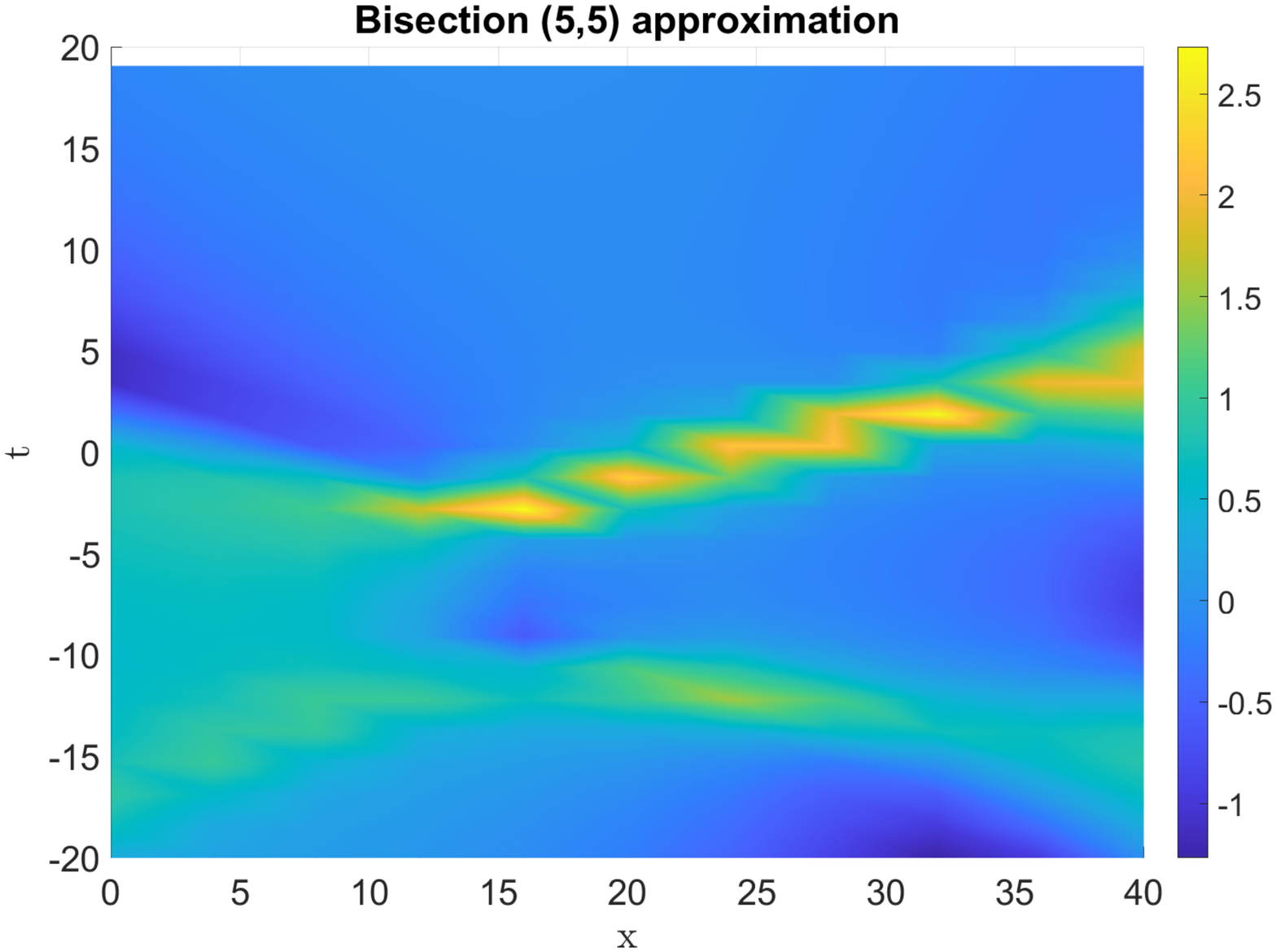}
		\includegraphics[width=60mm]{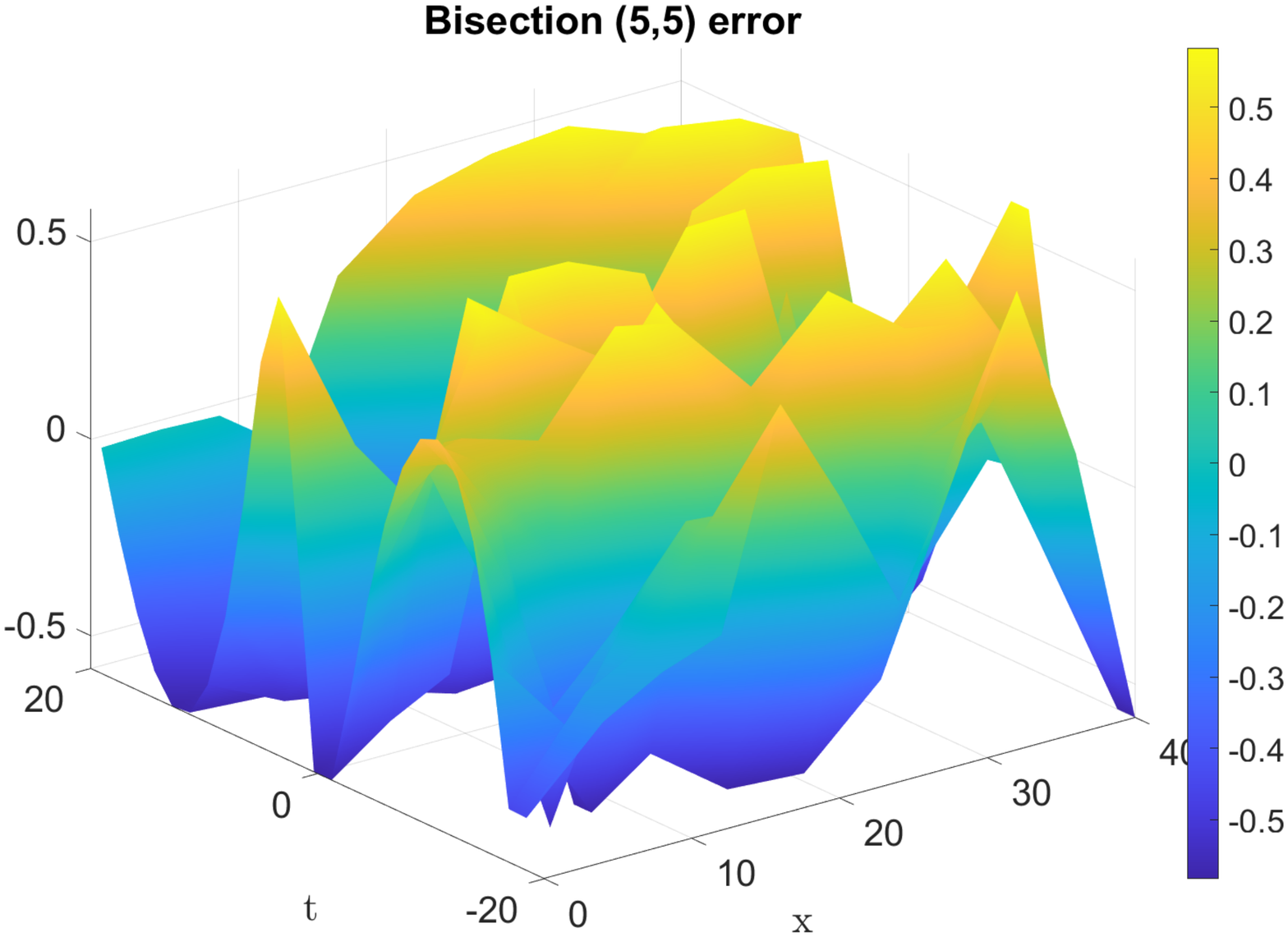}
		\includegraphics[width=60mm]{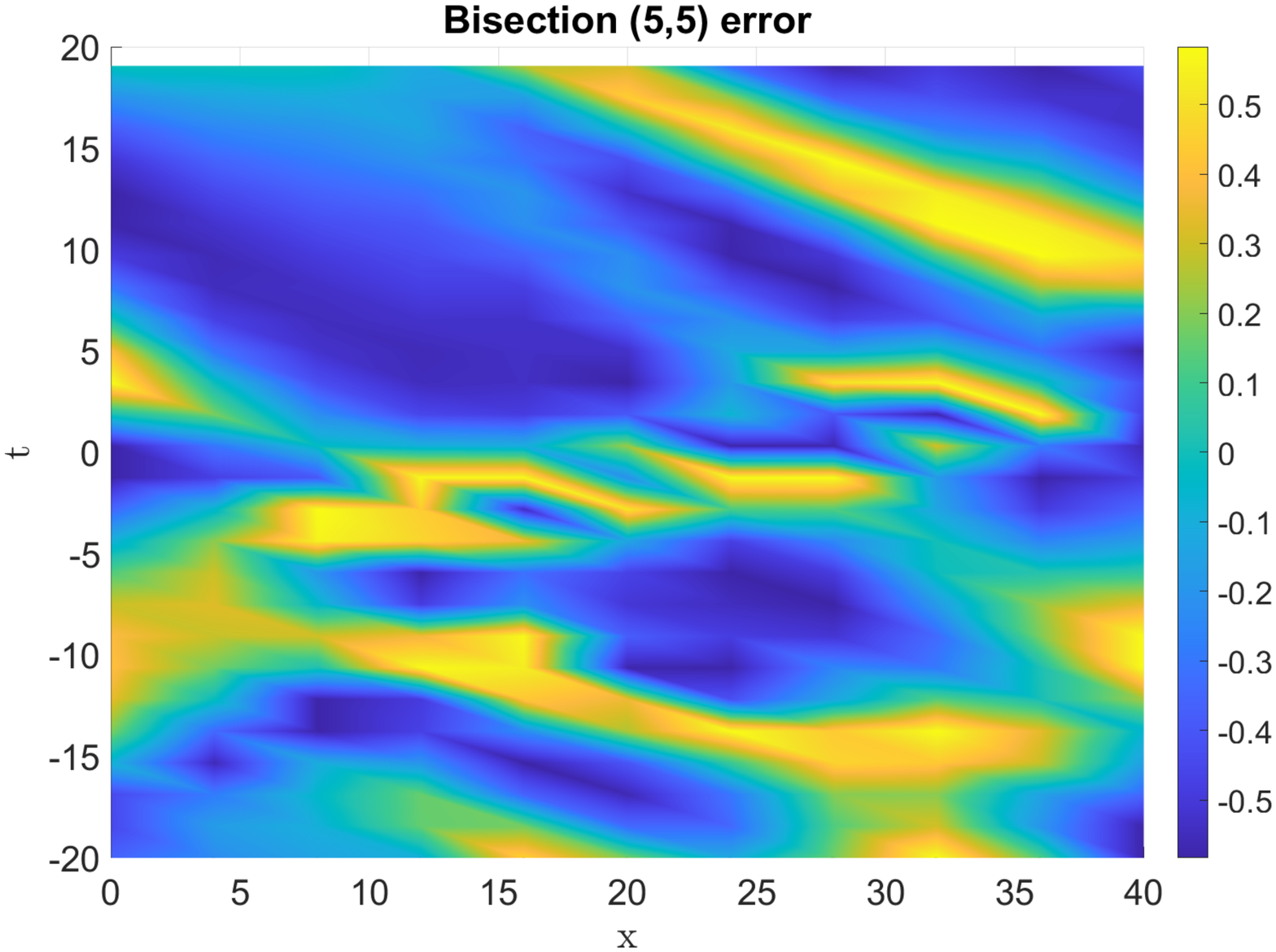}
	\caption{Approximations (above) and the error curves (below) are constructed on the domain which consists of pairs of every $20^{\text{th}}$ point of the original domain: 3D and 2D view.}
        \label{fig:Rational approximation of degree (5,5) multivariate domain of pairs of every 20th point}
\end{figure}

\begin{figure}
	\centering
		\includegraphics[width=60mm]{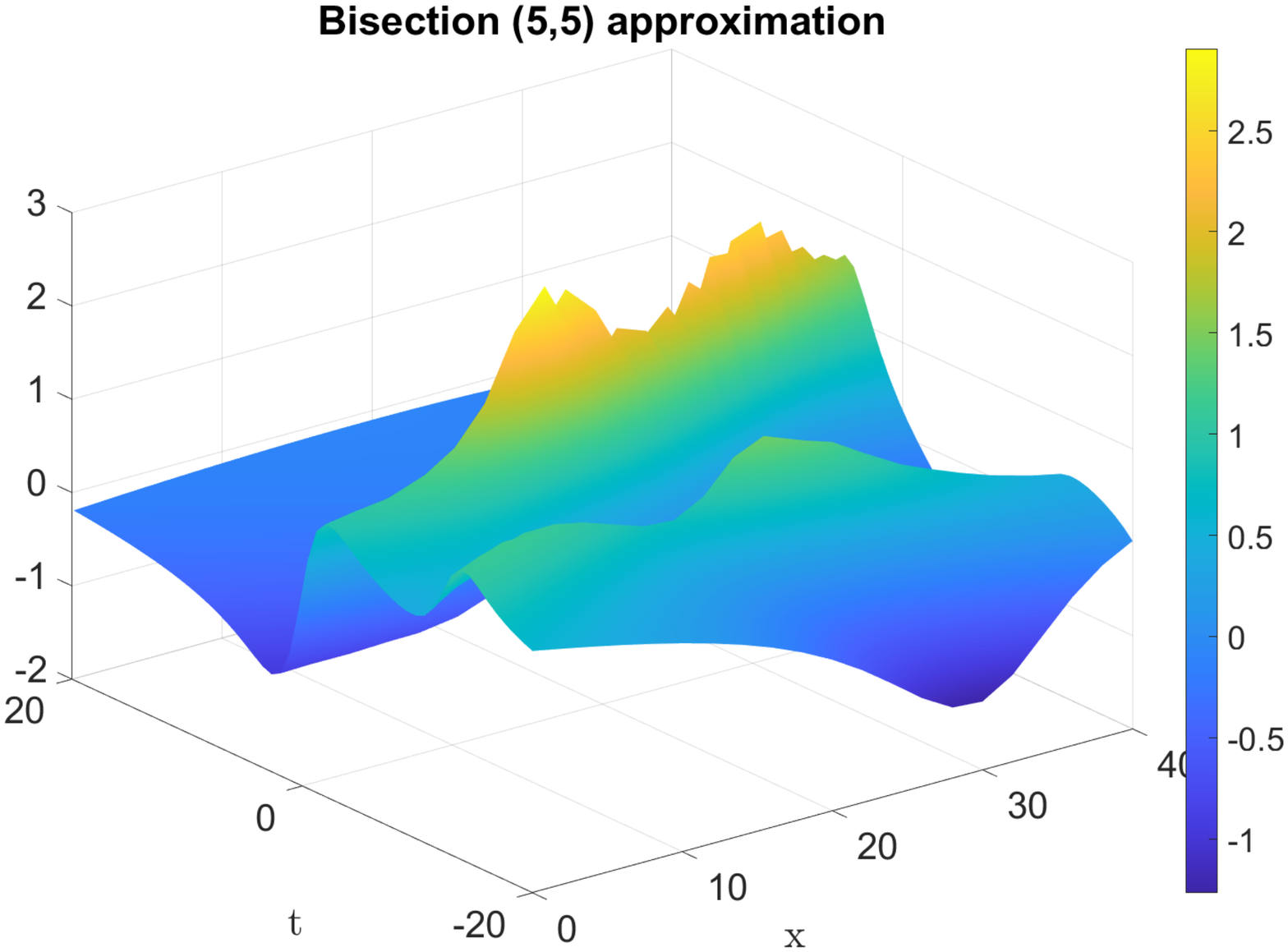}
		\includegraphics[width=60mm]{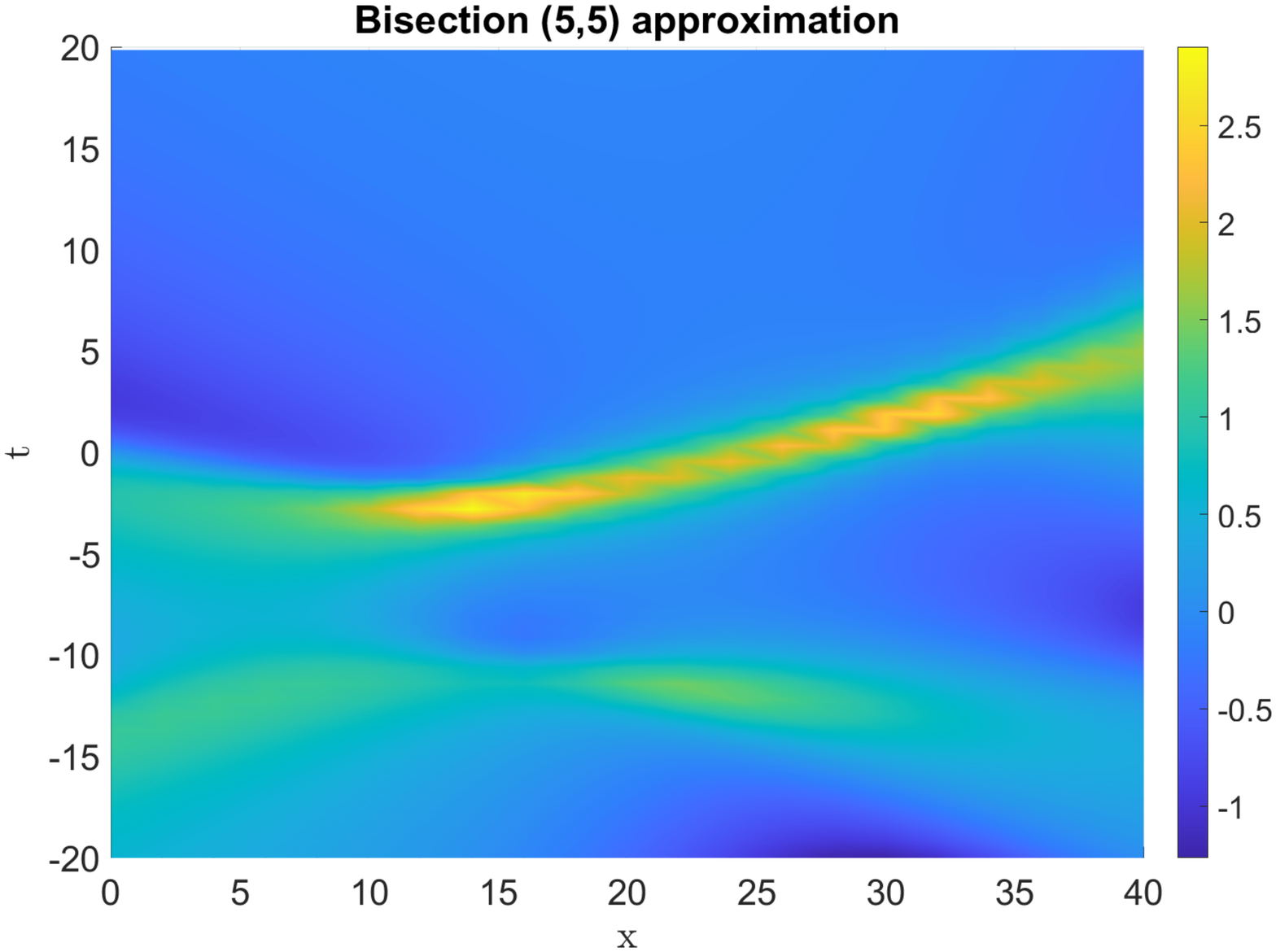}
		\includegraphics[width=60mm]{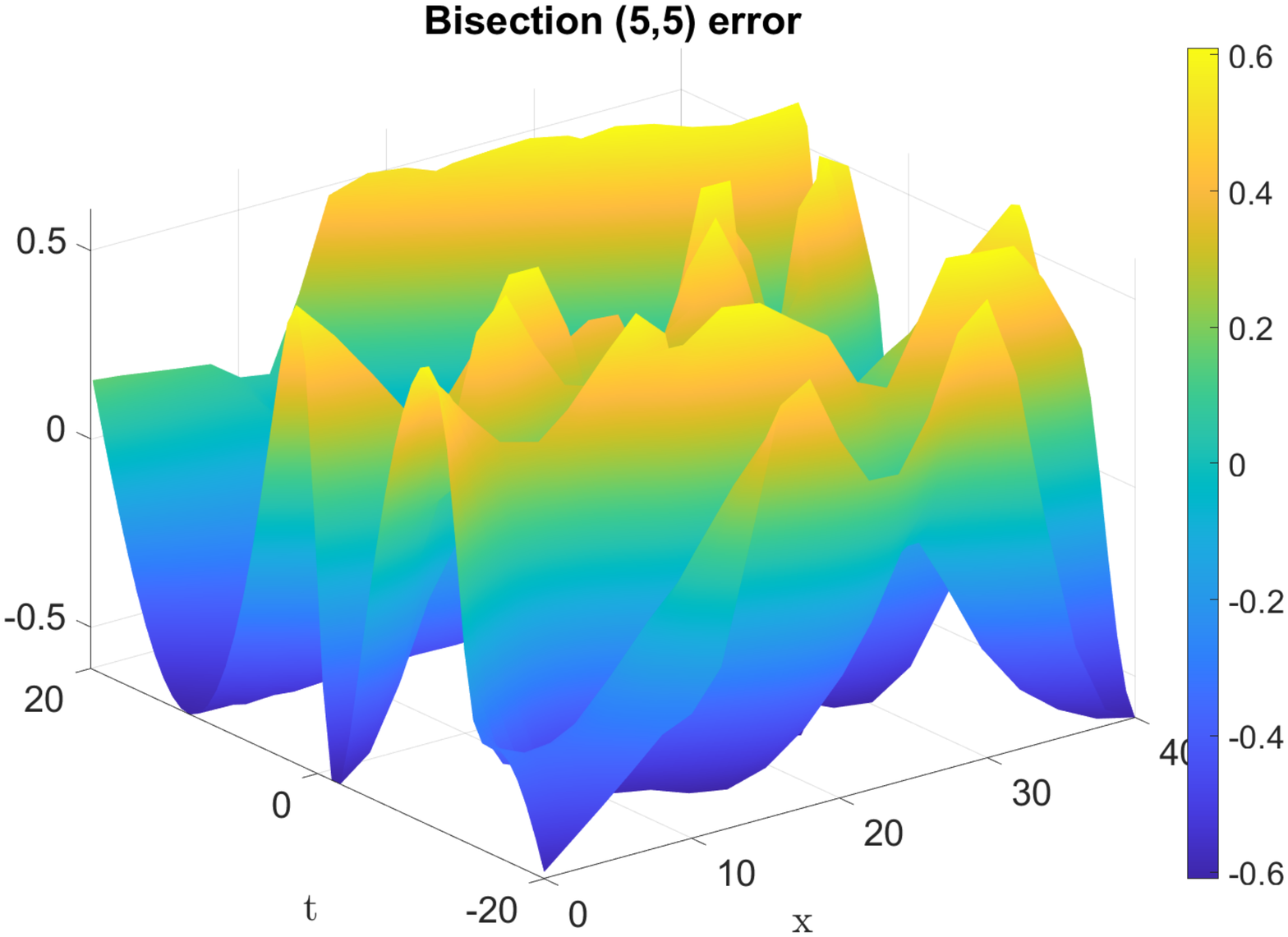}
		\includegraphics[width=60mm]{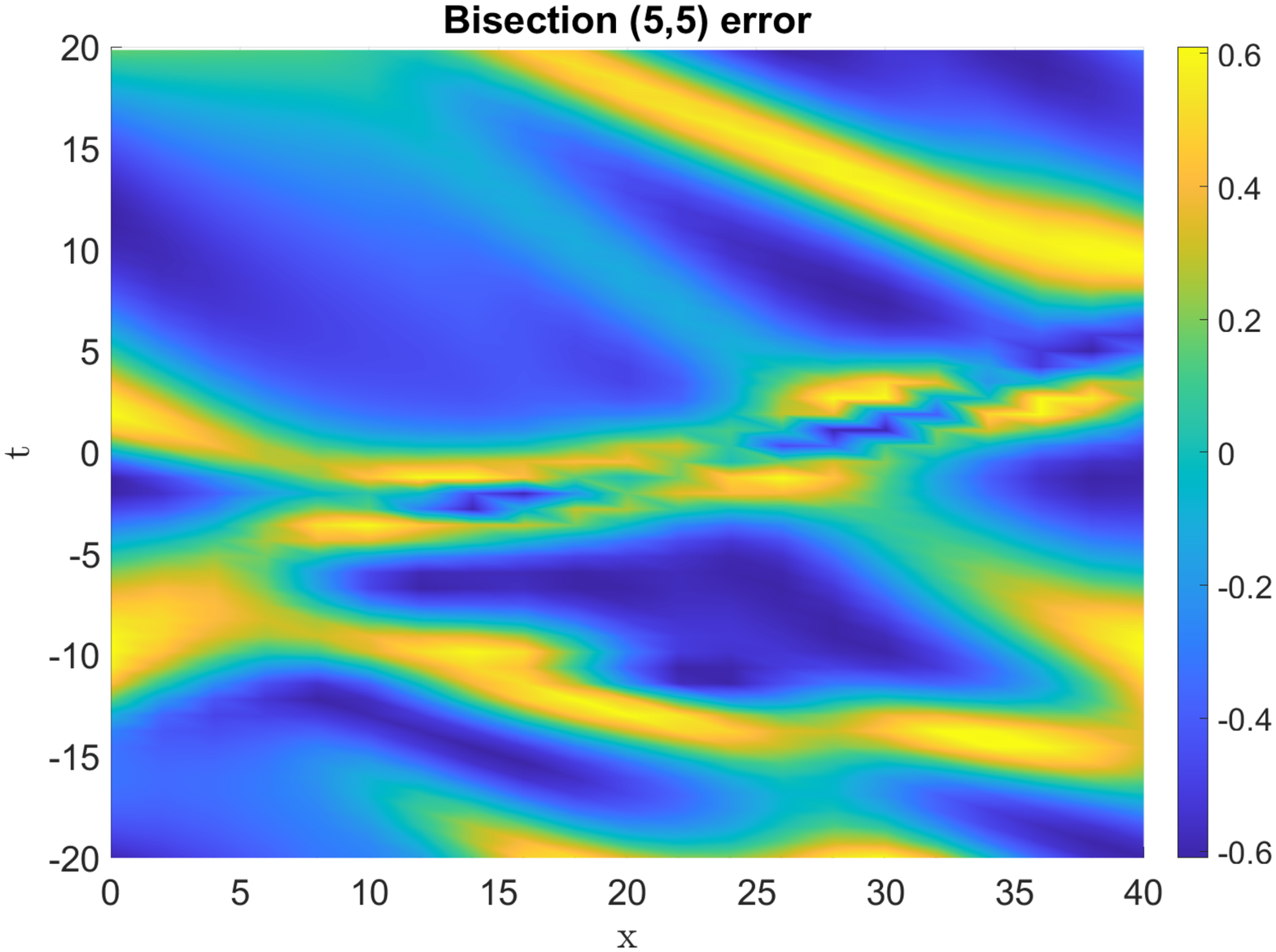}
	\caption{Approximations (above) and the error curves (below) are constructed on the domain which consists of pairs of every $10^{\text{th}}$ point of the original domain: 3D and 2D view.}
        \label{fig:Rational approximation of degree (5,5) multivariate domain of pairs of every 10th point}
\end{figure}

Figure~\ref{fig:Rational approximation of degree (5,5) multivariate domain of pairs of every 20th point} and Figure~\ref{fig:Rational approximation of degree (5,5) multivariate domain of pairs of every 10th point} consist of the approximations computed on the domains which contain pairs of every $20^{\text{th}}$ point of the domain and the $10^{\text{th}}$ point of the domain, respectively. The computational time and the uniform error is presented in Table~\ref{tab:Aproximation of degree (5,5)}. The approximations are better than the degree (2,2) case and the uniform error is drastically reduced.

\begin{table}
	\centering
	\begin{tabular}{|c|c|c|}
		\cline{2-3}
		\multicolumn{1}{c|}{}  &  Uniform error & Time (sec.)\\
		\hline
		Every $20^{\text{th}}$ point of the domain & 0.583136100002126 & 7.823996 \\
		\hline
		Every $10^{\text{th}}$ point of the domain & 0.609975243086065 & 10.455391 \\
		\hline
	\end{tabular}
	\caption{Uniform error and computational time for the degree (5,5) approximation}
        \label{tab:Aproximation of degree (5,5)}
\end{table}

\subsection{Rational approximation of degree (10,10)}
In this section, we approximate the original function by a rational approximation of degree (10,10). The approximations and the error curves are depicted in Figure~\ref{fig:Rational approximation of degree (10,10) multivariate domain of pairs of every 20th point} and Figure~\ref{fig:Rational approximation of degree (10,10) multivariate domain of pairs of every 10th point}. Approximations resemble the original functions well in both domains.

\begin{figure}
	\centering
		\includegraphics[width=60mm]{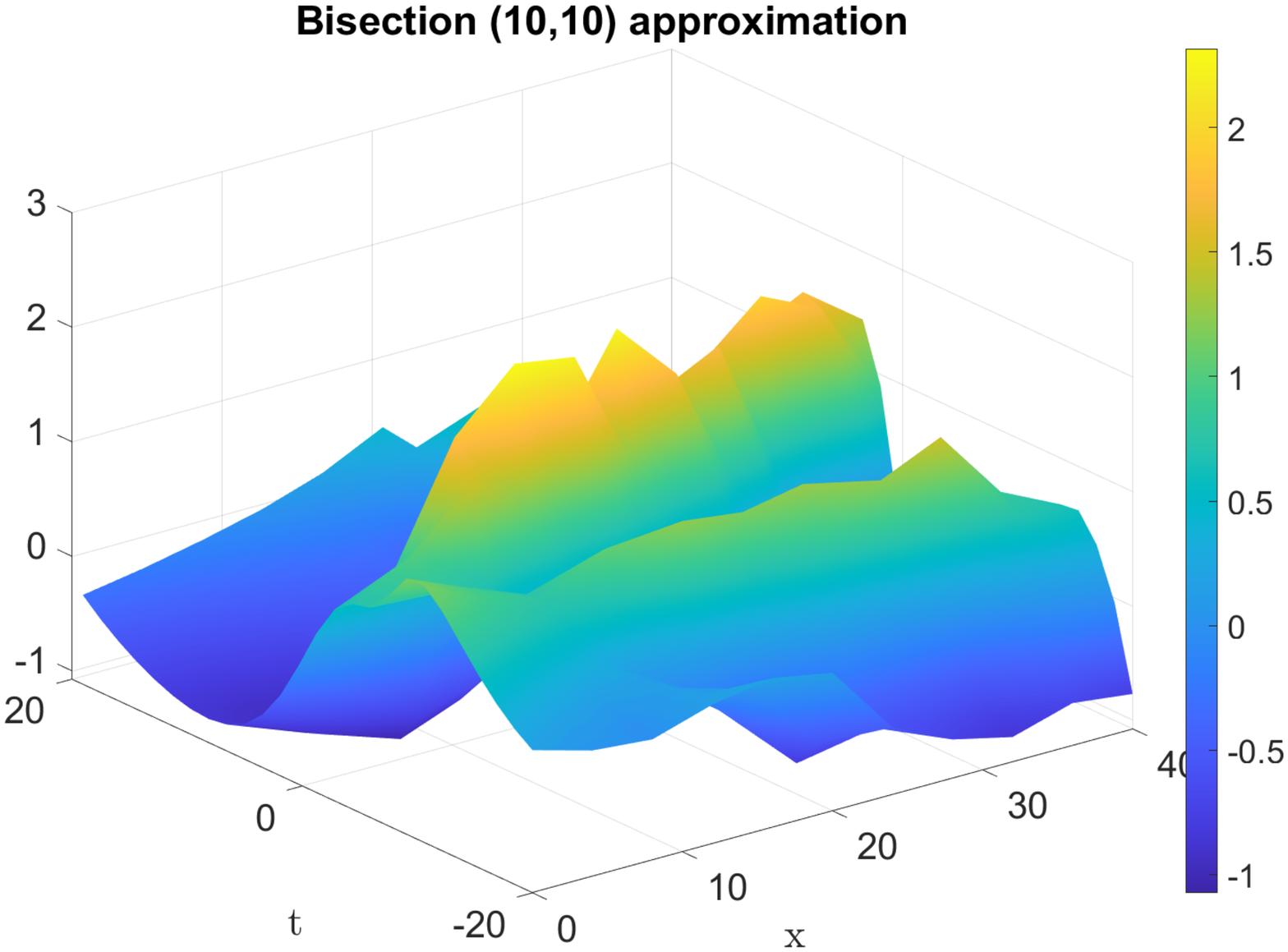}
		\includegraphics[width=60mm]{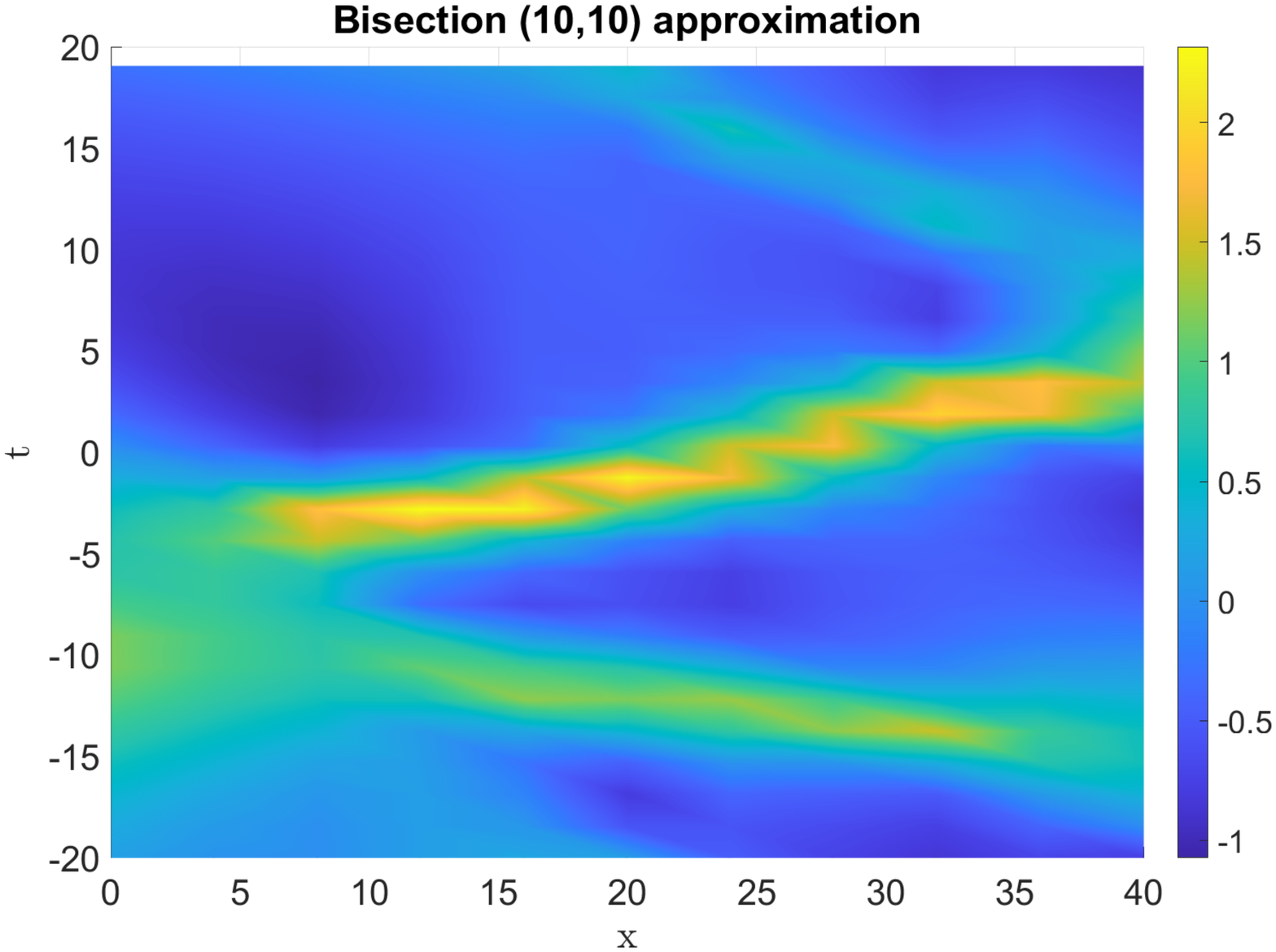}
		\includegraphics[width=60mm]{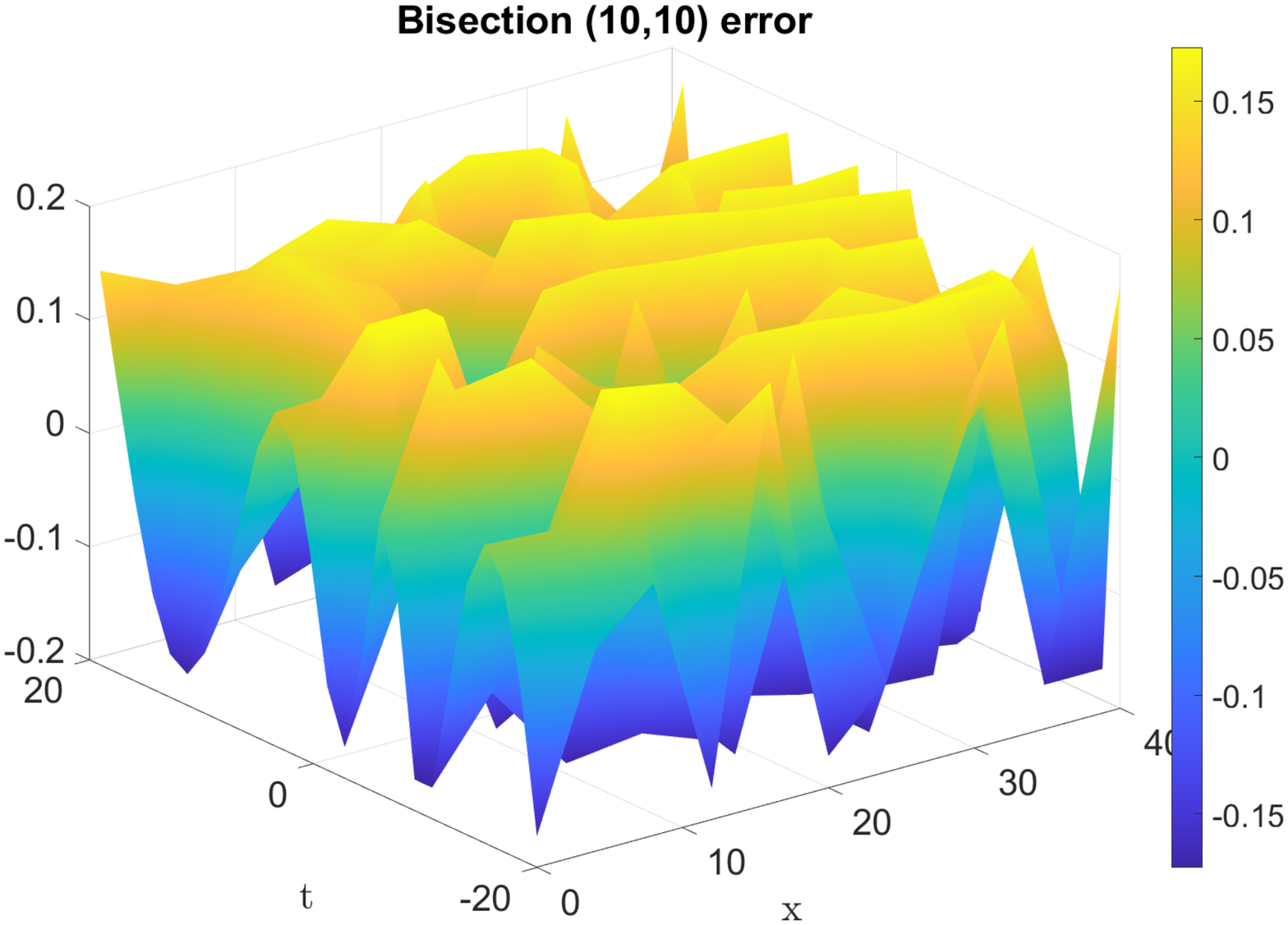}
		\includegraphics[width=60mm]{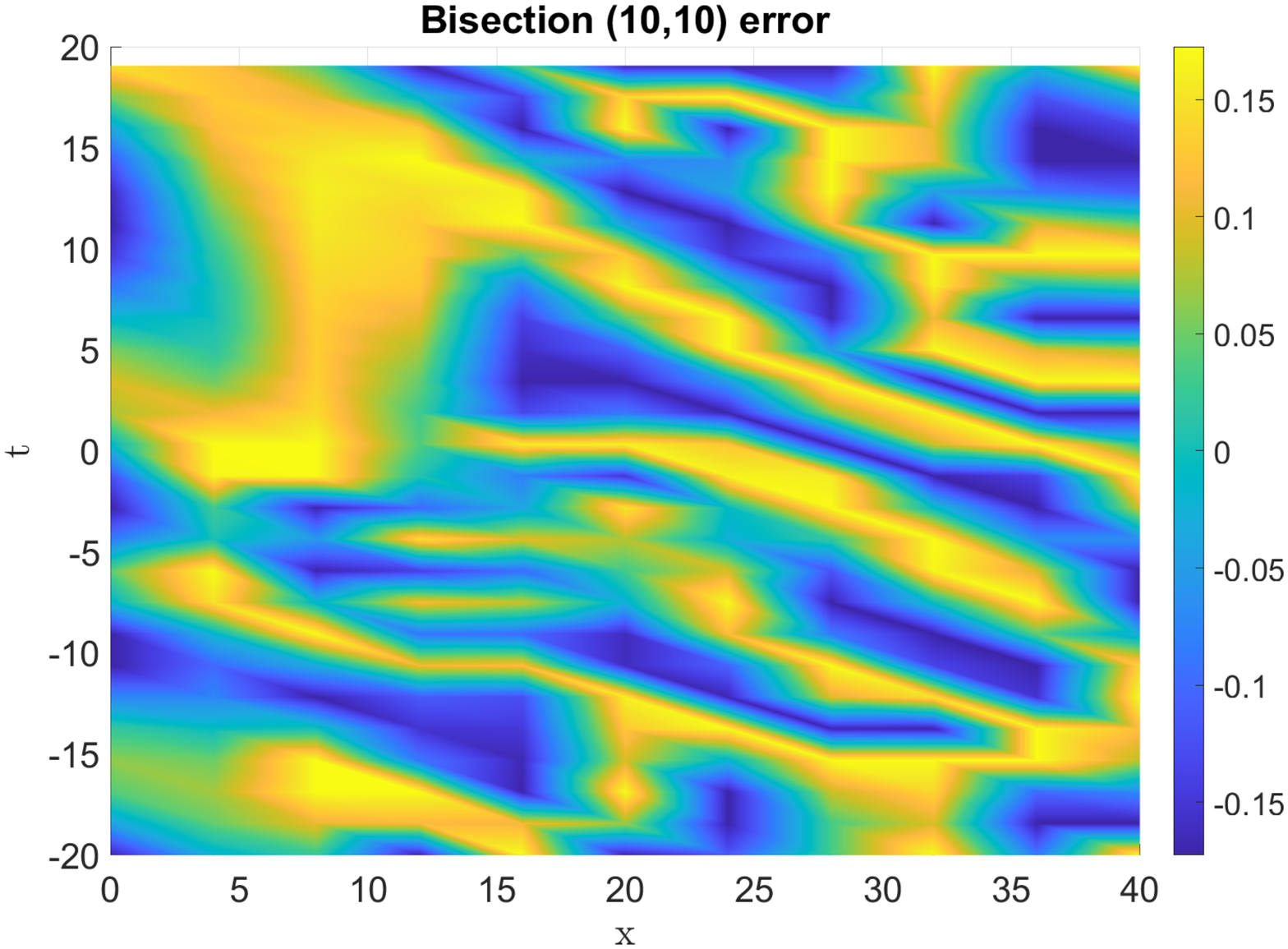}
	\caption{Approximations (above) and the error curves (below) are constructed on the domain which consists of pairs of every $20^{\text{th}}$ point of the original domain: 3D and 2D view.}
        \label{fig:Rational approximation of degree (10,10) multivariate domain of pairs of every 20th point}
\end{figure}

\begin{figure}
	\centering
		\includegraphics[width=60mm]{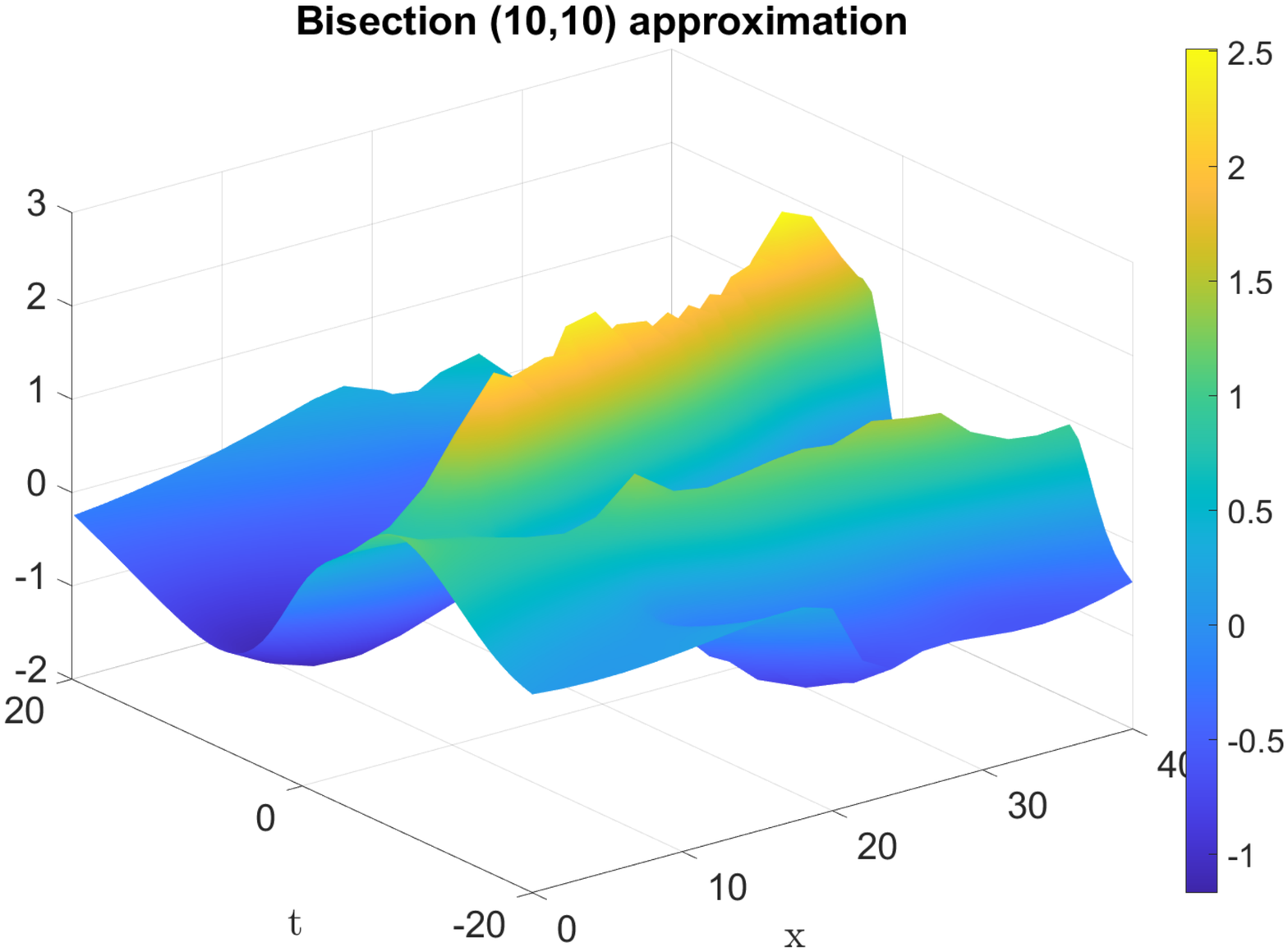}
		\includegraphics[width=60mm]{Rational_approximation_10_10_20th_2D_DNB.pdf}
		\includegraphics[width=60mm]{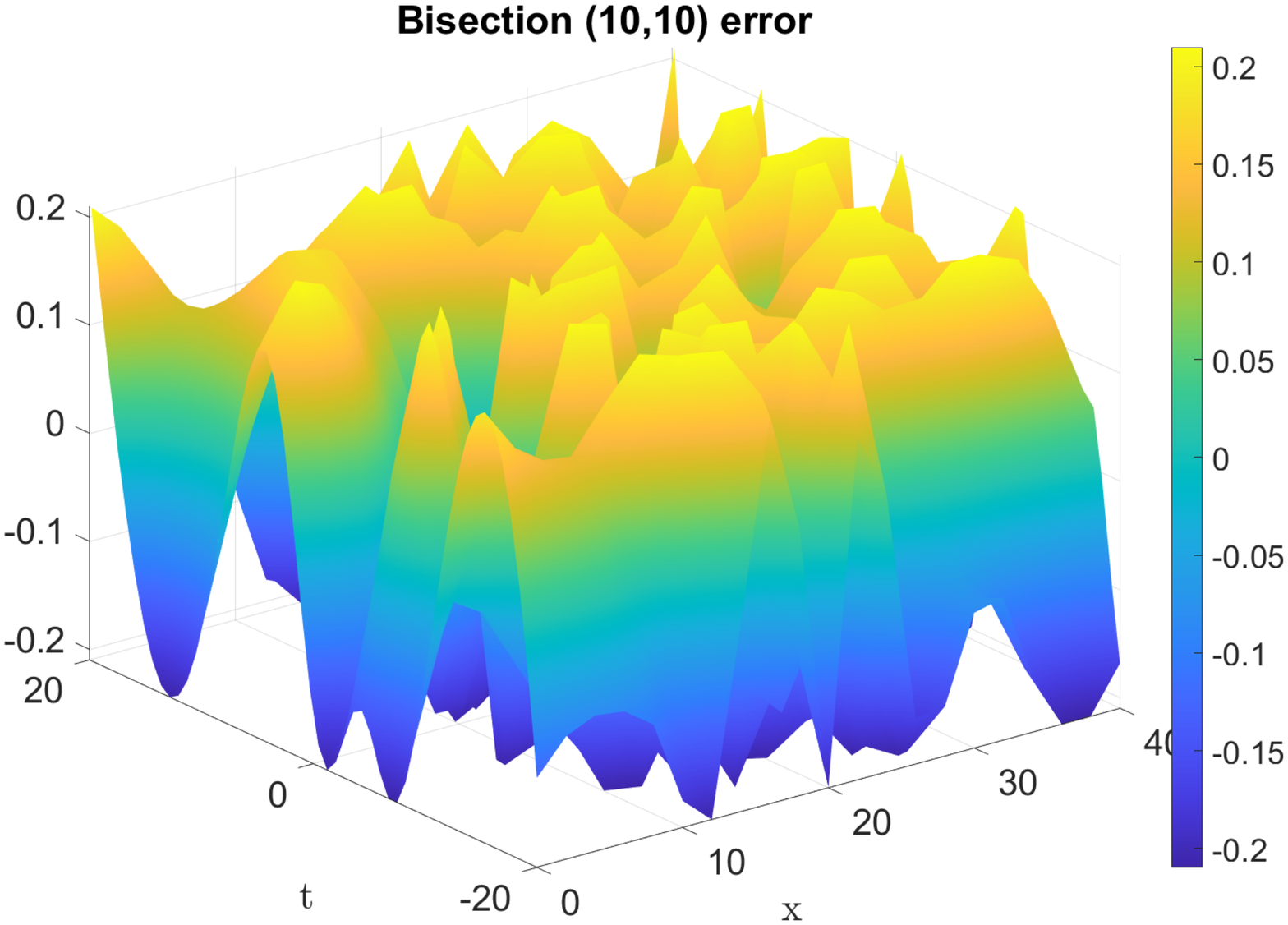}
		\includegraphics[width=60mm]{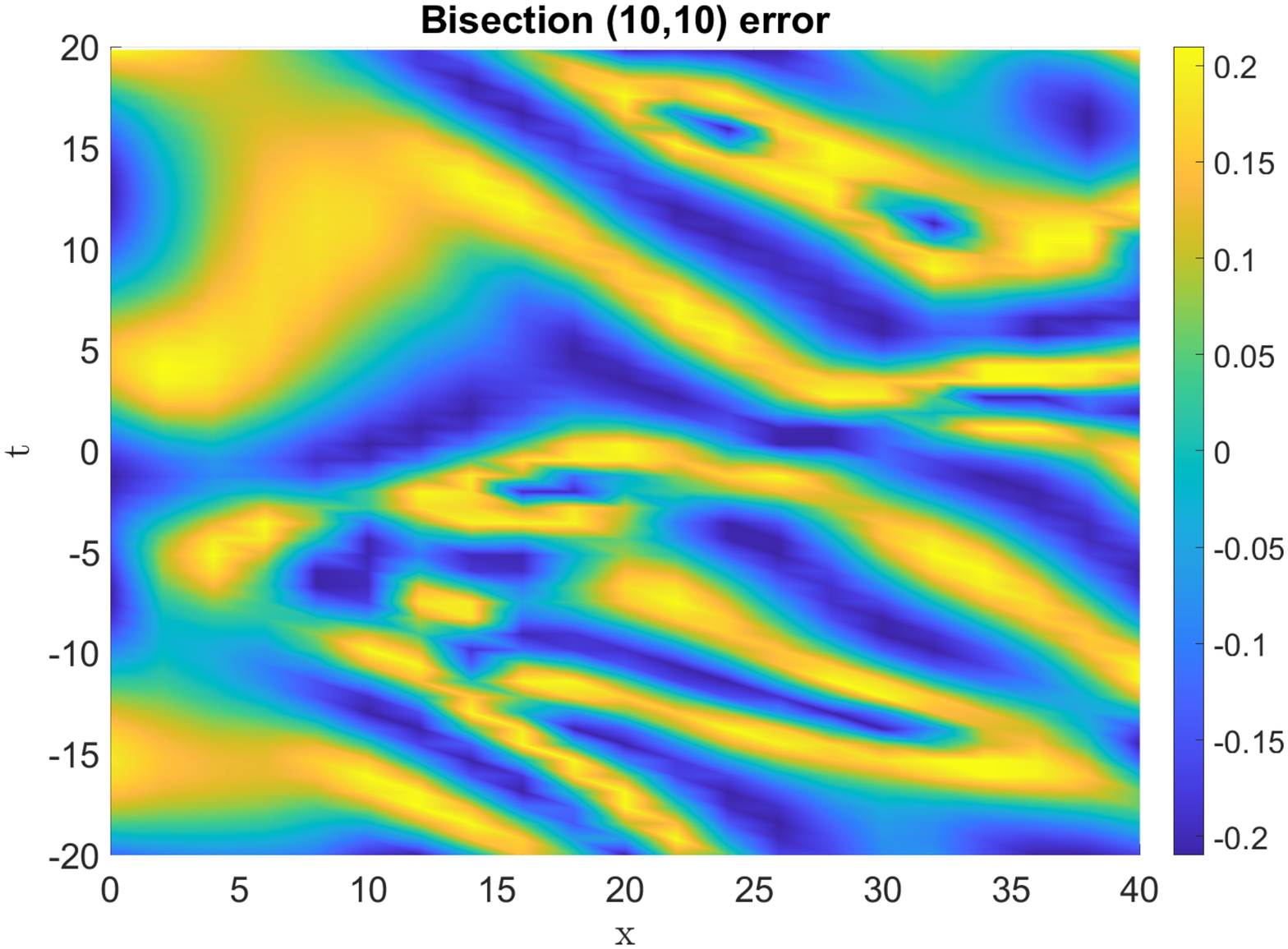}
	\caption{Approximations (above) and the error curves (below) are constructed on the domain which consists of pairs of every $10^{\text{th}}$ point of the original domain.}
        \label{fig:Rational approximation of degree (10,10) multivariate domain of pairs of every 10th point}
\end{figure}

The computational time and the uniform error are presented in Table~\ref{tab:Aproximation of degree (10,10)}. The uniform error term is much smaller than the case where the degree of the approximation is (5,5). On the other hand, the algorithm takes more time to compute the approximations as the degree increases,

\begin{table}
	\centering
	\begin{tabular}{|c|c|c|}
		\cline{2-3}
		\multicolumn{1}{c|}{}  &  Uniform error & Time (sec.)\\
		\hline
		Every $20^{\text{th}}$ point of the domain & 0.172559761245446 & 22.208227 \\
		\hline
		Every $10^{\text{th}}$ point of the domain & 0.209717782244555 & 42.384391 \\
		\hline
	\end{tabular}
	\caption{Uniform error and computational time for the degree (10,10) approximation}
        \label{tab:Aproximation of degree (10,10)}
\end{table}

\subsection{Rational approximation of degree (18,18)}
We now increase the degree of the approximation to (18,18). The approximations presented in Figure~\ref{fig:Rational approximation of degree (18,18) multivariate domain of pairs of every 20th point} and Figure~\ref{fig:Rational approximation of degree (18,18) multivariate domain of pairs of every 10th point} are very similar to the corresponding original functions presented in Figure~\ref{fig:original dataset with pairs of every 20th point of the domain} and Figure~\ref{fig:original dataset with pairs of every 10th point of the domain}. One can also see from table~\ref{tab:Aproximation of degree (18,18)} that the uniform error is much smaller than the previous cases.

\begin{figure}
	\centering
		\includegraphics[width=60mm]{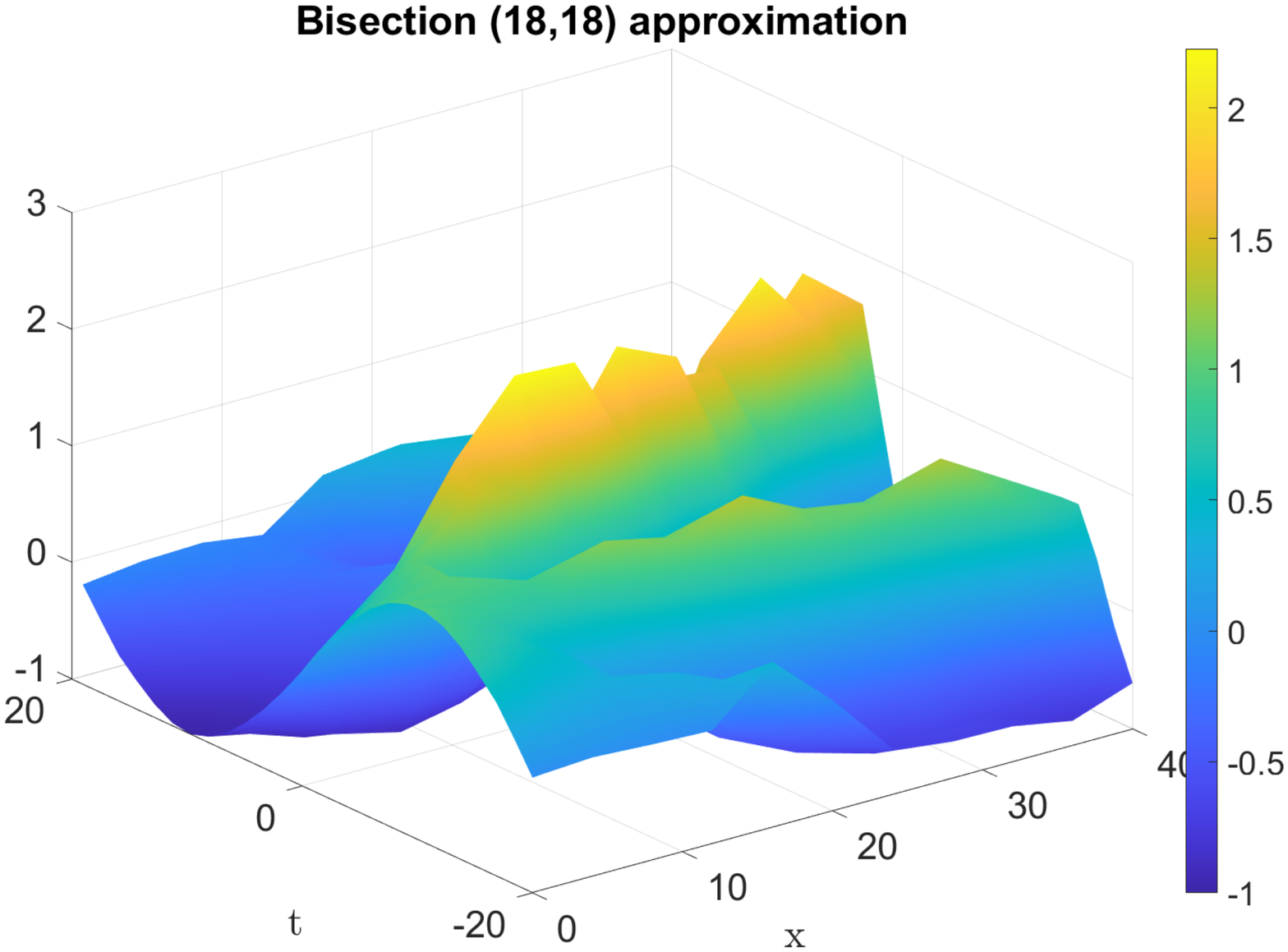}
		\includegraphics[width=60mm]{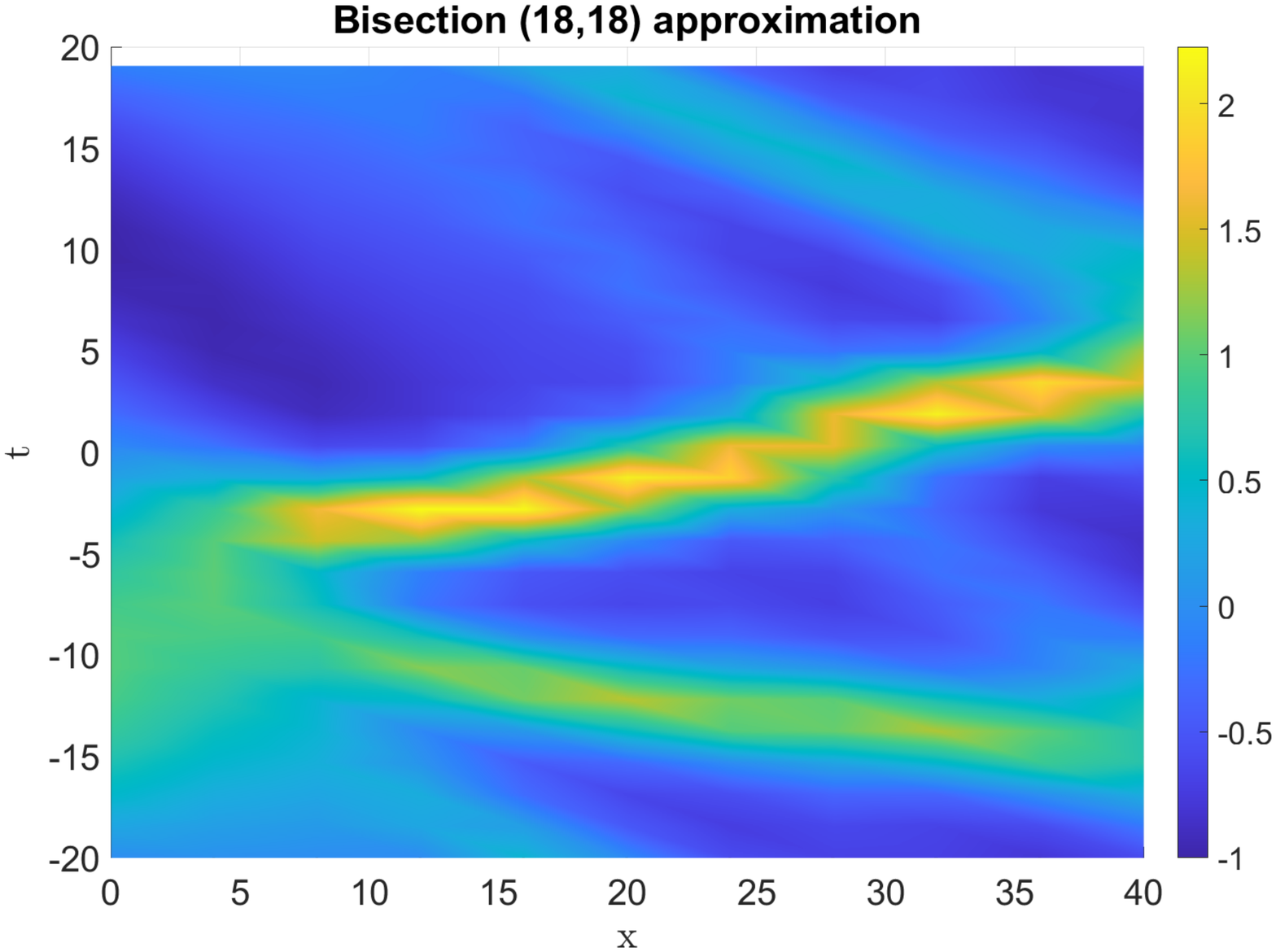}
		\includegraphics[width=60mm]{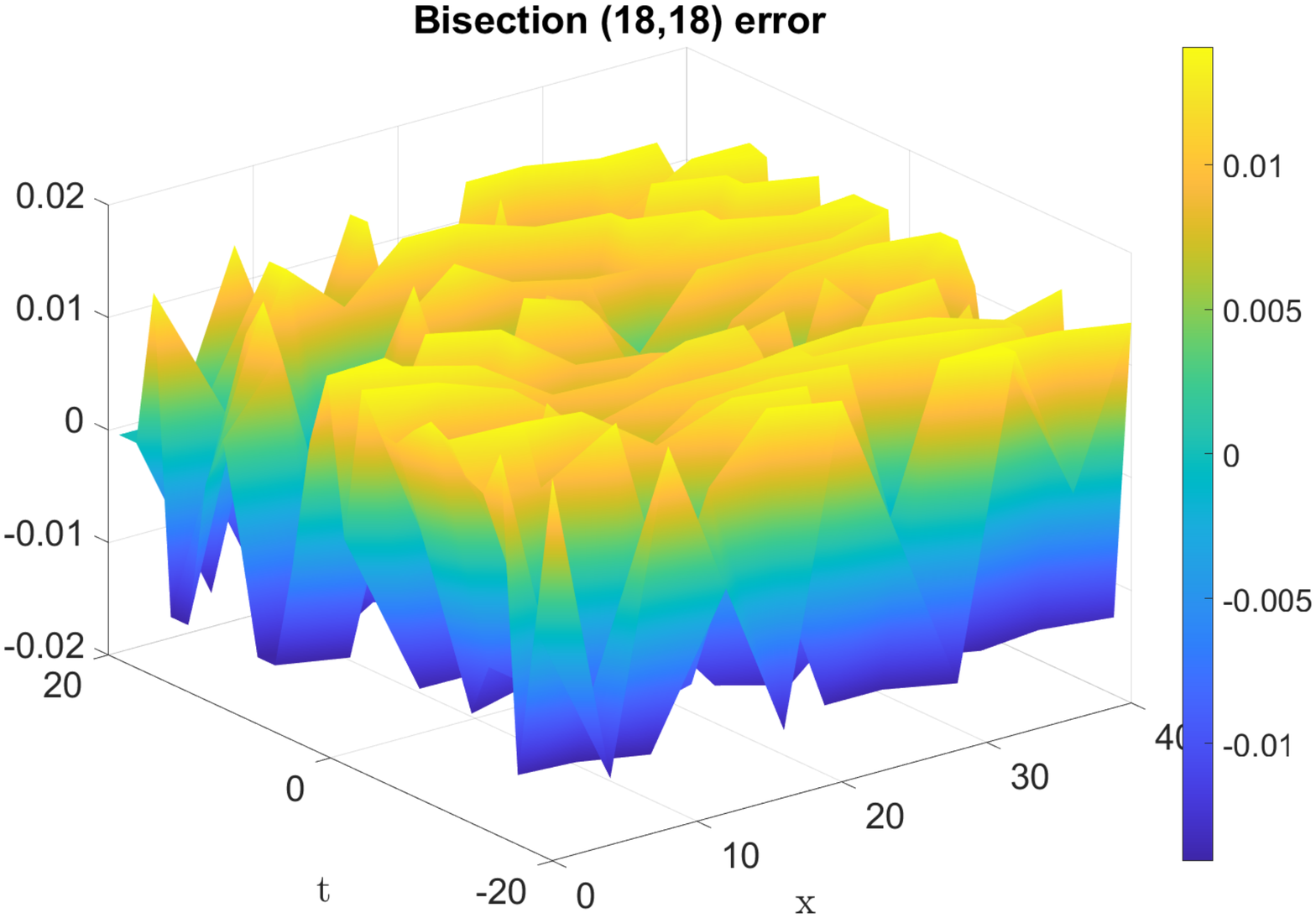}
		\includegraphics[width=60mm]{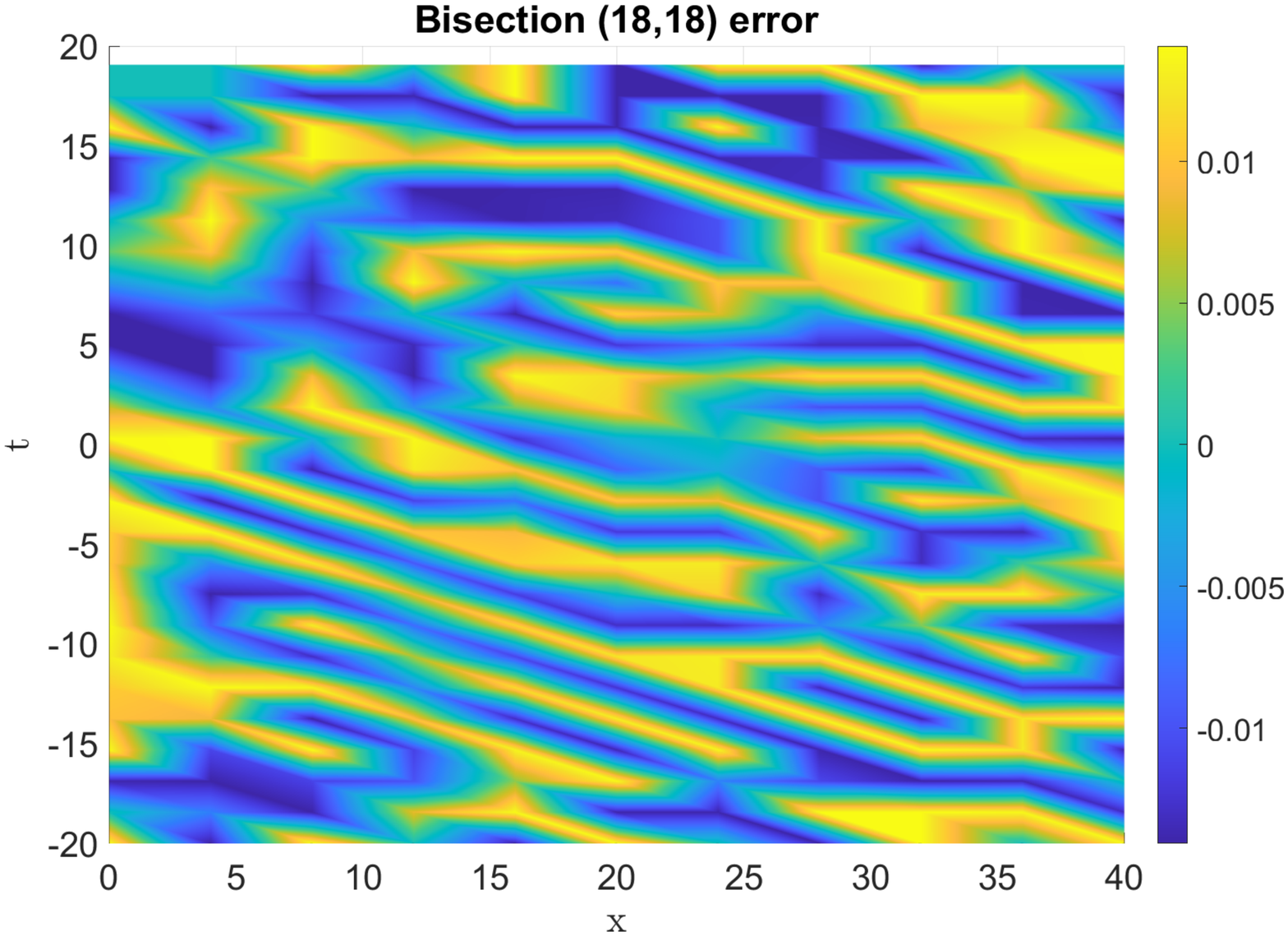}
	\caption{Approximations (above) and the error curves (below) are constructed on the domain which consists of pairs of every $20^{\text{th}}$ point of the original domain: 3D and 2D view.}
        \label{fig:Rational approximation of degree (18,18) multivariate domain of pairs of every 20th point}
\end{figure}

\begin{figure}
	\centering
		\includegraphics[width=60mm]{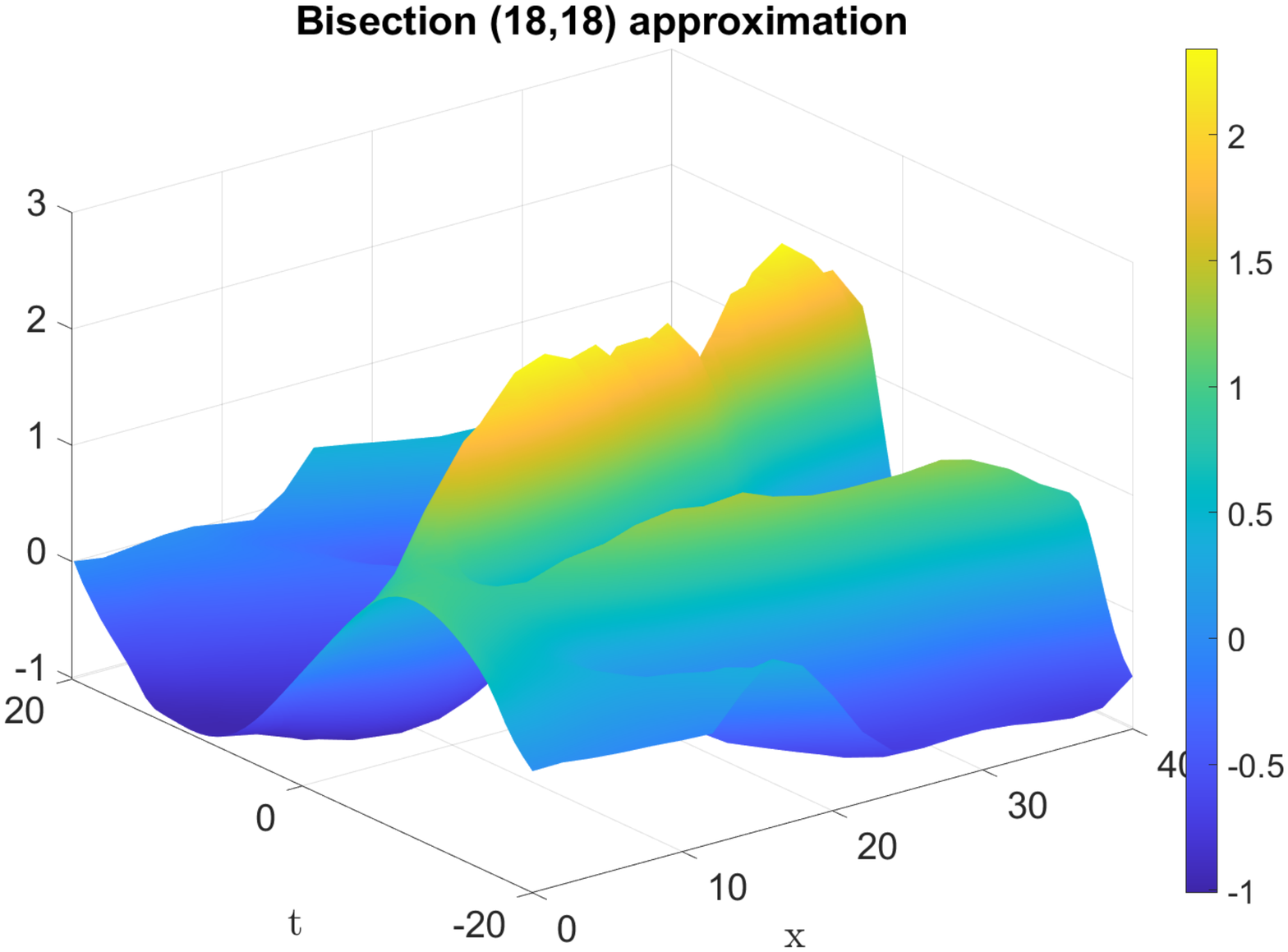}
		\includegraphics[width=60mm]{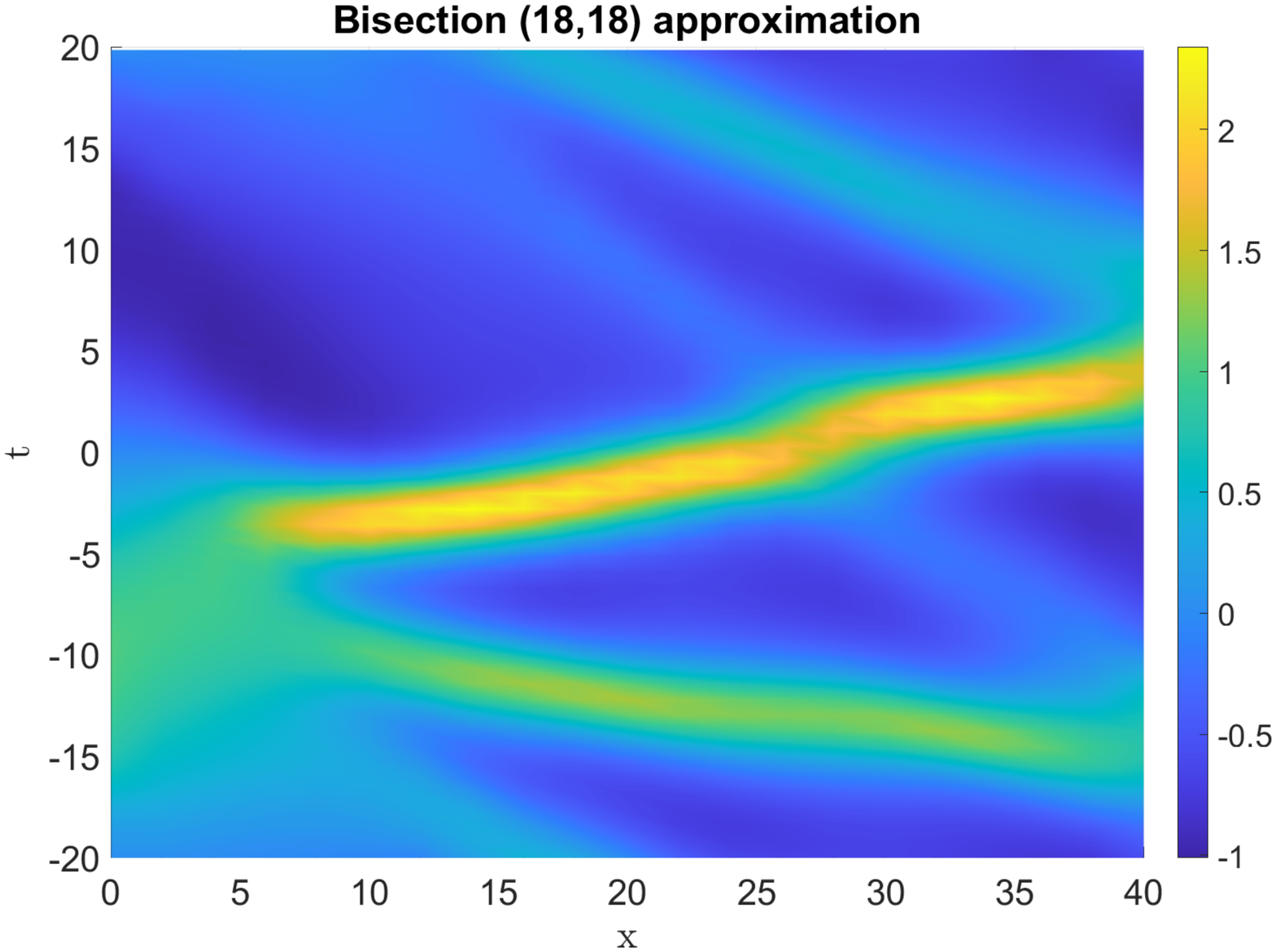}
		\includegraphics[width=60mm]{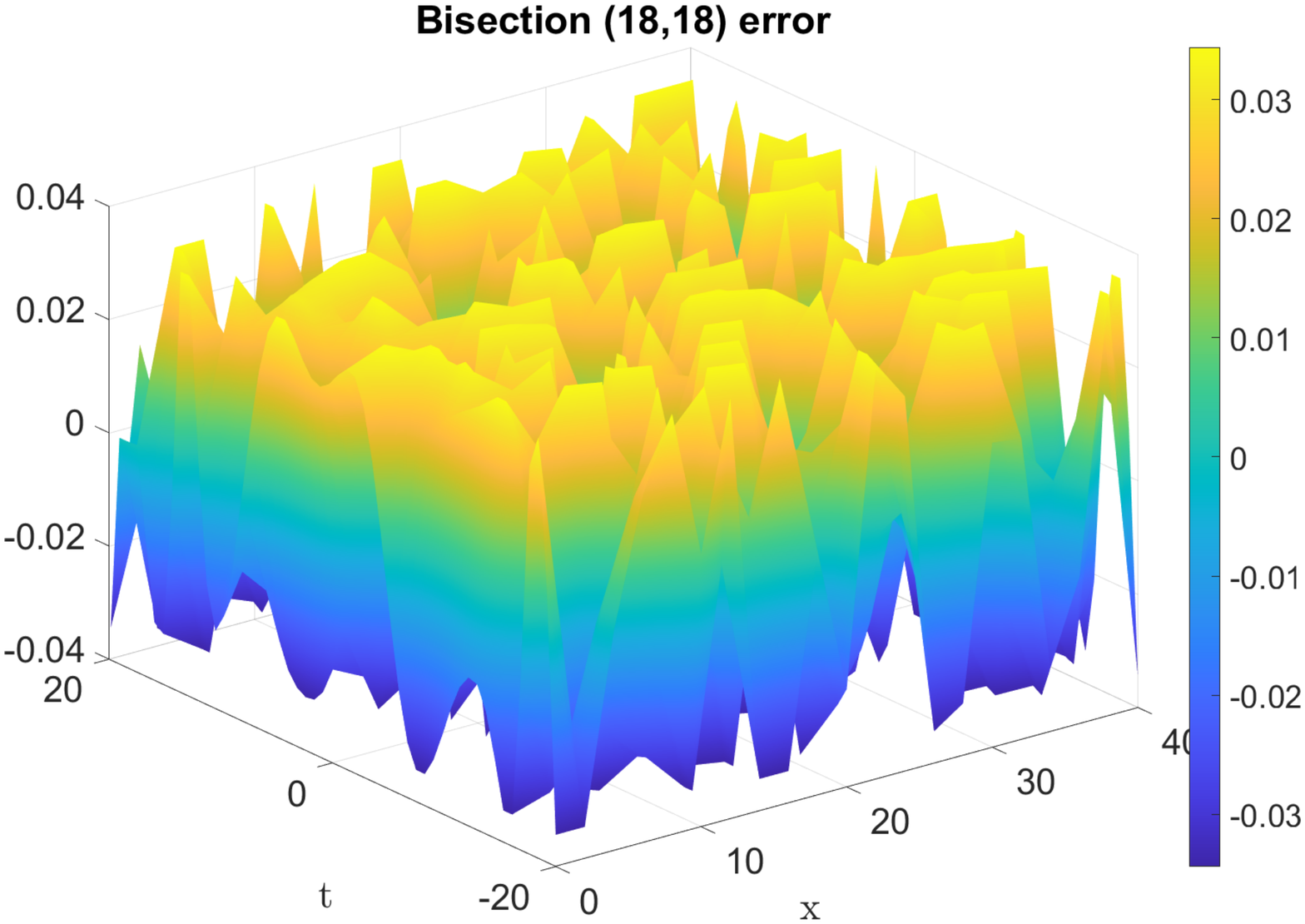}
		\includegraphics[width=60mm]{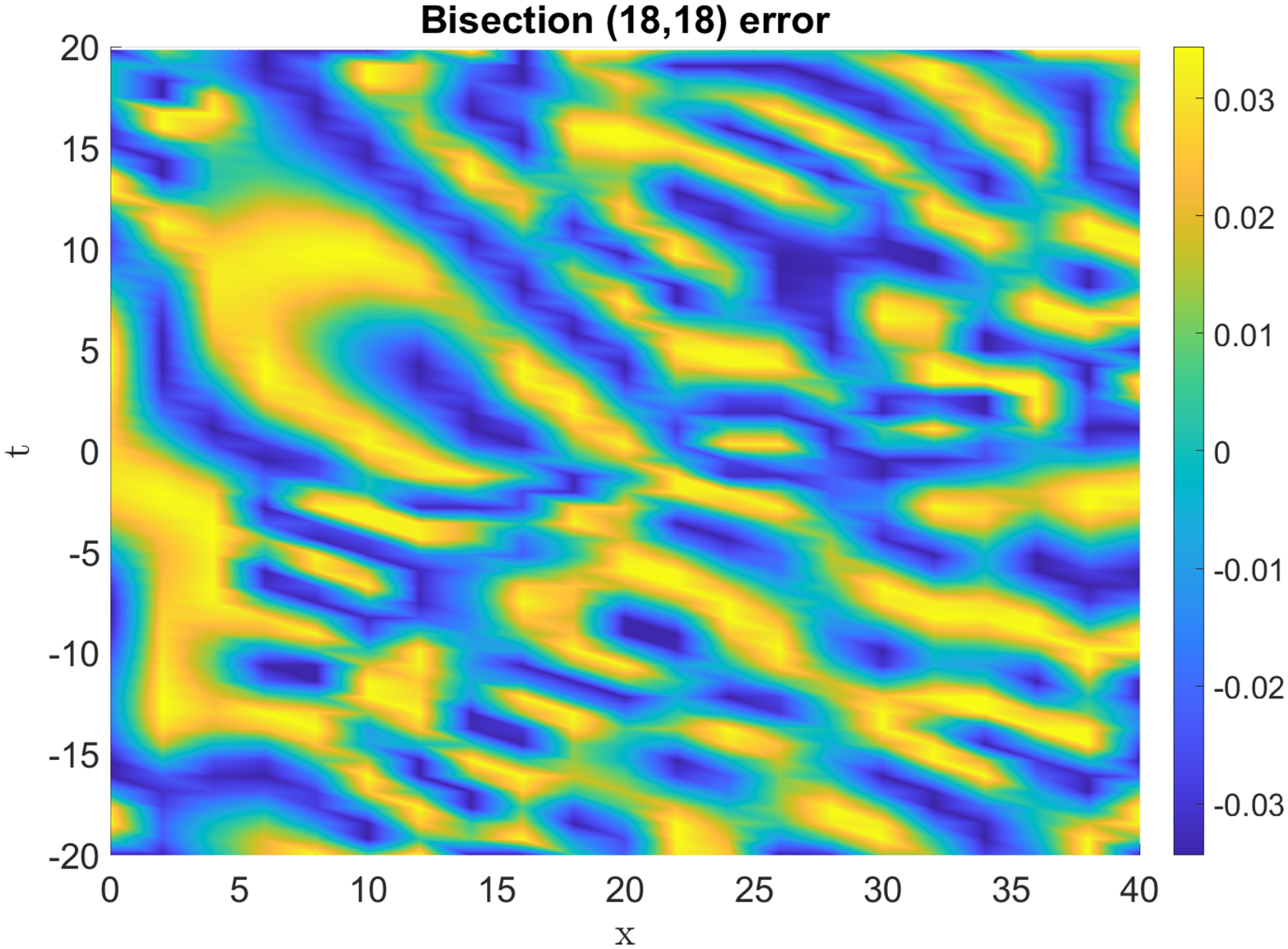}
	\caption{Approximations (above) and the error curves (below) are constructed on the domain which consists of pairs of every $10^{\text{th}}$ point of the original domain: 3D and 2D view.}
        \label{fig:Rational approximation of degree (18,18) multivariate domain of pairs of every 10th point}
\end{figure}

\begin{table}
	\centering
	\begin{tabular}{|c|c|c|}
		\cline{2-3}
		\multicolumn{1}{c|}{}   &  Uniform error & Time (sec.)\\
		\hline
		Every $20^{\text{th}}$ point of the domain & 0.014061551772703 & 66.460222 \\
		\hline
		Every $10^{\text{th}}$ point of the domain & 0.034389125232539 & 167.995076 \\
		\hline
	\end{tabular}
	\caption{Uniform error and computational time for the degree (18,18) approximation}
        \label{tab:Aproximation of degree (18,18)}
\end{table}

\subsection{Rational approximation of degree (20,20)}
In this section, we consider approximating the original functions in Figure~\ref{fig:original dataset with pairs of every 20th point of the domain} and Figure~\ref{fig:original dataset with pairs of every 10th point of the domain} by a rational function of degree (20,20). Even though the degree of the rational function is higher, the number of parameters (the coefficients of the polynomials in the numerator and the denominator) is smaller ($\approx 400$) compared to the number of parameters in the neural networks that the authors of~\cite{boulle2020rational} used in their experiments ($\approx 8035$). 

As the degree of the rational function increases, the condition number of the constraint matrices appearing in the auxiliary linear programming problem of the bisection method increases and in some cases, the MATLAB code tends to complain and crash. This issue can be solved by adding an extra linear constraint to the auxiliary problem of the bisection method which restricts the denominator polynomial from above. The bisection method has the flexibility of adding constraints and this additional adjustment allows us to compute higher degree rational approximations. Similar observation can be found in ~\cite{SharonPeirisSukhorukovaUgon} where authors use this property in matrix function evaluations for univariate approximation (lifting matrix functions). The experiments in this section where the degree of the approximation is $(20,20)$ were conducted with an additional constraint that restricts the denominator from above. The numerical experiments demonstrate that this extra improvement leads to more accurate approximations since the program can handle high degree polynomials for both the numerator and the denominator.

\begin{figure}
	\centering
		\includegraphics[width=60mm]{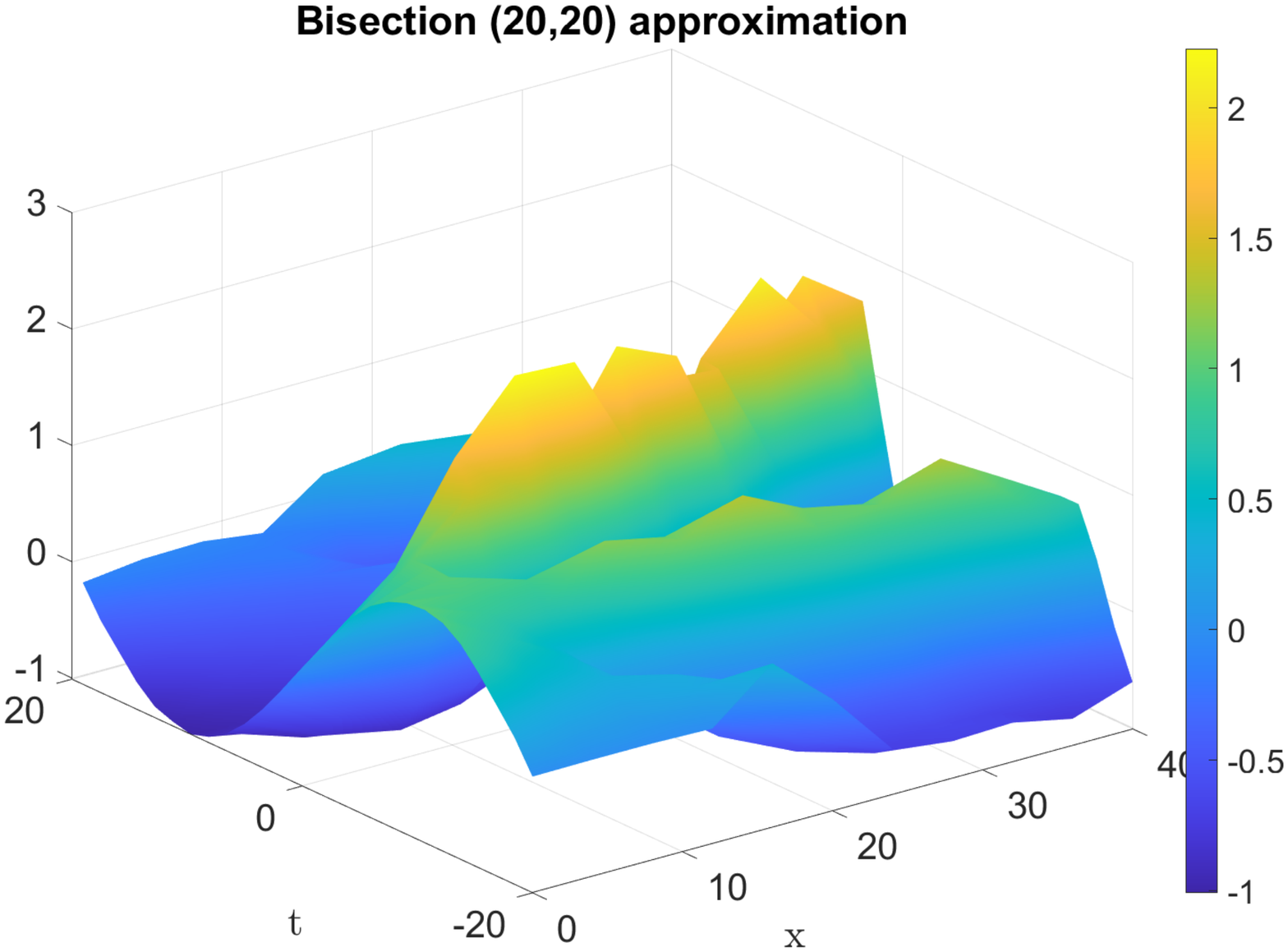}
		\includegraphics[width=60mm]{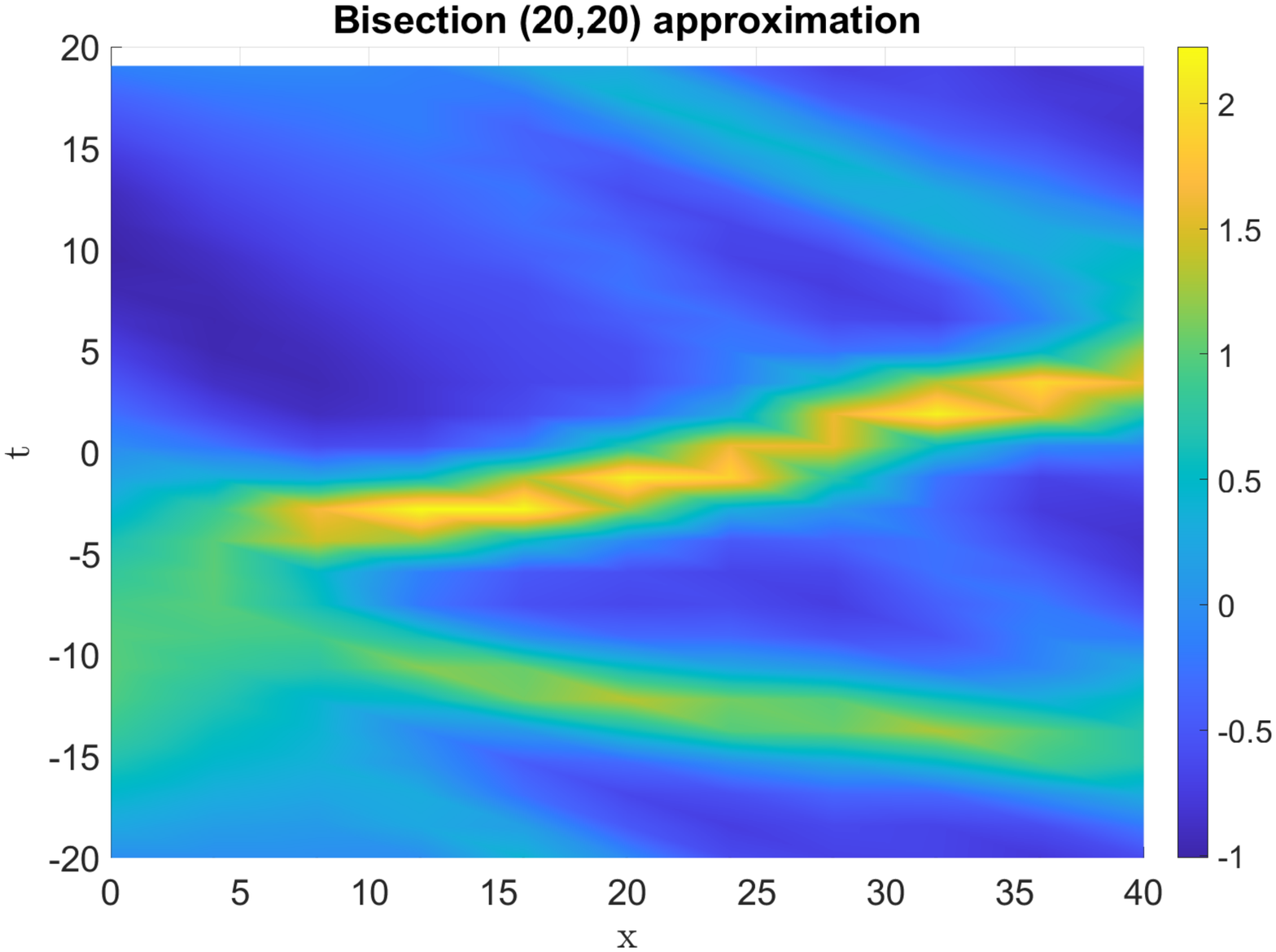}
		\includegraphics[width=60mm]{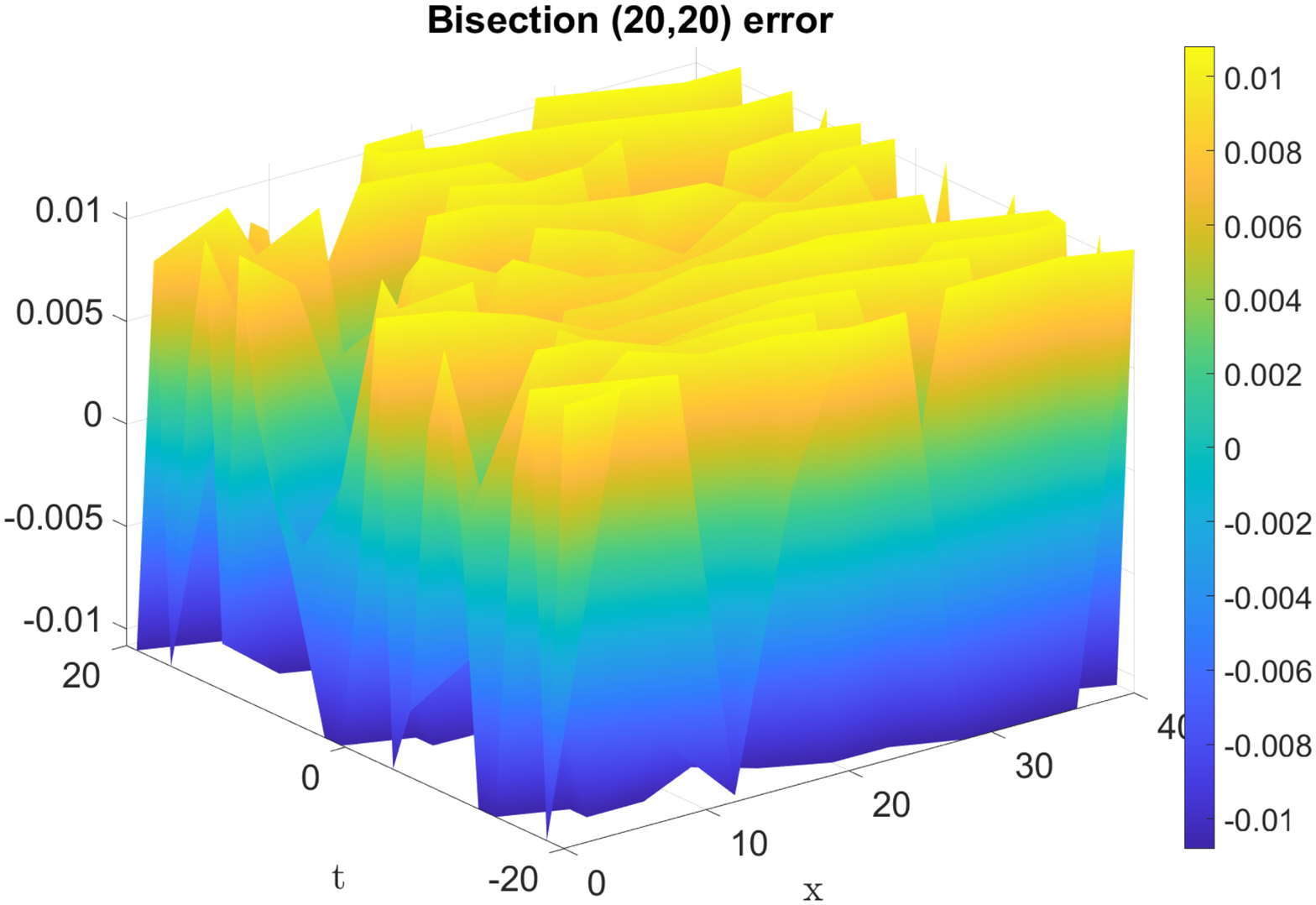}
		\includegraphics[width=60mm]{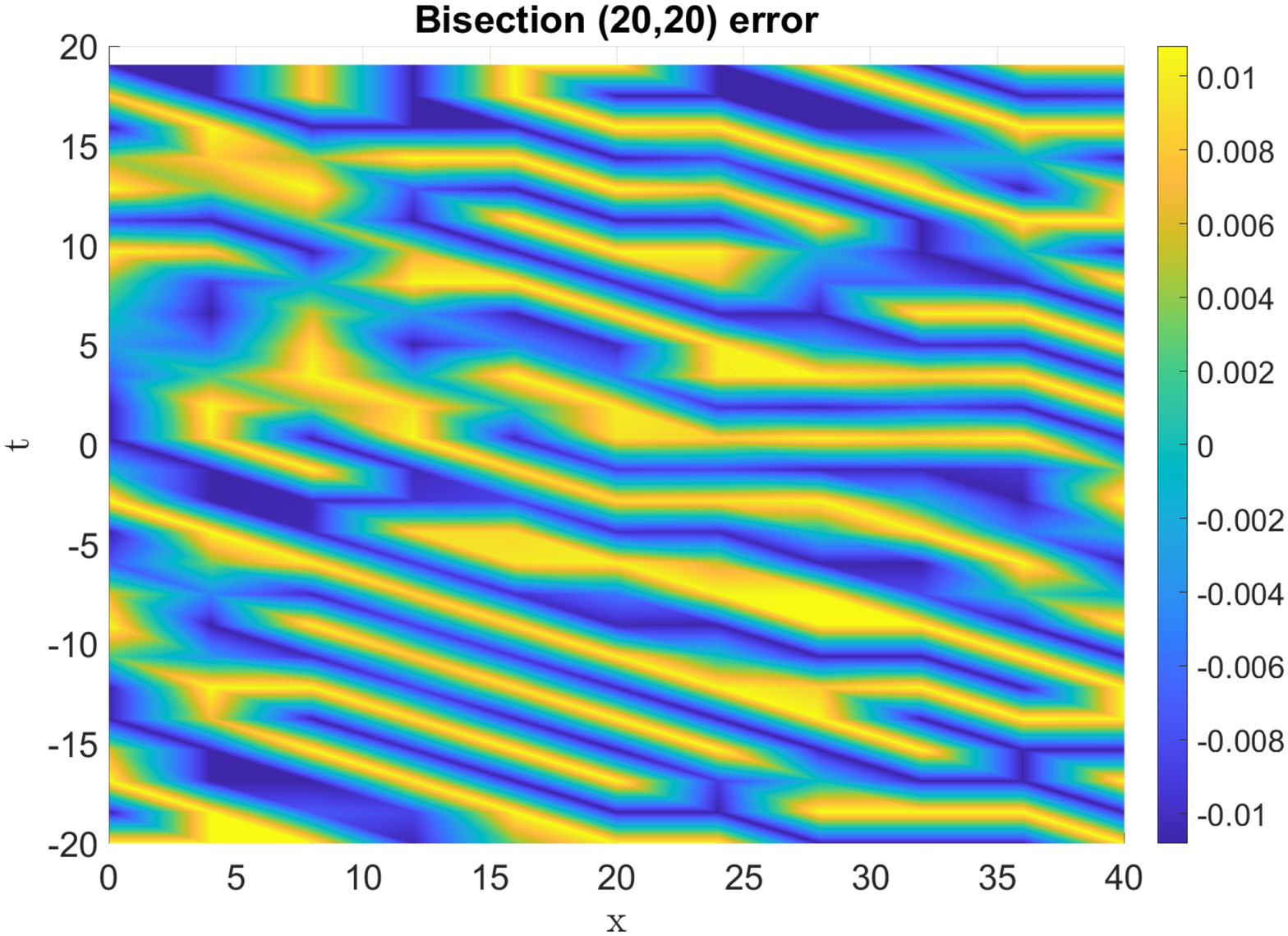}
	\caption{Approximations (above) and the error curves (below) are constructed on the domain which consists of pairs of every $20^{\text{th}}$ point of the original domain: 2D and 3D view.}
        \label{fig:Rational approximation of degree (20,20) multivariate domain of pairs of every 20th point}
\end{figure}

\begin{figure}
	\centering
		\includegraphics[width=60mm]{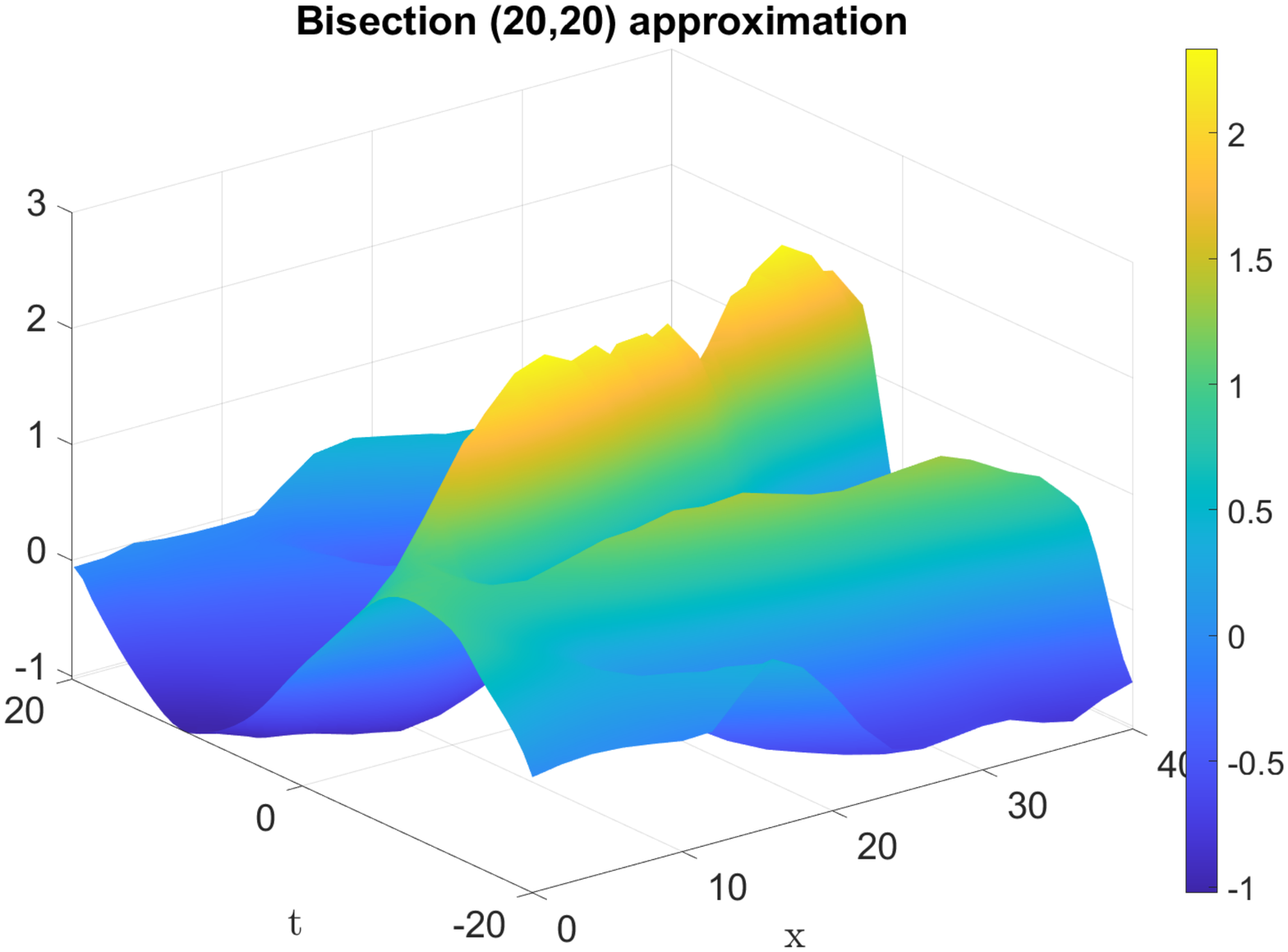}
		\includegraphics[width=60mm]{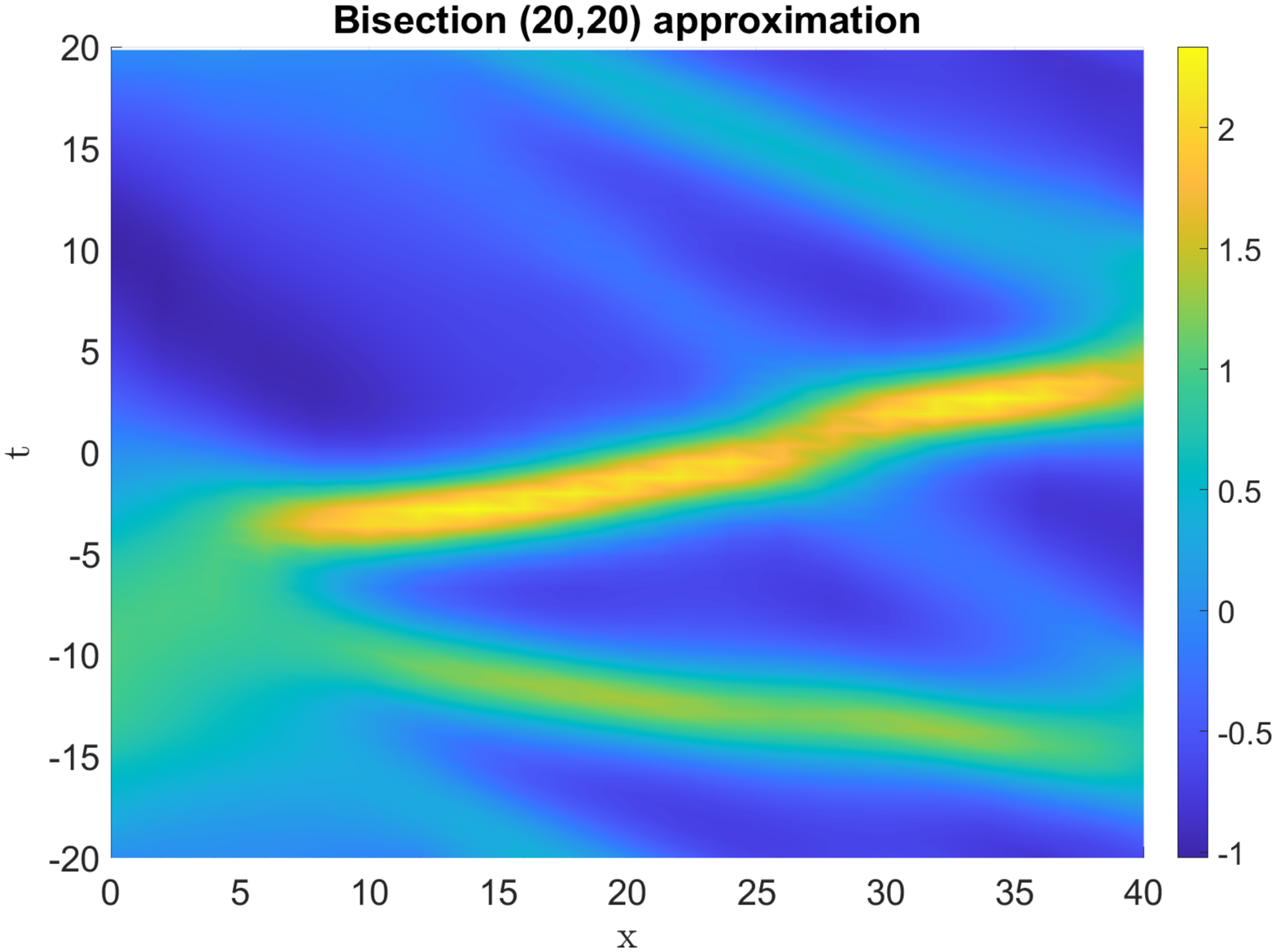}
		\includegraphics[width=60mm]{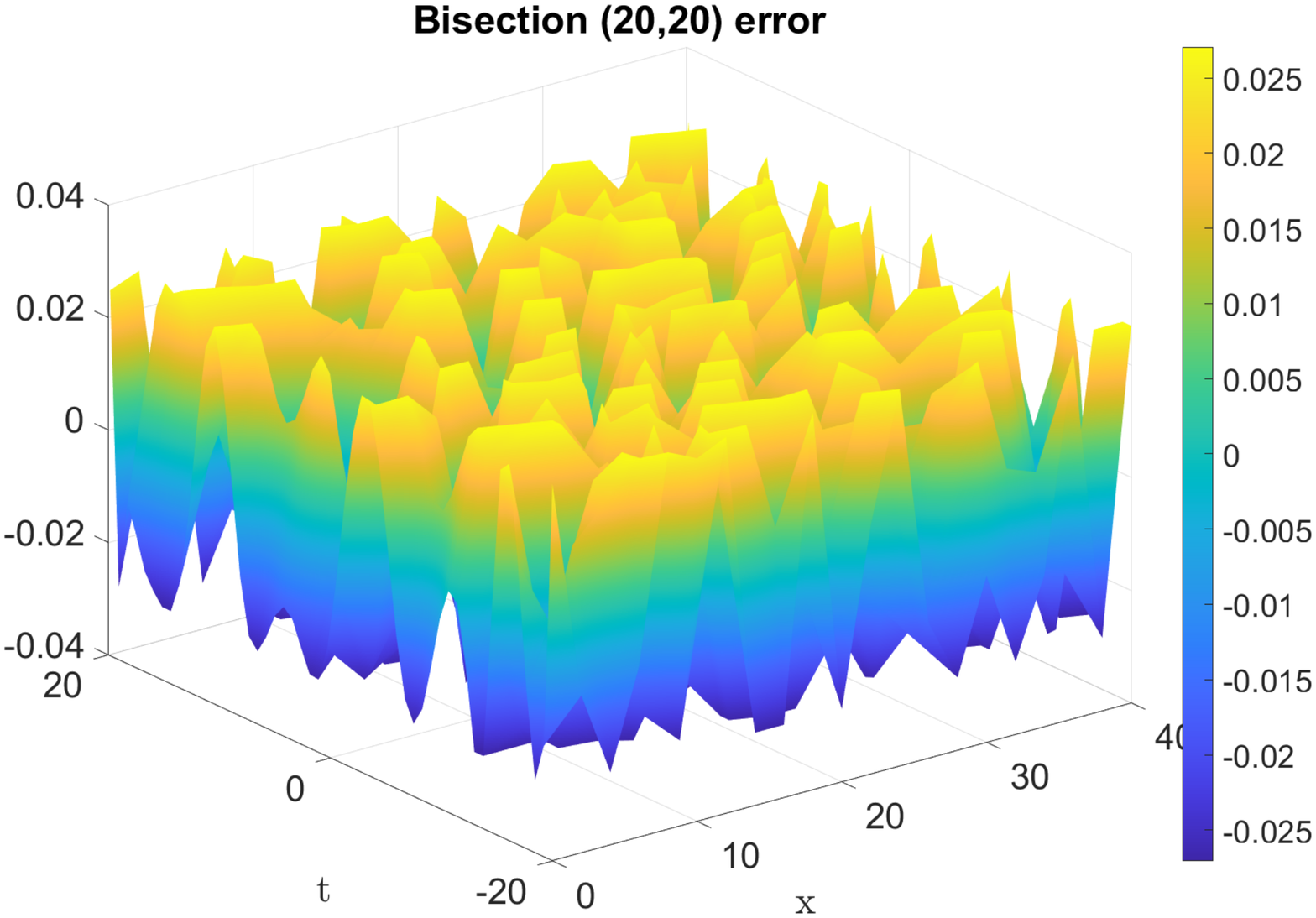}
		\includegraphics[width=60mm]{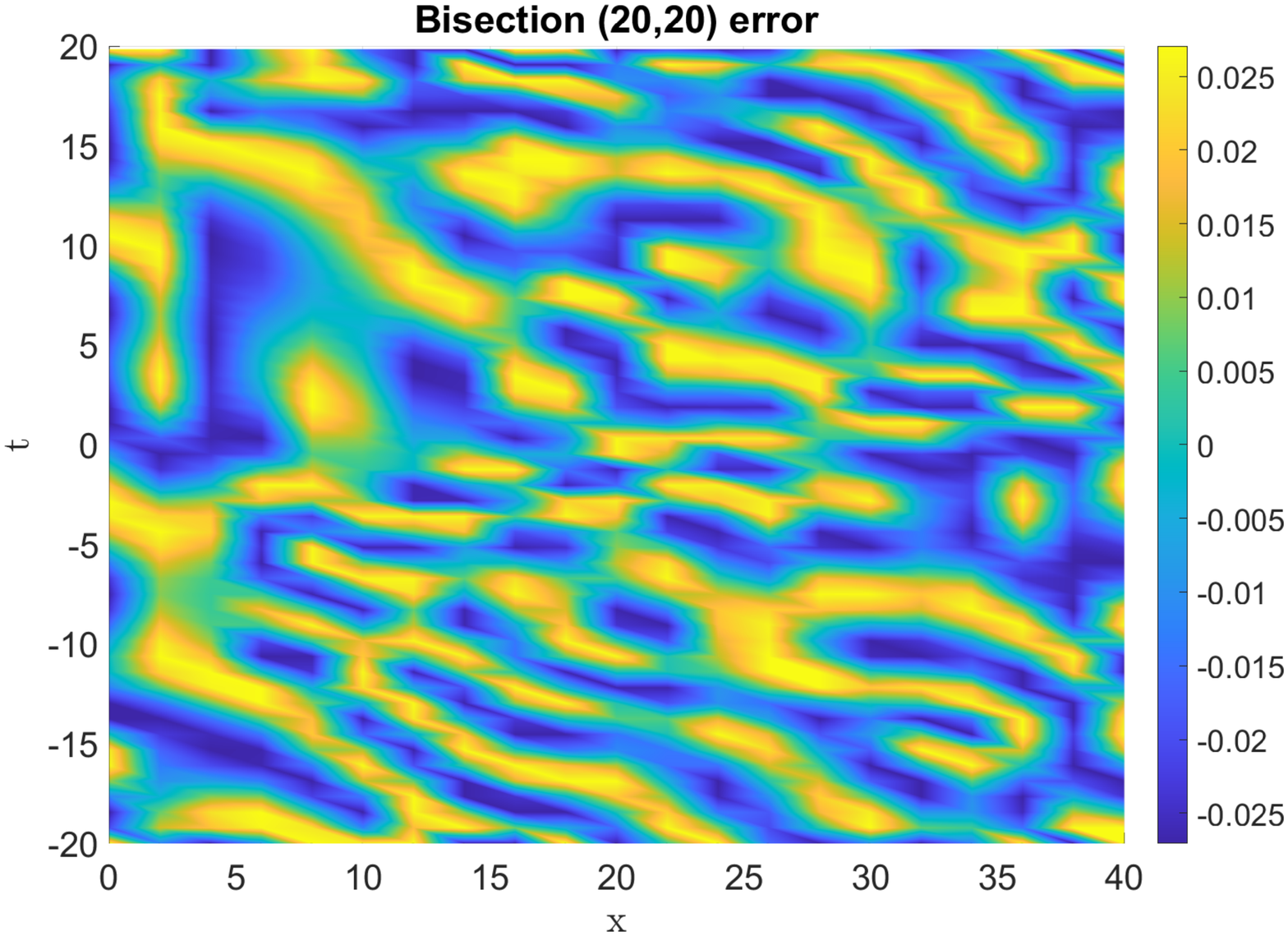}
	\caption{Approximations (above) and the error curves (below) are constructed on the domain which consists of pairs of every $10^{\text{th}}$ point of the original domain: 3D and 2D view.}
        \label{fig:Rational approximation of degree (20,20) multivariate domain of pairs of every 10th point}
\end{figure}

\begin{table}
	\centering
	\begin{tabular}{|c|c|c|}
		\cline{2-3}
		\multicolumn{1}{c|}{}   &  Uniform error & Time (sec.)\\
		\hline
		Every $20^{\text{th}}$ point of the domain & 0.010816799001175 & 158.933609 \\
		\hline
		Every $10^{\text{th}}$ point of the domain & 0.027058335297138 & 438.548139 \\
		\hline
	\end{tabular}
	\caption{Uniform error and computational time for the degree (20,20) approximation}
        \label{tab:Aproximation of degree (20,20)}
\end{table}

Clearly, the approximations in Figure~\ref{fig:Rational approximation of degree (20,20) multivariate domain of pairs of every 20th point} and Figure~\ref{fig:Rational approximation of degree (20,20) multivariate domain of pairs of every 10th point} coincide with the original functions and the uniform error is much smaller compared to the lower degree approximations.

Overall conclusion: Our results cannot be compared with the results in~\cite{boulle2020rational} since the errors are computed in different norms. However, we demonstrate that the rational approximation in a multivariate domain is as powerful as neural networks and it is capable of approximating complicated functions by using much fewer parameters than a neural network.

For each degree, we also computed the uniform error on the whole domain once the approximation is computed on a selected reduced domain. The results can be found in Table~\ref{tab:Aproximation errors computed on the whole domain}. Note that in the case of degree (20,20), the denominator is bounded from above.

\begin{table}
	\centering
	\begin{tabular}{|c|c|c|}
		\cline{2-3}
		\multicolumn{1}{c|}{} & Degree  &  Uniform error\\ 
		\hline
		   & (2,2) & 1.49062088756304\\ 
		\cline{2-3}
            Pairs of every & (5,5) & 1.57543434735584\\ 
            \cline{2-3}
		  $20^{\text{th}}$ point  & (10,10) & 245.34679368905\\ 
		\cline{2-3}
            of the domain & (18,18) & 11404.1284016\\ 
            \cline{2-3}
             & (20,20) & 88023.9109319522\\  
             \hline
		   & (2,2) & 1.40803015961408\\ 
		\cline{2-3}
            Pairs of every & (5,5) & 0.702071289628257\\ 
            \cline{2-3}
	    $10^{\text{th}}$ point  & (10,10) & 0.334154891477345\\ 
		\cline{2-3}
            of the domain & (18,18) & 15.3976686595557\\ 
            \cline{2-3}
             & (20,20) & 9.191409640164363 \\ 
		\hline
	\end{tabular}
	\caption{Uniform error computed on the whole domain}
        \label{tab:Aproximation errors computed on the whole domain}
\end{table}

One can clearly see that the uniform error computed on the whole domain is very similar to the uniform error computed on the domain with fewer selected points (reduced domain) when the degree of the approximation is small. This observation is very prominent in the case where the domain consists of pairs of every $10^{\text{th}}$ point of the domain. For higher degrees, the difference between the errors computed on the whole domain and the reduced domain gets much larger when the domain has much fewer points (when the domain has pairs of every $10^{\text{th}}$ point of the whole domain).

\section{Conclusions and future research directions}\label{sec:conclusions}

In this paper we compared several approaches to approximating functions, each being a slight modification of the other one. Starting from a classical neural network with a ReLU activation function, we gradually changed our set up, first by replacing the activation function with a rational function approximating ReLU, then including the coefficients of the rational activation function in the learning parameters, and finally replacing the neural network with a best rational approximation of the same number of parameters using traditional optimisation approaches.

Our finding is that as each step increases the flexibility and size of the feasible space, we observe an improvement in the approximation. In particular the optimisation-based approximation achieves better approximation than the neural network parts.

The training time is also in favour of the optimisation-based approximation, especially because it is often possible to achieve a good approximation with fewer parameters. As the dimension of the function to approximate increases, it is not clear whether this advantage remains.

Our study did not consider other aspects of the approximation, such as generalisability. It is well known in general that deep learning has good generalisability properties and avoid overfitting. However, it is not clear whether adding extra flexibility (for example including the coefficients of the activation function in the learning parameters) retain those generalisability properties. Further work needs to be done to gain better insights about this.

\section*{Acknowledgements, Statements and Declarations}

We are grateful to the Australian Research Council for supporting this work via Discovery Project DP180100602.



\appendix
\section{Appendix}  \label{sec:Appendix}

\subsection{Results: Neural network-based approximation.} \label{appendix:resultsNN}

\subsubsection{Neural Network with ReLU activation}\label{appendix:results1:Set1-4}

\paragraph{{\bf Set 1}: 2 nodes in the hidden layer, MSE loss, ADAM optimiser, ReLU activation}

Table~\ref{tab:Results: experiments set 1 of NN with ReLU activation} summarises the results for the case of 2~nodes in the hidden layer (NN, MSE loss, ADAM). One can see that the value of the loss function does not improve significantly when the number of epochs is increasing. Moreover,  it appeared that the minimum reported value of the loss function may take place before the final epoch. For example, in the case of 200~epochs, the minimal value (0.045839) was obtained at epoch~114, while the reported value (after 200~epochs) is~0.045934. Overall, the computational time (per epoch) does not increase significantly when the number of epochs increases from 50 to~200.  The last column of this table reports the running time per epoch and therefore the computational time is almost proportional to the number of epochs. Therefore, the large number of epochs increases the overall computational time proportionally without achieving significantly better accuracy.
\begin{table}  
            \centering
    \begin{tabular}{|c|c|c|c|}
    \hline
    Epoch & Final loss & Minimum loss & Run time (per epoch)\\
    \hline
    50   & 0.045868 &          & 2.51s $\pm$ 103ms\\
    44   &          & 0.045840 & \\
    \hline
    100  & 0.045864 &          & 2.7s $\pm$ 71.9ms\\
    63   &          & 0.045840 & \\
    \hline
    200  & 0.045934 &          &  2.77s $\pm$ 453ms\\
    114  &          & 0.045839 & \\
    \hline
    \end{tabular}
         \caption{Results: experiments set 1}
         \label{tab:Results: experiments set 1 of NN with ReLU activation}
\end{table}


Figure~\ref{fig:Relu activation - set 1, epoch is 50,100,200} depict the function~$f(x)$ and its approximations obtained in this experiment: 50, 100 and 200~epochs (Set~1). The pictures confirm that the approximation accuracy does not improve much when we increase the number of epochs, which is consistent with the results from table~\ref{tab:Results: experiments set 1 of NN with ReLU activation}, where the loss function values are similar for different numbers of epochs.
\begin{figure}
    \centering
    \includegraphics[width=40mm]{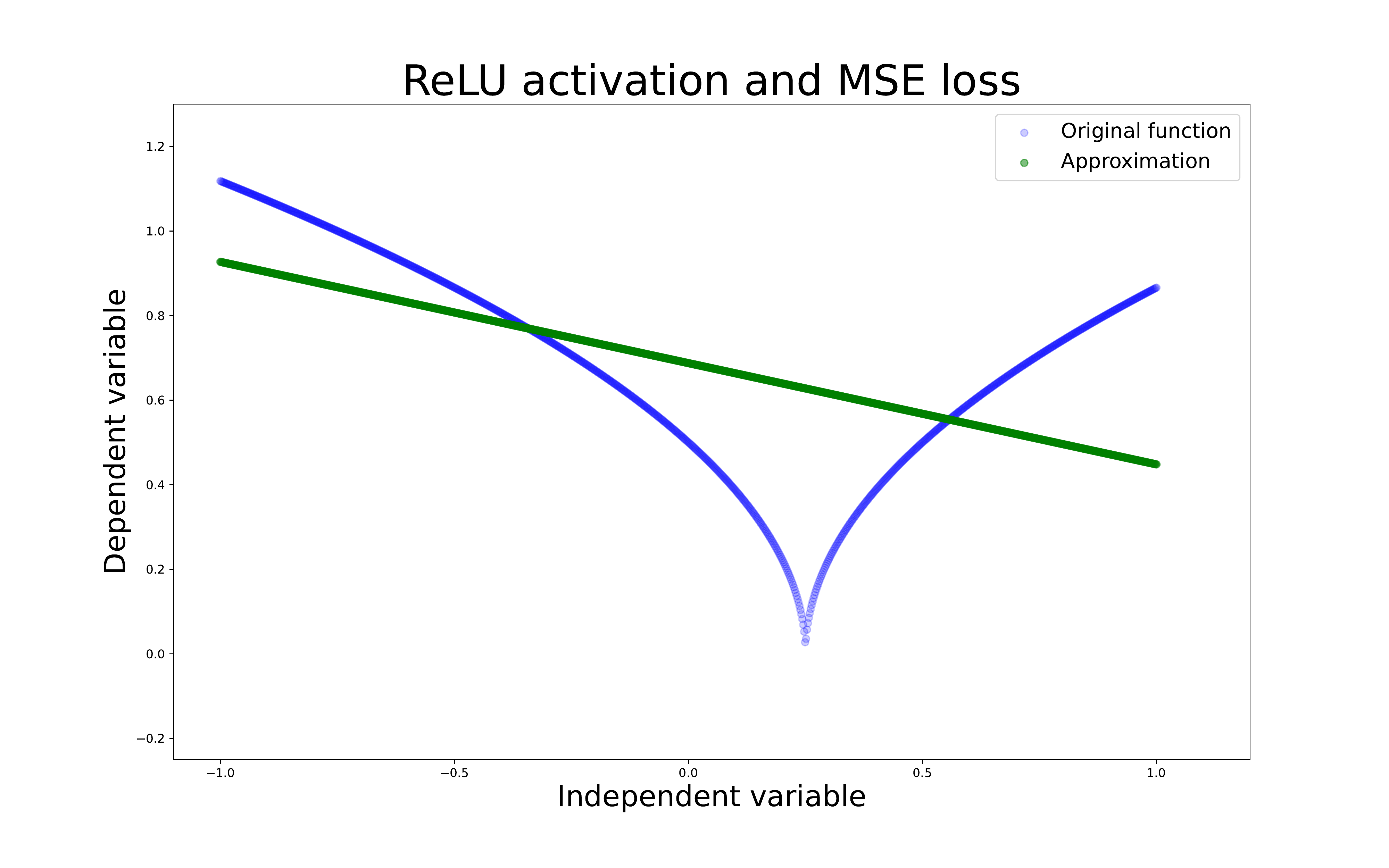}
    \includegraphics[width=40mm]{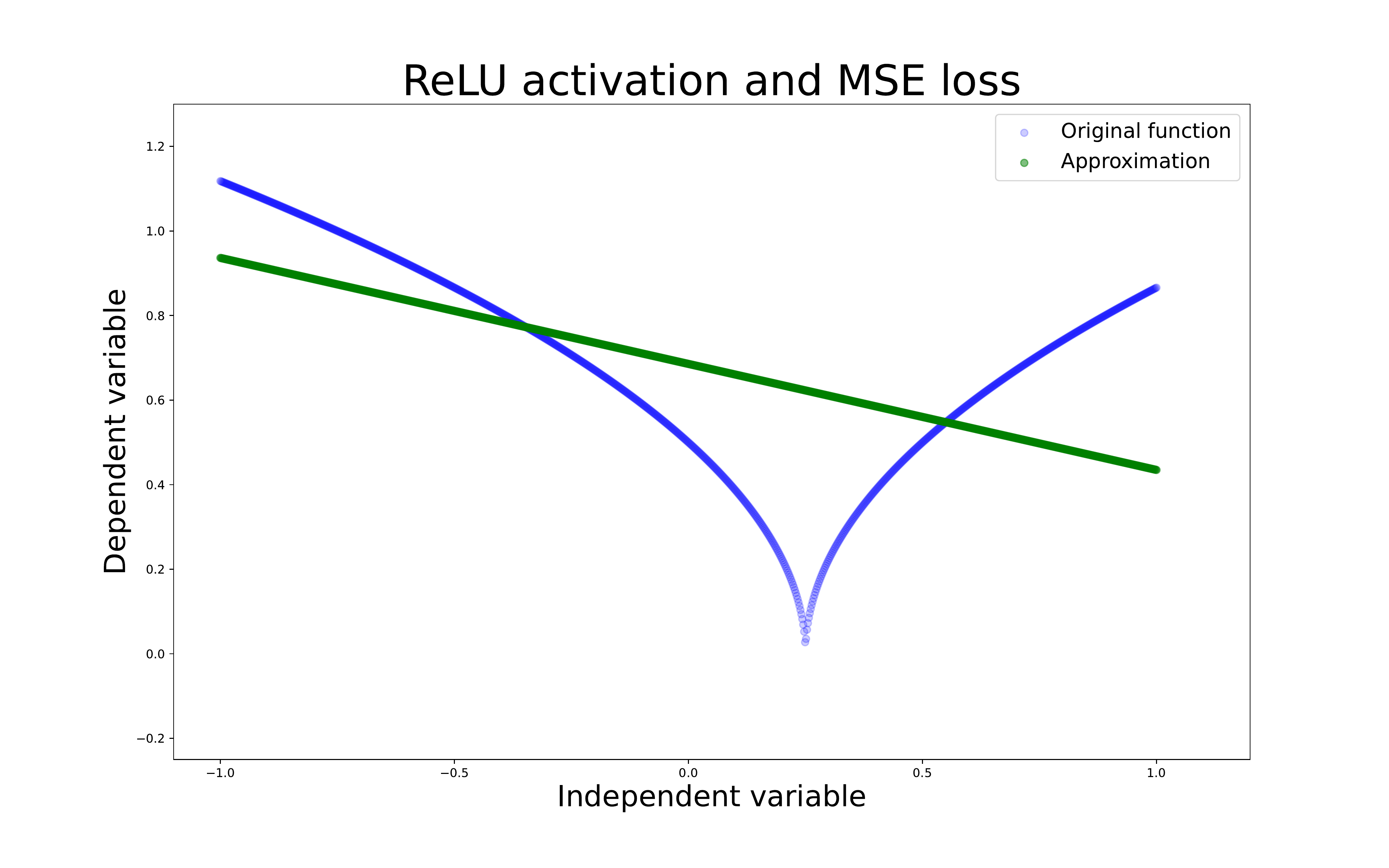}
       \includegraphics[width=40mm]{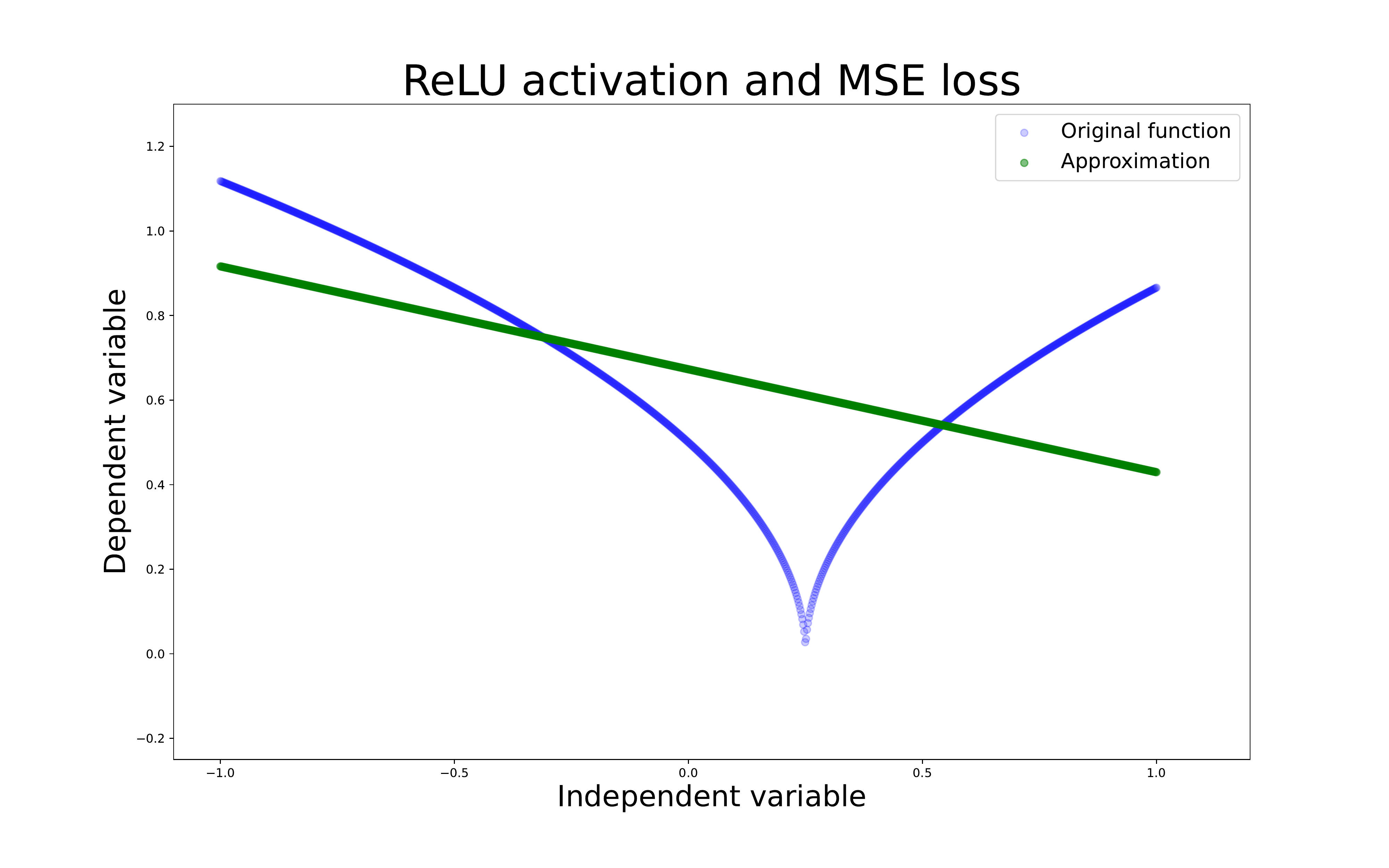}
    \caption{Set 1: Approximation is computed by a neural network with ReLU activation,  50, 100 and 200~epochs.}
    \label{fig:Relu activation - set 1, epoch is 50,100,200}
\end{figure}

\paragraph{{\bf Set 2}: 2 nodes in the hidden layer, Uniform loss, ADAMAX optimiser, ReLU activation}

In this set of experiments, we use the uniform loss function. This is the main difference with Set~1 experiments. We also use ADAMAX optimiser instead of ADAM, since it is specifically developed for the uniform loss function. 

It does not make sense to compare the loss function values in Set~1 and Set~2, since they correspond to very different inaccuracy measures. When we compare the results for different numbers of epochs in Set~2, we can see a similar pattern to Set~1 experiments:   the computational time is approximately proportional to the number of epochs, the value of the loss function does not improve significantly when the number of epochs is increasing and optimal loss is not decreasing monotonically with the epochs. Interestingly, the best recorded loss function value of~0.467289  is achieved after 83~epochs, which is not very different from~0.469656, achieved just after 13~epochs. Overall, the increase in the number of epochs leads to a higher computational time without a significant improvement in the optimal loss value.

\begin{table}
    \centering
    \begin{tabular}{|c|c|c|c|}
    \hline
    Epoch & Final loss & Minimum loss & Run time (per epoch) \\
    \hline
    50  & 0.492186 &          & 2.92s $\pm$ 707ms\\
    13  &          & 0.469656 & \\
    \hline
    100 & 0.505962 &          & 2.94s $\pm$ 535ms\\
    83  &          & 0.467289 & \\
    \hline
    200 & 0.503380 &          & 2.84s $\pm$ 314ms\\
    83  &          & 0.467289 & \\
    \hline
    \end{tabular}
    \caption{Results: experiments set 2}
    \label{tab:Results: experiments set 2 of NN with ReLU activation}
\end{table}



Figure~\ref{fig:Relu activation - set 2, epoch is 50,100,200} depict the function~$f(x)$ and its approximations obtained in this experiment: 50, 100 and 200~epochs (Set~2). The pictures confirm that the approximation accuracy does not improve much when we increase the number of epochs, which is consistent with the results from table~\ref{tab:Results: experiments set 2 of NN with ReLU activation}, since the loss function values are similar for different number of epochs.
\begin{figure}
    \centering
    \includegraphics[width=40mm]{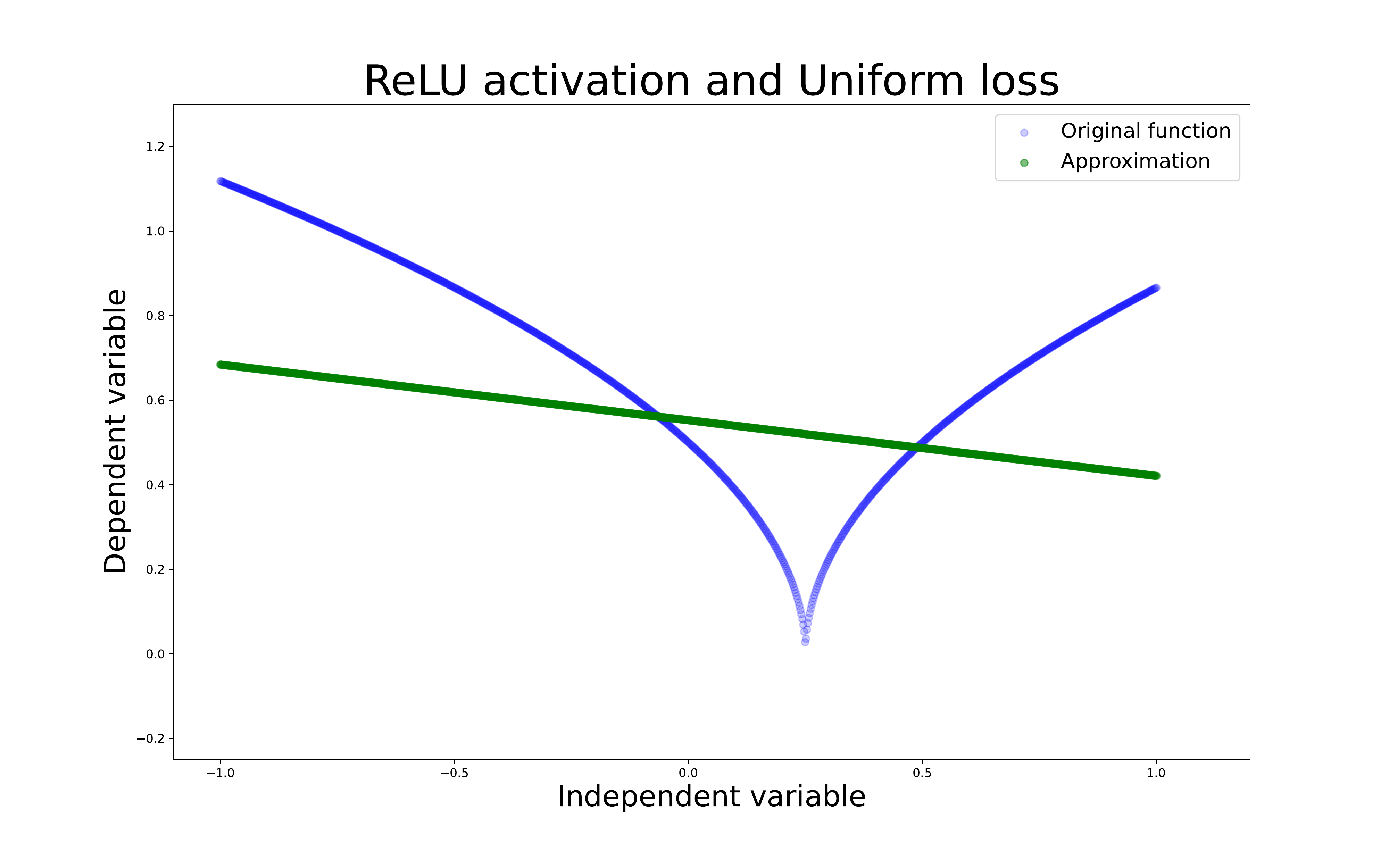}
     \includegraphics[width=40mm]{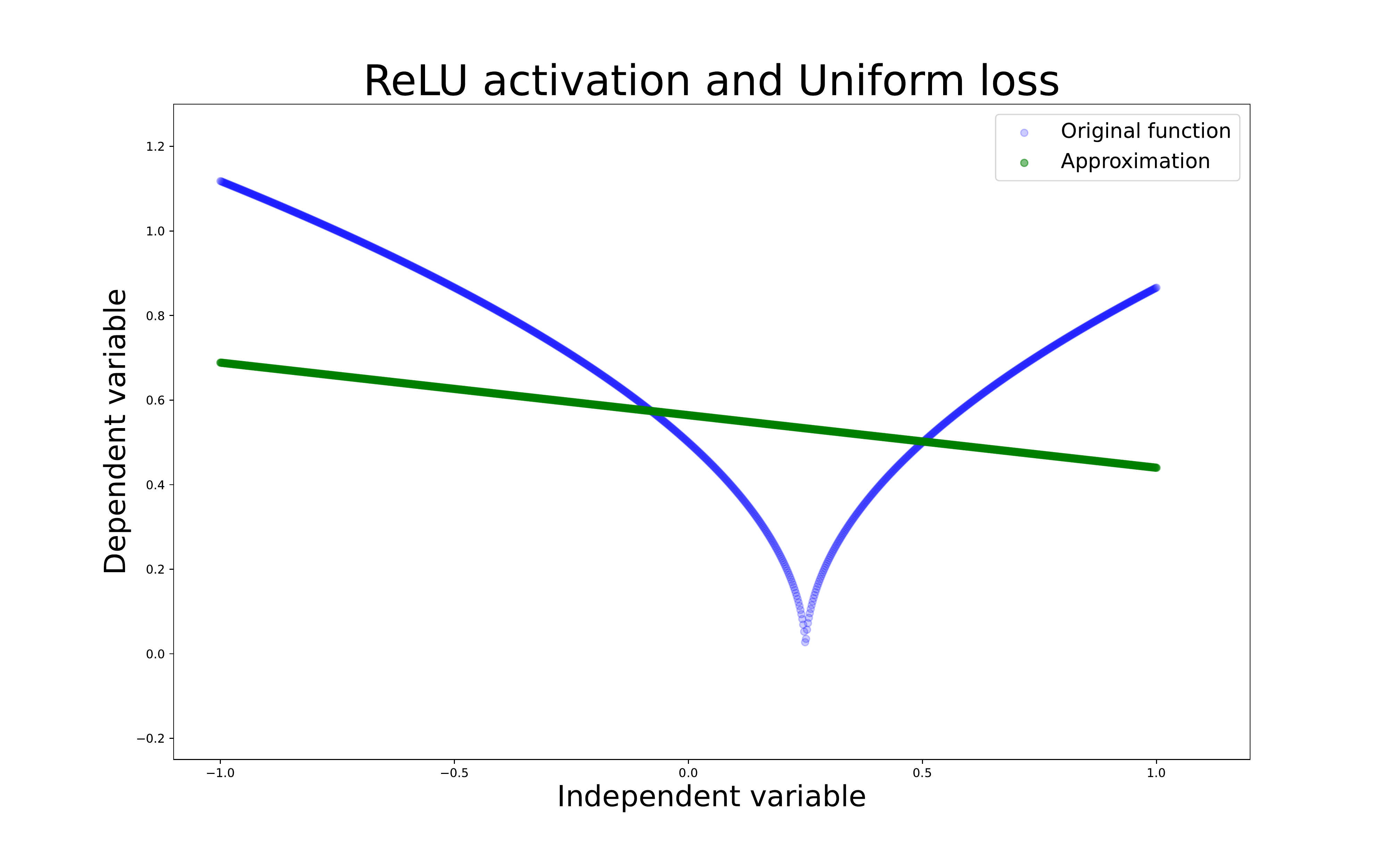}
      \includegraphics[width=40mm]{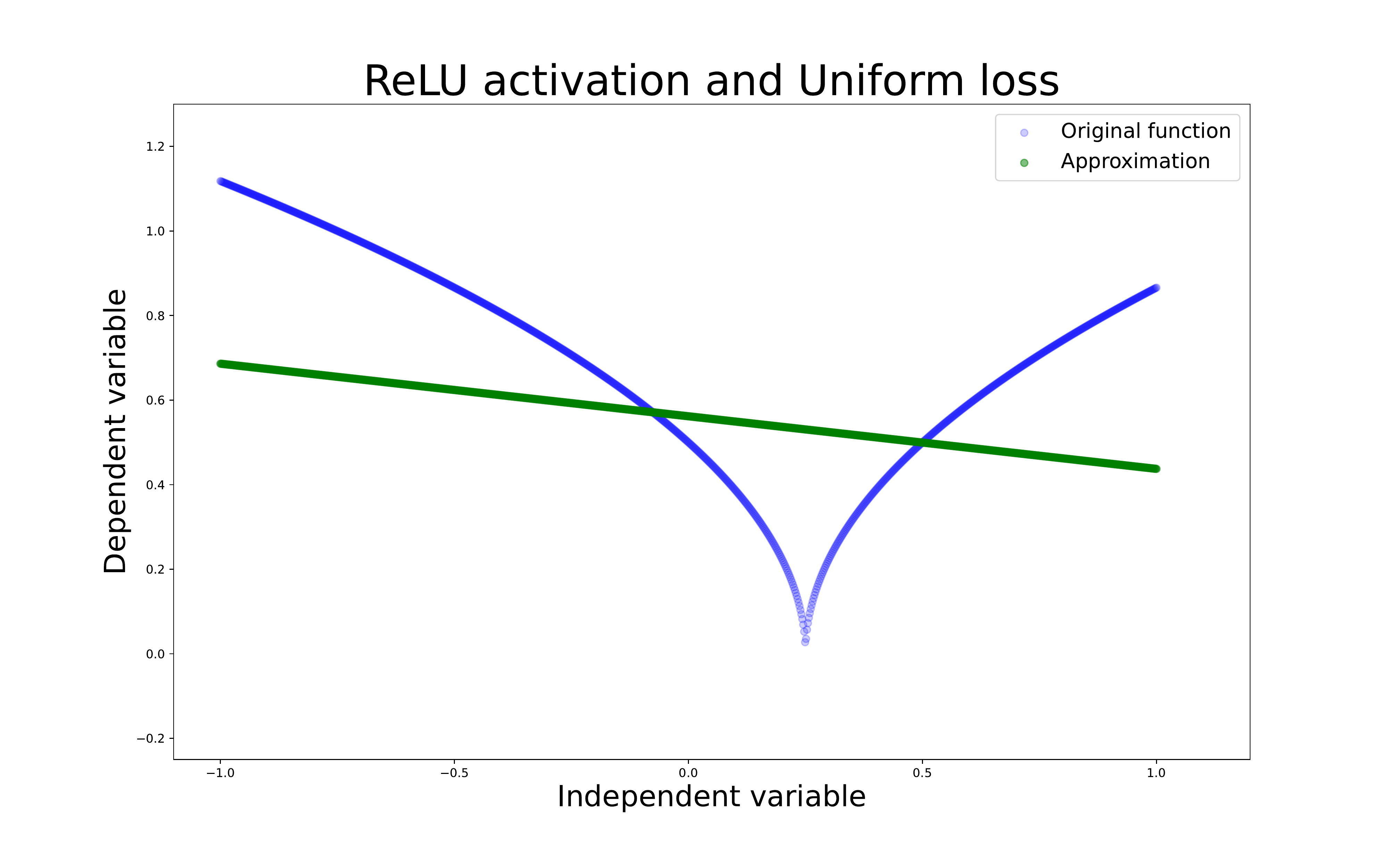}
    \caption{Set 2: Approximation is computed by a neural network with ReLU activation,  50, 100 and 200~epochs.}
    \label{fig:Relu activation - set 2, epoch is 50,100,200}
\end{figure}

Comparing the results of the experiments in Set~1 and Set~2, it is not clear which loss function (MSE or uniform) is more appropriate, but it is clear that nodes in the hidden layer is not enough for accurate results.

\paragraph{{\bf Set 3}: 10 nodes in the hidden layer, MSE loss, ADAM optimiser,  ReLU activation}

Our findings with 2~nodes in the hidden layer suggest that we need to increase the number of nodes. In this set of experiments, we have 10~nodes in the hidden layer.

Table~\ref{tab:Results: experiments set 3 of NN with ReLU activation} summarises the results of Set~3 experiments. The loss function values is significantly better than it was in the case of Set~1 experiments with only 2~nodes in the hidden layer. One can still observe some non-monotonicity of the optimal loss function value, but it is not as prominent as it was in the case of Set~1. Based on  the results, it is reasonable to stay within 100-150~epochs.

\begin{table}
    \centering
    \begin{tabular}{|c|c|c|c|}
    \hline
    Epoch & Final loss & Minimum loss & Run time (per epoch) \\
    \hline
    50  & 0.001848 &          & 2.89s $\pm$ 136ms\\
    50  &          & 0.001848 & \\
    \hline
    100 & 0.000915 &          & 2.5s $\pm$ 59.2ms\\
    99  &          & 0.000908 & \\
    \hline
    200 & 0.000914 &          & 2.53s $\pm$ 264ms\\
    141 &          & 0.000903 & \\
    \hline
    \end{tabular}
    \caption{Results: experiments set 3}
    \label{tab:Results: experiments set 3 of NN with ReLU activation}
\end{table}



Figure~\ref{fig:Relu activation - set 3, epoch is 50,100,200} depicts the results for Set~3 experiments. The picture clearly suggests that the increase in the number of nodes in the hidden layer improved the approximation accuracy and this accuracy is much better for 100 or 200~epochs. The inaccuracy near the ``difficult point'' $x=0.25$, where $f(x)$ is nonsmooth and non-Lipschitz is still significant, even for a larger number of epochs. It is also interesting to note that the computational time for Set~1 and Set~3 is very similar, despite an increase in the number of nodes in the hidden layer.

\begin{figure}
    \centering
    \includegraphics[width=40mm]{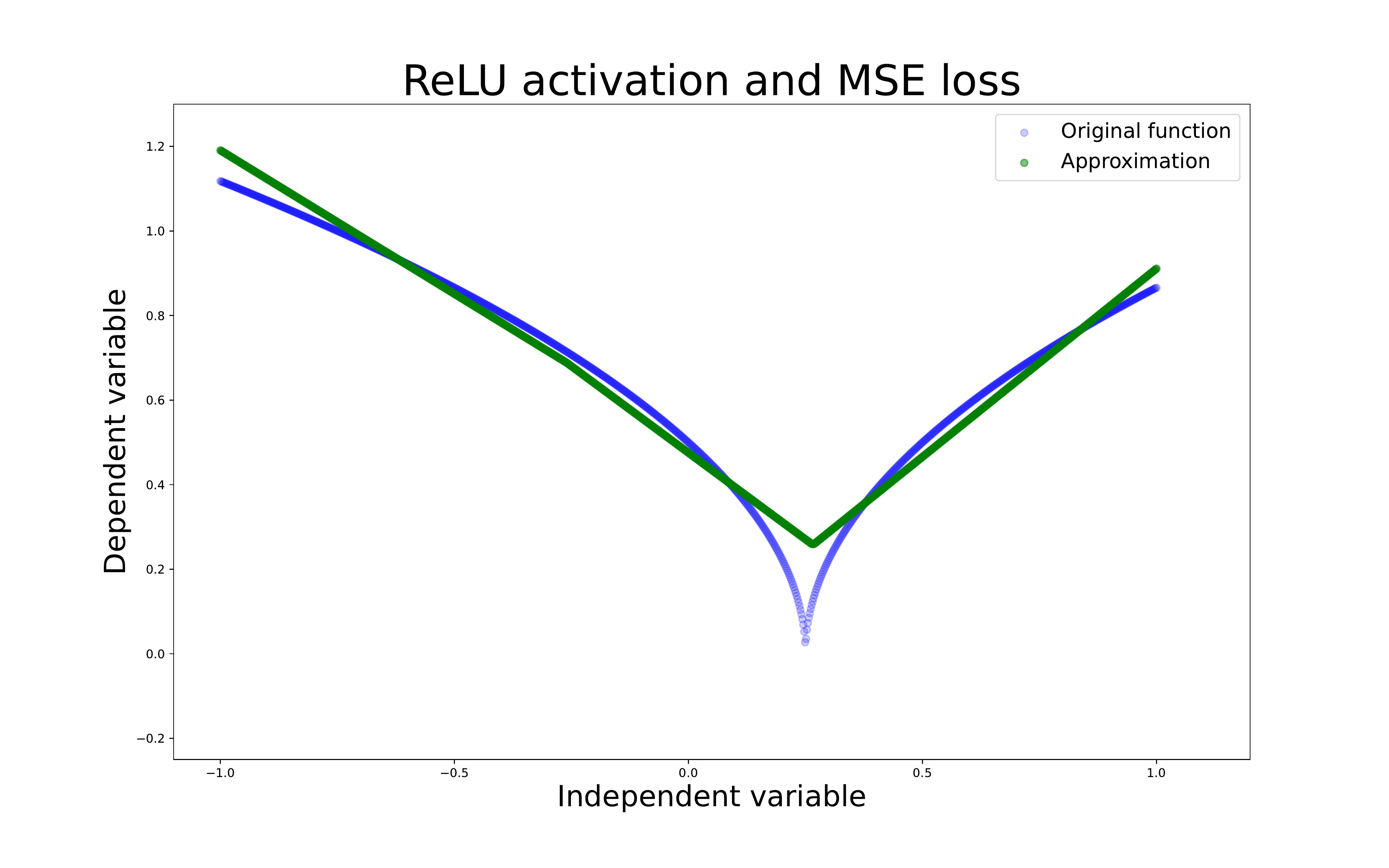}
    \includegraphics[width=40mm]{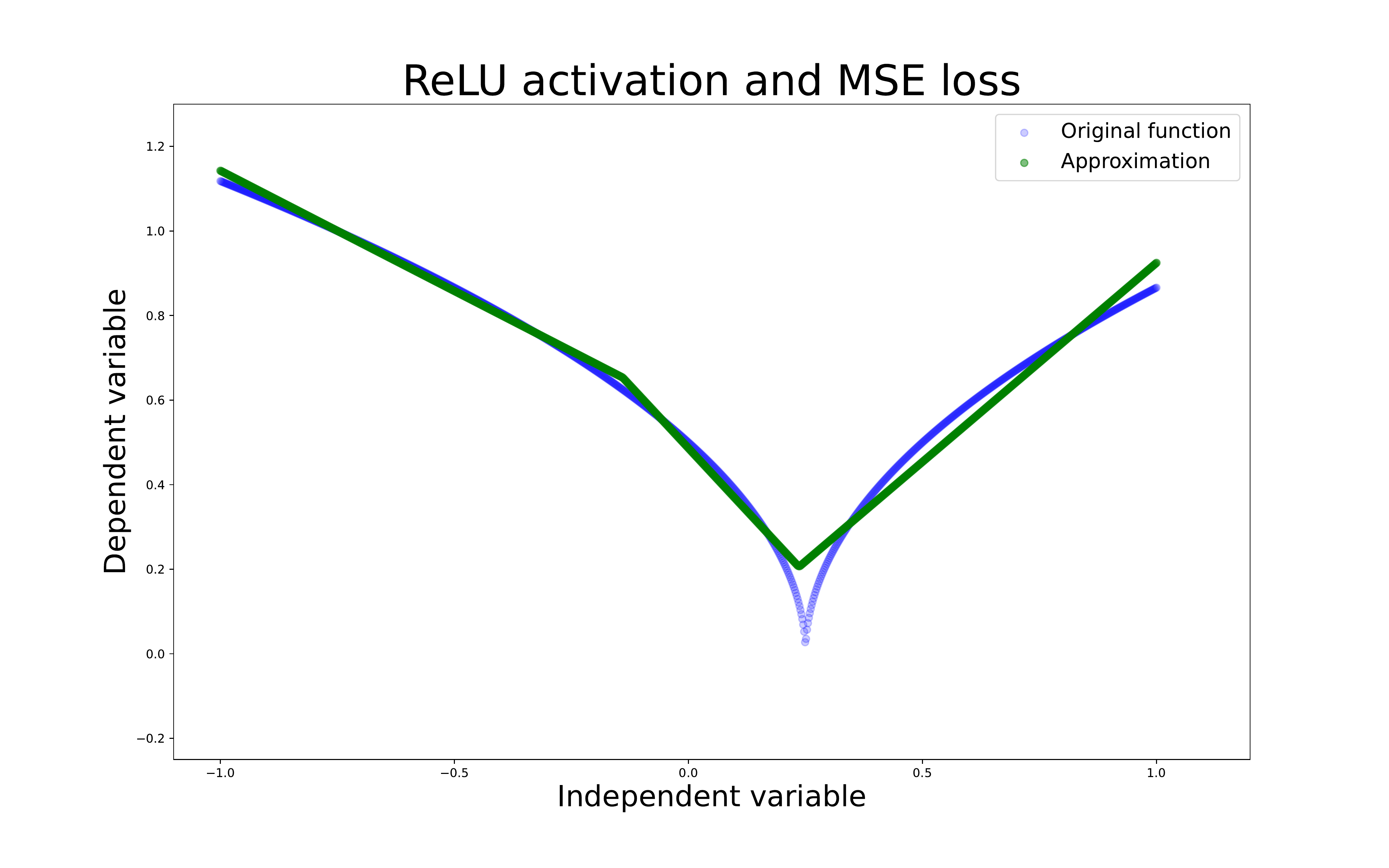}
    \includegraphics[width=40mm]{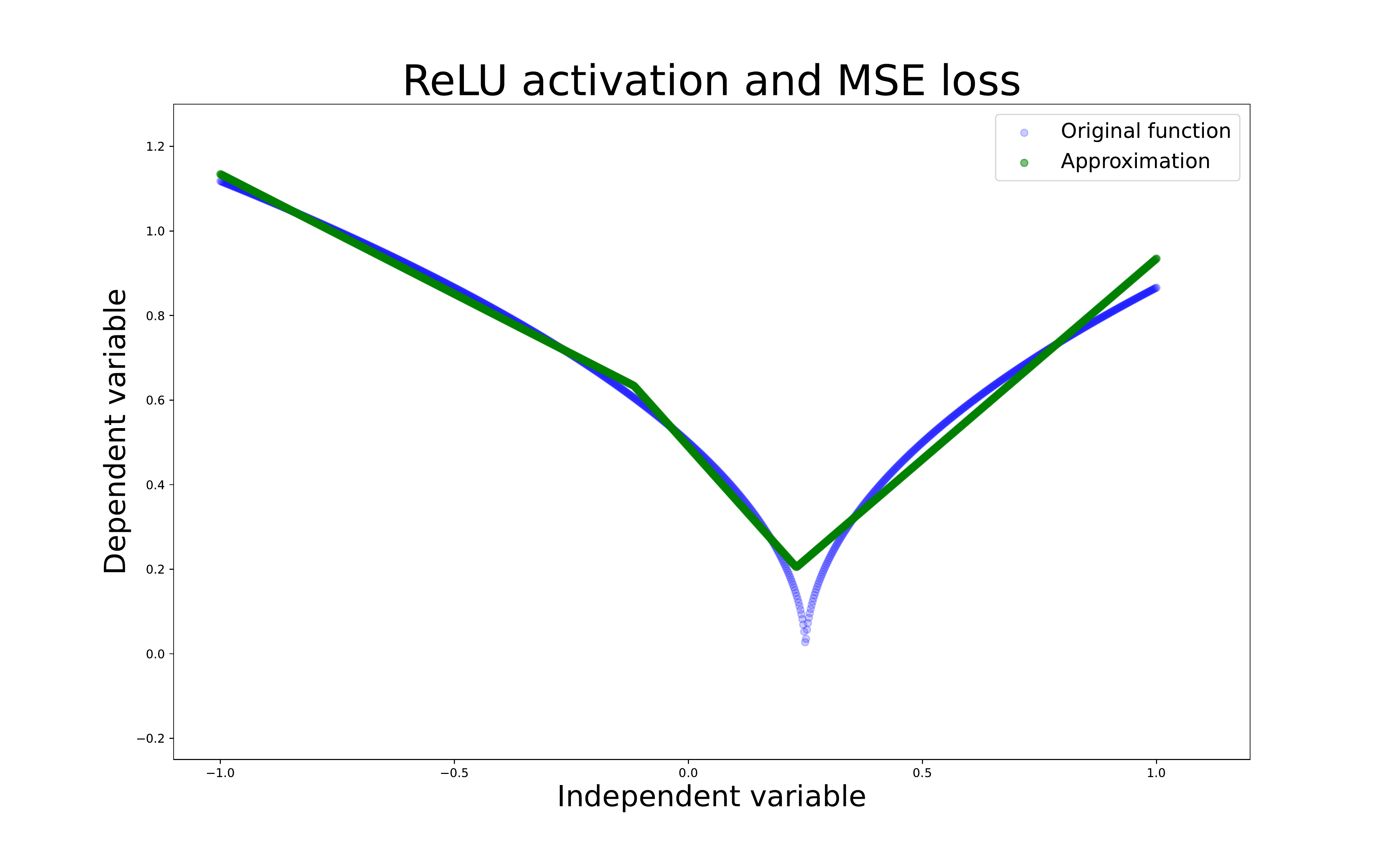}
    \caption{Set 3: Approximation is computed by a neural network with ReLU activation,  50, 100 and 200~epochs.}
    \label{fig:Relu activation - set 3, epoch is 50,100,200}
\end{figure}

\subsubsection{Neural Network with rational approximation to ReLU activation}\label{appendix:results2:Set1-4}

In this section, we present the results of numerical experiments, where the experiment settings are similar to those in section~\ref{appendix:results1:Set1-4}, but the activation function is the rational approximation to ReLU. This approximation was found by the differential correction method. 


Our activation function is a rational function of degree $(3,2)$: the rational approximation is the ratio of two polynomials, the degree of the numerator is~3 and the degree of the denominator is~2.  

The coefficients of the rational activation function is coming from the best rational $(3,2)$ approximation to the ReLU function. These coefficients are fixed throughout the whole training procedure. We do not learn the coefficients with the rest of the parameters of the network.

\paragraph{{\bf Set 1}: 2 nodes in the hidden layer, MSE loss, ADAM optimiser, rational approximation to ReLU is the activation function}

Table~\ref{tab:Results: experiments set 1 of NN with rational approximation to ReLU activation} summarises the results for the case of MSE loss, 2~nodes in the hidden layer and rational approximation to the activation function. The results clearly indicate that this approach is better than the direct use of ReLU as the activation function with 2~nodes in the hidden layer. The results are comparable with the results in the case of the direct use of ReLU as the activation function with 10~nodes in the hidden layer.
\begin{table}
    \centering
    \begin{tabular}{|c|c|c|c|}
    \hline
    Epoch & Final loss & Minimum loss & Run time (per epoch)\\
    \hline
    50   & 0.004621 &          & 2.85s $\pm$ 83.5ms\\
    50   &          & 0.004621 & \\
    \hline
    100  & 0.002096 &          & 2.94s $\pm$ 576ms\\
    100   &          & 0.002096 & \\
    \hline
    200  & 0.000578 &          &  2.87s $\pm$ 105ms\\
    199  &          & 0.000550 & \\
    \hline
    \end{tabular}
    \caption{Results: experiments set 1, rational approximation to ReLU}
    \label{tab:Results: experiments set 1 of NN with rational approximation to ReLU activation}
\end{table}



Figure~\ref{fig:Rational activation with ReLU coefficients - set 1, epoch is 50,100,200} confirms that in the case of only 2~nodes in the hidden layer, the rational approximation to ReLU is more accurate than ReLU itself as the activation function. The ``difficult point'' remains a challenge.

\begin{figure}
    \centering
    \includegraphics[width=40mm]{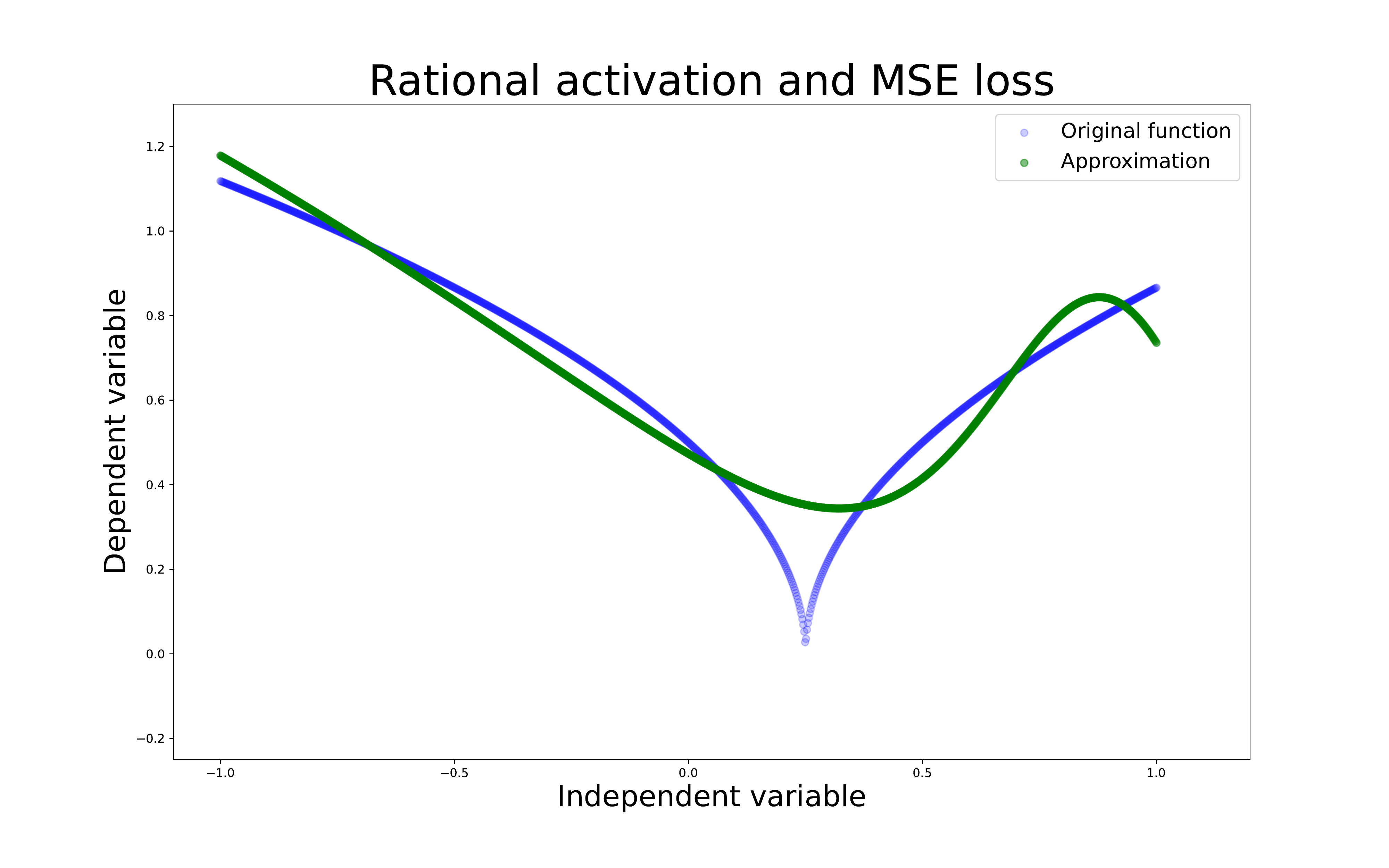}
        \includegraphics[width=40mm]{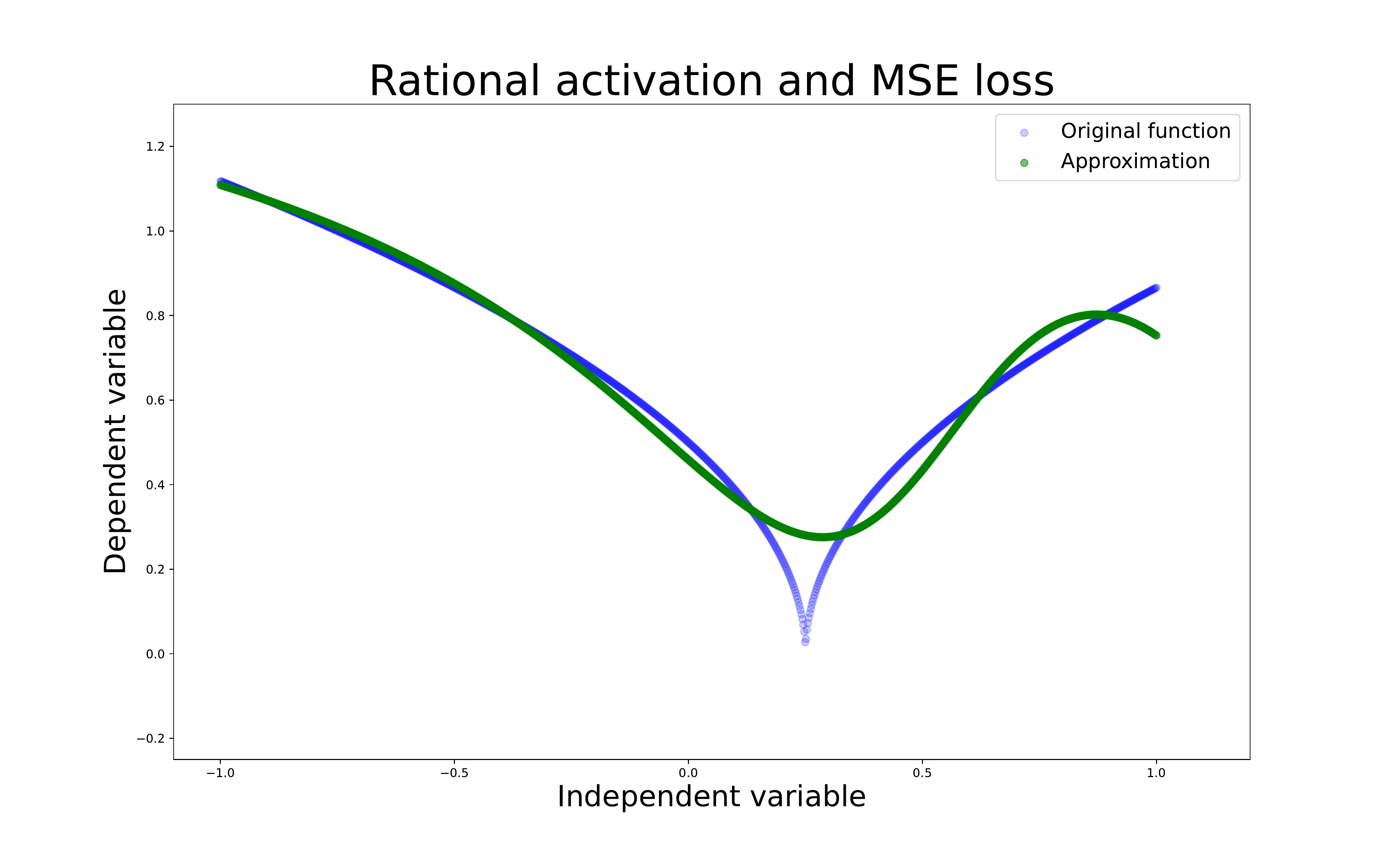}
            \includegraphics[width=40mm]{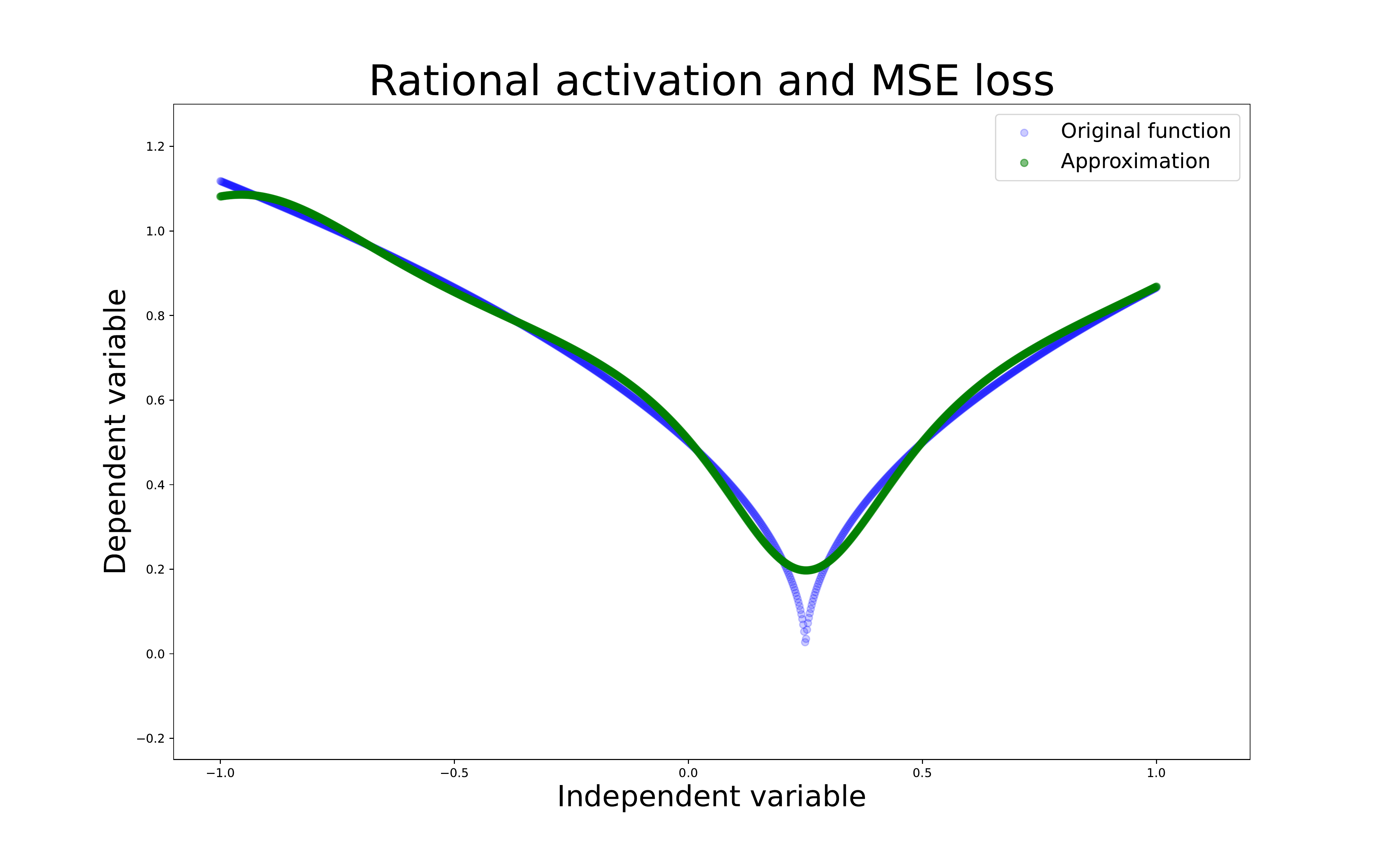}
    \caption{Set 1: Approximation is computed by a neural network with rational approximation to ReLU, 50, 100 and 200 epochs.}
    \label{fig:Rational activation with ReLU coefficients - set 1, epoch is 50,100,200}
\end{figure}

\paragraph{{\bf Set 2}: 2 nodes in the hidden layer, uniform loss, ADAMAX optimiser, rational approximation to ReLU is the activation function}

Table~\ref{tab:Results: experiments set 2 of NN with rational approximation to ReLU activation} summarises the results for the case of uniform loss, 2~nodes in the hidden layer and rational approximation to the activation function. Similarly to MSE loss, the results clearly indicate that this approach is better than the direct use of ReLU as the activation function with 2~nodes in the hidden layer. The results are comparable with the results in the case of the direct use of ReLU as the activation function with 10~nodes in the hidden layer.

\begin{table}
    \centering
    \begin{tabular}{|c|c|c|c|}
    \hline
    Epoch & Final loss & Minimum loss & Run time (per epoch) \\
    \hline
    50  & 0.211703 &          & 2.88s $\pm$ 726ms\\
    37  &          & 0.204005 & \\
    \hline
    100 & 0.212665 &          & 2.8s $\pm$ 125ms\\
    87  &          & 0.170603 & \\
    \hline
    200 & 0.166955 &          & 2.84s $\pm$ 337ms\\
    197  &          & 0.124133 & \\
    \hline
    \end{tabular}
    \caption{Results: experiments set 2, rational approximation to ReLU.}
    \label{tab:Results: experiments set 2 of NN with rational approximation to ReLU activation}
\end{table}



\begin{figure}
    \centering
    \includegraphics[width=40mm]{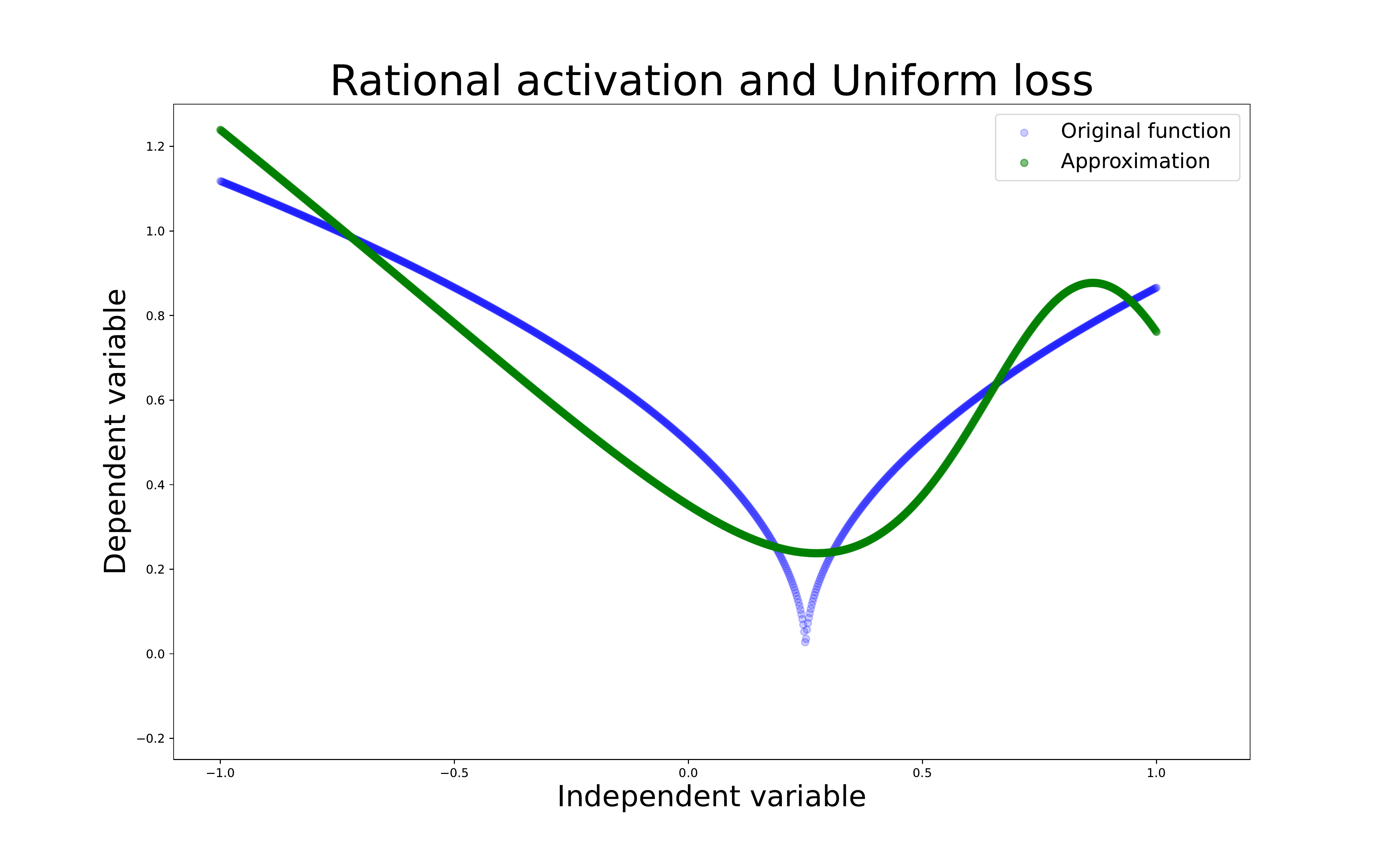}
     \includegraphics[width=40mm]{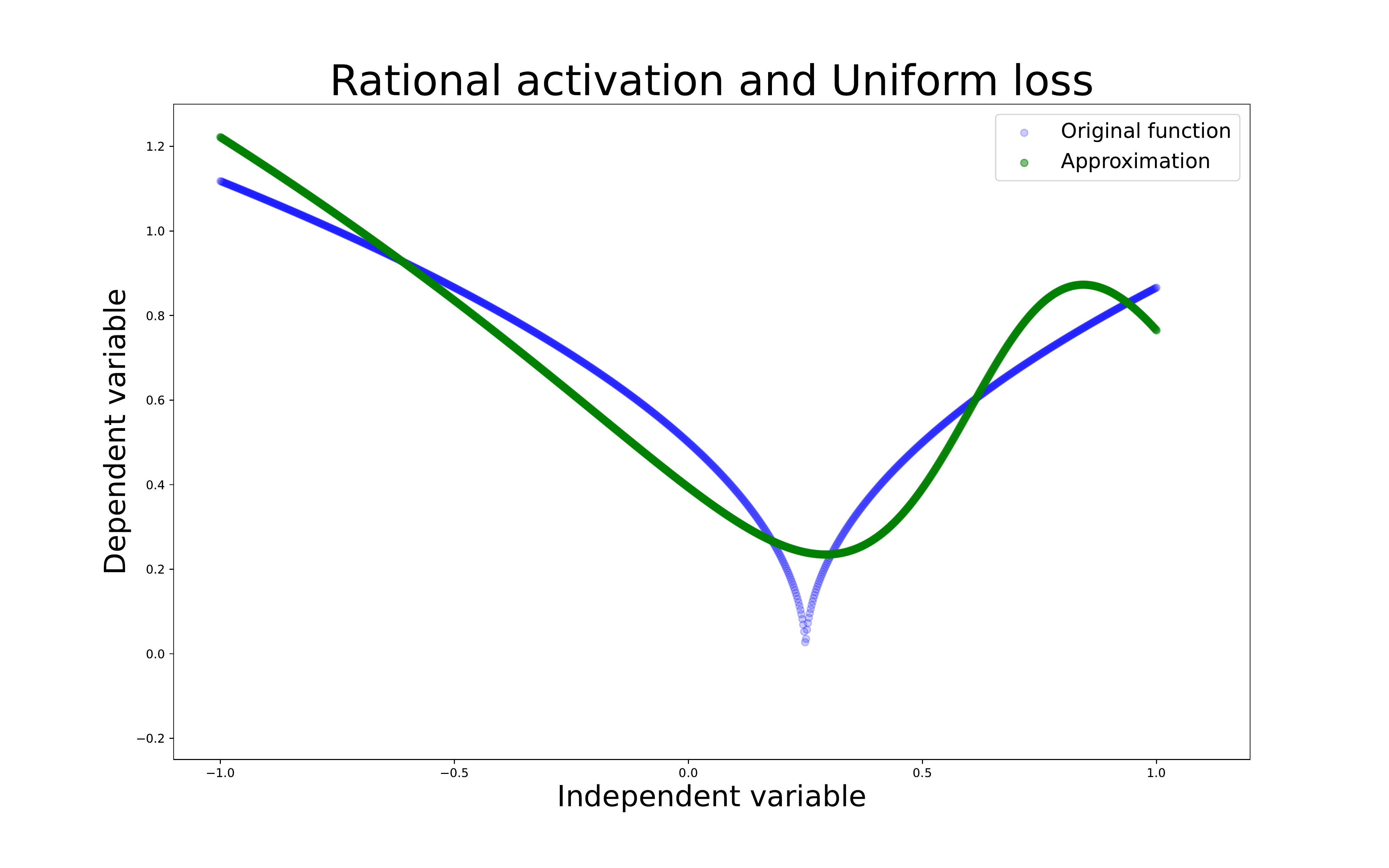}
      \includegraphics[width=40mm]{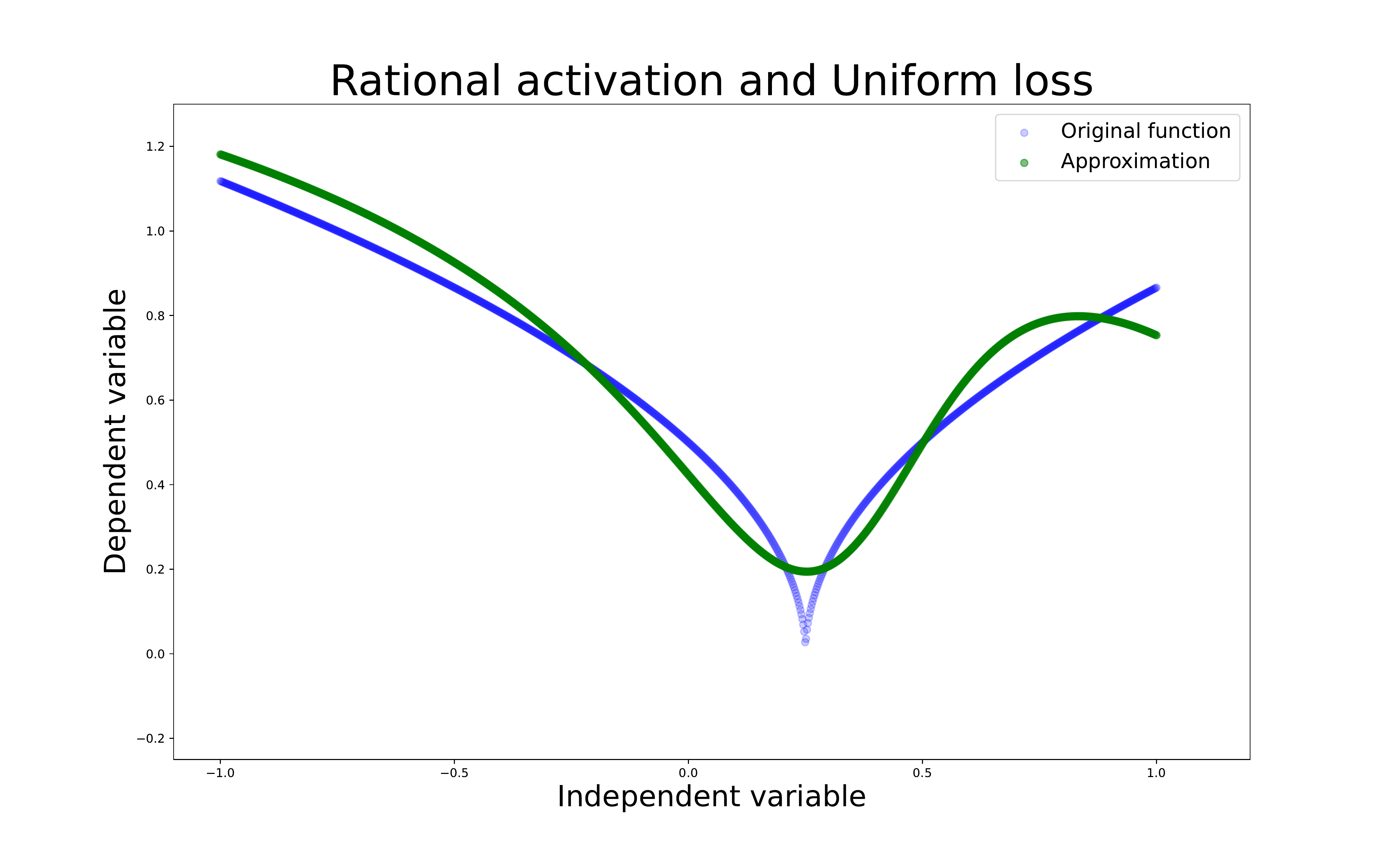}
    \caption{Set 2: Approximation is computed by a neural network with rational
activation to ReLU, 50, 100 and 200 epochs.}
    \label{fig:Rational activation with ReLU coefficients - set 2, epoch is 50}
\end{figure}


\paragraph{{\bf Set 3}: 10 nodes in the hidden layer, MSE loss, ADAM optimiser, rational approximation to ReLU is the activation function}

Table~\ref{tab:Results: experiments set 3 of NN with rational approximation to ReLU activation} shows further improvement in the accuracy when the number of nodes in the hidden layer is increased from~2 to~10. This observation is especially prominent when the number of epochs is~200. 

\begin{table}
    \centering
    \begin{tabular}{|c|c|c|c|}
    \hline
    Epoch & Final loss & Minimum loss & Run time (per epoch) \\
    \hline
    50  & 0.000885 &          & 2.78s $\pm$ 93.5ms\\
    50  &          & 0.000885 & \\
    \hline
    100 & 0.000153 &          & 2.85s $\pm$ 539ms\\
    100  &          & 0.000153 & \\
    \hline
    200 & 0.000046 &          & 2.79s $\pm$ 93.1ms\\
    190 &          & 0.000046 & \\
    \hline
    \end{tabular}
    \caption{Results: experiments set 3, rational approximation to ReLU}
    \label{tab:Results: experiments set 3 of NN with rational approximation to ReLU activation}
\end{table}



Figure~\ref{fig:Rational activation with ReLU coefficients - set 3, epoch is 50,100,200} illustrates that the accuracy of approximation is improving, even for the ``difficult point'', especially when the number of epochs is~200.

\begin{figure}
    \centering
    \includegraphics[width=40mm]{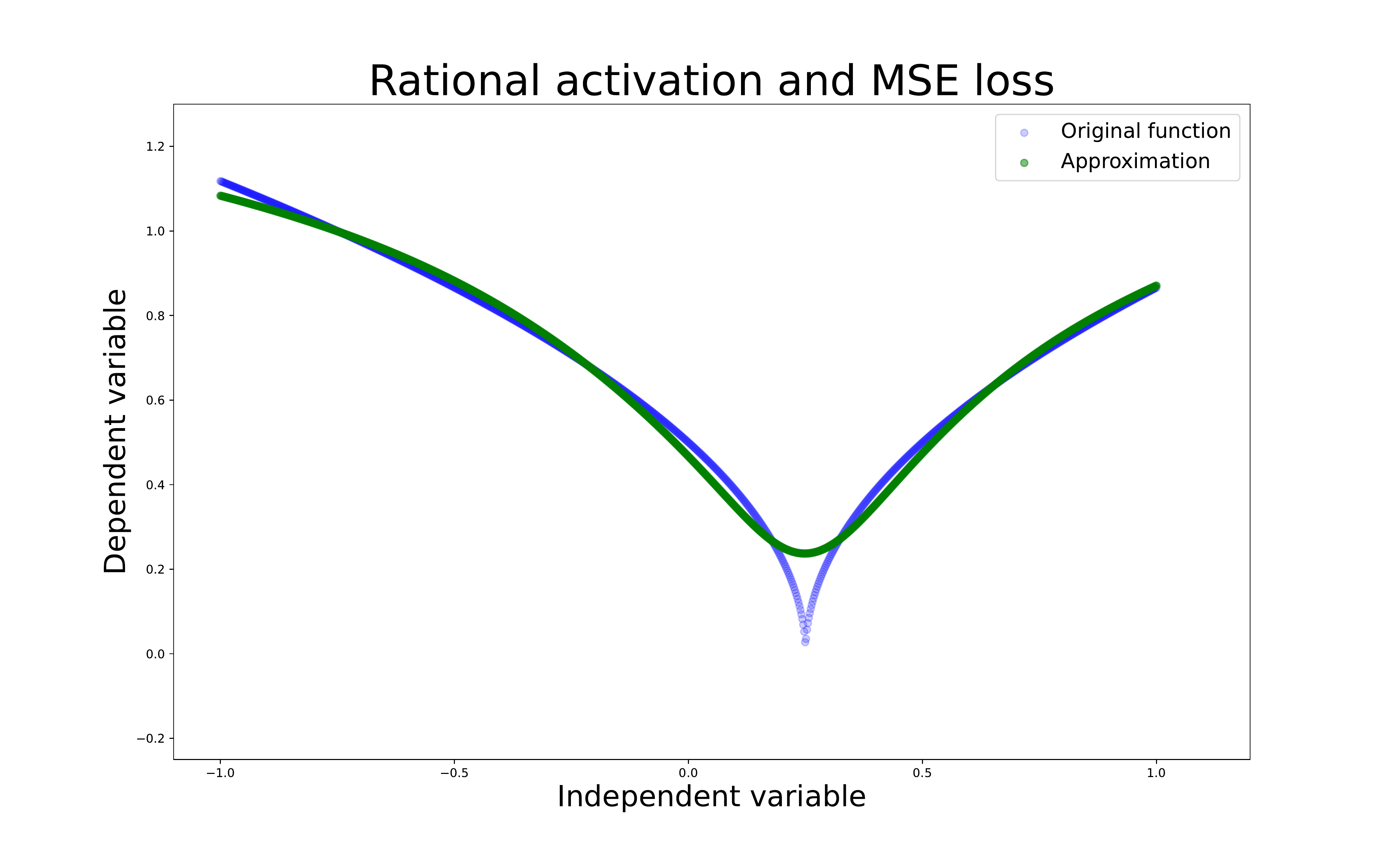}
        \includegraphics[width=40mm]{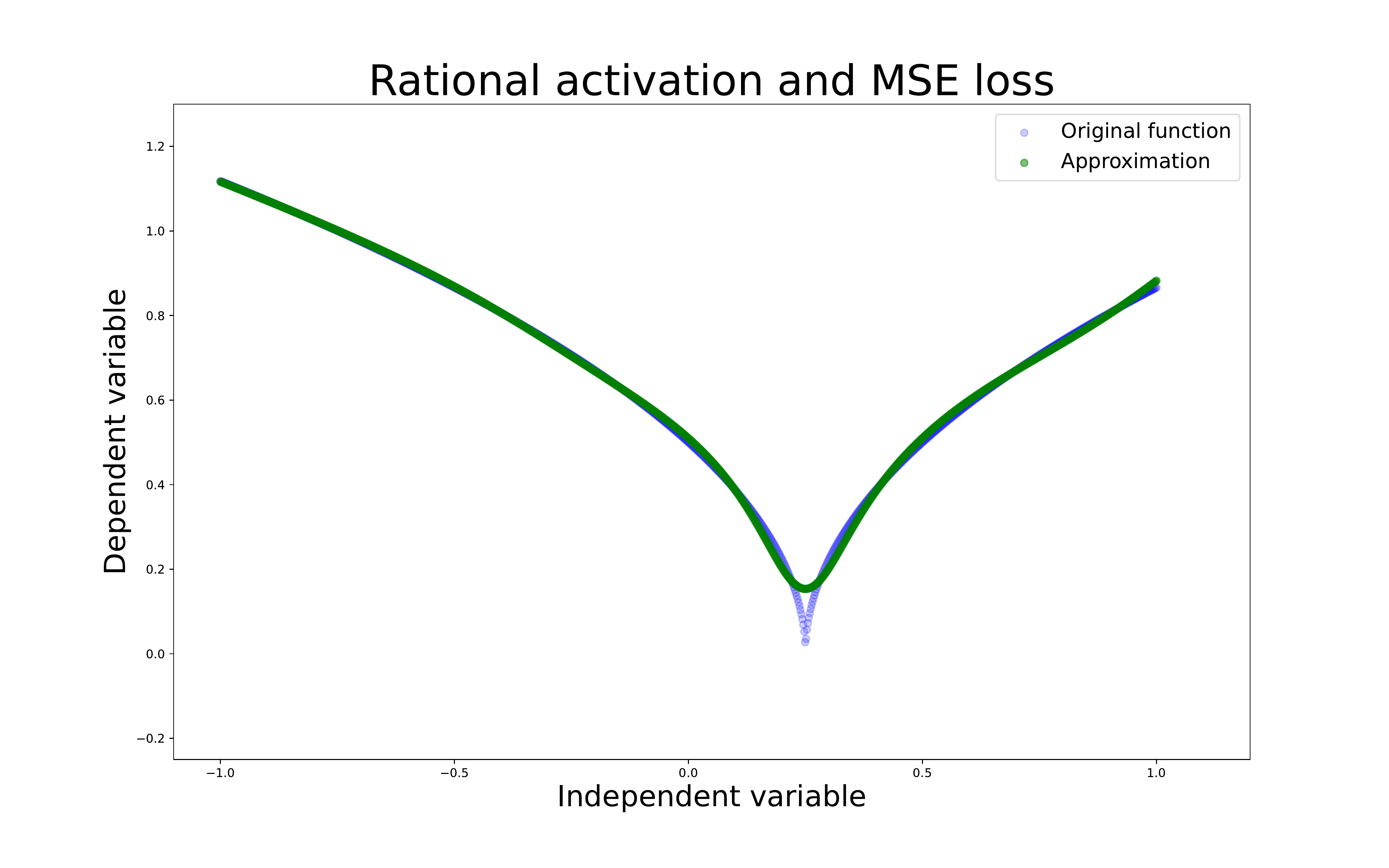}
            \includegraphics[width=40mm]{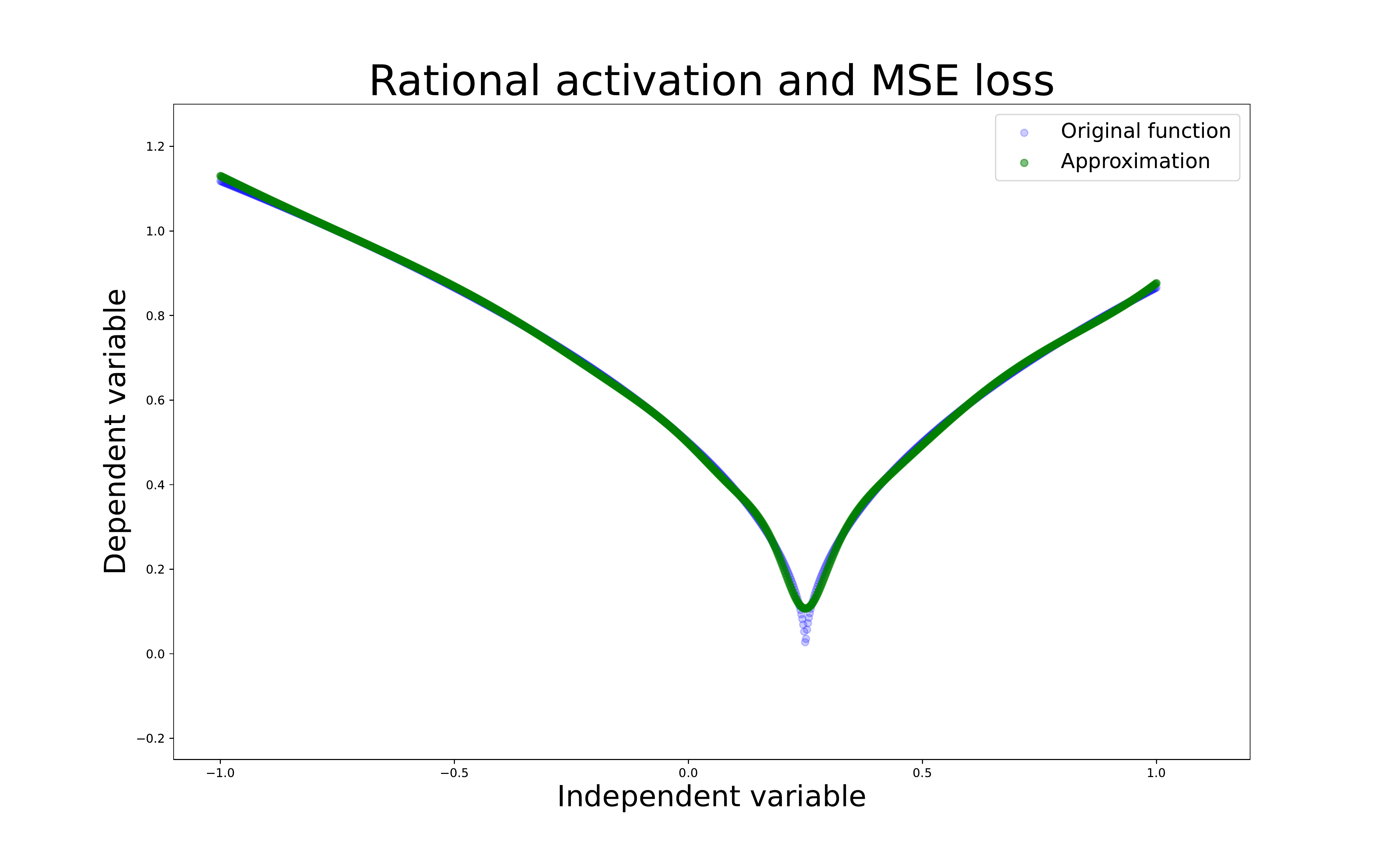}
    \caption{Set 3: Approximation is computed by a neural network with rational
activation to ReLU, 50, 100 and 200 epochs.}
    \label{fig:Rational activation with ReLU coefficients - set 3, epoch is 50,100,200}
\end{figure}

\subsubsection{Neural Network with rational activation}\label{appendix:results3:Set1-4}

In this case, our activation function is a rational function of degree $(3,2)$. 
The coefficients of the rational activation function are now a part of the parameter set. We learn these coefficients as we learn other parameters during the training procedure.

\paragraph{{\bf Set 1}: 2 nodes in the hidden layer, MSE loss, ADAM optimiser, rational activation function}

Table~\ref{tab:Results: experiments set 1 of NN with rational approximation and usual training} shows that the values of the loss function are better than they are in the case of ReLU and comparable with the results in the case of rational approximation to ReLU (rational approximation to ReLU is slightly more accurate).  Figure~\ref{fig:Rational activation with usual training - set 1, epoch is 50,100,200} confirms this observation.

Interestingly, the computational time per epoch is increasing, but the overall running time is similar to ReLU and ReLU approximated by rational functions.

\begin{table}
    \centering
    \begin{tabular}{|c|c|c|c|}
    \hline
    Epoch & Final loss & Minimum loss & Run time (per epoch)\\
    \hline
    50   & 0.005893 &          & 2.54s $\pm$ 111ms\\
    50   &          & 0.005893 & \\
    \hline
    100  & 0.001330 &          & 2.81s $\pm$ 339ms\\
    100   &          & 0.001330 & \\
    \hline
    200  & 0.000773 &          &  3.05s $\pm$ 433ms\\
    170  &          & 0.000745 & \\
    \hline
    \end{tabular}
    \caption{Results: experiments set 1, rational activation function}
    \label{tab:Results: experiments set 1 of NN with rational approximation and usual training}
\end{table}



\begin{figure}
    \centering
    \includegraphics[width=40mm]{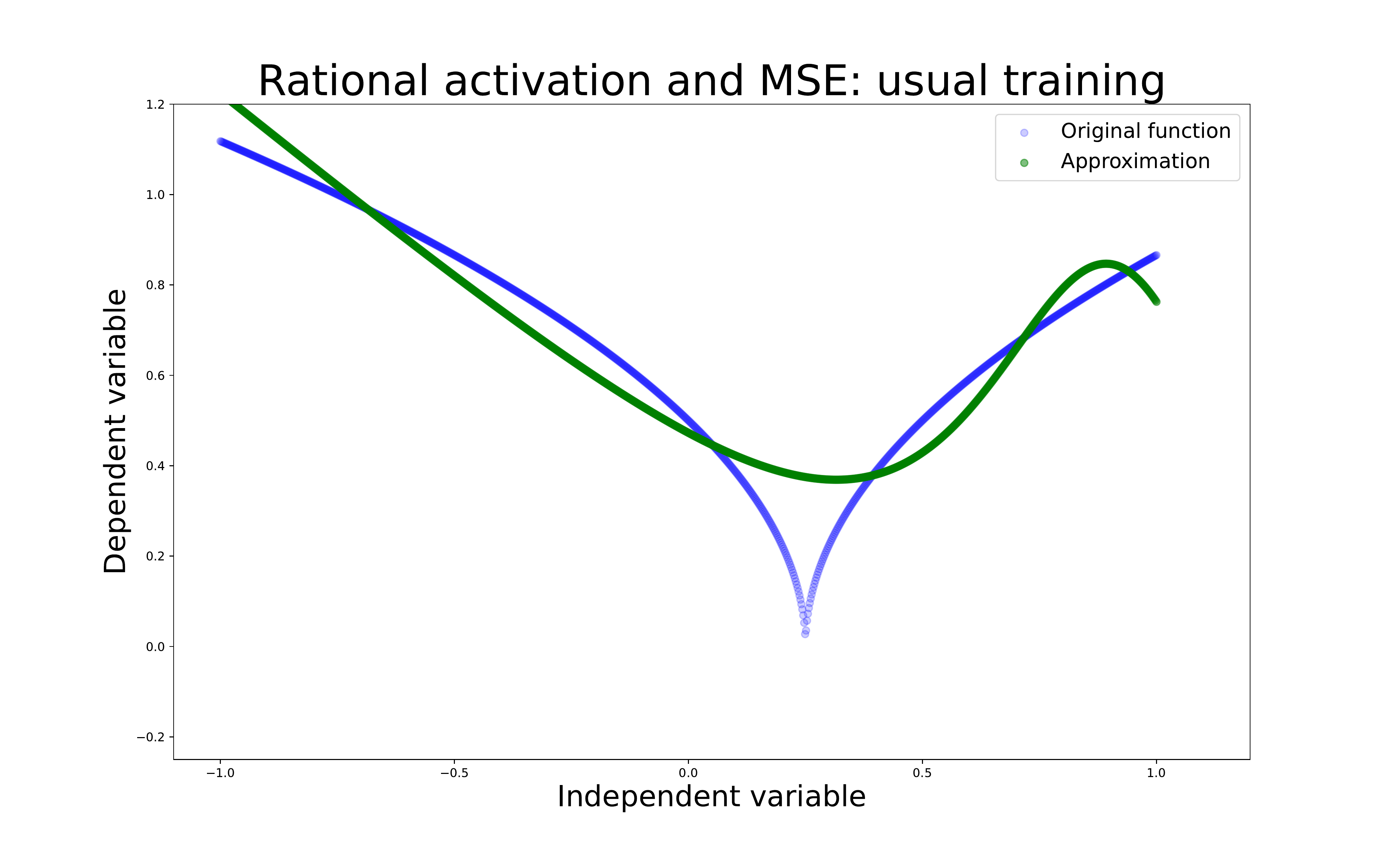}
    \includegraphics[width=40mm]{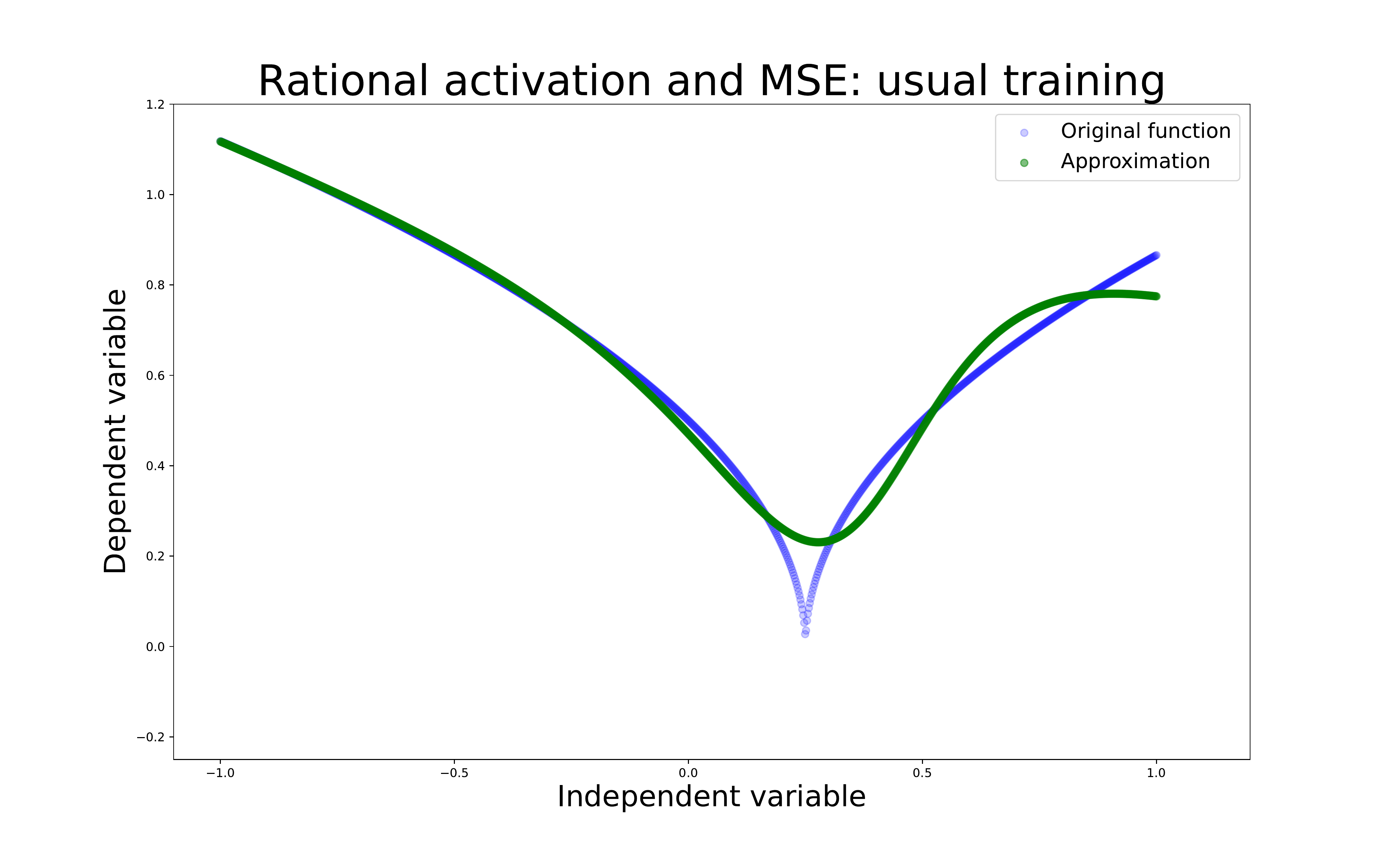}
    \includegraphics[width=40mm]{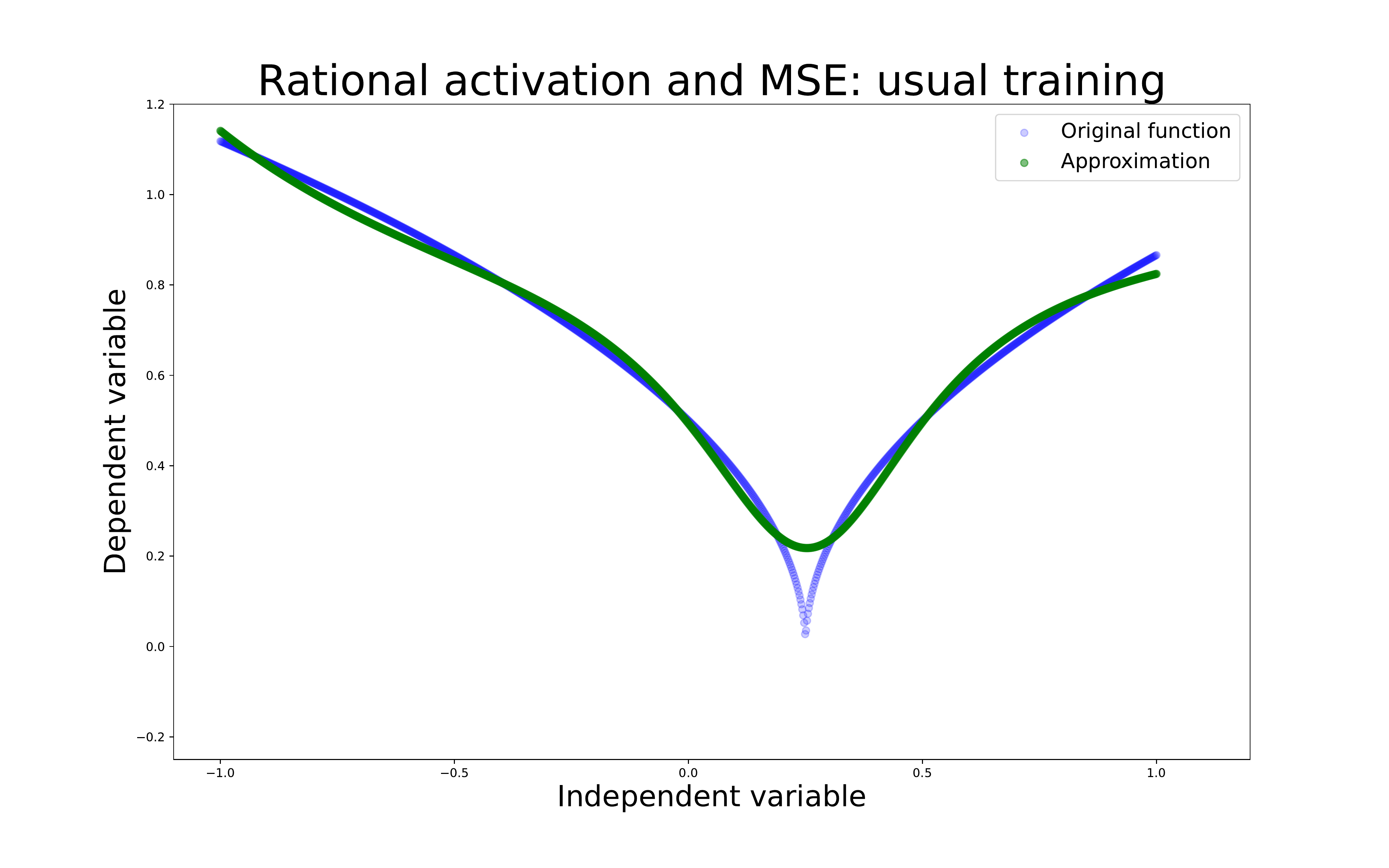}
    \caption{Approximation is computed by the usual training process,  50, 100 and 200~epochs.}
    \label{fig:Rational activation with usual training - set 1, epoch is 50,100,200}
\end{figure}


\paragraph{{\bf Set 2}: 2 nodes in the hidden layer, Uniform loss, ADAMAX optimiser, rational activation function}

The situation is similar to MSE loss. 
Table~\ref{tab:Results: experiments set 2 of NN with rational approximation and usual training} shows that the values of the loss function are better than they are in the case of ReLU and comparable with the results in the case of rational approximation to ReLU.  Figure~\ref{fig:Rational activation with usual training - set 2, epoch is 50} confirms this observation. 

\begin{table}
    \centering
    \begin{tabular}{|c|c|c|c|}
    \hline
    Epoch & Final loss & Minimum loss & Run time (per epoch) \\
    \hline
    50  & 0.219781 &          & 2.57s $\pm$ 616ms\\
    43  &          & 0.206314 & \\
    \hline
    100 & 0.183087 &          & 2.84s $\pm$ 82.3ms\\
    100  &          & 0.183087 & \\
    \hline
    200 & 0.151559 &          & 2.87s $\pm$ 85.4ms\\
    198  &          & 0.123221 & \\
    \hline
    \end{tabular}
    \caption{Results: experiments set 2, rational activation function}
    \label{tab:Results: experiments set 2 of NN with rational approximation and usual training}
\end{table}



\begin{figure}
    \centering
    \includegraphics[width=40mm]{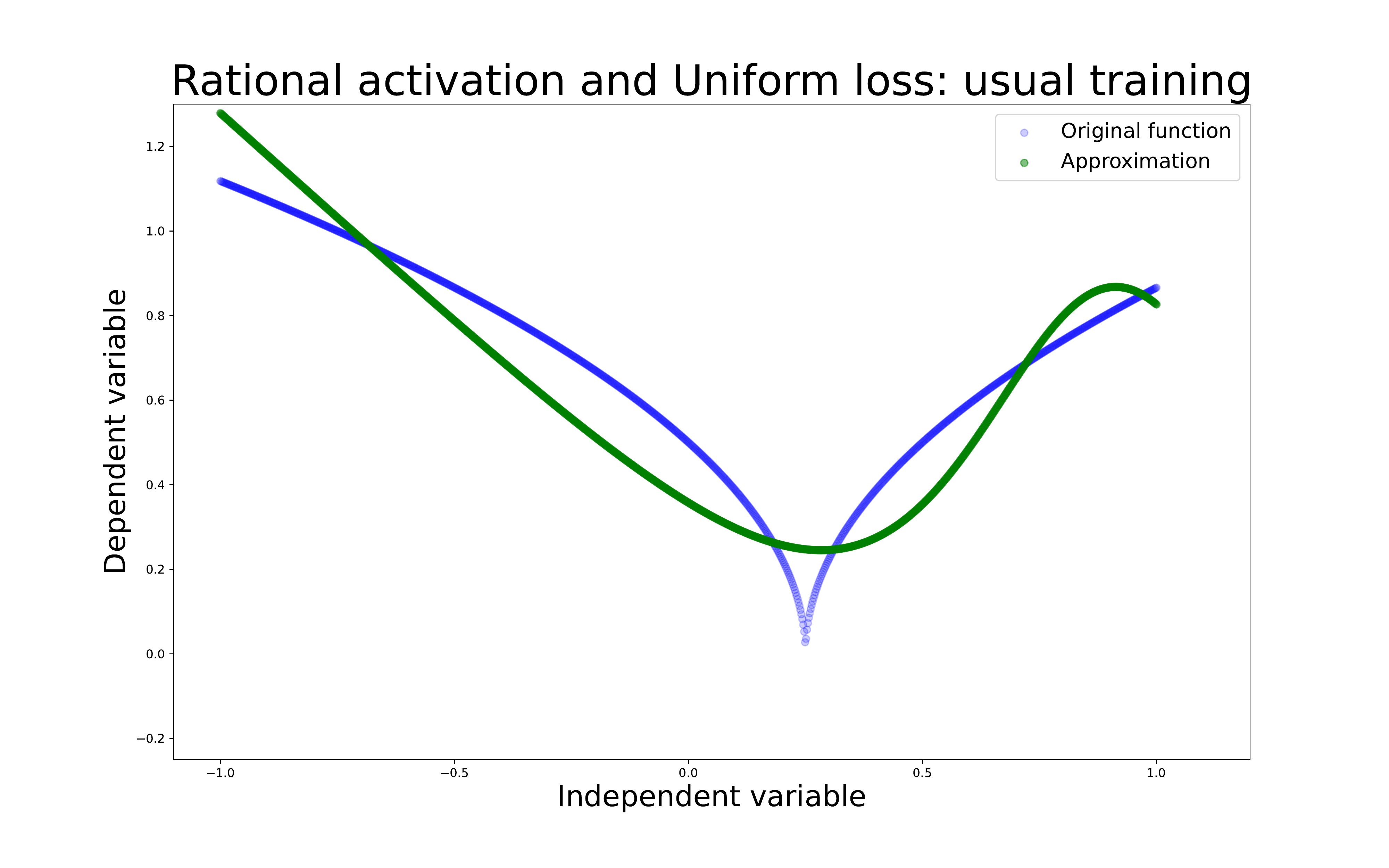}
        \includegraphics[width=40mm]{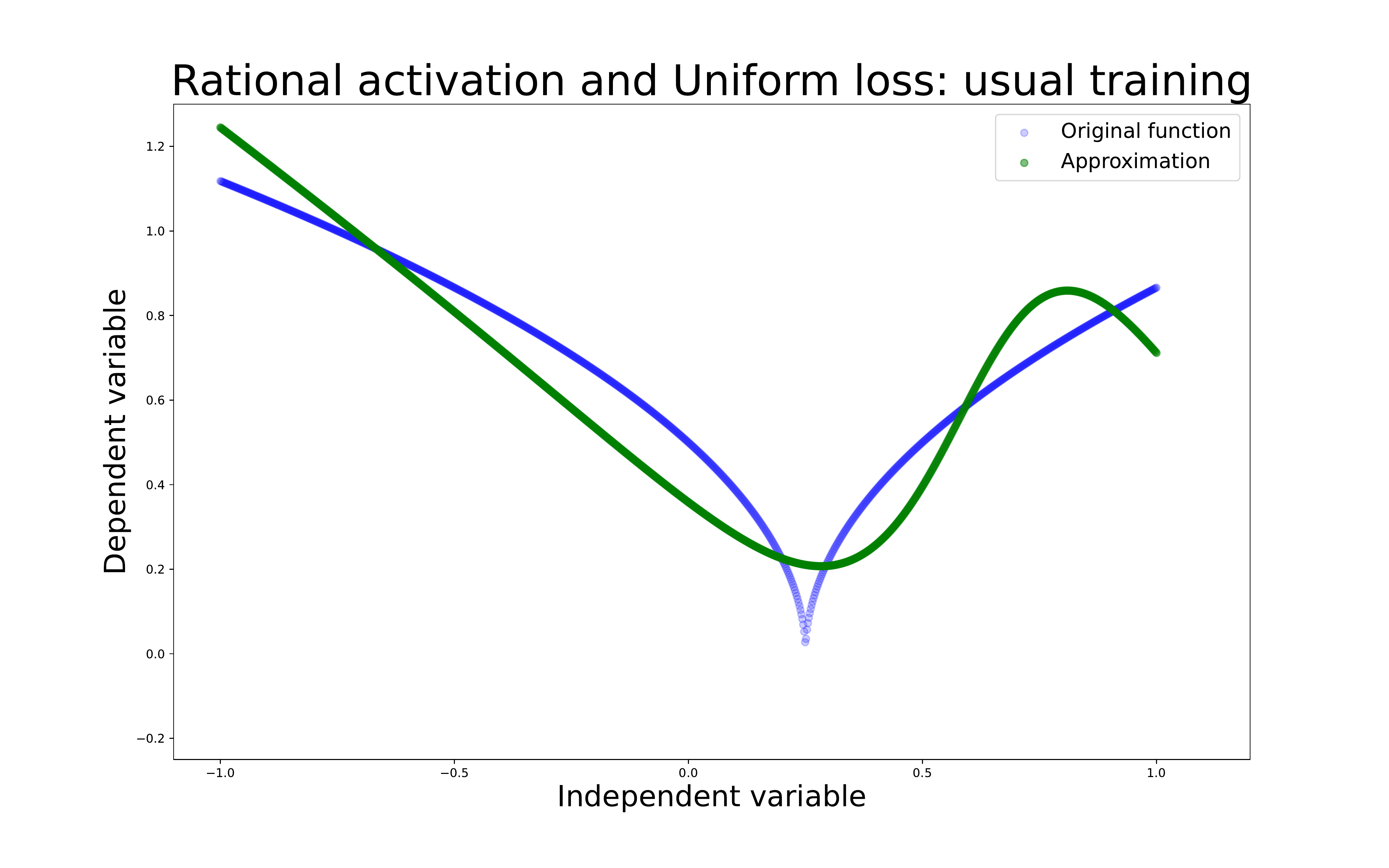}
            \includegraphics[width=40mm]{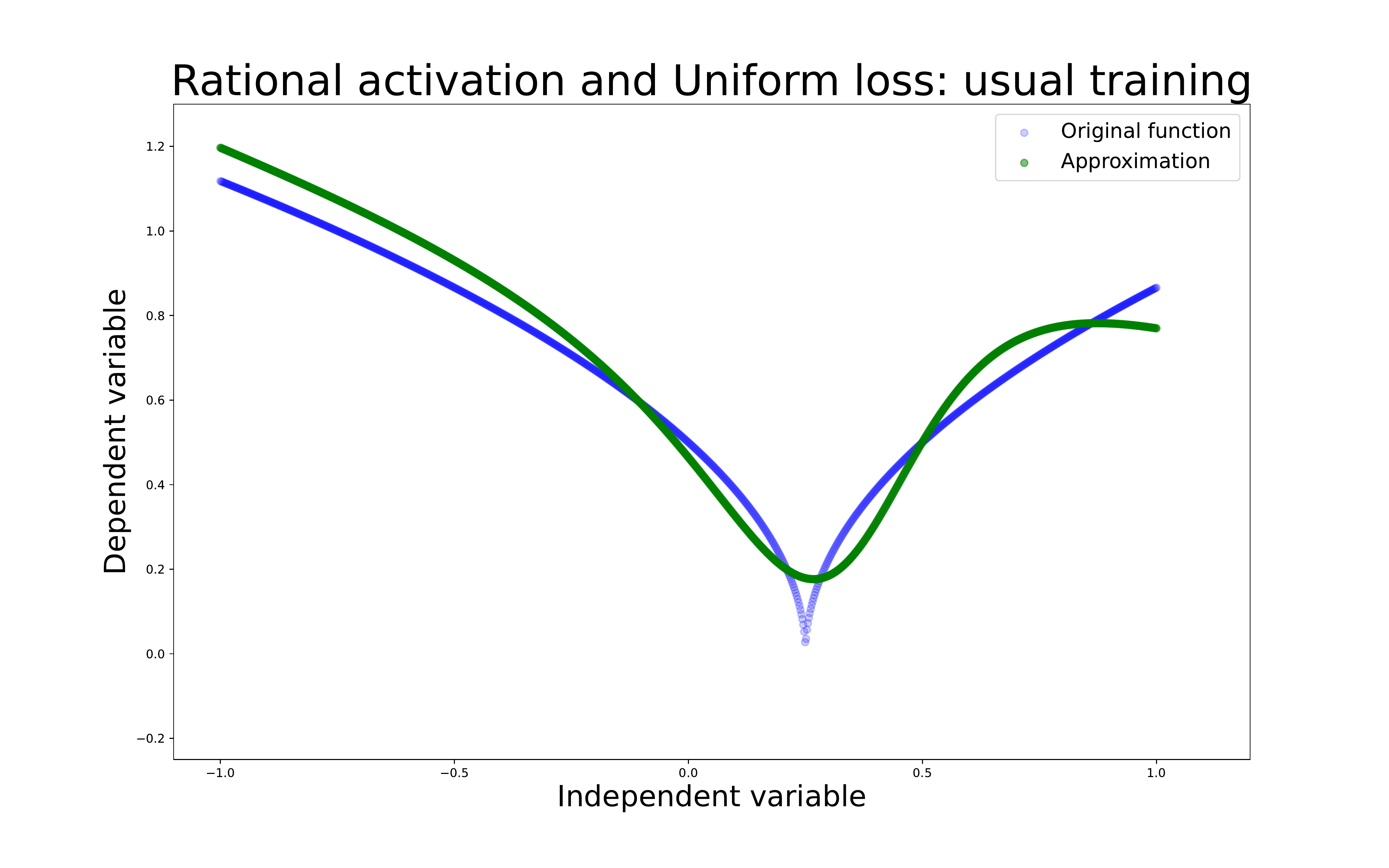}
    \caption{Approximation is computed by the usual training process,  50, 100 and 200~epochs.}
    \label{fig:Rational activation with usual training - set 2, epoch is 50}
\end{figure}


\paragraph{{\bf Set 3}: 10 nodes in the hidden layer, MSE loss, ADAM optimiser, rational activation function}

Table~\ref{tab:Results: experiments set 3 of NN with rational approximation and usual training} shows that when the number of nodes in the hidden layer is increased to~10, the optimal loss function values are very close to the case when ReLU was approximated by the rational function. Figure~\ref{fig:Rational activation with usual training - set 3, epoch is 50,100,200} confirms these results. It is important to note that the increase in the number of nodes in the hidden layer does not increase significantly the computational time.

\begin{table}
    \centering
    \begin{tabular}{|c|c|c|c|}
    \hline
    Epoch & Final loss & Minimum loss & Run time (per epoch) \\
    \hline
    50  & 0.000935 &          & 2.91s $\pm$ 711ms\\
    50  &          & 0.000935 & \\
    \hline
    100 & 0.000190 &          & 2.76s $\pm$ 91.6ms\\
    100  &          & 0.000190 & \\
    \hline
    200 & 0.000070 &          & 2.52s $\pm$ 295ms\\
    198 &          & 0.000043 & \\
    \hline
    \end{tabular}
    \caption{Results: experiments set 3, rational activation function}
    \label{tab:Results: experiments set 3 of NN with rational approximation and usual training}
\end{table}



The approximation appears to be accurate even around the ``difficult point'' when the number of epochs is~200.

\begin{figure}
    \centering
    \includegraphics[width=40mm]{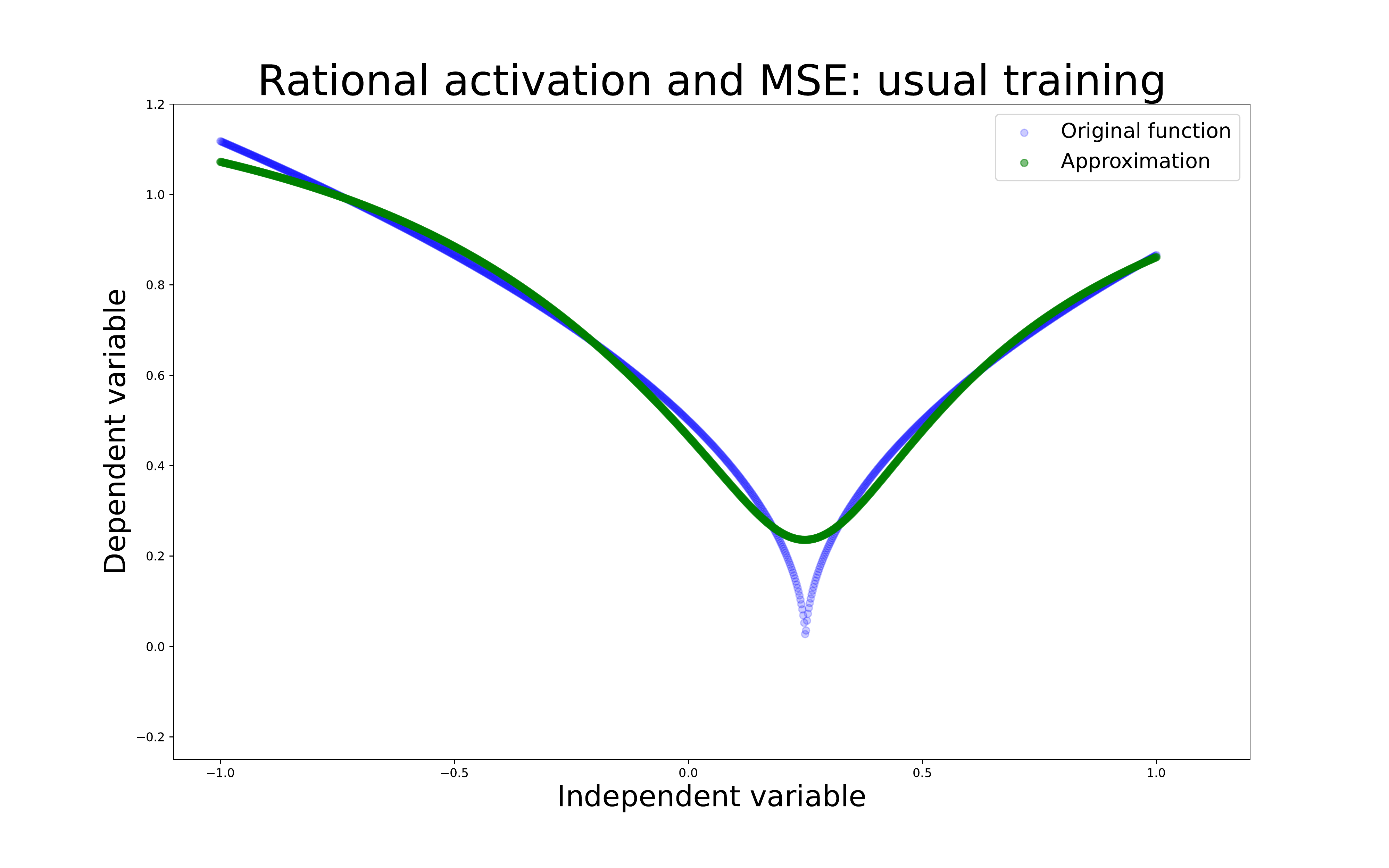}
        \includegraphics[width=40mm]{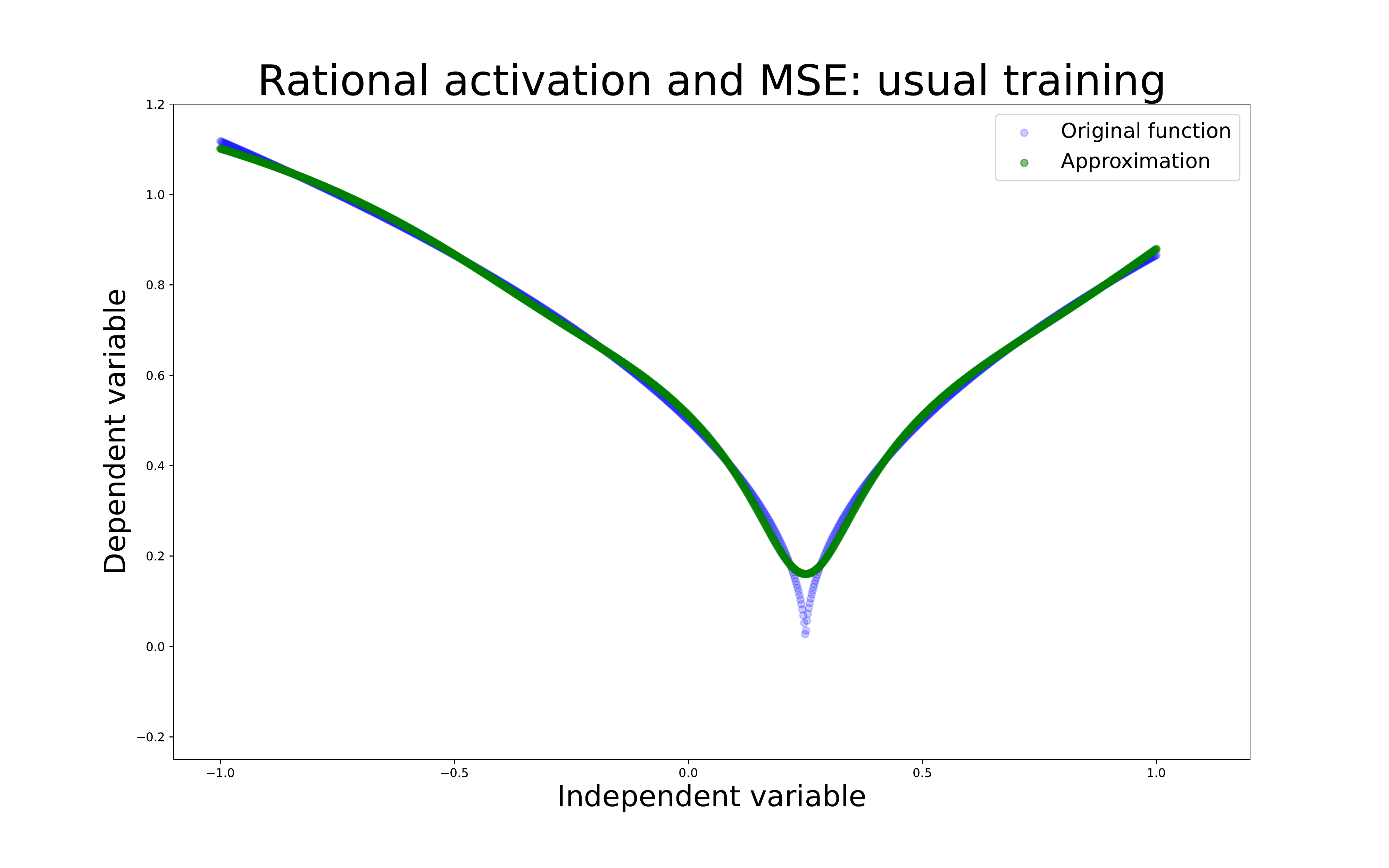}
            \includegraphics[width=40mm]{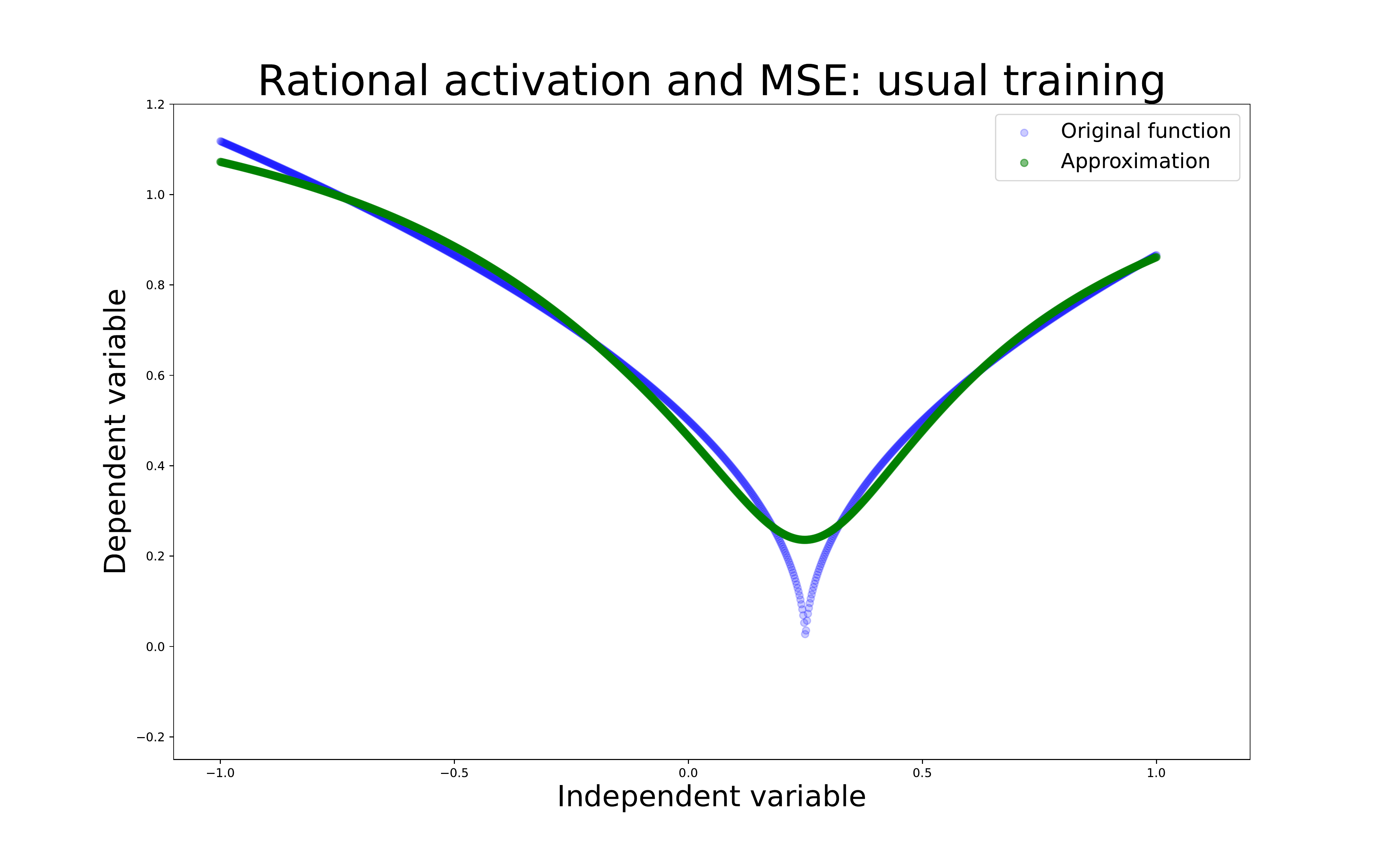}
    \caption{Approximation is computed by the usual training process,  50, 100 and 200~epochs.}
    \label{fig:Rational activation with usual training - set 3, epoch is 50,100,200}
\end{figure}

\subsubsection{Neural Network with rational activation, learn with split training method} \label{appendix:results4:Set1-4}

In this section, our activation function is a rational function of degree $(3,2)$. The coefficients are now a part of the parameter set. We learn these coefficients as we learn other parameters during the training procedure using a split method. 

\paragraph{{\bf Set 1}: 2 nodes in the hidden layer, MSE loss, ADAM optimiser, split method}

Table~\ref{tab:Results: experiments set 1 of NN with rational approximation and split training} shows that the results are similar to those with rational approximation to ReLU from the previous section in the sense of optimal loss function values and the corresponding computational time. The computational time is also similar. The results are improving with the increase in the number of epochs. Figure~\ref{fig:Rational activation with split training - set 1,  50,100,200} confirms this similarity.

\begin{table}
    \centering
    \begin{tabular}{|c|c|c|c|}
    \hline
    Epoch & Final loss & Minimum loss & Run time (per epoch)\\
    \hline
    50   & 0.005291 &          & 2.51s $\pm$ 70.5ms\\
    50   &          & 0.005291 & \\
    \hline
    100  & 0.002445 &          & 2.65s $\pm$ 448ms\\
    100   &          & 0.002445 & \\
    \hline
    200  & 0.000602 &          & 3.39s $\pm$ 585ms\\
    197  &          & 0.000577 & \\
    \hline
    \end{tabular}
    \caption{Results: experiments set 1, split method}
    \label{tab:Results: experiments set 1 of NN with rational approximation and split training}
\end{table}



\begin{figure}
    \centering
    \includegraphics[width=40mm]{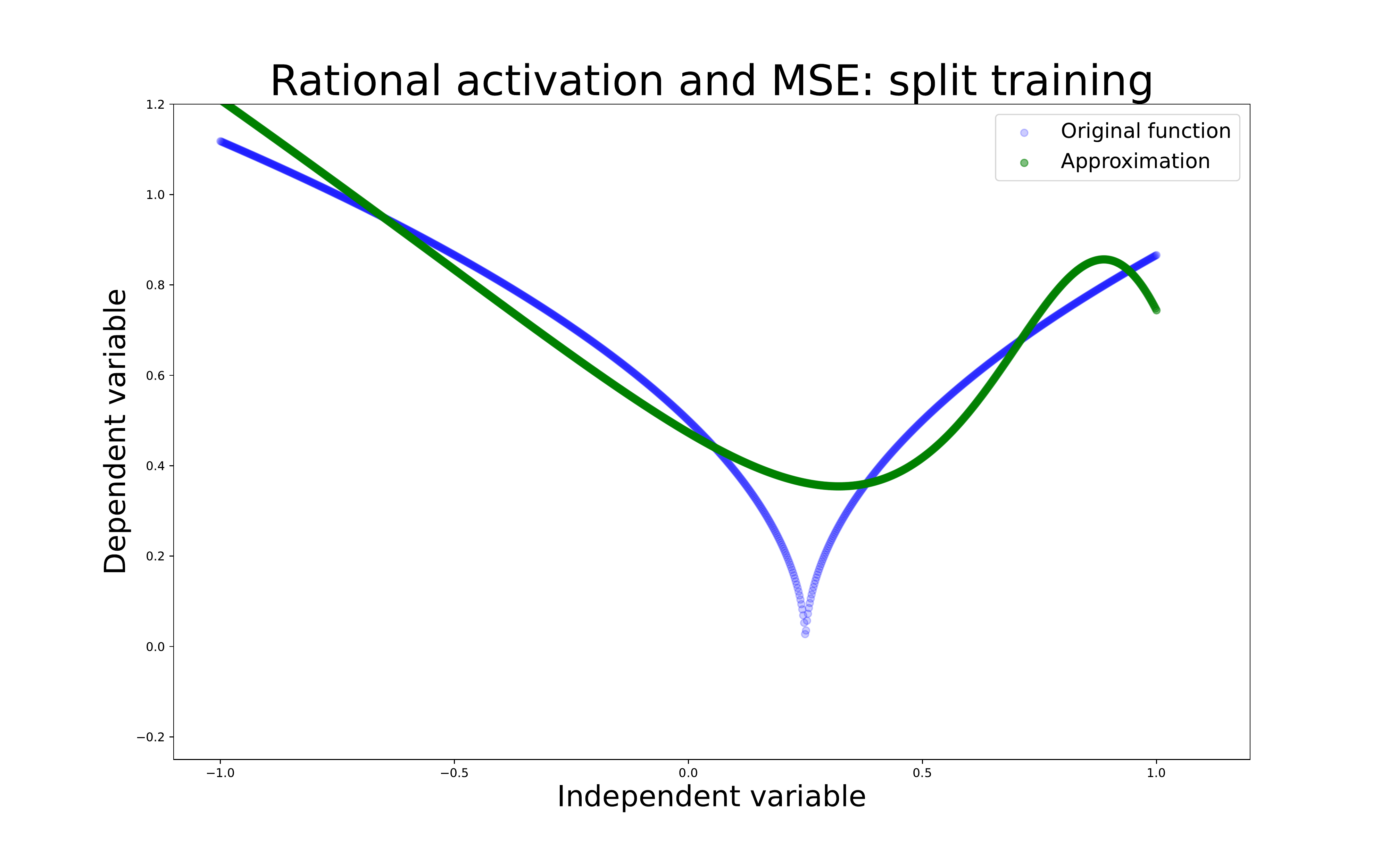}
 \includegraphics[width=40mm]{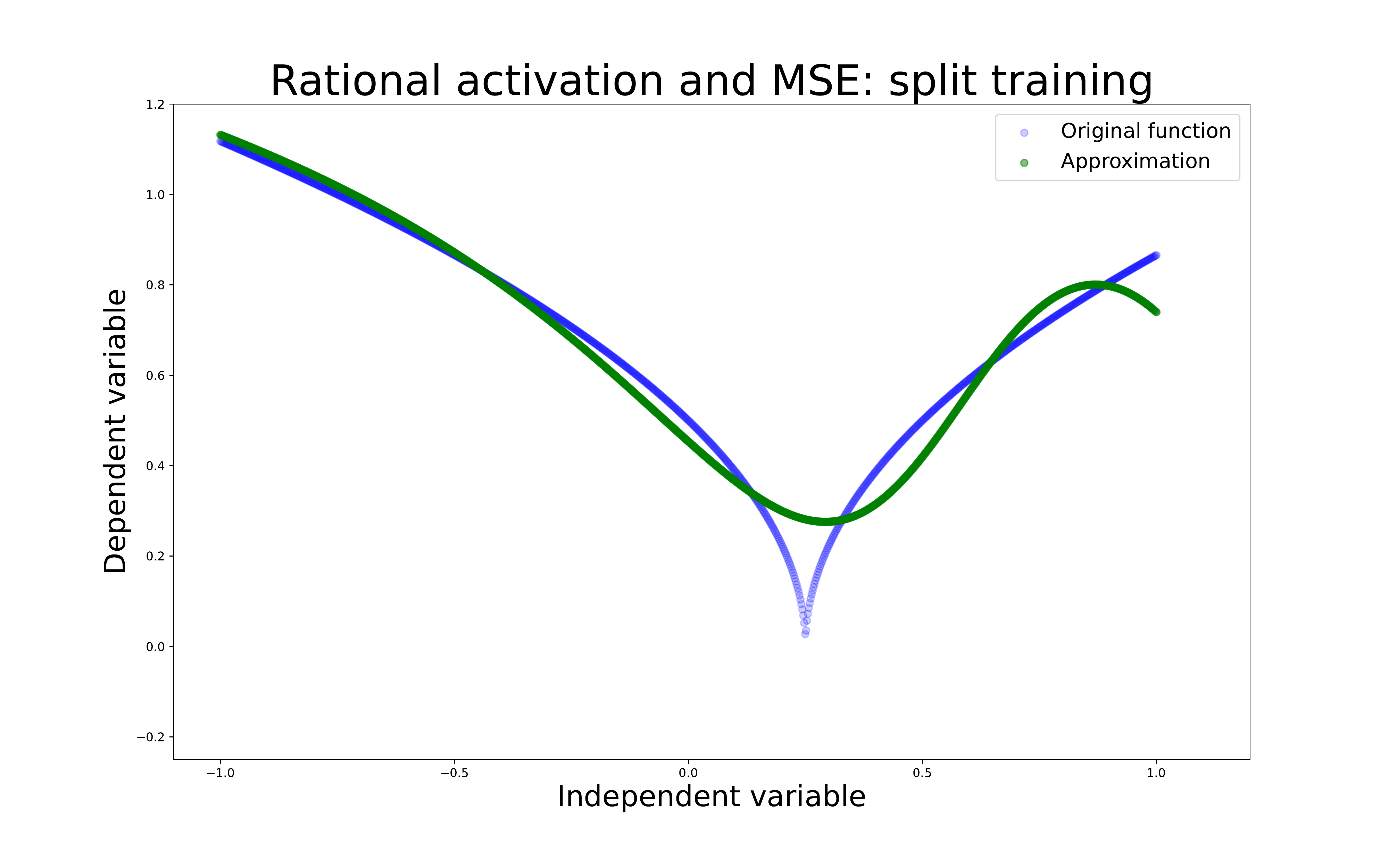}
      \includegraphics[width=40mm]{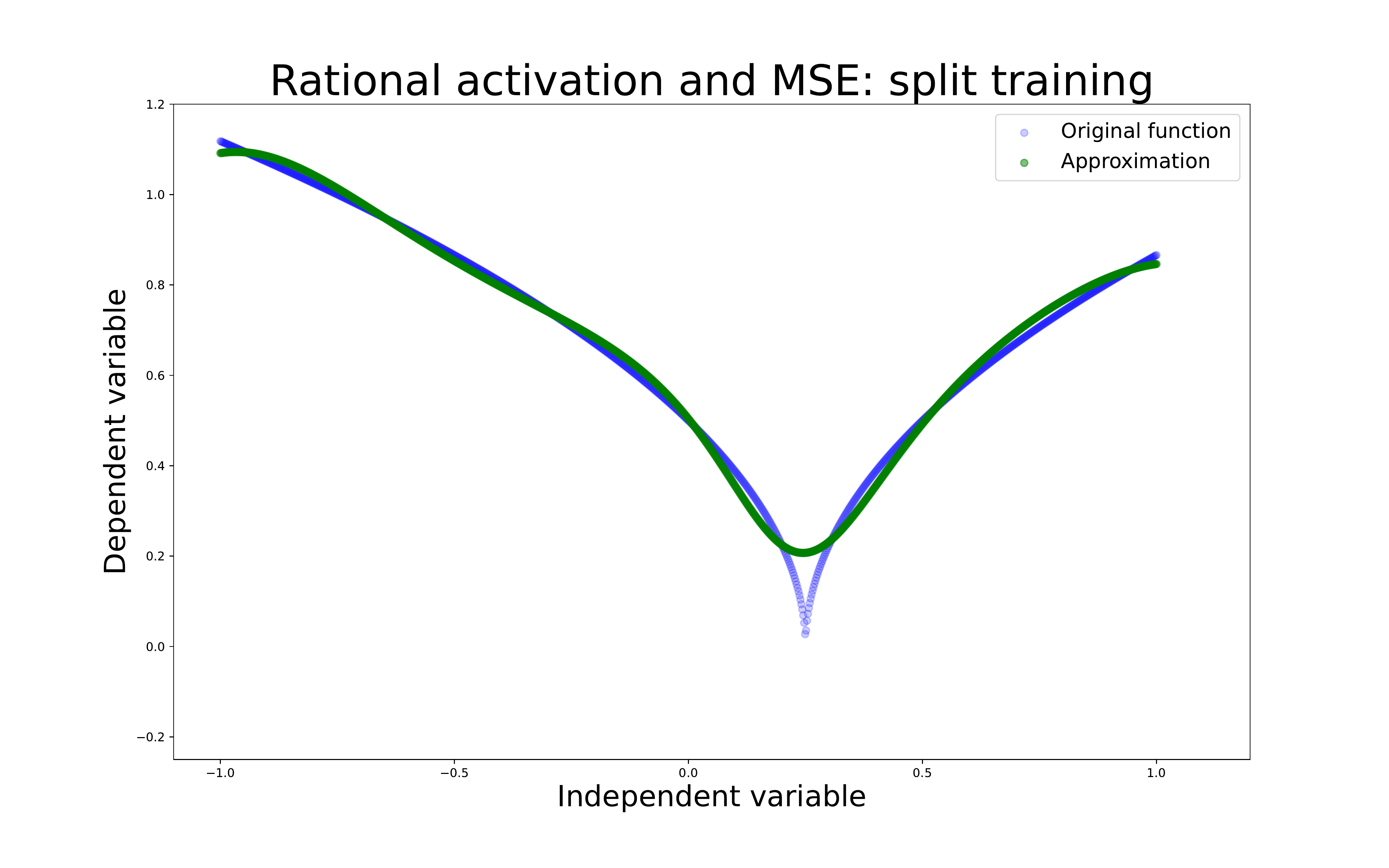}
    \caption{Approximation is computed by the split training process, the epoch is 50, 100 and 200.}
    \label{fig:Rational activation with split training - set 1,  50,100,200}
\end{figure}


\paragraph{{\bf Set 2}: 2 nodes in the hidden layer, Uniform loss, ADAMAX optimiser, split method.}

The results with only 2~nodes in the hidden layer are better for the split method than they are in the case of rational activation function or rational approximation to ReLU activation in terms of the optimal loss function value (table~\ref{tab:Results: experiments set 2 of NN with rational approximation and split training}). The computational time is slightly higher for the split method, but the difference is very insignificant. The average computational time per epoch remains the same when the number of epochs is increasing.

\begin{table}
    \centering
    \begin{tabular}{|c|c|c|c|}
    \hline
    Epoch & Final loss & Minimum loss & Run time (per epoch) \\
    \hline
    50  & 0.196001 &          & 3.02s $\pm$ 700ms\\
    43  &          & 0.144596 & \\
    \hline
    100 & 0.144665 &          & 3.07s $\pm$ 549ms\\
    78  &          & 0.106687 & \\
    \hline
    200 & 0.122333 &          & 3.03s $\pm$ 110ms\\
    145  &          & 0.085358 & \\
    \hline
    \end{tabular}
    \caption{Results: experiments set 2, split method.}
    \label{tab:Results: experiments set 2 of NN with rational approximation and split training}
\end{table}



\begin{figure}
    \centering
    \includegraphics[width=40mm]{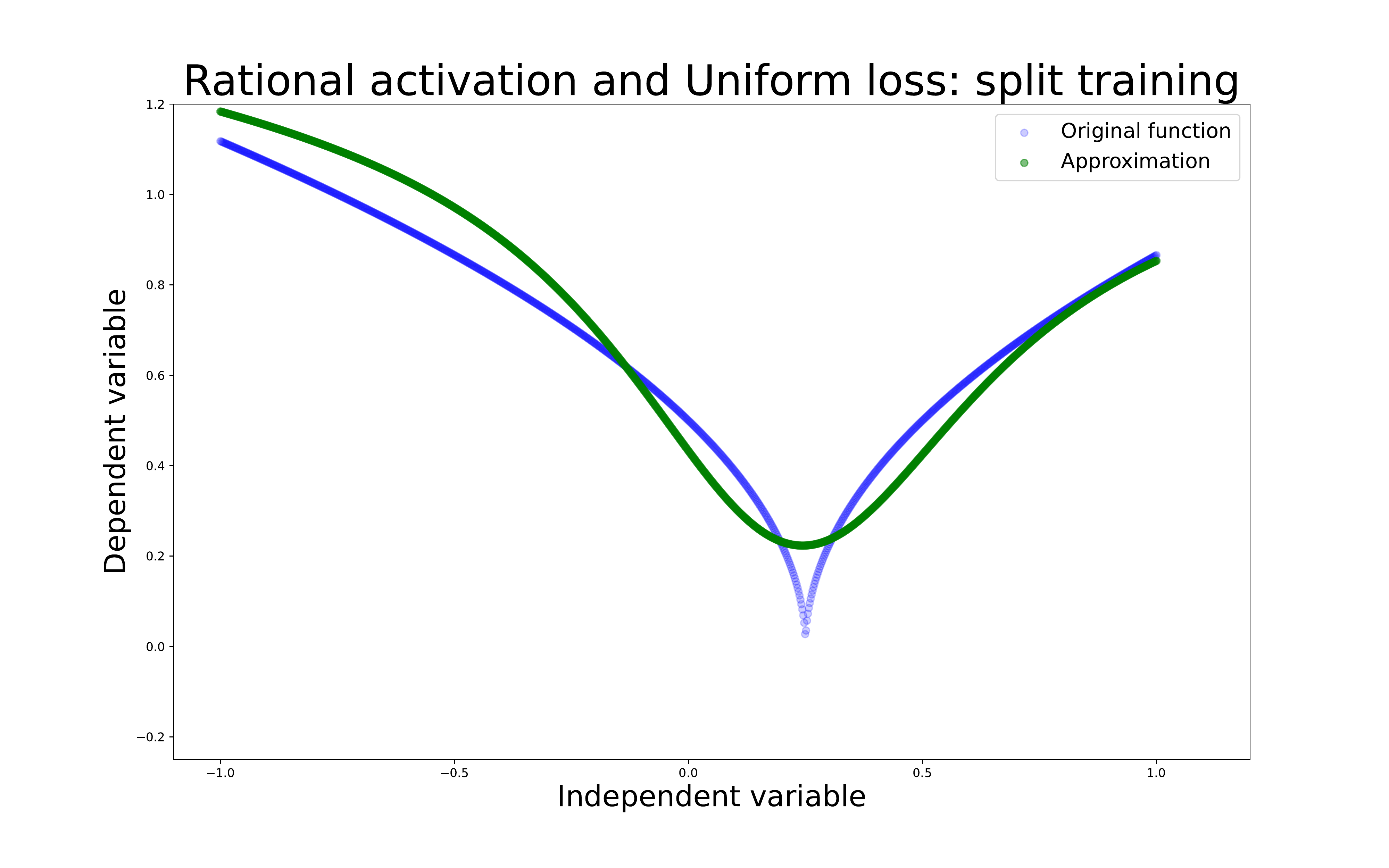}
        \includegraphics[width=40mm]{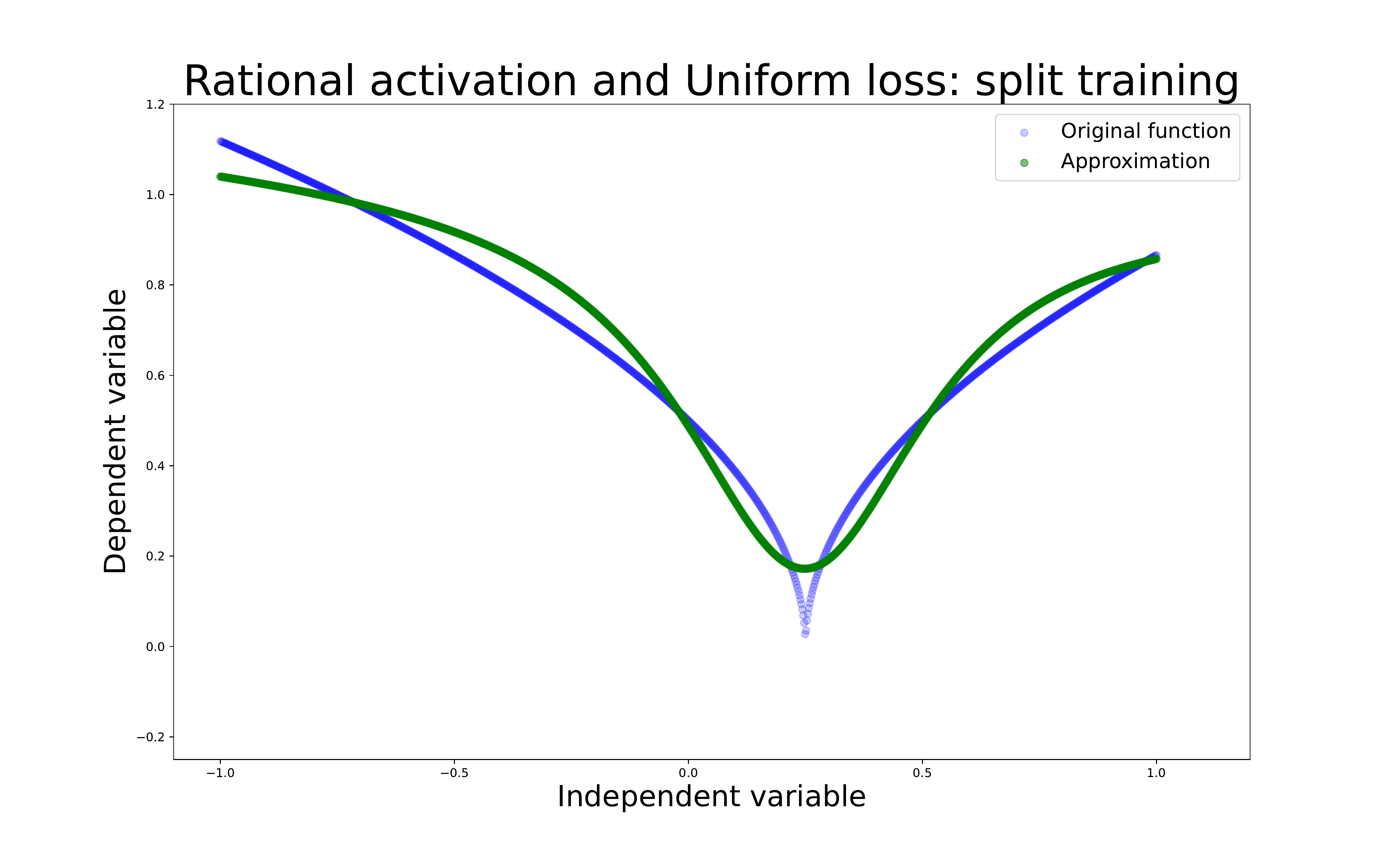}
            \includegraphics[width=40mm]{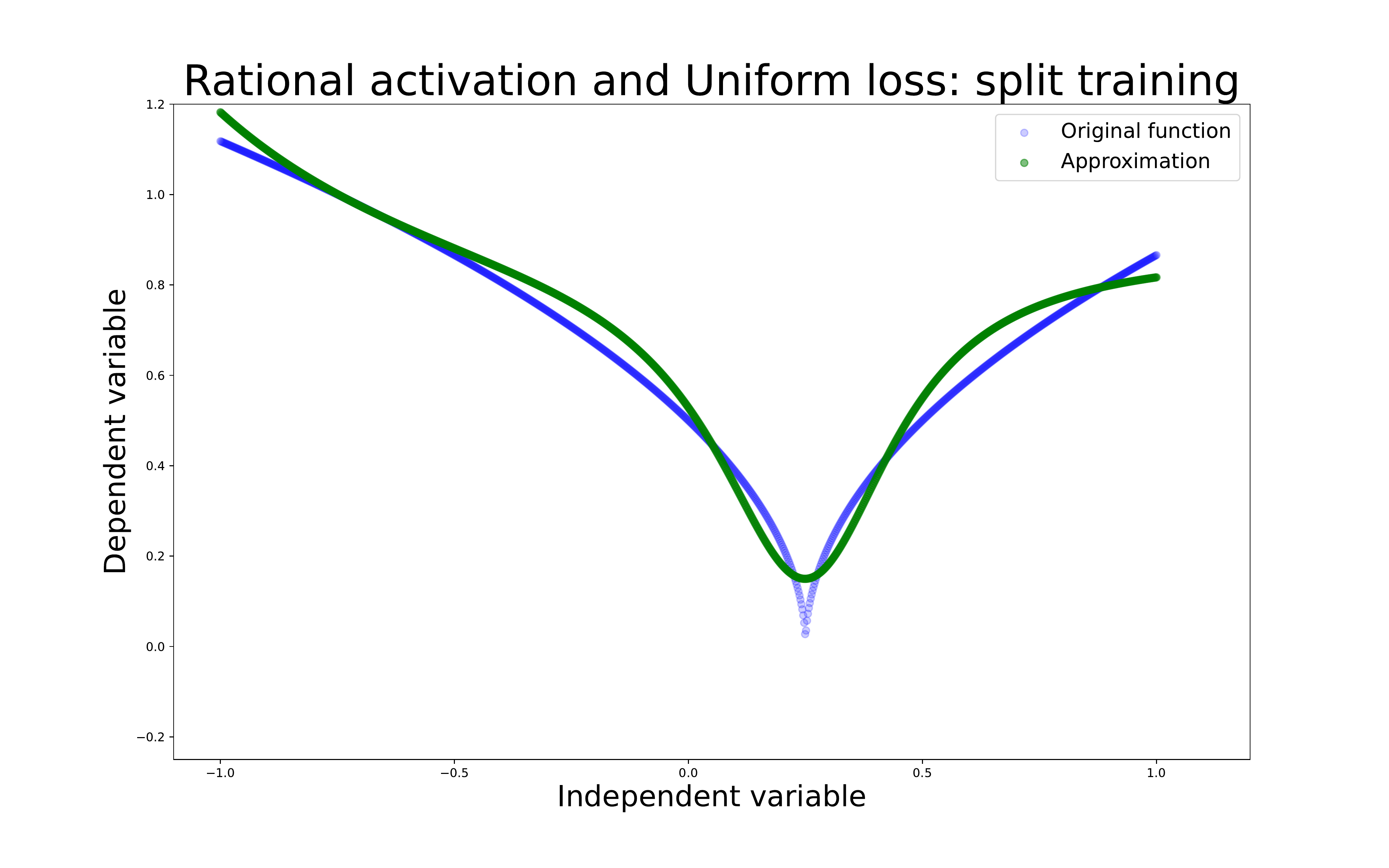}
    \caption{Approximation is computed by the split training process, 50, 100 and 200 epochs.}
    \label{fig:Rational activation with split training - set 2, epoch is 50}
\end{figure}

Figure~\ref{fig:Rational activation with split training - set 2, epoch is 50} shows that the approximation inaccuracy mostly appears at the ``difficult point''.


\paragraph{{\bf Set 3}: 10 nodes in the hidden layer, MSE loss, ADAM optimiser, split method}

When the number of nodes in the hidden layer increases to~10, the approximation accuracy in terms of the optimal loss function value is improving, but the results are not as good as they were in the case of rational activation and rational approximation to ReLU activation (table~\ref{tab:Results: experiments set 3 of NN with rational approximation and split training}).

\begin{table}
    \centering
    \begin{tabular}{|c|c|c|c|}
    \hline
    Epoch & Final loss & Minimum loss & Run time (per epoch) \\
    \hline
    50  & 0.000318 &          & 2.68s $\pm$ 131ms\\
    50  &          & 0.000318 & \\
    \hline
    100 & 0.000107 &          & 2.59s $\pm$ 79ms\\
    100  &          & 0.000107 & \\
    \hline
    200 & 0.011278 &          & 2.65s $\pm$ 305ms\\
    143 &          & 0.000064 & \\
    \hline
    \end{tabular}
    \caption{Results: experiments set 3, split method}
    \label{tab:Results: experiments set 3 of NN with rational approximation and split training}
\end{table}



\begin{figure}
    \centering
    \includegraphics[width=40mm]{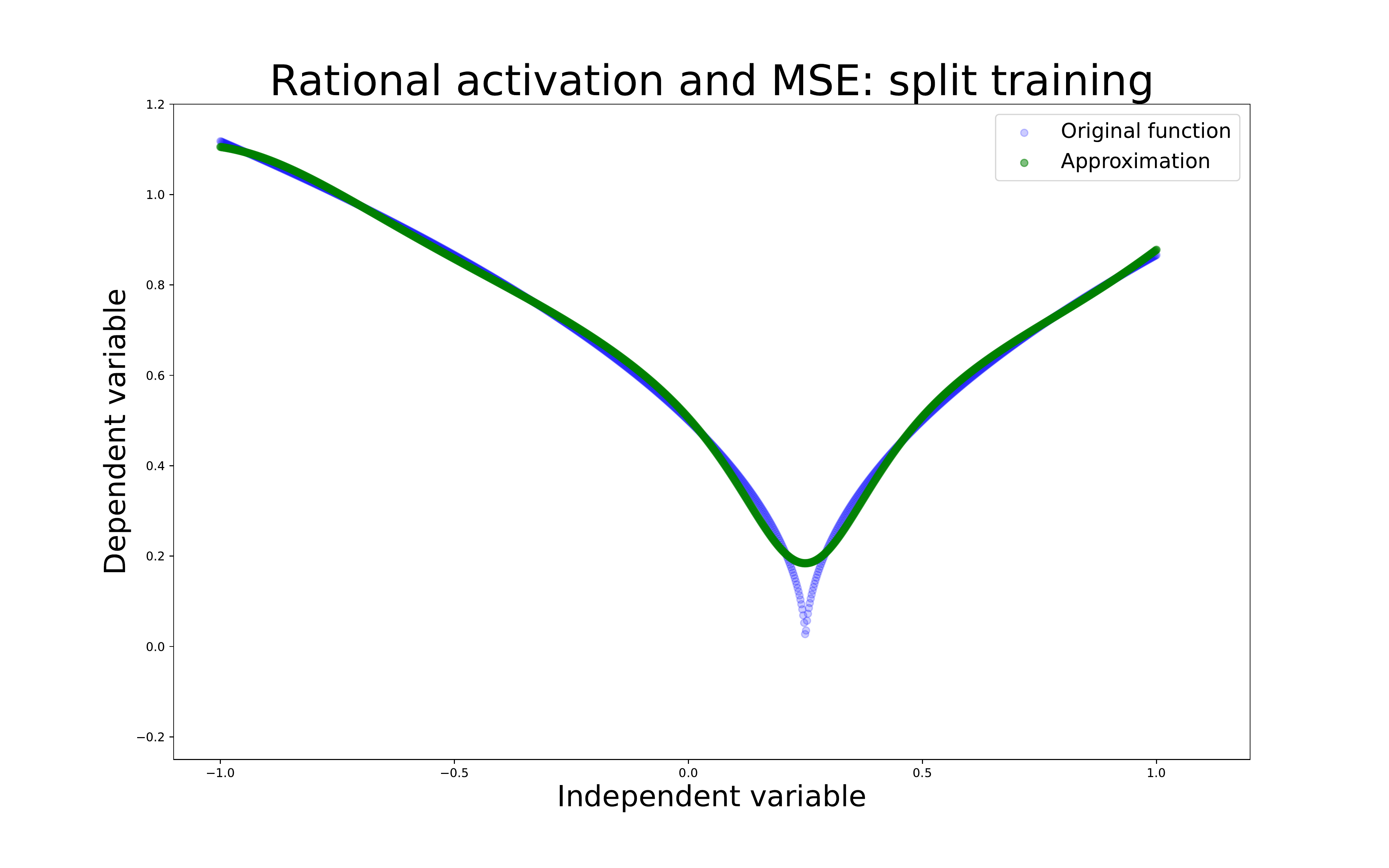}
    \includegraphics[width=40mm]{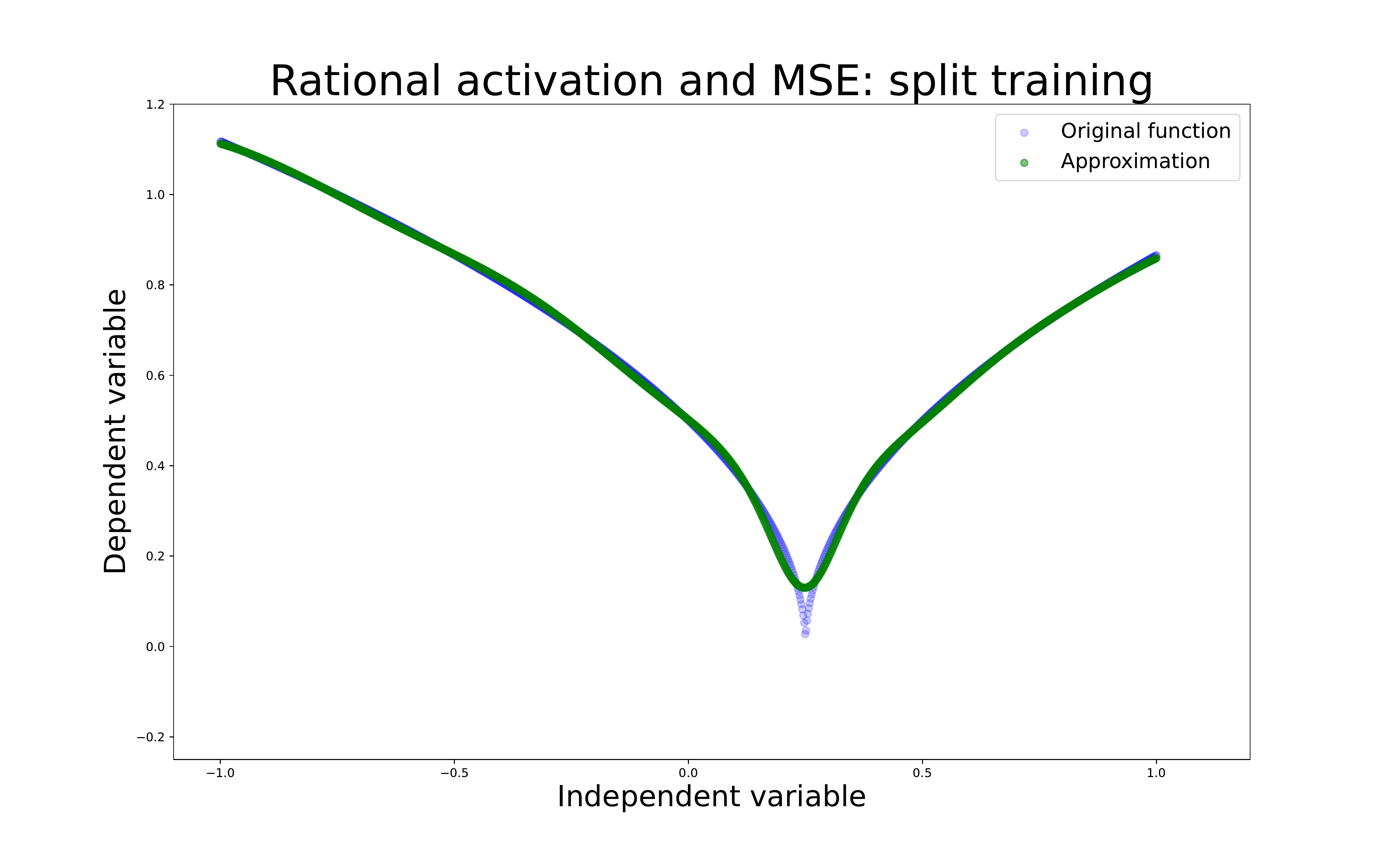}
    \includegraphics[width=40mm]{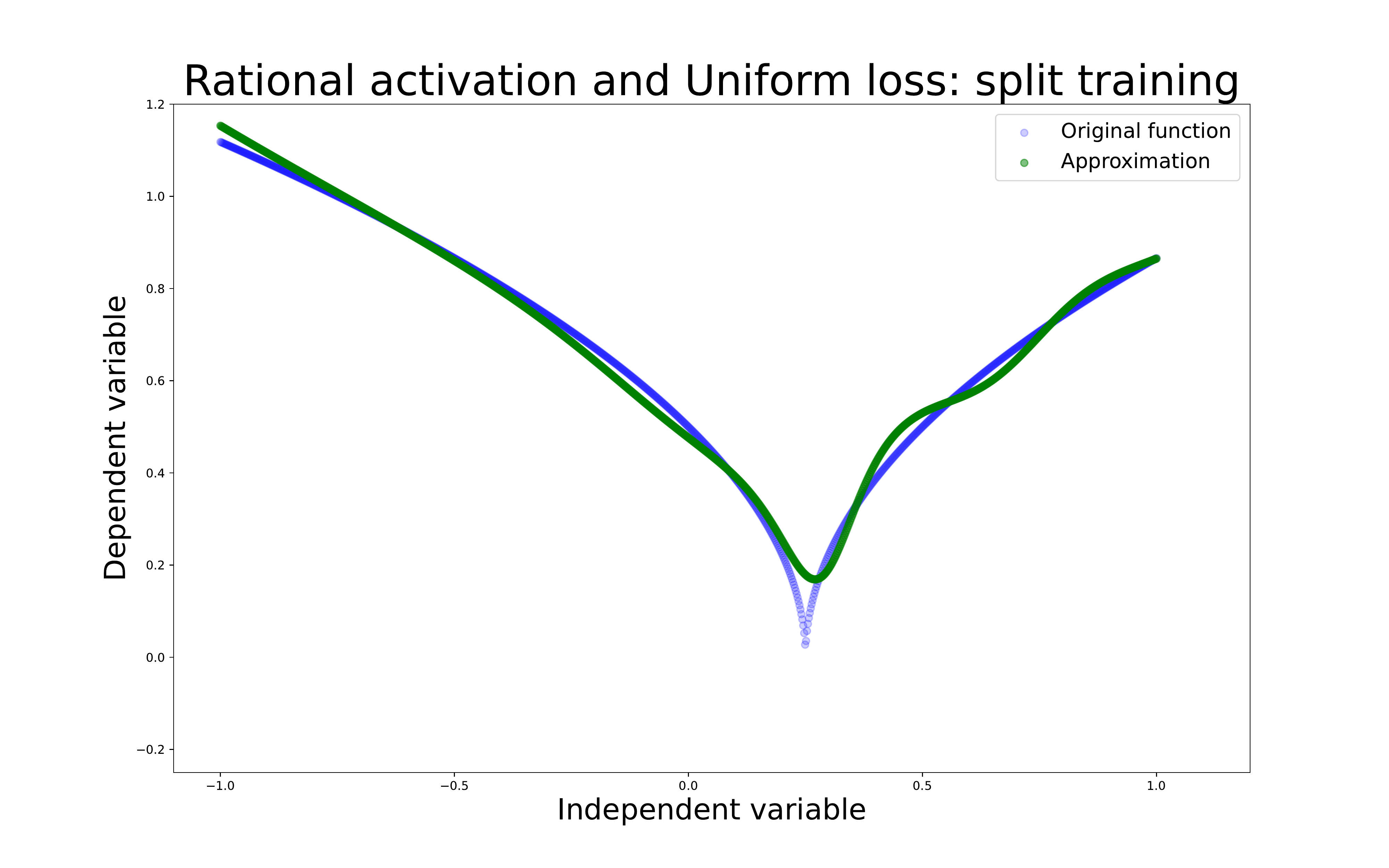}
    \caption{Approximation is computed by the split training process,  50, 100 and 200~epochs.}
    \label{fig:Rational activation with split training - set 3, epoch is 50,100,200}
\end{figure}

Figure~\ref{fig:Rational activation with split training - set 3, epoch is 50,100,200} shows that there are some oscillations when the number of epochs is~200. However, if the training process is stopped exactly at epoch 143, a better approximation can be achieved (Figure~\ref{fig:Rational activation with split training - set 3, epoch is 143}).

\begin{figure}
    \centering
    \includegraphics[width=50mm]{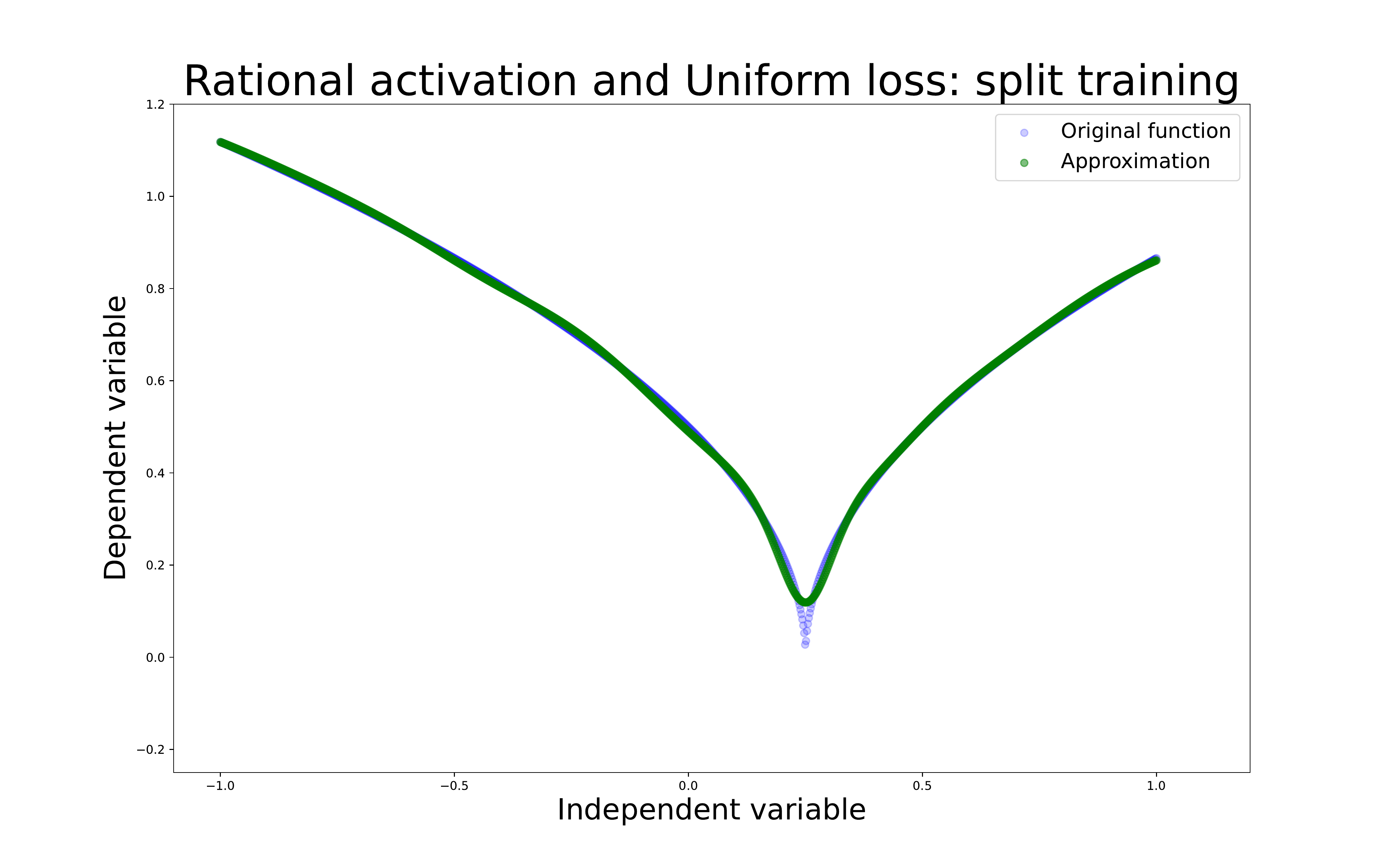}
    \caption{Approximation is computed by the split training process for 143 epochs.}
    \label{fig:Rational activation with split training - set 3, epoch is 143}
\end{figure}

\subsection{Direct rational approximation approach: the differential correction and AAA method.}

In this section, all the experiments correspond to uniform loss. We compare two main methods for constructing rational approximations: the differential correction method and the AAA method.

\subsubsection{The differential correction method}

\paragraph{{\bf Rational approximation of degree (5,4), the differential correction method.}}

When one considers a neural network with three layers and just 2~nodes in the hidden layer and the rational activation function of degree~$(3,2)$, the output is a rational function of degree $(5,4)$. Therefore, we can approximate the function $f(x) = \sqrt{|x-0.25|}$ by a rational function of degree~$(5,4)$  using the differential correction method and compare this approximation with the results obtained for uniform loss function and 2~nodes in the hidden layer. In theory, the direct approximation by a rational function of degree~$(5,4)$ is more flexible. We use Python  in our experiments. 
We also compute the rational approximations of degree~$(4,4)$ and~$(5,5)$ for comparison. We will also use these additional results to compare with the AAA method, which is expecting the same degree in the numerator and the denominator. Table~\ref{tab:Results: experiments set 1 of rational approximation through differential correction method} summarises the computational time and the optimal loss. The results demonstrate that the computational time is very small for all three settings and corresponding errors are very similar. Interestingly, the increase in the degree does not always improve the optimal loss function value. This can be explained as a result of computational instability when the dimension of the corresponding optimisation problems is increasing. 

\begin{table}
    \centering
    \begin{tabular}{|c|c|c|}
    \hline
    Degrees & Error (uniform norm) & Run time \\
    \hline
    (5,4) & 0.05598587273362793 & 2.22s\\
    \hline
    (4,4) & 0.08595288703758228 & 1.73s\\
    \hline
    (5,5) & 0.03475461848979777 & 2.59s\\
    \hline
    \end{tabular}
    \caption{Results: direct rational approximation for degrees $(5,4)$, $(4,4)$ and $(5,5)$}
    \label{tab:Results: experiments set 1 of rational approximation through differential correction method}
\end{table}
Figure~\ref{fig:Rational approximation of degree (4,4), (5,4) and (5,5) via the differential correction method} depicts the function and the corresponding rational approximation of degree~(4,4), (5,4) and (5,5).

\begin{figure}
    \centering
    \includegraphics[width=40mm]{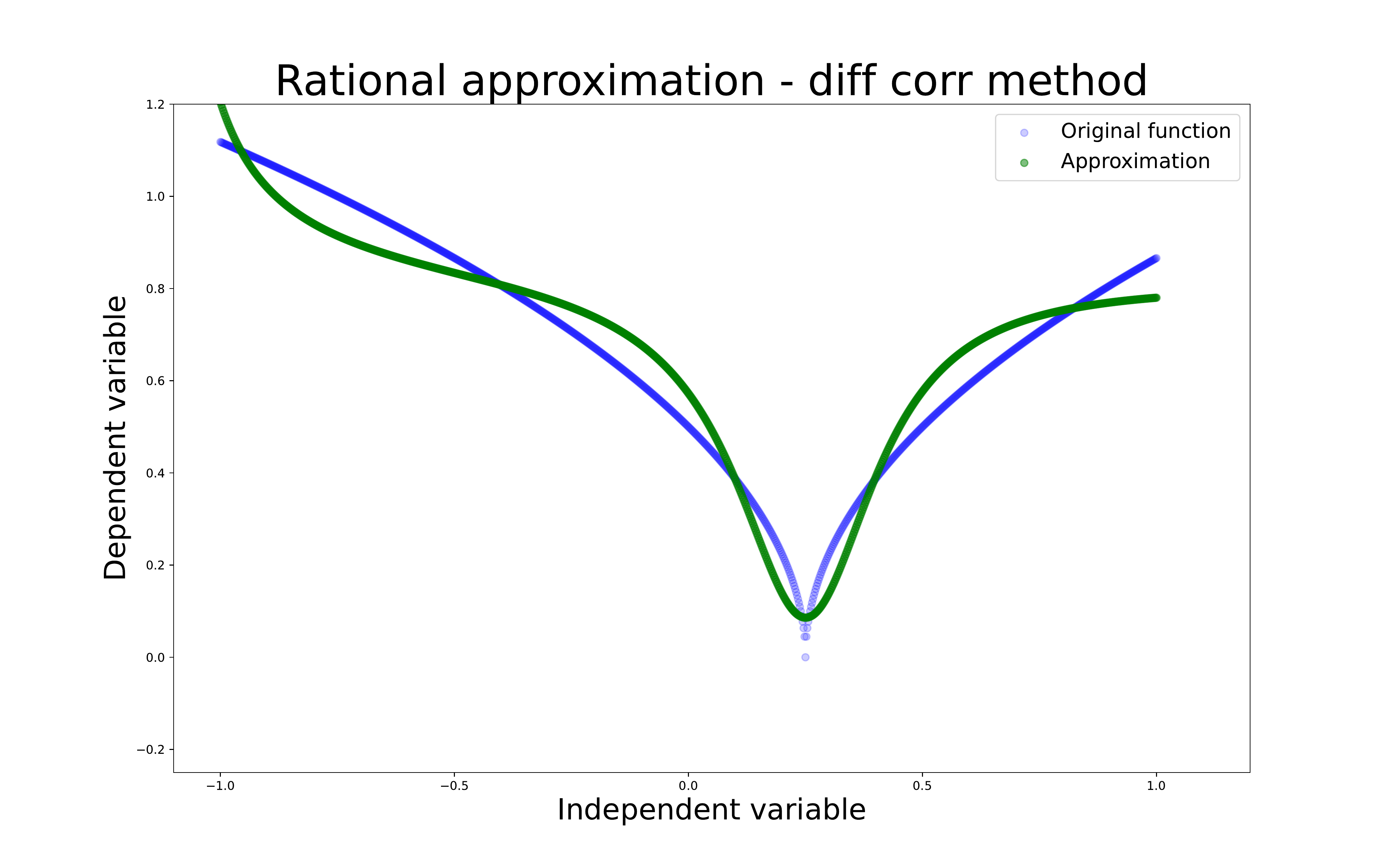}
    \includegraphics[width=40mm]{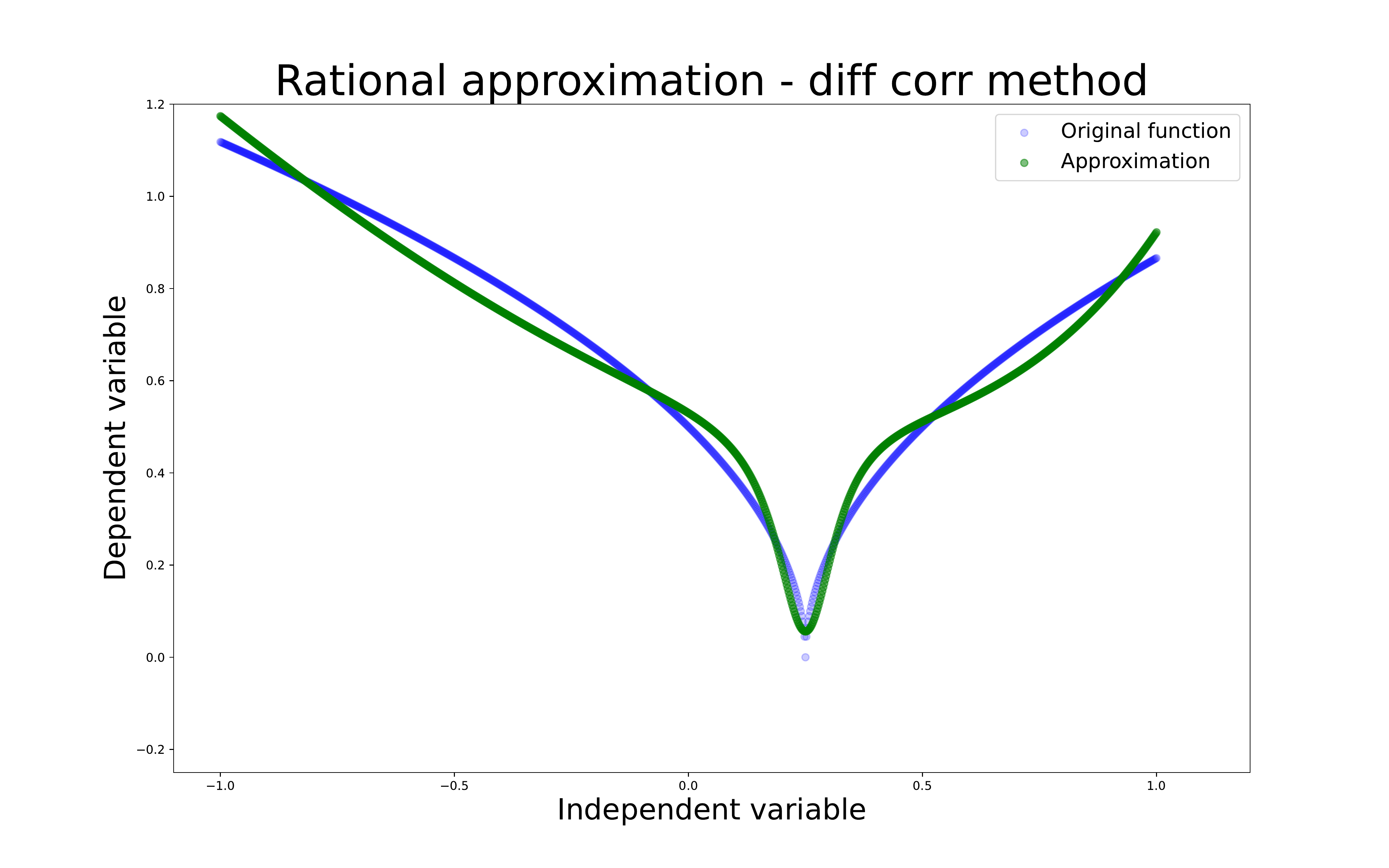}
    \includegraphics[width=40mm]{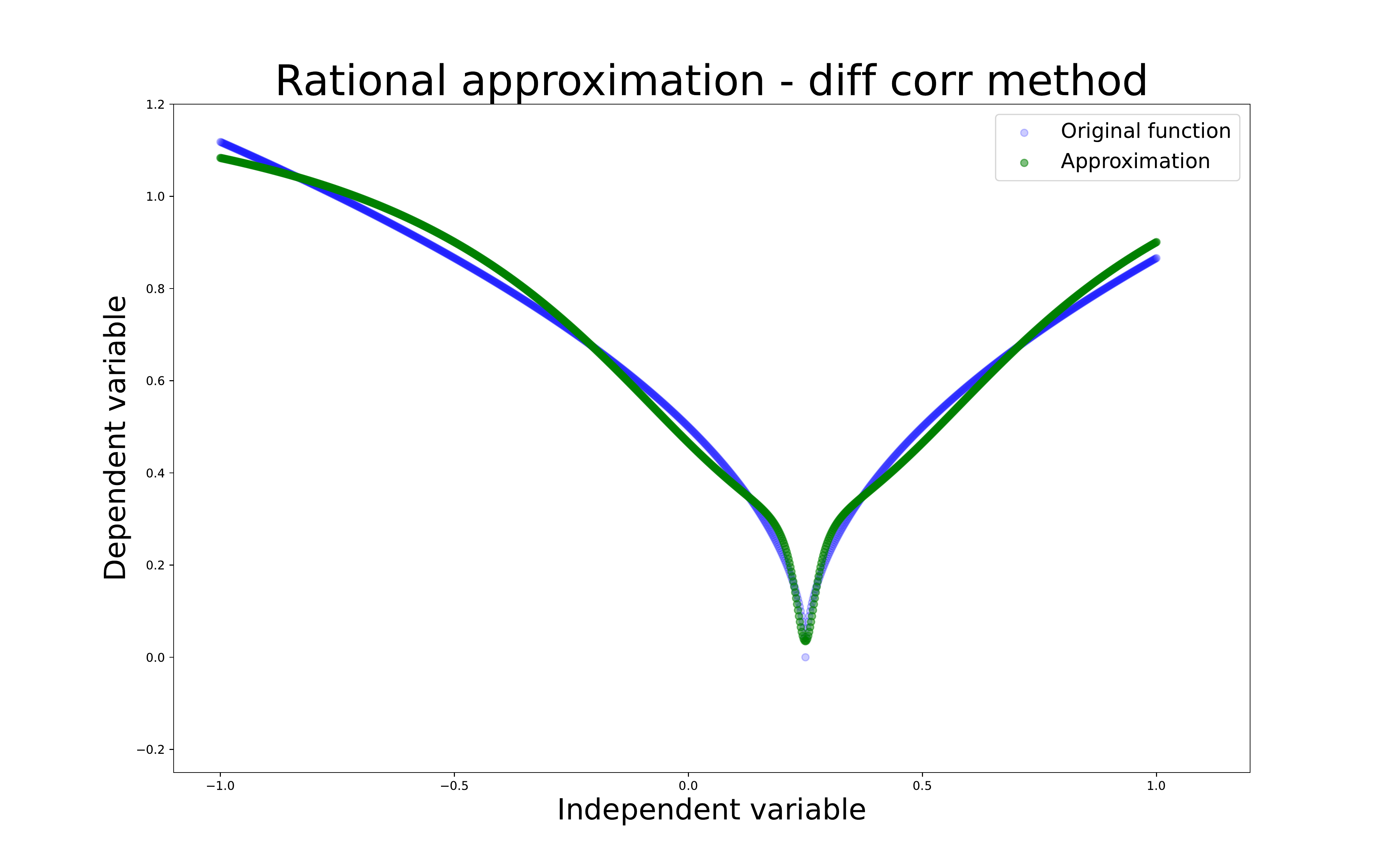}
    \caption{Approximation by the differential correction method, rational approximation degree is $(4,4)$, $(5,4)$ and $(5,5)$ respectively.}
    \label{fig:Rational approximation of degree (4,4), (5,4) and (5,5) via the differential correction method}
\end{figure}

\subsubsection{The AAA method}

\paragraph{{\bf Rational approximation of degree (5,4) with the AAA method}}
In order to compare the differential correction method with AAA, we apply the AAA method in the same experimental settings as they are in the case of the differential correction method. Due to the limitations of the AAA method, we cannot use different degrees for the numerator and the denominator. Hence, for these experiments, we use the rational approximation of degree (4,4) and the rational approximation of degree (5,5). The codes are developed in Python by C. Hofreither~\cite{Hofreither21}. 


Table~\ref{tab:Results: experiments set 1 of rational approximation through AAA method} summarises the results for the AAA method. The AAA method is designed for uniform loss, but we report the corresponding MSE loss as well for comparison. One can see that the optimal loss is high. This is due to instability. Unfortunately, this instability  cannot be fixed by limiting the condition number (differential correction).  


\begin{table}
    \centering
    \begin{tabular}{|c|c|c|c|}
    \hline
    Degrees & Error type & Error & Run time \\
    \hline
    \multirow{2}{4em}{(4,4)} & Uniform error &
    1707.081229434771   & \multirow{2}{4em}{53.1ms}\\
    \cline{2-3}
     & MSE & 2658.9489333999118 & \\
    \hline
    \multirow{2}{4em}{(5,5)} & Uniform error &
    0.14229129790010392   & \multirow{2}{4em}{85.5ms}\\
    \cline{2-3}
     & MSE & 1.1112525485550502 & \\
    \hline
    \end{tabular}
    \caption{Results: AAA method, degrees (4,4) and (5,5)}
    \label{tab:Results: experiments set 1 of rational approximation through AAA method}
\end{table}

Figure~\ref{fig:Rational approximation of degree (4,4) and (5,5) via AAA method} highlights the difficulties of approximations. These difficulties are especially prominent in the case of degree~(4,4).


\begin{figure}
    \centering
    \includegraphics[width=50mm]{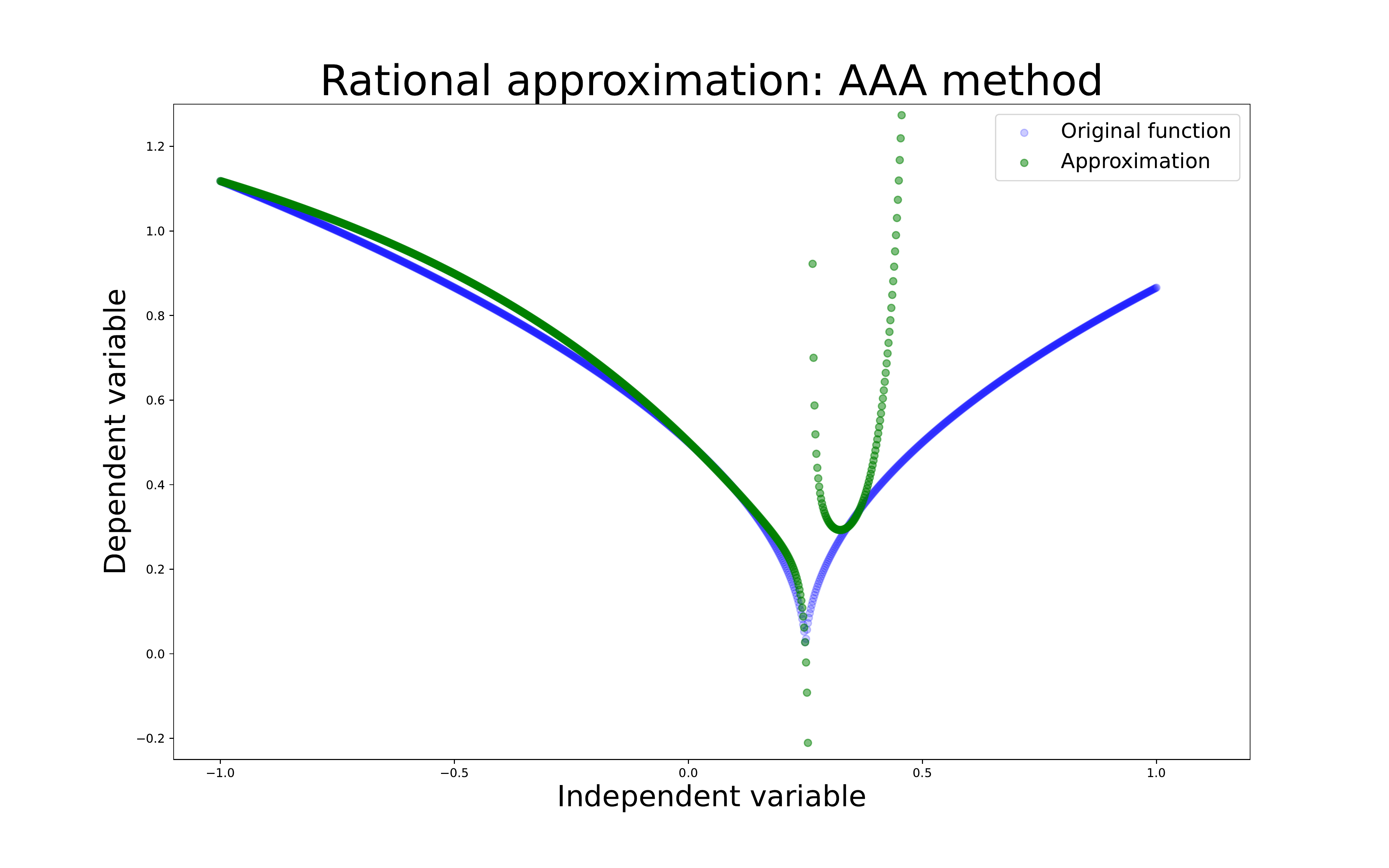}
    \includegraphics[width=50mm]{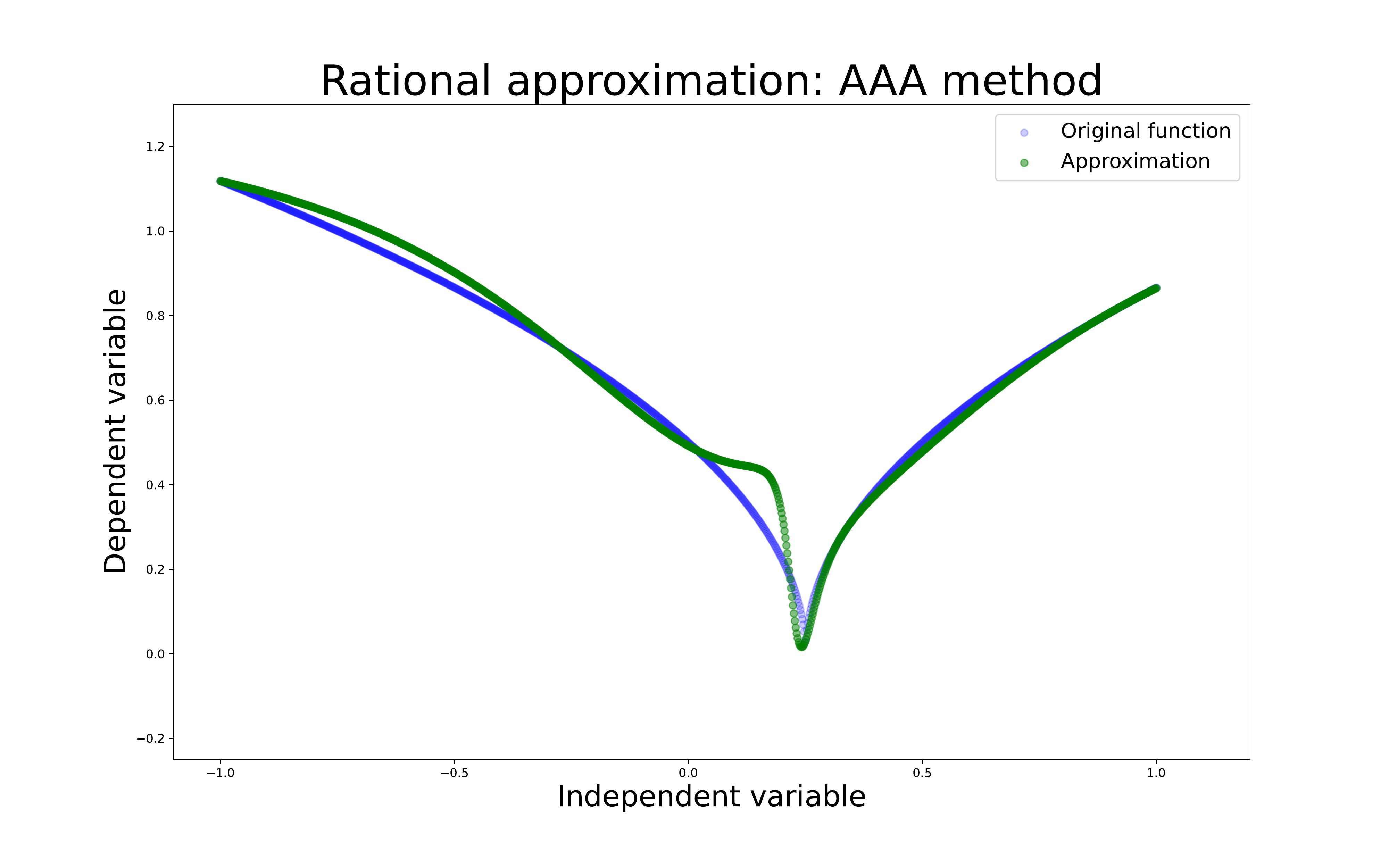}
    \caption{Approximation is computed by the AAA method, degree (4,4) and (5,5).}
    \label{fig:Rational approximation of degree (4,4) and (5,5) via AAA method}
\end{figure}

\end{document}